\documentclass[twoside]{article}

\usepackage[accepted]{aistats2024}

\usepackage{amsthm}
\usepackage{amsmath}
\usepackage{amssymb}
\usepackage{graphicx}
\usepackage{mathtools}
\usepackage{enumerate}
\usepackage{enumitem}
\usepackage{footnote}
\usepackage{float}
\usepackage{xspace}
\usepackage{multirow}
\usepackage{nicefrac}
\usepackage{wrapfig}
\usepackage{framed}
\usepackage{url}
\usepackage[colorlinks=true, linkcolor=blue, citecolor=blue]{hyperref}
\usepackage{blkarray}
\PassOptionsToPackage{round}{natbib}
\usepackage{algorithmic}

\usepackage{accents}

\usepackage{tikz}
\usetikzlibrary{arrows,chains,matrix,positioning,scopes}

\let\hat\widehat

\newtheorem{lemma}{Lemma}
\newtheorem{theorem}{Theorem}
\newtheorem{cor}{Corollary}

\newtheorem{definition}{Definition}
\newtheorem{assumption}{Assumption}

\DeclareMathOperator*{\minimize}{minimize}
\DeclareMathOperator*{\maximize}{maximize}
\DeclareMathOperator*{\subject}{subject~to}

\newcommand{\calA}{{\mathcal{A}}}

\newcommand{\calS}{{\mathcal{S}}}

\DeclareMathOperator*{\argmin}{argmin}
\DeclareMathOperator*{\argmax}{argmax}

\newcommand{\norm}[1]{\left\|{#1}\right\|}

\DeclareFontFamily{OMX}{MnSymbolE}{}
\DeclareFontShape{OMX}{MnSymbolE}{m}{n}{
    <-6>  MnSymbolE5
   <6-7>  MnSymbolE6
   <7-8>  MnSymbolE7
   <8-9>  MnSymbolE8
   <9-10> MnSymbolE9
  <10-12> MnSymbolE10
  <12->   MnSymbolE12}{}
\DeclareSymbolFont{mnlargesymbols}{OMX}{MnSymbolE}{m}{n}
\SetSymbolFont{mnlargesymbols}{bold}{OMX}{MnSymbolE}{b}{n}
\DeclareMathDelimiter{\llangle}{\mathopen}{mnlargesymbols}{'164}{mnlargesymbols}{'164}
\DeclareMathDelimiter{\rrangle}{\mathclose}{mnlargesymbols}{'171}{mnlargesymbols}{'171}

\usepackage{prettyref}

\newcommand{\savehyperref}[2]{\texorpdfstring{\hyperref[#1]{#2}}{#2}}
\newrefformat{eq}{\savehyperref{#1}{\textup{(\ref*{#1})}}}
\newrefformat{eqn}{\savehyperref{#1}{Equation~\eqref{#1}}}
\newrefformat{lem}{\savehyperref{#1}{Lemma~\ref*{#1}}}
\newrefformat{def}{\savehyperref{#1}{Definition~\ref*{#1}}}
\newrefformat{line}{\savehyperref{#1}{line~\ref*{#1}}}
\newrefformat{thm}{\savehyperref{#1}{Theorem~\ref*{#1}}}
\newrefformat{corr}{\savehyperref{#1}{Corollary~\ref*{#1}}}
\newrefformat{cor}{\savehyperref{#1}{Corollary~\ref*{#1}}}
\newrefformat{sec}{\savehyperref{#1}{Section~\ref*{#1}}}
\newrefformat{app}{\savehyperref{#1}{Appendix~\ref*{#1}}}
\newrefformat{assum}{\savehyperref{#1}{Assumption~\ref*{#1}}}
\newrefformat{ex}{\savehyperref{#1}{Example~\ref*{#1}}}
\newrefformat{fig}{\savehyperref{#1}{Figure~\ref*{#1}}}
\newrefformat{alg}{\savehyperref{#1}{Algorithm~\ref*{#1}}}
\newrefformat{rem}{\savehyperref{#1}{Remark~\ref*{#1}}}
\newrefformat{conj}{\savehyperref{#1}{Conjecture~\ref*{#1}}}
\newrefformat{prop}{\savehyperref{#1}{Proposition~\ref*{#1}}}
\newrefformat{proto}{\savehyperref{#1}{Protocol~\ref*{#1}}}
\newrefformat{prob}{\savehyperref{#1}{Problem~\ref*{#1}}}
\newrefformat{claim}{\savehyperref{#1}{Claim~\ref*{#1}}}
\newrefformat{que}{\savehyperref{#1}{Question~\ref*{#1}}}
\usepackage{microtype}
\usepackage{graphicx}
\usepackage{subfigure}
\usepackage{booktabs} 
\usepackage[colorinlistoftodos]{todonotes}
\usepackage{multicol}
\usepackage{algorithm}
\usepackage{float}
\usepackage{xcolor}
\usepackage{tikz,pgfplots}
\usetikzlibrary{calc}
\newcommand{\DefinedAs}[0]{\mathrel{\mathop:}=}
\begin{document}

%

%

\twocolumn[

\aistatstitle{Resilient Constrained Reinforcement Learning
}

\aistatsauthor{ Dongsheng~Ding \And Zhengyan~Huan \And  Alejandro Ribeiro}

\aistatsaddress{ \{dongshed, zhhuan,  aribeiro\}@seas.upenn.edu \\University of Pennsylvania } 
]

\begin{abstract}
      We study a class of constrained reinforcement learning (RL) problems in which multiple constraint specifications are not identified before training. It is challenging to identify appropriate constraint specifications due to the undefined trade-off between the reward maximization objective and the constraint satisfaction, which is ubiquitous in constrained decision-making. To tackle this issue, we propose a new constrained RL approach that searches for policy and constraint specifications together. This method features the adaptation of relaxing the constraint according to a relaxation cost introduced in the learning objective. Since this feature mimics how ecological systems adapt to disruptions by altering operation, our approach is termed as resilient constrained RL. Specifically, we provide a set of sufficient conditions that balance the constraint satisfaction and the reward maximization in notion of resilient equilibrium, propose a tractable formulation of resilient constrained policy optimization that takes this equilibrium as an optimal solution, and advocate two resilient constrained policy search algorithms with non-asymptotic convergence guarantees on the optimality gap and constraint satisfaction. Furthermore, we demonstrate the merits and the effectiveness of our approach in computational experiments.      
\end{abstract}

\section{INTRODUCTION}

Constrained reinforcement learning (RL) is a constrained control problem in which an agent aims to maximize its expected cumulative reward while satisfying a given constraint by interacting with an environment over time. 
Multiple requirements are of growing interest in constrained RL towards designing an agent to meet more than one constraint, e.g., resource allocation for many users~\cite{de2021constrained} and safe learning in robotics~\cite{brunke2022safe}. Real-world constrained RL often engages practical problems with unwell-specified requirements, e.g., human-satisfaction in human-robot interaction~\cite{el2020towards} and safety level of robotic agents~\cite{zhang2020cautious}. Hence, it is challenging to determine reasonable constraint specifications for making trade-off between reward maximization and constraint satisfaction.  

Although reward shaping has been widely used to aggregate multiple requirements into a single reward, e.g.,~\cite{perez2021robot}, it doesn't guarantee constraint satisfaction for each requirement. A known reason for this single-reward failure is that the solutions generated by standard RL algorithms do not necessarily satisfy required constraints, which is known as ``scalarization fallacy"~\cite{szepesvari2020,zahavy2021reward,calvo2023state}. Therefore, it is crucial to directly impose the constraints that result from multiple requirements~\cite{roy2022direct}, which has been studied by a lot of recent efforts, e.g.,~\cite{chow2017risk,paternain2019constrained,ding2020natural}. However, such results are based on known feasible constraints, not applicable in the situations with unknown constraint specifications. 

To fill this gap, we aim to automate the constraint specifications during constrained RL training by facilitating the trade-off between reward maximization and constraint satisfaction. The focal RL environment of this paper is the constrained Markov decision process (MDP) that constrains expected cumulative utilities~\cite{altman1999constrained}, which has been widely-used in many constraint-rich domains, e.g.,  resource allocation, robotic planning, and financial management; see more in~\cite{garcia2015comprehensive,de2021constrained,gu2022review,brunke2022safe}. 

\paragraph{Contribution.} 

Our contributions are threefold. 

\begin{itemize}
    \item We first introduce nominal constraints that are possibly infeasible, so that they can be relaxed to compromise reward maximization for constraint satisfaction. We provide the sensitivity analysis of the optimal reward value function to the perturbations in constraints. Since the compromise mimics how ecological systems \emph{adapt to disruptions} by changing operating conditions, we term this as \emph{resilient} constrained policy optimization, and broadly as \emph{resilient} constrained RL. 
    \item To specify the levels or thresholds of constraints, we introduce a user-defined cost function that establishes a price for relaxing nominal constraints, and exploit the relative difficulty of relaxing different constraints to define a trade-off solution: resilient equilibrium. We provide a tractable formulation of resilient constrained policy optimization that takes this equilibrium as an optimal solution, and establish its duality theory under less restrictive feasibility assumption. 
    \item To find policy and constraint specification, we extend two non-resilient policy gradient algorithms for our resilient constrained policy optimization problem, and prove that they converge to a optimal solution with non-asymptotic convergence guarantees on the optimality gap and constraint satisfaction.
    To the best of our knowledge, for the first time we establish provably resilient constrained policy search algorithms against uncertain constraints.
    Moreover, we provide computational experiments to show the merits and the effectiveness of our approach. 
\end{itemize}

\paragraph{Related Work.} Our problem formulation is based on the constrained MDP framework~\cite{altman1999constrained}. Constrained MDPs with well-specified constraints are relatively well-studied in constrained RL literature, e.g., policy gradient methods~\cite{ding2022convergence} and model-based algorithms~\cite{efroni2020exploration}, under the strict feasibility assumption on constraints; see more related works in this line in~\cite{gu2022review}. However, it is intractable to determine the feasibility of constraints in many scenarios, e.g., budget distribution for many users~\cite{boutilier2016budget,vora2023welfare} and online budget level~\cite{diaz2023flexible} are unknown for feasibility-checking, and safety constraints in training are different from those for real robotics, which are expensive to model~\cite{kaspar2020sim2real}. Although this issue can be alleviated by some heuristic methods in the references aforementioned, their optimality and constraint satisfaction are not well-understood. We note that the reward and constraint trade-off essentially reduces to the sensitivity of the optimal reward value function, which was only studied using the parameter  perturbations of a constrained MDP~\cite{altman1991sensitivity,altman1993stability}. Compared with this line of work, in this paper we exploit the sensitivity analysis of the optimal reward value function against the perturbations in constraints to strike a balance between reward maximization and constraint satisfaction. We also establish two constrained policy search algorithms for finding optimal policy and constraint specification simultaneously, with theoretical guarantees.

Our work is also pertinent to recent efforts of augmenting a RL agent with the adaptation to the interference that is potentially catastrophic to system~\cite{yang2021causal,huang2022reinforcement}. This capability is often termed as ``resilience'' that draws the ability of ecological systems to adapt to disrupted environment~\cite{holling1973resilience,holling1996engineering}. Resilience to perturbations in agent-environment interaction has been studied in several prior works~\cite{yang2021causal,phan2021resilient,gao2022resilient}; yet the resilience to corrupted constraints on system or performance was not studied. Recently, the adaptation of trained policy to unknown constraint specifications is investigated in control~\cite{chamon2020resilient}, constrained learning~\cite{hounie2023resilient}, and constrained offline RL~\cite{liu2023constrained,zhang2023saformer}. In contrast, in this paper we investigate ``resilience'' for constrained policy optimization and provably convergent constrained policy search algorithms, with a focus on the adaptation of trained policy to unknown constraint specifications.   


    

\section{CONSTRAINED MDP}\label{sec:CMDPs}

We consider an infinite-horizon constrained MDP, 
\[
    \text{CMDP}
    \left(\, 
    S,\, A,\, P,\, r, \, \{u_i\}_{i\,=\,1}^m, \, \{b_i\}_{i\,=\,1}^m, \, \gamma, \, \rho 
    \,\right)
\]
where $S$ and $A$ are finite state/action spaces, $P$ is a transition kernel that specifies the probability $P(s'\,\vert\,s,a)$ from state $s$ to next state $s'$ under action $a\in A$, $r$, $u_i$: $S\times A\to [0,1]$ are reward/utility functions, $b_i$ is a constraint threshold for the $i$th utility, $\gamma\in [0,1)$ is a discount factor, and $\rho$ is an initial distribution. A stochastic policy $\pi$: $S\to \Delta(A)$ determines a probability distribution $\Delta(A)$ over the action space $A$ based on the current state, i.e., $a_t\sim\pi(\cdot\,\vert\,s_t)$ at time $t$. Let $\Pi$ be the set of all possible stochastic policies. A policy $\pi \in \Pi$, together with the initial state distribution $\rho$, induces a distribution over trajectories $\tau = \{(s_t, a_t, r_t, \{u_{i,t}\}_{i\,=\,1}^m)\}_{t\,=\,0}^\infty$, where $s_0\sim\rho$, $a_t\sim\pi(\cdot\,\vert\,s_t)$, $r_t = r(s_t,a_t)$, $u_{i,t} = u_i(s_t,a_t)$, and $s_{t+1}\sim P(\cdot\,\vert\,s_t, a_t)$ for all $t\geq 0$.

Given a policy $\pi$, the value functions $V_r^\pi$, $V_{u_i}^\pi$: $S\to \mathbb{R}$ associated with the reward $r$ or the utility $u_i$ are given by the expected sums of discounted rewards or utilities received under policy $\pi$, respectively,
\[
    \begin{array}{rcl}
         V_r^\pi (s)
         & \DefinedAs &
         \displaystyle
         \mathbb{E}\left[ \sum_{t\,=\,0}^\infty \gamma^t r(s_t, a_t) \,\vert\, \pi, s_0 = s \right] 
    \end{array}
\]
where $\mathbb{E}$ is expected over the randomness in the trajectory $\tau$ induced by $\pi$; similarly, we define $V_{u_i}^\pi (s)$ for the utility $u_i$. The expected values over the initial distribution $\rho$ are given by $V_r^\pi(\rho) = \mathbb{E}_{s\,\sim\,\rho} \left[\, V_r^\pi(s) \,\right]$ and $V_{u_i}^\pi(\rho) = \mathbb{E}_{s\,\sim\,\rho} \left[\, V_{u_i}^\pi(s) \,\right]$. It is useful to introduce the discounted state visitation distribution, $d_{s_0}^\pi(s) = (1-\gamma)\sum_{t\,=\,0}^\infty \gamma^t \text{Pr}(s_t = s\,\vert\,\pi,s_0)$ which adds up discounted probabilities of visiting $s$ in the execution of $\pi$ starting from $s_0$. Denote $d_\rho^\pi(s) \DefinedAs \mathbb{E}_{s_0\,\sim\,\rho}[\, d_{s_0}^\pi(s) \,]$ and thus $d_\rho^\pi(s) \geq (1-\gamma)\rho(s)$ for any $\rho$ and $s$. Furthermore, for the reward $r$, we introduce the state-action value function $Q_r^\pi$: $S\times A\to \mathbb{R}$ when the agent begins with a state-action pair $(s,a)$ and follows a policy $\pi$, together with its advantage function $A_r^\pi$: $S\times A\to\mathbb{R}$,
\[
    \begin{array}{rcl}
         Q_r^\pi (s,a)
         & \DefinedAs &
         \displaystyle
         \mathbb{E}\left[ \sum_{t\,=\,0}^\infty \gamma^t r(s_t, a_t) \,\vert\, \pi, s_0 = s, a_0 = a\right] 
         \\[0.4cm]
         A_r^\pi (s,a)
         & \DefinedAs & 
         \displaystyle
         Q_r^\pi (s,a) - V_r^\pi(s). 
    \end{array}
\]
Similarly, we define $Q_{u_i}^\pi$, $A_{u_i}^\pi$: $S\times A\to\mathbb{R}$ for $u_i$.

The constrained MDP aims to find a policy that maximizes the reward value function $V_r^\pi(\rho)$ while the utility value function $V_{u_i}^\pi(\rho)$ is above some threshold $b_i$,
\begin{equation}\label{eq:CMDP}
    \begin{array}{rl}
         \displaystyle \maximize_{\pi\,\in\,\Pi} & V_r^\pi(\rho)
         \\
         \subject&  
         V_{u_i}^\pi(\rho) \; \geq \; b_i, \; i = 1,\ldots,m
    \end{array}
\end{equation}
where $b_i$ is the specified threshold \emph{priori} for the $i$th utility value function. Since $V_r^\pi(\rho)$, $V_{u_i}^\pi(\rho)\in [0,1/(1-\gamma)]$, we assume $b_i\in(\,0,1/(1-\gamma)\,]$. Thus, by taking $g_i$: $S\times A\to [-1,1]$ with $g_i = u_i - (1-\gamma)b_i$, equivalently we translate the constraint $V_{u_i}^\pi(\rho) \geq b_i$ into $V_{g_i}^\pi(\rho) \geq 0$, which is our focal constraint. Those utility constraints often result from additional requirements on the system operation, e.g., budget or safety constraints~\cite{boutilier2016budget,paternain2022safe}. We denote the optimal value for Problem~\eqref{eq:CMDP} by $V^\star$ which takes $V^\star = V_r^\star(\rho)$ at an optimal policy $\pi^\star$ if it is feasible; $V^\star = -\infty$ otherwise. 

Although an optimal policy always exists in the unconstrained case, i.e., $m=0$ in Problem~\eqref{eq:CMDP}, it is not necessarily true in the constrained case because of potentially infeasible constraints. Thus, it is important to specify relevant constraints for Problem~\eqref{eq:CMDP}. In many scenarios, how to specify the constraint thresholds $\{b_i\}_{i\,=\,1,\ldots,m}$ is not known \emph{priori}. For instance, the threshold $b_i$ means the $i$th user's negative budget that is often time-varying in resource allocation~\cite{boutilier2016budget,vora2023welfare}; to trade-off many rewards in preference-based RL~\cite{eysenbach2019reinforcement,liang2022reward}, the threshold $b_i$ is often unknown for preference $i$; see their details and more examples in Section~\ref{subsec:CMDP_without_prior}. In practice, only a nominal constraint specification is given with unknown feasibility, and we have to relax (or tighten) the nominal constraints for guarding the feasibility. With a slight abuse of notation, we use notation $\{b_i\}_{i\,=\,1}^m$ to denote the nominal constraint specifications that might be too conservative (or loose).

To study the effect of constraint specifications, we form a variant of constrained MDP with flexible constraints,
\begin{equation}\label{eq:CMDP_relaxed}
    \begin{array}{rl}
         \displaystyle  
         \!\!\!\!
         \!\!\!\!
         V^\star(\xi)
         \;\DefinedAs\;
         \maximize_{\pi\,\in\,\Pi} & V_r^\pi(\rho)
         \\
         \!\!\!\!
         \!\!\!\!
         \subject&  
         V_{g_i}^\pi(\rho) \; \geq \; \xi_i, \; i = 1,\ldots,m
    \end{array}
\end{equation}
where $\xi \in \mathbb{R}^m$ is the \emph{unknown} perturbation that relaxes the constraint when $\xi_i<0$ (or tightens the $i$th inequality constraint when $\xi_i>0$), and $V^\star(\xi)$ is the primal value function: $V^\star(\xi) = V_r^\star(\rho)$ at an optimal policy $\pi^\star(\xi)$ if it is feasible; $V^\star(\xi) = -\infty$ otherwise. Since $V^\star (0) = V^\star$ for $\xi = 0$, Problem~\eqref{eq:CMDP} is our nominal problem, and Problem~\eqref{eq:CMDP_relaxed} is our perturbed problem.

Since $|V_{g_i}^\pi(\rho)| \leq 1/(1-\gamma)$, it is natural to restrict $|\xi_i|\leq 1/(1-\gamma)$ since Problem~\eqref{eq:CMDP_relaxed} is infeasible when $\xi_i>1/(1-\gamma)$ for all $i=1,\ldots,m$, and it is unconstrained when $\xi_i<-1/(1-\gamma)$ for all $i=1,\ldots,m$. Denote $\mathbb{R}_{\gamma}^m \DefinedAs\{ \xi\in\mathbb{R}^m\,\vert\, |\xi_i|\leq 1/(1-\gamma)\}$. 
Let the domain of $V^\star(\xi)$ be $\Xi$, which is the set of all $\xi$ for which the constraint set $\{ \pi\in\Pi\,\vert\, V_{g_i}^\pi(\rho)\geq \xi_i, i = 1,\ldots\,m\}$ is non-empty, or equivalently, $\Xi \DefinedAs \{ \xi\in\mathbb{R}_\gamma^m\,\vert\, V^\star(\xi)>-\infty \}$.  



It's known from the non-convexity of value functions in policy~\cite{agarwal2021theory} that Problem~\eqref{eq:CMDP_relaxed} is non-convex. Nevertheless, we prove that the primal function $V^\star(\xi)$ has several nice properties inherited from the duality analysis. Lemma~\ref{lem:primal_function} shows that the primal function $V^\star(\xi)$ enjoys monotonicity and it is concave over the domain $\Xi$; see Appendix~\ref{app:CMDPs} for proof.

\begin{lemma}[Coordinate-Wise Monotonicity and Concavity]
    \label{lem:primal_function}
    For Problem~\eqref{eq:CMDP_relaxed},
    (i) the primal function $V^\star(\xi)$ is monotonically non-increasing with respect to the coordinates of $\xi \in \Xi$, i.e., $V^\star(\xi) \leq V^\star(\xi')$ when $\xi_j > \xi_j'$ for some $j$ and $\xi_i = \xi_i'$ for $i\neq j$;
    (ii) the primal function $V^\star(\xi)$ is a concave function over $\xi\in\Xi$.
\end{lemma}

Lemma~\ref{lem:primal_function} shows that relaxing the constraints more (or decreasing $\xi$) may yield a larger optimal value for Problem~\eqref{eq:CMDP_relaxed}; however this relaxed problem becomes more far away from the nominal problem~\eqref{eq:CMDP}. Similarly, tightening the constraints can decrease the optimal value, even nullifying the constraints. 

To stay close to the nominal problem~\eqref{eq:CMDP} while specifying constraints efficiently, we will exploit
the properties of $V^\star(\xi)$ together with a relaxation cost to introduce a trade-off solution for Problem~\eqref{eq:CMDP_relaxed} in Section~\ref{sec:resilient CRL}. We feature this solution with \emph{resilient} due to its adaptation of primal value function to the changing constraint specifications, and we term the problem of finding optimal policy and constraint specification together in   Problem~\eqref{eq:CMDP_relaxed} as resilient constrained constrained policy optimization.

\subsection{Examples with Unspecified Constraints}\label{subsec:CMDP_without_prior}

We showcase that constraint specifications are often not \emph{a priori} knowledge for Problem~\eqref{eq:CMDP}. 

\noindent\textbf{Resource allocation.}
In a system of $m$ users sharing a transition kernel~\cite{boutilier2016budget,vora2023welfare}, each user $i$ has a reward function $r_i$ and a cost function $c_i$ (or consumption), and the budget $B>0$ is given. The resource allocation is to decide a budget assignment $\{ \bar c_i \}_{i\,=\,1}^m$ for $m$ users by maximizing the average reward $r = \frac{1}{m}\sum_{i\,=\,1}^m r_i$ and restraining the total cost  $\sum_{i\,=\,1}^m\bar c_i \leq B$. Formulation~\eqref{eq:CMDP} applies when we take $u_i = - c_i$ and  $b_i = - \bar c_i$, and evaluate their discounted value functions. One application scenario is the robot monitoring problem experimented in Section~\ref{sec:experiments}.   
However, the budget assignment $\{b_i\}_{i\,=\,1}^m$ is a decision variable to be determined. Moreover, the budget can be uncertain, e.g., multi-arm bandits with limited resources~\cite{diaz2023flexible}.
 
\noindent\textbf{Many rewards trade-off.} Multiple rewards appear in real RL applications~\cite{shelton2000balancing,liu2014multiobjective,eysenbach2019reinforcement,liang2022reward}. For instance, in preference-based RL~\cite{liang2022reward}, extrinsic rewards $\{r_i\}_{i\,=\,1}^m$ are preferences of human feedback while an intrinsic reward $r$ captures the uncertainty in the disagreement among humans.
Take $u_i = r_i$ in Problem~\eqref{eq:CMDP}. To encourage exploration under  alignment with preferences, Problem~\eqref{eq:CMDP} aims to maximize the intrinsic reward value function while constraining extrinsic reward value functions above some thresholds $\{b_i\}_{i\,=\,1}^m$. However, such thresholds are unknown due to varying human's preferences.


\section{RESILIENT CONSTRAINED RL}\label{sec:resilient CRL}

To specify appropriate constraints, we introduce a cost function of relaxing the constraints, and a resilient equilibrium that balances the relaxation and the primal value function in Section~\ref{sec: RE}. We provide a tractable formulation based on regularization to find a resilient equilibrium in Section~\ref{sec: RE regularization}. 

\subsection{Resilient Equilibrium}\label{sec: RE}


We characterize the change of the primal value function $V^\star(\xi)$ to relaxation $\xi$ via the subgradient and geometric multiplier in nonlinear programming~\cite{bertsekas2016nonlinear}. Let the dual function for Problem~\eqref{eq:CMDP} be $D(\lambda) \DefinedAs \sup_{\pi\,\in\,\Pi} V_{r+\lambda^\top g}^\pi(\rho)$ and its domain be $\Lambda \DefinedAs \{ \lambda\in\mathbb{R}_+^m\,\vert\, D(\lambda) > -\infty\}$. A relation between the primal value function and the dual function is, for any $\lambda\geq 0$, 
\begin{equation}\label{eq:conjugate_primal_function}
    D(\lambda) 
    \; = \;
    \sup_{ \xi \,\in\,\Xi}  \;  \left\{ \lambda^\top \xi - (- V^\star(\xi)) \right\}
\end{equation}
which can be derived in Appendix~\ref{app:dual_function}. In other words, $D(\lambda)$ is the conjugate convex function of $-V^\star(\xi)$ for $\xi\in\Xi$ and the domain $\Lambda$ is a convex set. By the concavity of the primal function $V^\star(\xi)$ in Lemma~\ref{lem:primal_function} and $\Xi$, there exists a subgradient for $V^\star(\xi)$ at any interior point $\xi\in\Xi$. This subgradient naturally connects to the geometric multiplier $\lambda$, i.e.,  $V^\star(\xi) = \sup_{\pi\,\in\,\Pi} \{V_{r+\lambda^\top g}^\pi(\rho)-\lambda^\top \xi\}$. Thus, we can interpret
the subgradient of the negative primal function as a geometric multiplier for Problem~\eqref{eq:CMDP_relaxed}; 
see Appendix~\ref{app:geometric_multiplier} for proof.

\begin{lemma}[Subgradient and Geometric Multiplier]\label{lem:geometric_multiplier}
    In Problem~\eqref{eq:CMDP_relaxed} with any $\xi \in \Xi$, these are equivalent:
    (i) $\lambda$ is a subgradient of $-V^\star(\xi)$ at $\xi$; (ii) $\lambda$ is a geometric multiplier for Problem~\eqref{eq:CMDP_relaxed}.
\end{lemma}

Having described the sensitivity of the primal function, we next introduce a notion of resilient equilibrium in supervised learning~\cite{hounie2023resilient} to determine the constraint specifications according to the difficulty in solving the nominal problem~\eqref{eq:CMDP}. The non-increasing property in  Lemma~\ref{lem:primal_function} is that the primal function $V^\star(\xi)$ would be increased by decreasing $\xi$ coordinate-wise (relaxing the constraints). However, more relaxation $\xi$ yields looser constraints. We introduce a convex cost function $h(\xi)$: $\Xi\to\mathbb{R}$ to measure the relaxation cost. Without loss of generality, we use the case: $\xi=0$ to match our nominal problem and set the cost to be zero, i.e., $h(0)=0$. To govern the change in $V^\star(\xi)$, we use the marginal price of relaxing constraint: $\nabla h(\xi)$, to define a resilient equilibrium $\xi^\star$.
\begin{definition}[Resilient Equilibrium]\label{def:RE}
    For any cost function $h$ that is continuously differentiable, concave, and non-increasing coordinate-wise,
    a resilient equilibrium for $V^\star(\xi)$ is a relaxation $\xi^\star \in \Xi$,
    \[
    \nabla h(\xi^\star)  
    \;\in\;
    \partial V^\star(\xi^\star).
    \]
\end{definition}

At a resilient equilibrium $\xi^\star$, for some $\epsilon>0$, relaxing it to $\xi^\star-\epsilon$ would increase the cost by $-\epsilon \nabla h(\xi^\star)$ and the primal function may get larger; similarly, tightening it to $\xi^\star+\epsilon$ would decrease the cost by $\epsilon \nabla h(\xi^\star)$ and the primal function may become smaller. The `resilient' captures how much we can relax (or tighten) the constraints before we observe significant improvement (or degradation) in the primal value function. Thus, a resilient constrained policy optimization problem reduces to solving Problem~\eqref{eq:CMDP_relaxed} for some $\xi^\star$ that is a resilient equilibrium for an user-defined cost function $h$.

\begin{figure}[tbh]
	\begin{center}
	   {\includegraphics[scale=0.14]{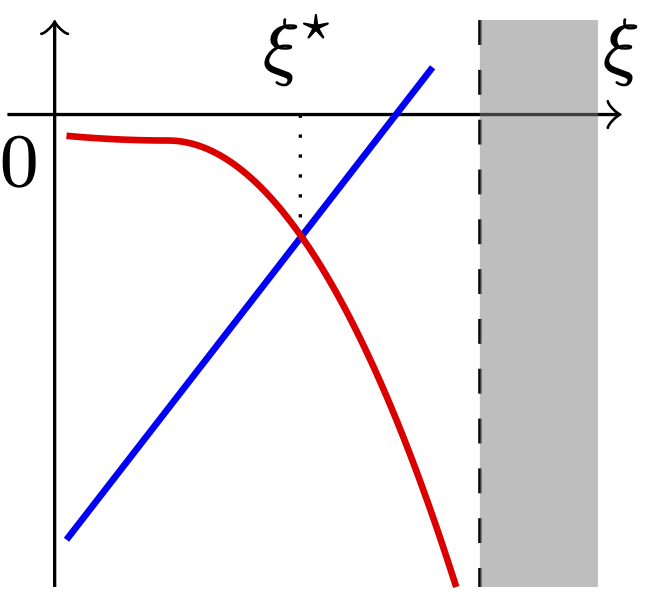}}
	\end{center}
        \vspace{-0.4cm}
	\caption{ Resilient equilibrium for Problem~\eqref{eq:CMDP_relaxed} with $m=1$ and a quadratic function $h(\xi)$ for $\xi \in \mathbb{R}$. The horizontal axis is the relaxation $\xi$, and the vertical axis is the (sub)gradient values: 
 $\nabla h(\xi)$~\mbox{(\textbf{\color{blue}---})} and 
 $\partial V(\xi)$~\mbox{(\textbf{\color{red}---})}. The shaded area means the infeasibility when $\xi$ is large.  
	}
	\label{fig: resilient equilibrium}
\end{figure}

We note that $V^\star(\xi) - h(\xi)$ is a concave function. The existence of resilient equilibrium can be easily obtained in Lemma~\ref{lem:resilient_eq}, and we delay it to  Appendix~\ref{app:resilient_eq}.
\begin{lemma}[Equilibrium Existence]\label{lem:resilient_eq}
    There exists a resilient equilibrium $\xi^\star \in\Xi$, and it is unique when $h$ is strictly convex. 
\end{lemma}

Behind the existence, Lemma~\ref{lem:resilient_mono} shows that relaxing the constraint (or decrease $\xi$) increases the cost sensitivity $|\nabla h(\xi)|$ while the effect on the primal value function is decreased; see Appendix~\ref{app:resilient_mono} for proof. Thus, they must cross at a resilient equilibrium as shown in Figure~\ref{fig: resilient equilibrium}.

\begin{lemma}[Coordinate-Wise Monotonicity]\label{lem:resilient_mono}
    Let $\xi$, $\xi' \in \Xi$ satisfy  $\xi_i'<\xi_i$ and $\xi_j' = \xi_j$ for $j\neq i$. 
    Then, $(\nabla h(\xi'))_i \leq (\nabla h(\xi))_i$, and $(\partial V^\star(\xi))_i \leq (\partial V^\star(\xi'))_i$, where $(\cdot)_i$ is the $i$th entry. 
\end{lemma}

However, the primal value function is unavailable. Theorem~\ref{thm:resilient_geometric_multiplier} gives a sufficient condition on a resilient equilibrium via the geometric multiplier, and we relate it to the duality next; see Appendices~\ref{app:resilient_geometric_multiplier}--\ref{app:resilient_geometric_multiplier_cor} for proofs. 

\begin{theorem}[Geometric Multiplier Condition]\label{thm:resilient_geometric_multiplier}
    For Problem~\eqref{eq:CMDP_relaxed} with $\bar\xi\in\Xi$, if $\lambda$ is a geometric multiplier and 
    $\nabla  h(\bar\xi)+\lambda = 0$, then $\bar\xi$ is a resilient equilibrium.
\end{theorem}

As a corollary of Theorem~\ref{thm:resilient_geometric_multiplier}, we relate the resilient equilibrium to an optimal Lagrange multiplier. Let the standard Lagrangian for Problem~\eqref{eq:CMDP_relaxed} for $\xi\in\mathbb{R}^m$ be,
\[
    \begin{array}{rcl}
        L(\pi,\lambda; \xi) 
    & = &\displaystyle
    V_r^\pi(\rho) \,+\, \sum_{i\,=\,1}^m \lambda_i (V_{g_i}^\pi(\rho) - \xi_i)
    \\[0.2cm]
        & \DefinedAs & \displaystyle 
    V_{r+\lambda^\top g}^\pi(\rho) - \lambda^\top \xi
    \end{array}
\]
and the associated dual function be,
\[
    D(\lambda; \xi) 
    \; = \;
    \max_{\pi\,\in\,\Pi} \;
    L(\pi,\lambda; \xi) \; \text{ for any } \lambda\geq 0.
\]
The optimal dual function  $D^\star(\xi) = \min_{\lambda\,\geq\,0} D(\lambda;\xi)$ is achieved at an optimal Lagrange multiplier $\lambda^\star(\xi)$. By the weak duality, $D^\star(\xi) \geq V^\star(\xi)$ for any $\xi\in\mathbb{R}^m$.

\begin{cor}\label{cor:resilient_geometric_multiplier}
    Let the strong duality hold for Problem~\eqref{eq:CMDP_relaxed} with some $\bar\xi\in\Xi$, i.e., $V^\star(\bar\xi) = D^\star(\bar\xi)$.
    If  
    $\nabla h(\bar\xi) + \lambda^\star(\bar\xi) = 0$, then $\bar\xi$ is a resilient equilibrium.
\end{cor}

The strong duality in Corollary~\ref{cor:resilient_geometric_multiplier} only concerns the relaxed problem~\eqref{eq:CMDP_relaxed}, which is much weaker than the the strong duality for the nominal problem (e.g.,~\cite{paternain2019constrained,ding2020natural}). We also remark that the strong duality is stronger than the geometric multiplier condition~\cite{bertsekas2016nonlinear}, which is more general than the study~\cite{hounie2023resilient}.   

\subsection{Resilience via Regularization }\label{sec: RE regularization}

Although a resilient equilibrium always exists under mild regularity conditions, it is not straightforward to design efficient policy learning algorithms for achieving such an equilibrium. To address this issue, we introduce the cost function into Problem~\eqref{eq:CMDP} as regularization,
\begin{equation}\label{eq:CMDP_regularized}
    \begin{array}{rl}
         \displaystyle \maximize_{\pi\,\in\,\Pi,\, \xi\,\in\,\Xi} & V_r^\pi(\rho) - h(\xi)
         \\
         \subject&  
         V_{g_i}^\pi(\rho) \; \geq \; \xi_i \; \text{ for } \; i = 1,\ldots,m
    \end{array}
\end{equation}
where $h(\xi)$ is a regularizer that is monotonically non-increasing coordinate-wise. Relaxing the constraint, i.e., decreasing $\xi$ would increase $h(\xi)$; tightening the constraint is the opposite. Lemma~\ref{lem:regularized_optimal_policy} states that Problem~\eqref{eq:CMDP_regularized} provides a resilient equilibrium and an associated resilient policy; see  Appendix~\ref{app:regularized_optimal_policy} for proof.

\begin{lemma}[Regularized Solution]\label{lem:regularized_optimal_policy}
    Let $(\bar\pi^\star,\bar\xi^\star)$ an optimal solution to Problem~\eqref{eq:CMDP_regularized}. Then, $\bar\xi^\star$ is a resilient equilibrium and $\bar\pi^\star$ is an associated resilient policy.
\end{lemma}

Problem~\eqref{eq:CMDP_regularized} is a practical extension of constrained policy optimization to jointly optimizing over policy and relaxation. Naturally, we can enable the extension of existing constrained policy search algorithms to being resilient; see two of them in Section~\ref{sec: resilient CPL}. Before that, we first show some important properties of Problem~\eqref{eq:CMDP_regularized}. 

We denote the optimal value for Problem~\eqref{eq:CMDP_regularized} by $V_h^\star \DefinedAs V_r^\star(\rho) - h(\bar\xi^\star)$ that is evaluated at an optimal solution $(\bar\pi^\star,\bar\xi^\star)$ if it is feasible; $V_h^\star = -\infty$ otherwise. Let the dual function for Problem~\eqref{eq:CMDP_regularized} be $D_h(\lambda) \DefinedAs \sup_{\pi\,\in\,\Pi, \xi\,\in\,\Xi} \{ V_{r+\lambda^\top g}^\pi(\rho) - h(\xi) - \lambda^\top \xi\}$ and the optimal dual function be $D_h^\star \DefinedAs D_h(\bar\lambda^\star)$ that is achieved at an optimal dual variable $\bar\lambda^\star$.

\begin{assumption}[Strict feasibility]
    \label{as:feasibility_regularized}
    There exist a pair of $(\bar\pi,\bar\xi) \in \Pi\times\Xi$ and a constant $c>0$ such that $V_{g_i}^{\bar\pi}(\rho)-\bar\xi_i\geq c$ for all $i=1,\ldots,m$.
\end{assumption}

Due to the flexibility of selecting $\bar\xi$, Assumption~\ref{as:feasibility_regularized} is weaker than the usual Slater condition~\cite{altman1999constrained}.



Thus, the strong duality and dual boundedness hold for Problem~\eqref{eq:CMDP_regularized}; see Appendices~\ref{app:strong_duality}--\ref{app:bounded_dual} for proofs.

\begin{theorem}[Strong Duality for Regularized Problem]\label{thm:strong_duality}
     Let Assumption~\ref{as:feasibility_regularized} hold. Then, the strong duality holds for Problem~\eqref{eq:CMDP_regularized}, i.e., $V_h^\star = D_h^\star$.
\end{theorem}

\begin{cor}[Dual Boundedness]\label{cor:bounded_dual}
    Let Assumption~\ref{as:feasibility_regularized} hold. Then, the optimal dual is bounded, i.e., 
    \[
        0 \; \leq \;
        \bar\lambda_i^\star 
        \;\leq\; 
       \frac{V_r^\star(\rho)  - h(\bar\xi^\star)  - (V_r^{\bar\pi} - h(\bar\xi))}{c}
       \; \DefinedAs \;
       C_h.
    \]
\end{cor}

The strict feasibility of $(\bar\pi,\bar\xi)$ and the optimality of $(\bar{\pi}^\star,\bar{\xi}^\star)$ leads to $C_h>0$. 
Corollary~\ref{cor:bounded_dual} restricts dual variables in $\Lambda \DefinedAs \{  \lambda\in\mathbb{R}_+^m \,\vert\, \lambda_i \leq C_h, i = 1,\ldots,m \}$. Let the standard Lagrangian for Problem~\eqref{eq:CMDP_regularized} be
\[
    \begin{array}{rcl}
         L_h(\pi,\xi; \lambda) 
    & \DefinedAs &
    V_{r+\lambda^\top g}^\pi(\rho) - h(\xi)- \lambda^\top \xi. 
    \end{array}
\]
By the strong duality, Problem~\eqref{eq:CMDP_regularized} is equivalent to the following constrained saddle-point problem,
\[
    \begin{array}{rcl}
         &&  \displaystyle
         \!\!\!\!  \!\!\!\!  \!\!\!\!  
         \!\!\!\!
         \!\!\!\!
         \!\!\!\!
         \!\!\!\!
         \!\!\!\!
         \maximize_{\pi\,\in\,\Pi,\, \xi\,\in\,\Xi} \;\minimize_{\lambda\,\in\,\Lambda} \;
        L_h(\pi,\xi; \lambda)
        \\[0.2cm]
        & = & \displaystyle
        \minimize_{\lambda\,\in\,\Lambda} \;
        \maximize_{\pi\,\in\,\Pi,\, \xi\,\in\,\Xi} \; L_h(\pi,\xi; \lambda).
    \end{array}
\]
Let $\Pi^\star\times \Xi^\star\times \Lambda^\star$ be a set of saddle points of $L_h(\pi,\xi;\lambda)$ over $\Pi\times \Xi\times \Lambda$.
From Theorem~\ref{thm:strong_duality}, there always exists such a saddle point, i.e., $\Pi^\star\times \Xi^\star\times \Lambda^\star \neq \emptyset$. From the definition of $\Xi^\star$, $|\xi_i|\leq 1/(1-\gamma)$ for any $\xi\in\Xi^\star$. Aided by these nice properties, we next introduce two constrained policy search algorithms to find optimal policy and constraint specification. 


\section{RESILIENT CONSTRAINED POLICY LEARNING}\label{sec: resilient CPL}

We provide two constrained policy gradient algorithms for searching for policy and constraint specification together, in Section~\ref{sec: PGPD} and Section~\ref{sec: OPGPD}, respectively. 

\subsection{Resilient Policy Gradient Primal-Dual (ResPG-PD)  Method}\label{sec: PGPD}

We generalize the policy gradient primal-dual mirror descent~\cite{ding2022policy} for our resilient problem~\eqref{eq:CMDP_regularized} by adding a relaxation update. The resilient policy gradient primal-dual (ResPG-PD) method in Algorithm~\ref{alg: resilient PG} maintains three sequences for primal and dual variables via Primal update~\eqref{eq:policy_gradient_primal} and Dual update~\eqref{eq:policy_gradient_dual}: two primal sequences $(\{ \pi_t\}_{t\,\geq\,1}, \{\xi_t\}_{t\,\geq\,1})$ for policy and relaxation, 
and a dual sequence $\{ \lambda_t \}_{t\,\geq\,1}$, where $\eta$ is the stepsize, $\pi_0$ is the uniform distribution over the action space, $\xi_0 = 0$, and $\lambda_0 = 0$. In Primal update~\eqref{eq:policy_gradient_primal}, the policy update works as the projected $Q$-ascent~\cite{bhandari2021linear,xiao2022convergence} and the relaxation update performs the projected gradient ascent. Dual update~\eqref{eq:policy_gradient_dual} is the standard projected gradient descent.  When the relaxation is fixed, i.e., $\xi_t=\xi$, the relaxation update does not impact the dual update, and thus Algorithm~\ref{alg: resilient PG} reduces to the policy gradient primal-dual method in Euclidean space. By viewing this, we next extend the average-value convergence analysis for our resilient problem by incorporating the additional relaxation update $\{\xi_t\}_{t\,\geq\,1}$ in Theorem~\ref{thm: average-value convergence} and delay its proof to Appendix~\ref{app: average-value convergence}. 

We measure the performance of Algorithm~\ref{alg: resilient PG} by comparing the sequences $\{\pi_t,\xi_t,\lambda_t\}_{t\,\geq\,1}$ with the optimal solution $(\bar\pi^\star,\bar\xi^\star)$ in the standard notion of regret,
\[
\begin{array}{rcl}
     R_{\text{opt}} 
     & = &  \displaystyle
     \frac{1}{T} \sum_{t \, = \,0}^{T-1} (V_r^\star(\rho) -h(\bar\xi^\star) - (V_r^{\pi_t}(\rho) -h(\xi_t))) 
     \\[0.2cm]
     R_{\text{vio}}
     & =  & \displaystyle
    \sum_{i\,=\,1}^m \left[ \frac{1}{T} \sum_{t \, = \,0}^{T-1}(\xi_{i,t}- V_{g_i}^{\pi_t}(\rho)) \right]_+
\end{array}
\]
where $R_{\text{opt}}$ is the average of the sub-optimal gaps and $R_{\text{vio}}$ is the sum of the averaged constraint violations.  

\begin{algorithm*}[h]
	\caption{Resilient policy gradient primal-dual (ResPG-PD) method }
	\label{alg: resilient PG}
	\begin{algorithmic}[1]
		\STATE
		\textbf{Parameters:} $\eta>0$. 
  \\
  \textbf{Initialization}: Let $\pi_0 (a\,\vert\,s) = {1}/{A}$ for $s\in\calS$, $a\in\calA$, and $\xi_0 =  0$, and $\lambda_0  = 0$.
		\FOR{step $t=0,\ldots,T-1$} 
		\STATE Primal-dual update 
            \begin{subequations} \label{eq:policy_gradient}
                    \begin{equation}\label{eq:policy_gradient_primal}
    \begin{array}{rcl}
         \pi_{t+1}(\cdot\,\vert\,s)
         &  =  & 
         \displaystyle\argmax_{\pi(\cdot\,\vert\,s)\,\in\,\Pi}\;
         \left\{
                \sum_{a} \pi(a\,\vert\,s) Q_{r+\lambda_{t}^\top g}^{\pi_{t}}(s,a)
                - 
                \frac{1}{2\eta} \norm{\pi(\cdot\,\vert\,s) - \pi_t(\cdot\,\vert\,s)}^2
         \right\}
         \\[0.4cm]
         \xi_{t+1}
         &  =  &
         \displaystyle
         \argmax_{\xi\,\in\,\Xi} \;
         \left\{
                \xi^\top\left( -\nabla h(\xi_{t}) - \lambda_{t}\right) -
                \frac{1}{2\eta} \norm{\xi - \xi_t}^2
         \right\}
    \end{array}
\end{equation}
\begin{equation}
    \label{eq:policy_gradient_dual}
    \begin{array}{rcl}
         \lambda_{t+1}
         &  =  &
         \displaystyle
         \argmin_{\lambda\,\in\,\Lambda} \;
         \left\{
                \lambda^\top\left( V_{g}^{\pi_{t}}(\rho)- \xi_{t}\right) +
                \frac{1}{2\eta} \norm{\lambda - \lambda_t}^2
         \right\}
    \end{array}
\end{equation}
\end{subequations}
		\ENDFOR
	\end{algorithmic}
\end{algorithm*}


\begin{theorem}[Regret-Type Performance]\label{thm: average-value convergence}
    Let Assumption~\ref{as:feasibility_regularized} hold. Suppose $\Lambda = \left[0, 2C_h\right]$, and $h(\xi)$ has Lipschitz continuous gradient with parameter $L_h$ over $\xi\in\Xi$. If $\eta = 1/\sqrt{T}$ for Algorithm~\ref{alg: resilient PG}, then,
    \[
    \begin{array}{rcl}
         \displaystyle
          R_{\normalfont\text{opt}} 
         &\leq& \displaystyle
         \frac{m(7+(L_h+1)^2)}{\sqrt{T}}
         \\[0.2cm]
         \displaystyle
            R_{\normalfont\text{vio}}
          &\leq& \displaystyle
          \frac{(8+(L_h+1)^2) m /C_h + mC_h}{(1-\gamma)^2 \sqrt{T}}.
    \end{array}
    \]
\end{theorem}

Theorem~\ref{thm: average-value convergence} states that the average sub-optimal gaps and constraint violations of the primal-dual iterates of Algorithm~\ref{alg: resilient PG} decay to zero with rate $1/\sqrt{T}$. This rate matches the rate of non-resilient algorithms~\cite{ding2022policy} and is independent of MDP's dimension. Due to the regularization and relaxation in Problem~\eqref{eq:CMDP_regularized}, our proof  handles regularized reward value and relaxed utility value together, generalizing the prior art for a broader class of problems. Since each primal-dual iteration involves projections to a probability simplex and intervals, requiring linear complexity, Algorithm~\ref{alg: resilient PG} has polynomial computational complexity. Denote $\xi_i' \DefinedAs \frac{1}{T}\sum_{t\,=\,0}^{T-1}\xi_{i,t}$. After $T=O(1/\epsilon^2)$ iterations of Algorithm~\ref{alg: resilient PG}, we can select the best policy $\pi'$ from $T$ steps,
\[
    \begin{array}{rcl}    
        V_r^{\star} - h(\bar\xi^\star) - (V_r^{\pi'} - h(\xi')) & = & O(\epsilon)
        \\[0.2cm]
        \displaystyle\left[
        \xi_i'-
V_{g_i}^{\pi'}
\right]_+ & = & O(\epsilon)
    \end{array}
\]
However, safety-critical systems demand training stability of policy iterates, which can't be guaranteed by regret performance. This issue was addressed by the state-augmentation method~\cite{calvo2023state} and the regularized or optimistic policy gradient methods~\cite{ding2023last} in non-resilient problems. We next address this issue in the resilient context by offering a resilient optimistic policy gradient method.


%

\subsection{Resilient Optimistic Policy Gradient Primal-Dual (ResOPG-PD) Method}\label{sec: OPGPD}

We extend Algorithm~\ref{alg: resilient PG} to an optimistic variant via the optimistic method~\cite{rakhlin2013online}, which is detailed in Algorithm~\ref{alg: resilient OPG} in Appendix~\ref{app: resilient OPG}. The resilient optimistic policy gradient primal-dual (ResOPG-PD) method maintains two sets of primal-dual sequences $\{\pi_t,\xi_t,\lambda_t\}_{t\,\geq\,1}$ and $\{\hat\pi_t,\hat\xi_t,\hat\lambda_t\}_{t\,\geq\,1}$. The update for $\{\hat\pi_t,\hat\xi_t,\hat\lambda_t\}_{t\,\geq\,1}$ is similar as~\eqref{eq:policy_gradient} that can be viewed as a real update, except that their gradients are computed at some intermediate iterates $\{\pi_t,\xi_t,\lambda_t\}_{t\,\geq\,1}$ that serve as predictions, instead of previous iterates. Thus, the real step is optimistic about the predictions, which is used to stabilize the dynamics of gradient-based algorithms~\cite{popov1980modification}. We can also view Algorithm~\ref{alg: resilient OPG} as a resilient version of the optimistic policy gradient primal-dual method~\cite{ding2023last} with the introduction of the relaxation update $\{\xi_t,\hat\xi_t\}_{t\,\geq\,1}$. By accounting for the relaxation update, we establish convergence guarantee on the primal-dual iterates in Theorem~\ref{thm: last-iterate convergence}; see Appendix~\ref{app: last-iterate convergence} for proof. 

We first state a few notations. The distribution mismatch coefficient over $\rho$ is $\kappa \DefinedAs \sup_{\pi} \Vert d_\rho^\pi / \rho \Vert_\infty$, where the division is component-wise. Clearly, $\kappa \leq {1}/{\rho_{\text{min}}}$, where $\rho_{\text{min}} \DefinedAs \min_s \rho(s)$. The projection operator $\mathcal{P}_X(\cdot)$ is given by $\mathcal{P}_X(x) \DefinedAs \argmin_{x'\,\in\, X} \norm{x-x'}$.


\begin{theorem}[Last-Iterate Convergence]\label{thm: last-iterate convergence}
    Let Assumption~\ref{as:feasibility_regularized} hold. Suppose $\Lambda = \left[0, 2C_h\right]$, $\rho_{\normalfont\text{min}}>0$, $h(\xi)$ is strongly convex and has Lipschitz continuous gradient with parameter $L_h$, and the optimal state visitation distribution is unique, i.e., $d_\rho^{\pi^\star} = d_\rho^{\pi}$ for $\pi\in\Pi^\star$. If we set stepsize $\eta \leq \eta_{\max}$, where~$\eta_{\max}$~is given in Appendix~\ref{app: last-iterate convergence} for Algorithm~\ref{alg: resilient OPG} in Appendix~\ref{app: resilient OPG}, then for any $t$,
    \[
    \begin{array}{rcl}
          & &
         \displaystyle
         \!\!\!\!  \!\!\!\!
         \!\!\!\!
         \frac{1}{2(1-\gamma)} \sum_{s}d_\rho^{\pi^\star}(s) \norm{\mathcal{P}_{\Pi^\star}( \hat\pi_{t}(\cdot\,\vert\,s)) - \hat\pi_{t}(\cdot\,\vert\,s)}^2
         \\[0.2cm]
        &  & \displaystyle
        \!\!\!\!  \!\!\!\!
         \!\!\!\!
         +\,
         \frac{1}{2}\norm{\mathcal{P}_{\Xi^\star}(\hat\xi_{t})-\hat\xi_{t}}^2 
         +
         \frac{1}{2}  
            \norm{\mathcal{P}_{\Lambda^\star}(\hat\lambda_{t})-\hat\lambda_{t}}^2
        \;=\; O\left(\frac{1}{t}\right)
    \end{array}
\]
where $O(\cdot)$ hides a problem-dependent constant $C_{\rho,\gamma,\sigma}$ in Lemma~\ref{lem:problem constant} in Appendix~\ref{app: last-iterate convergence}.
\end{theorem}

Theorem~\ref{thm: last-iterate convergence} states that the primal-dual iterates of Algorithm~\ref{alg: resilient OPG} converge to a set of $(\bar\pi^\star,\bar\xi^\star,\bar\lambda^\star)$ in a sublinear rate. Due to the introduction of regularization and relaxation, our proof distinguishes itself from the linear rate~\cite{ding2023last}, e.g., a new quadratic term enters into the lower bound in Lemma~\ref{lem:problem constant}. An immediate implication of Theorem~\ref{thm: last-iterate convergence} is that the primal iterates $(\hat\pi_t,\hat\xi_t)$ are $\epsilon$-near optimal after $O(1/\epsilon^2)$ iterations; see Appendix~\ref{app:policy and relaxation} for proof.

\begin{cor}[Near-Optimal Policy and Relaxation]\label{cor:policy and relaxation}
    Let assumptions in Theorem~\ref{thm: last-iterate convergence} hold. For a desired level of accuracy $\epsilon>0$, if the stepsize $\eta$ is provided by Theorem~\ref{thm: last-iterate convergence}, then for $t = \Omega(1/\epsilon^2)$,
    \[
    \begin{array}{rcl}
         V_r^{\star}(\rho) - h(\bar\xi^\star) - (V_r^{\hat\pi_t}(\rho) - h(\hat\xi_t)) & = & O(\epsilon)
         \\[0.2cm]
         \norm{\hat\xi_{t} - V_{g}^{\hat\pi_t}(\rho)} & = & O(\epsilon)
    \end{array}
    \]
    where $\Omega(\cdot)$ hides some problem-dependent constant.
\end{cor}

Corollary~\ref{cor:policy and relaxation} states that the last primal iterate $(\hat\pi_t,\hat\xi_t)$ is $\epsilon$-near optimal after $\Omega(1/\epsilon^2)$ iterations. This iteration complexity is similar as the one for Algorithm~\ref{alg: resilient PG}, as well as the computational complexity. However, policy convergence in Corollary~\ref{cor:policy and relaxation} is stated per iterate, which is stronger than the one in expectation. 

\definecolor{mblue}{HTML}{1f77b4}
\definecolor{mgreen}{HTML}{2ca02c}
\scalebox{0}{%
\begin{tikzpicture}
    \begin{axis}[hide axis]
        \addplot [
        color=red,
        dashed,
        line width=0.5pt,
        forget plot
        ]
        (0,0);\label{legend:red}
                \addplot [
        color=red,
        solid,
        line width=0.5pt,
        forget plot
        ]
        (0,0);\label{legend:redsolid}
                \addplot [
        color=mgreen,
        solid,
        dashed,
        line width=0.5pt,
        forget plot
        ]
        (0,0);\label{legend:greendash}
        \addplot [
        color=blue,
        line width=0.5pt,
        forget plot
        ]
        (0,0);\label{legend:blue}
        \addplot [
        color=black,
        dotted,
        line width=0.5pt,
        forget plot
        ]
        (0,0);\label{legend:black}
                \addplot [
        color=blue,
        dashed,
        line width=0.5pt,
        forget plot
        ]
        (0,0);\label{legend:bluedash}
        \addplot [
        color=orange,
        line width=1.2pt,
        forget plot
        ]
        (0,0);\label{legend:orange}
        \addplot [
        color=mgreen,
        line width=1.2pt,
        forget plot
        ]
        (0,0);\label{legend:green}
        \addplot[mark=x, color = red] [
    only marks,
    error bars/.cd,
      y dir=both, y explicit,
      error bar style={color=mblue},
    ] coordinates {
    (0, 0) +- (0, .5)
    };    \addlegendentry{xmarker}\label{legend:x}
    \addlegendimage{color=blue,mark=|,yshift = -0.15cm,
    mark options={mark repeat=2,mark phase=1},yscale= 0.5,rotate= 90,mark size=2.5pt }
    \addlegendentry{errorbar}\label{legend:errorbar}

            \addplot[mark=x, color = black] [
    only marks,
    error bars/.cd,
      y dir=both, y explicit,
      error bar style={color=mblue},
    ] coordinates {
    (0, 0) +- (0, .5)
    };    \addlegendentry{xmarker}\label{legend:blackx}
    \addlegendimage{color=mgreen,mark=|,yshift = -0.15cm,
    mark options={mark repeat=2,mark phase=1},yscale= 0.5,rotate= 90,mark size=2.5pt }
    \addlegendentry{errorbar}\label{legend:greenerrorbar}
    \end{axis}
\end{tikzpicture}
}

\begin{figure}[h]
    \centering
        \begin{tabular}{ccc}
             \rotatebox{90}{optimality gap}
             & 
             \includegraphics[width = 0.2\textwidth]{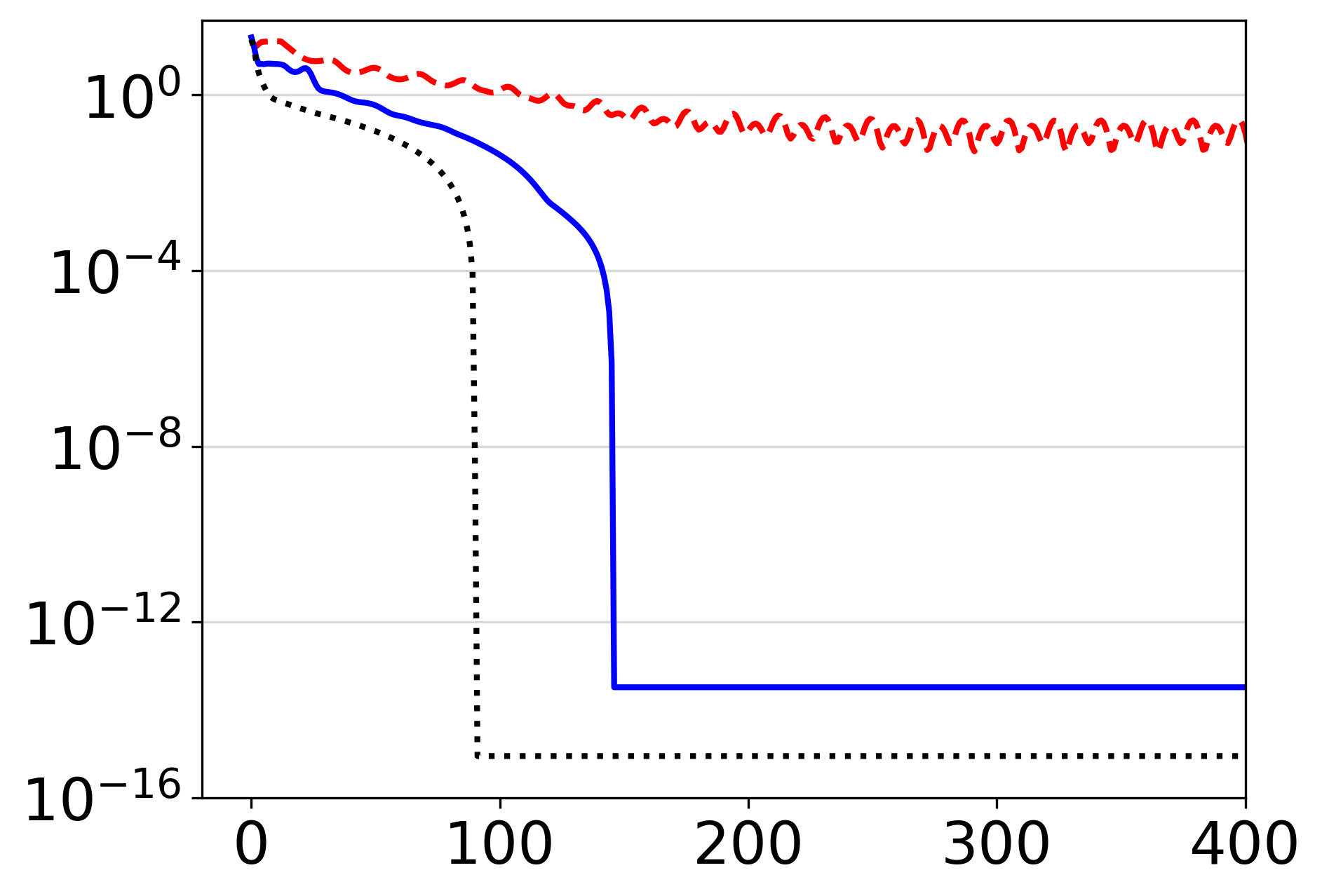} \!\!\!\!&\!\!\!\!
             \includegraphics[width = 0.2\textwidth]{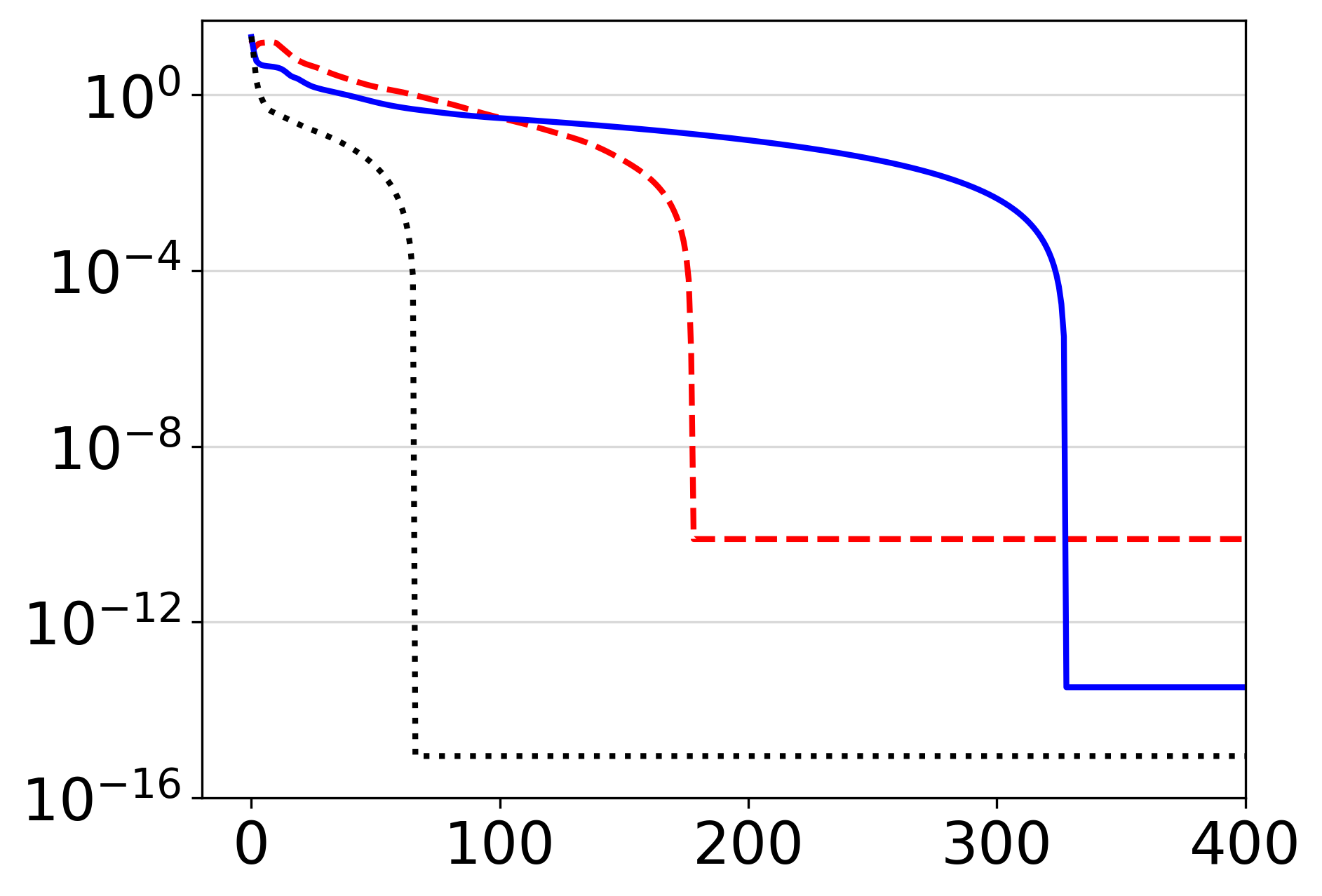}
             \\
             & iteration &  iteration
        \end{tabular}
    \\[-0.2cm]
    \caption{Policy optimality gaps of ResPG-PD (Algorithm~\ref{alg: resilient PG}, left) and ResOPG-PD (Algorithm~\ref{alg: resilient OPG}, right), with three cost functions $h(\xi) = \alpha \xi^2$ for $\alpha = 0.03$ (\ref{legend:red}) $\alpha = 0.2$ (\ref{legend:blue}), $\alpha = 1$ (\ref{legend:black}), and stepsize $\eta=0.2$.
    }
    \label{fig:DistanceVsIteropt}
\end{figure}

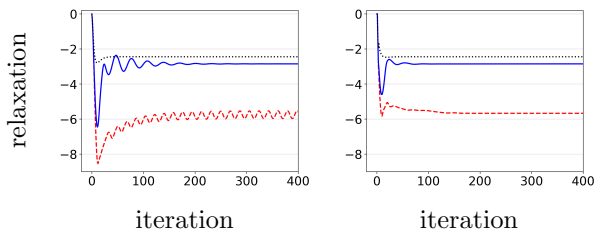
\begin{figure}[h]
    \centering
        \begin{tikzpicture}
        \node[](img1) at(0,0) {\includegraphics[width = 0.2\textwidth]{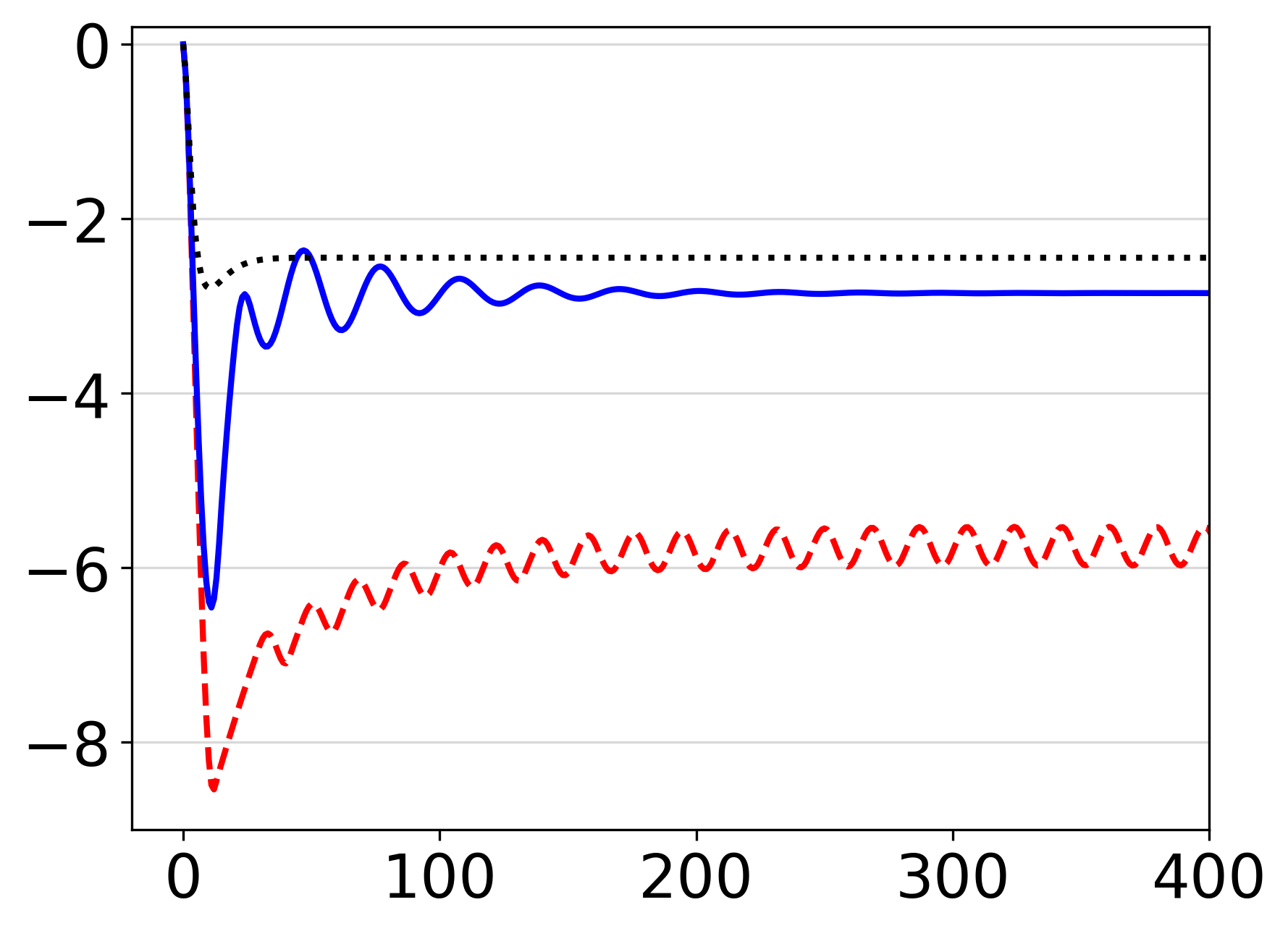}};
        \node[below = 0cm of img1]{iteration};
         \node[left = 0.4cm of img1,yshift=1.0cm, rotate=90]{relaxation};
    \end{tikzpicture}
    \begin{tikzpicture}
        \node[](img1) at(0,0) {\includegraphics[width = 0.2\textwidth]{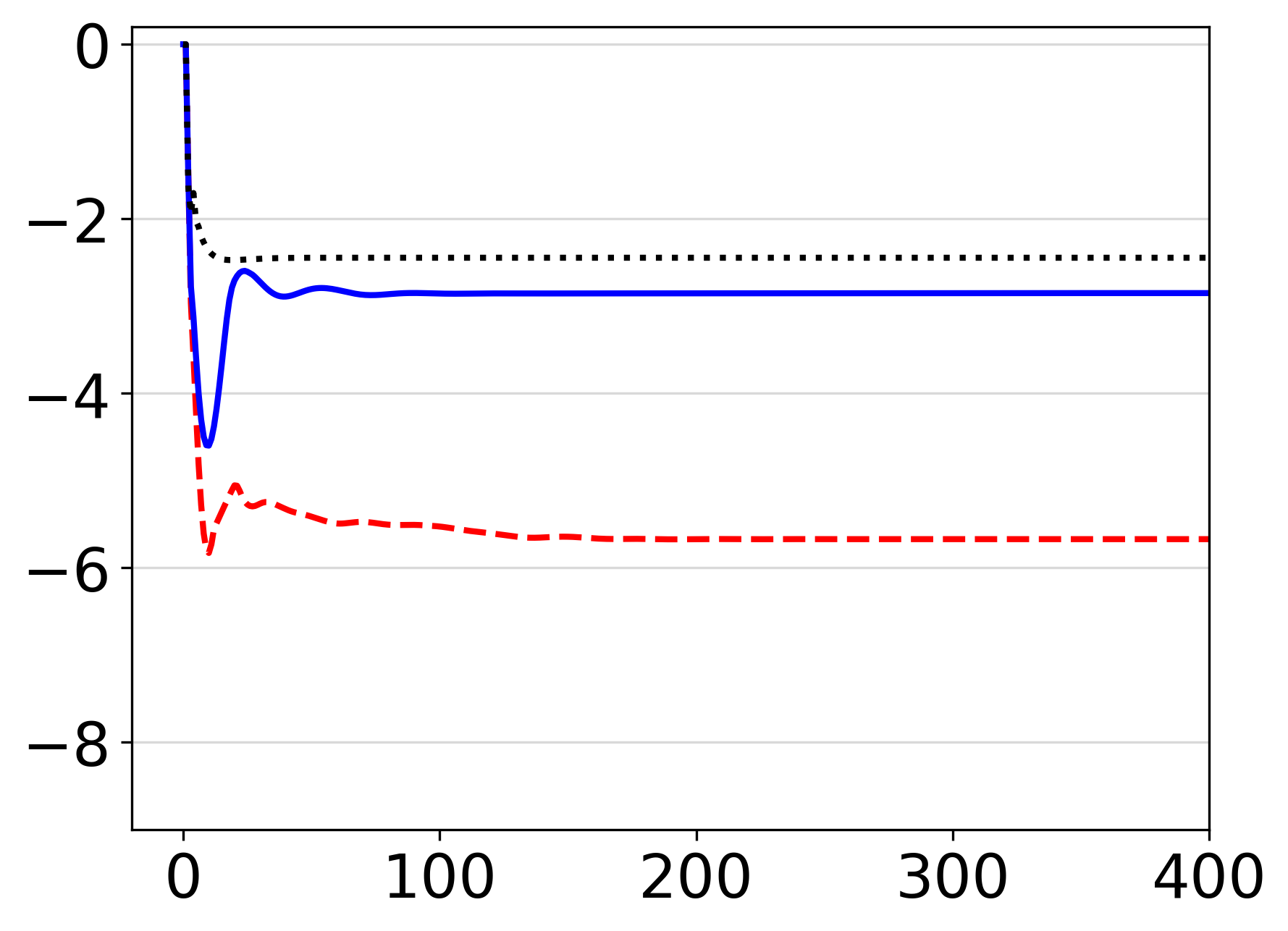}};
        \node[below = 0cm of img1]{iteration};
    \end{tikzpicture}
    \\[-0.2cm]
    \caption{Relaxation of ResPG-PD (Algorithm~\ref{alg: resilient PG}, left) and ResOPG-PD (Algorithm~\ref{alg: resilient OPG}, right), with three cost functions $h(\xi) = \alpha \xi^2$ for $\alpha = 0.03$ (\ref{legend:red}), $\alpha = 0.2$ (\ref{legend:blue}), $\alpha = 1$ (\ref{legend:black}) and stepsize $\eta=0.2$.
    }
    \label{fig:E1XivsIter}
    \end{figure}

\begin{figure}[h]
    \centering
        \begin{tikzpicture}
        \node[](img1) at(0,0) {
        \includegraphics[width = 0.3\textwidth]{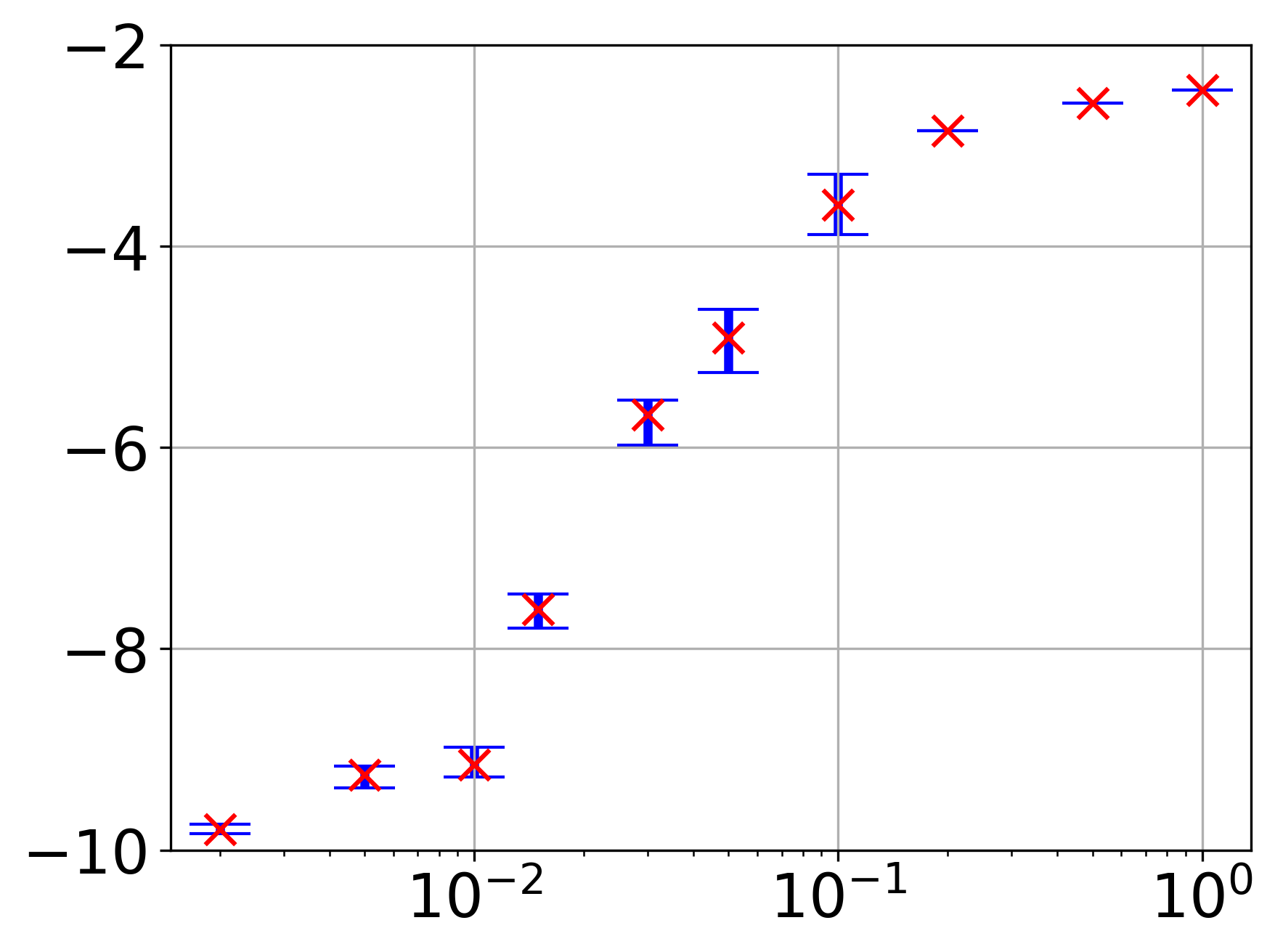}};
        \node[below = 0cm of img1]{relaxation cost};
         \node[left = 0.4cm of img1,yshift=1.0cm, rotate=90]{relaxation};
    \end{tikzpicture}
    \\[-0.2cm]
    \caption{Constraint specifications under different relaxation costs for Algorithm~\ref{alg: resilient PG} (ResPG-PD, \ref{legend:errorbar}\,) and Algorithm~\ref{alg: resilient OPG} (ResOPG-PD, \ref{legend:x}\,). The relaxation cost function is $h(\xi) = \alpha\xi^2$.
    The horizontal axis is the value of $\alpha$ and the vertical axis is the relaxation $\xi$. 
    The height of~\ref{legend:errorbar} is the oscillation magnitude of ResPG-PD. We run algorithms for $2000$ iterations with stepsize $\eta = 0.2$ and uniform initial distribution $\rho$.
    }
    \label{fig:E1XivsAlpha}
\end{figure}

\section{EXPERIMENTS}\label{sec:experiments}

We show the merits and the effectiveness of our resilient policy search algorithms: ResPG-PD (Algorithm~\ref{alg: resilient PG}) and ResOPG-PD (Algorithm~\ref{alg: resilient OPG}) in three experiments.

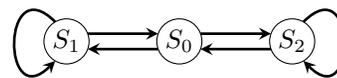
\begin{figure}[h]
    \centering
    \vspace{-4ex}
\scalebox{1}{\begin{tikzpicture}
  \coordinate (R0) at (0,0);
  \coordinate (R1) at (-1.5,0);
  \coordinate (R2) at (1.5,0);
  
  \draw (R0) circle (0.3cm);
  \draw (R1) circle (0.3cm);
  \draw (R2) circle (0.3cm);

  \node at (R0) {$S_0$};
  \node at (R1) {$S_1$};
  \node at (R2) {$S_2$};

    \draw[->, line width=1pt, >=stealth] ($(R0)+(200:0.3cm)$) -- ($(R1)+(-20:0.3cm)$) node[midway, above] {};
    \draw[->, line width=1pt, >=stealth] ($(R1)+(20:0.3cm)$) -- ($(R0)+(160:0.3cm)$) node[midway, above] {};

    \draw[->, line width=1pt, >=stealth] ($(R0)+(20:0.3cm)$) -- ($(R2)+(160:0.3cm)$) node[midway, above] {};
    \draw[->, line width=1pt, >=stealth] ($(R2)+(200:0.3cm)$) -- ($(R0)+(-20:0.3cm)$) node[midway, above] {};

    \draw[->,line width=1pt, >=stealth] ($(R1)+(120:0.3cm)$) to[out=135, in=-135, looseness=5] ($(R1)+(-120:0.3cm)$);
        \draw[->,line width=1pt, >=stealth] ($(R2)+(60:0.3cm)$) to[out=45, in=-45, looseness=5] ($(R2)+(-60:0.3cm)$);
\end{tikzpicture}}
    \\[-0.5cm]
    \caption{Robot monitoring of three locations. }\label{fig:MointoringProblem}
\end{figure}

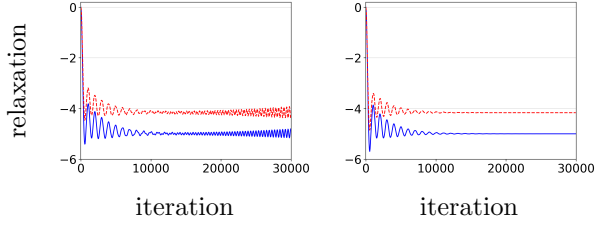
\begin{figure}[h]
    \centering
        \begin{tikzpicture}
        \node[](img1) at(0,0) {\includegraphics[width = 0.2\textwidth]{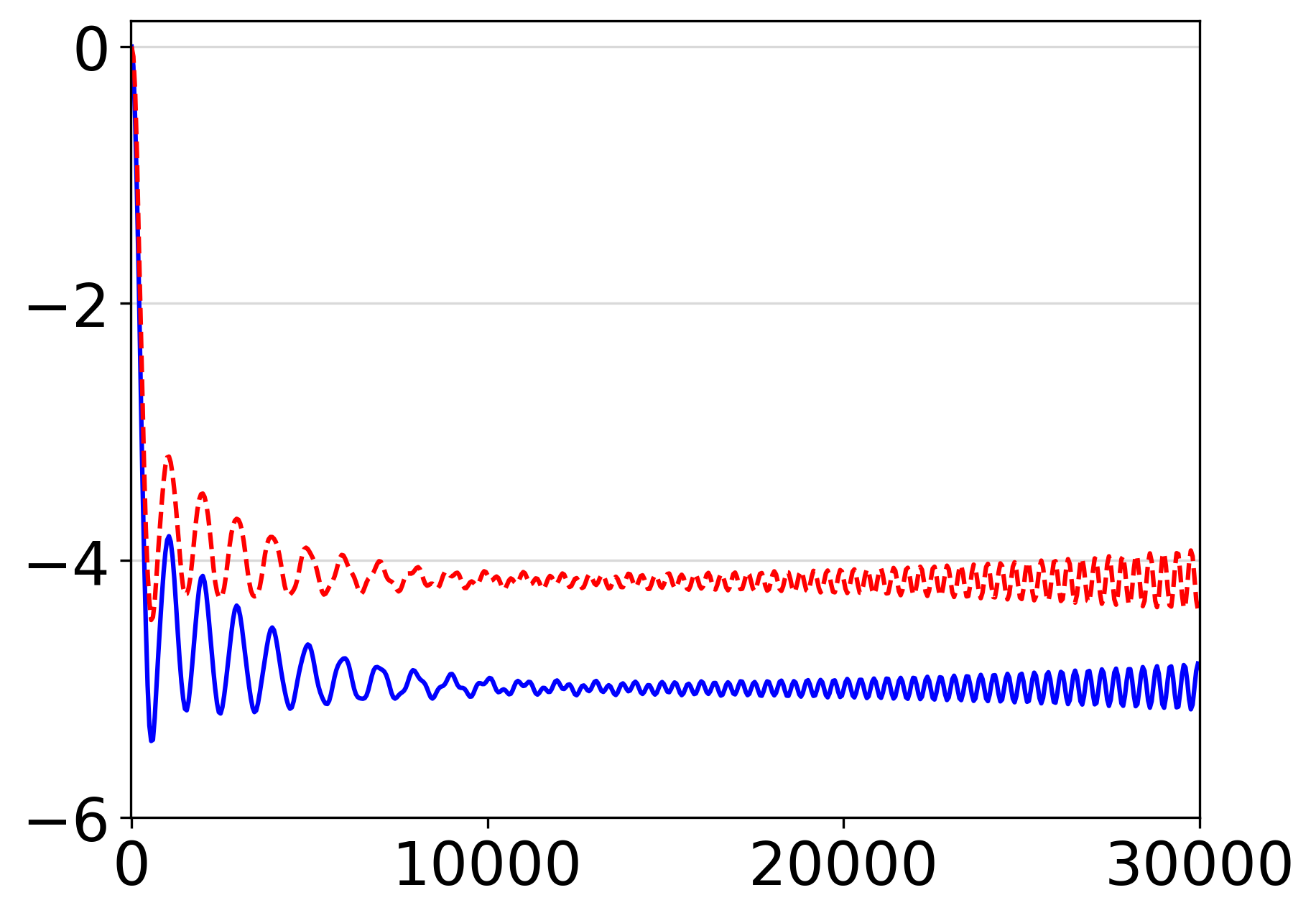}};
        \node[below = 0cm of img1]{iteration};
         \node[left = 0.4cm of img1,yshift=1.0cm, rotate=90]{relaxation};
    \end{tikzpicture}
\begin{tikzpicture}
        \node[](img1) at(0,0) {\includegraphics[width = 0.2\textwidth]{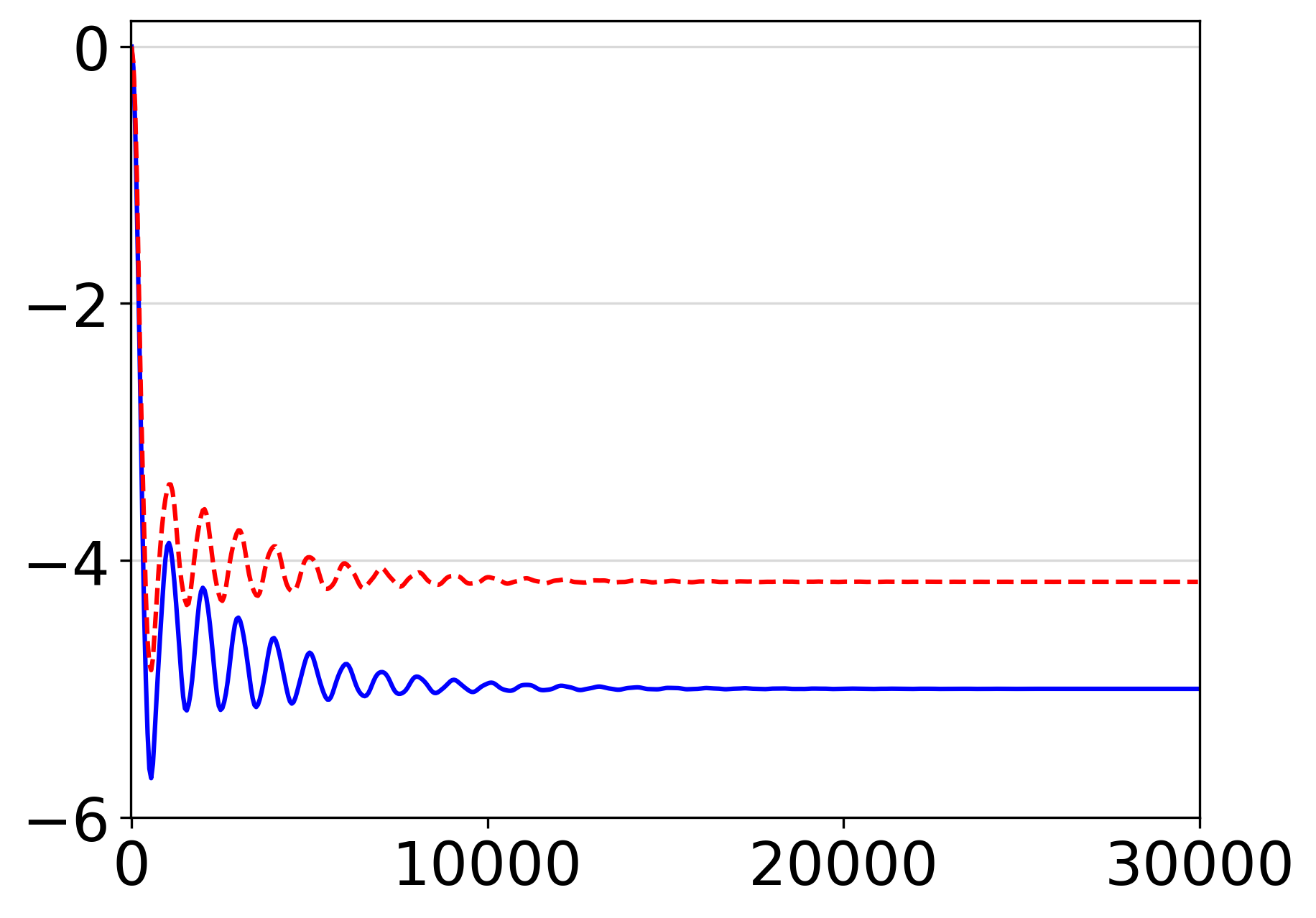}};
        \node[below = 0cm of img1]{iteration};
    \end{tikzpicture}
    \caption{Relaxations ($\xi_1$:~\ref{legend:blue}, $\xi_2$:~\ref{legend:red}) of ResPG-PD (Algorithm~\ref{alg: resilient PG}, left) and ResOPG-PD (Algorithm~\ref{alg: resilient OPG}, right), with a cost functions $h(\xi) = \alpha \norm{\xi}^2$ for $\alpha = 0.1$, and stepsize $\eta = 0.005$. 
    }
    \label{fig:ResE1}
\end{figure}

\begin{figure}[h]
    \centering
        \begin{tikzpicture}
        \node[](img1) at(0,0) {\includegraphics[width = 0.3\textwidth]{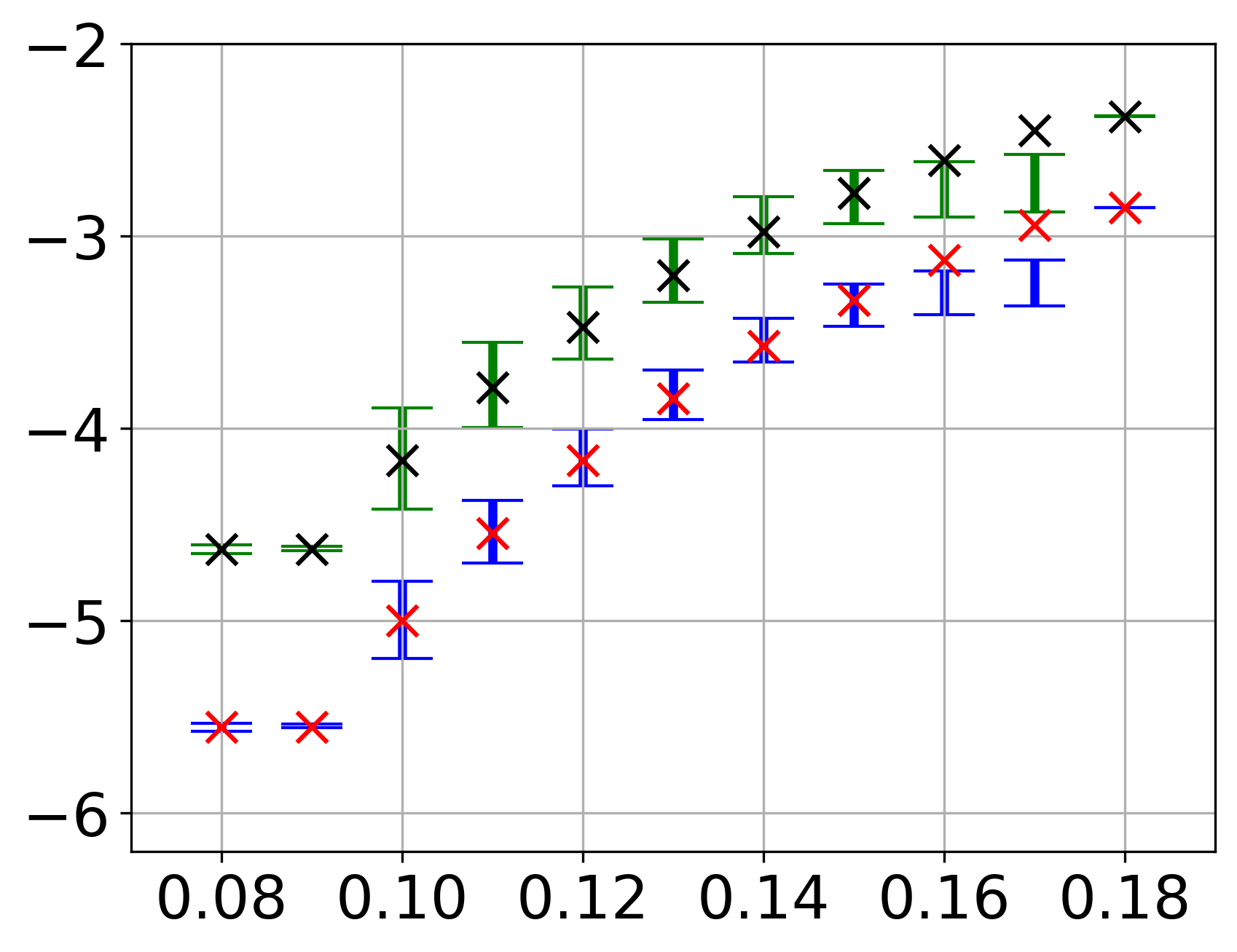}};
        \node[below = 0cm of img1]{relaxation cost};
         \node[left = 0.4cm of img1,yshift=1.0cm, rotate=90]{relaxation};
    \end{tikzpicture}
    \caption{Constraint specifications under different relaxation costs for Algorithm~\ref{alg: resilient PG} (ResPG-PD, $\xi_1$:~\ref{legend:errorbar}, $\xi_2$:~\ref{legend:greenerrorbar}\,) and Algorithm~\ref{alg: resilient OPG} (ResOPG-PD, $\xi_1$:~\ref{legend:x},  $\xi_2$:~\ref{legend:blackx}\,). The relaxation cost function is $h(\xi) = \alpha \norm{\xi}^2$.
    The horizontal axis is the value of $\alpha$ and the vertical axes are relaxations $\xi_1$ and $\xi_2$. 
    The height of~\ref{legend:errorbar} is the oscillation magnitude of ResPG-PD. We run algorithms for $100000$ iterations with stepsize $\eta = 0.005$ and uniform initial distribution $\rho$.
    }
    \label{fig:ResE2}
\end{figure}

We first use a randomly-generated constrained MDP with state/action size $(20,5)$ and calculate its optimal policy and relaxation for comparison; see Appendix~\ref{app: experiments} for the detail. Figure~\ref{fig:DistanceVsIteropt} shows that ResOPG-PD's policy iterates converge for different cost functions and so does ResPG-PD except for a small cost function, which verifies the average and the last-iterate performance in Theorems~\ref{thm: average-value convergence}--\ref{thm: last-iterate convergence}; see the convergence of relaxation in Figure~\ref{fig:E1XivsIter}. This indicates the capability of ResPG-PD and ResOPG-PD to find an optimal policy associated with a proper constraint specification. Figure~\ref{fig:E1XivsAlpha} shows that increasing relaxing cost leads to less relaxation, providing relaxation in accord with the relaxing cost. 

Second, 
we consider a monitoring problem in Figure~\ref{fig:MointoringProblem} in which an agent needs to spend as much time as possible at $S_0$ and must stay $S_1$ and $S_2$ with some time~\cite{calvo2023state}; see Appendix~\ref{app: experiments} for the detail. Due to the unknown feasibility of the time for $S_1$ or $S_2$, it warrants a resilient approach to relax (or tighten) either constraints. We choose the initial constraints to be \emph{infeasible}. Figure~\ref{fig:ResE1} shows that both algorithms can adapt two relaxations to the difficulty of constraints: one is relaxed more than the other, and Figure~\ref{fig:ResE2} shows two relaxation curves for two algorithms against different relaxation cost functions.

\begin{figure}[b]
    \centering
        \begin{tikzpicture}
        \node[](img1) at(0,0) {\includegraphics[width = 0.2\textwidth]{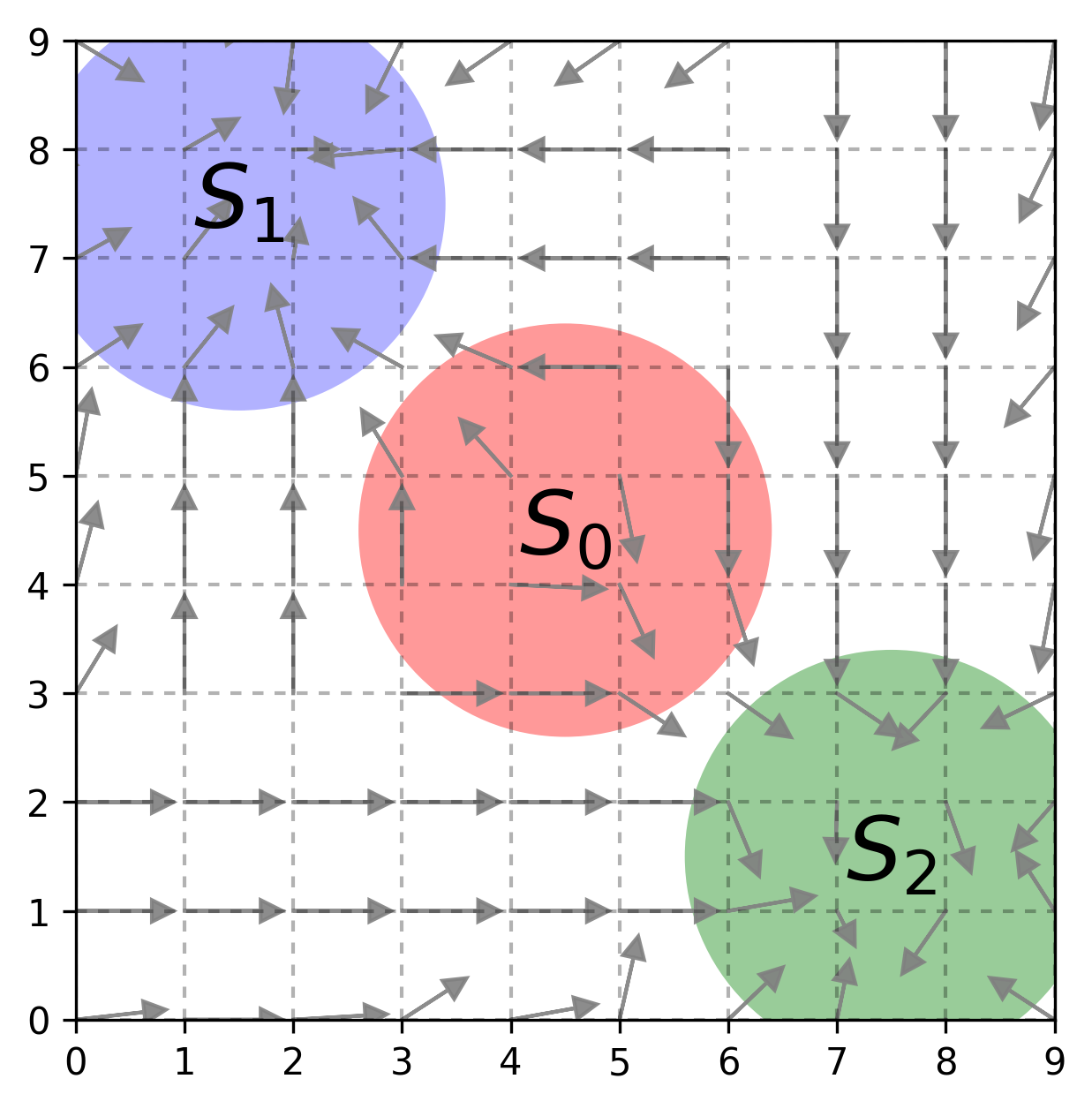}};
        \node[below = 0cm of img1]{(a) non-resilient policy};
        \end{tikzpicture}
        \begin{tikzpicture}
        \node[](img2) at(0,0) {\includegraphics[width = 0.2\textwidth]{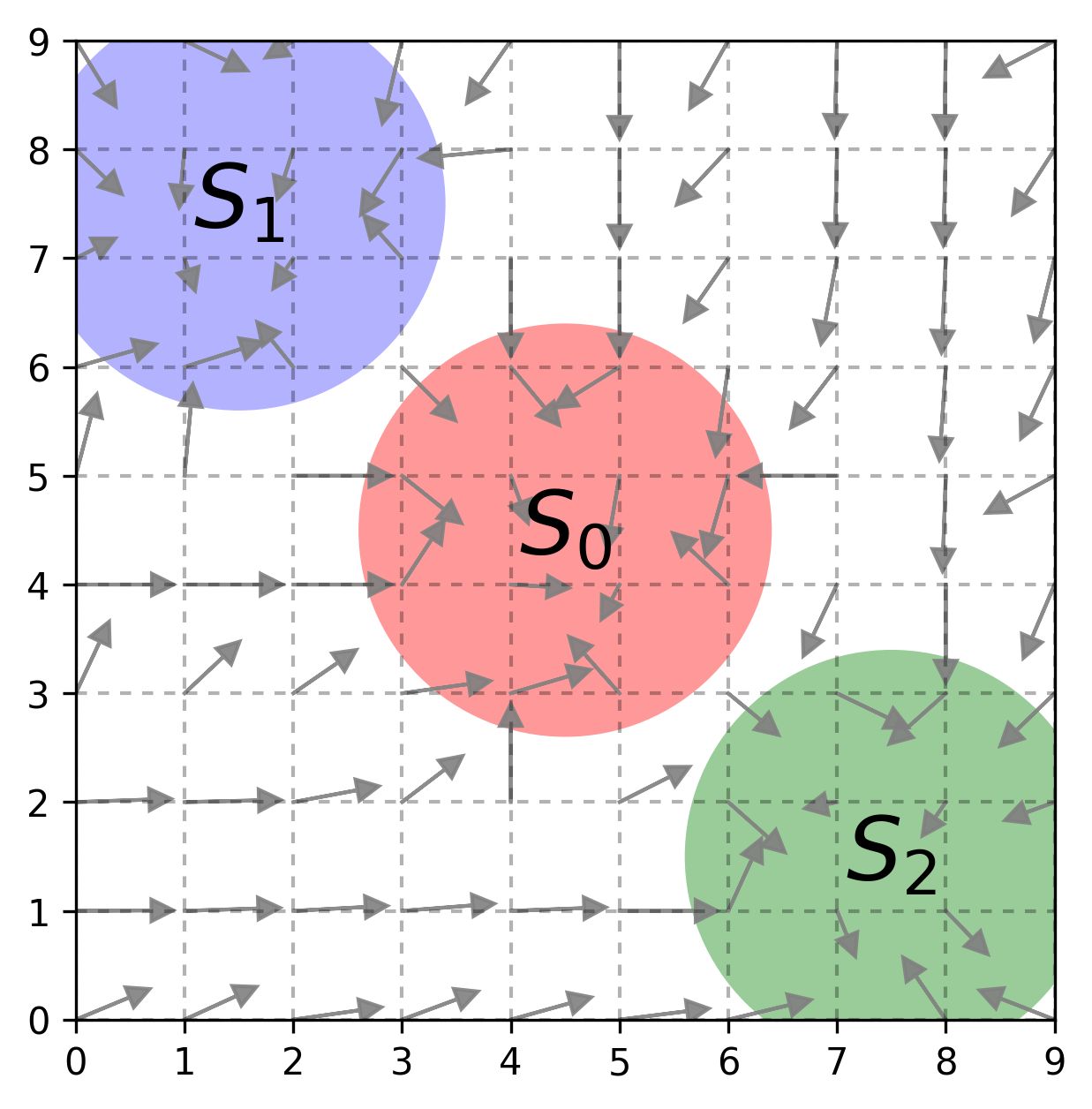}};
        \node[below = 0cm of img2]{(b) resilient policy};
    \end{tikzpicture}
    \caption{Robot monitoring of three areas. Arrows mean the moving directions of a policy generated by (a) a non-resilient algorithm: Algorithm~\ref{alg: resilient OPG} (ResOPG-PD) without relaxation updates; (b) a resilient algorithm: Algorithm~\ref{alg: resilient OPG} (ResOPG-PD), with a cost functions $h(\xi) = \alpha \norm{\xi}^2$ for $\alpha = 0.08$, and stepsize $\eta=0.05$.
    }
    \label{fig:Monitoring_nonres}
\end{figure}

\begin{figure}[tbh]
    \centering
    \begin{tikzpicture}
        \node[](img1) at(0,0) {\includegraphics[width = 0.18\textwidth]{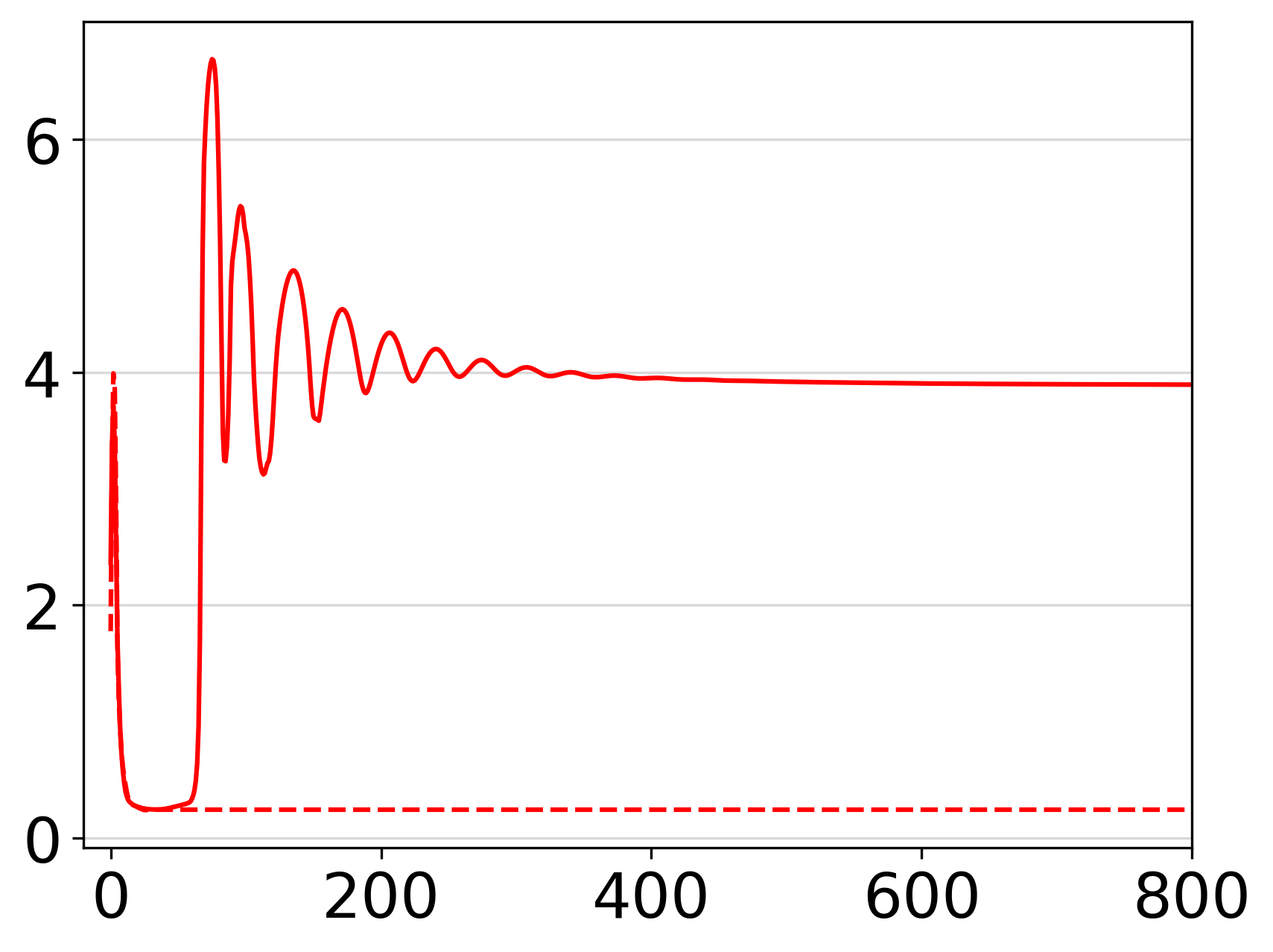}};
        \node[left = 0.4cm of img1,yshift=1.2cm, rotate=90]{reward value};
        \node[below = 0cm of img1]{iteration};
    \end{tikzpicture}
    \begin{tikzpicture}
        \node[](img2) at(0,0) {\includegraphics[width = 0.18\textwidth]{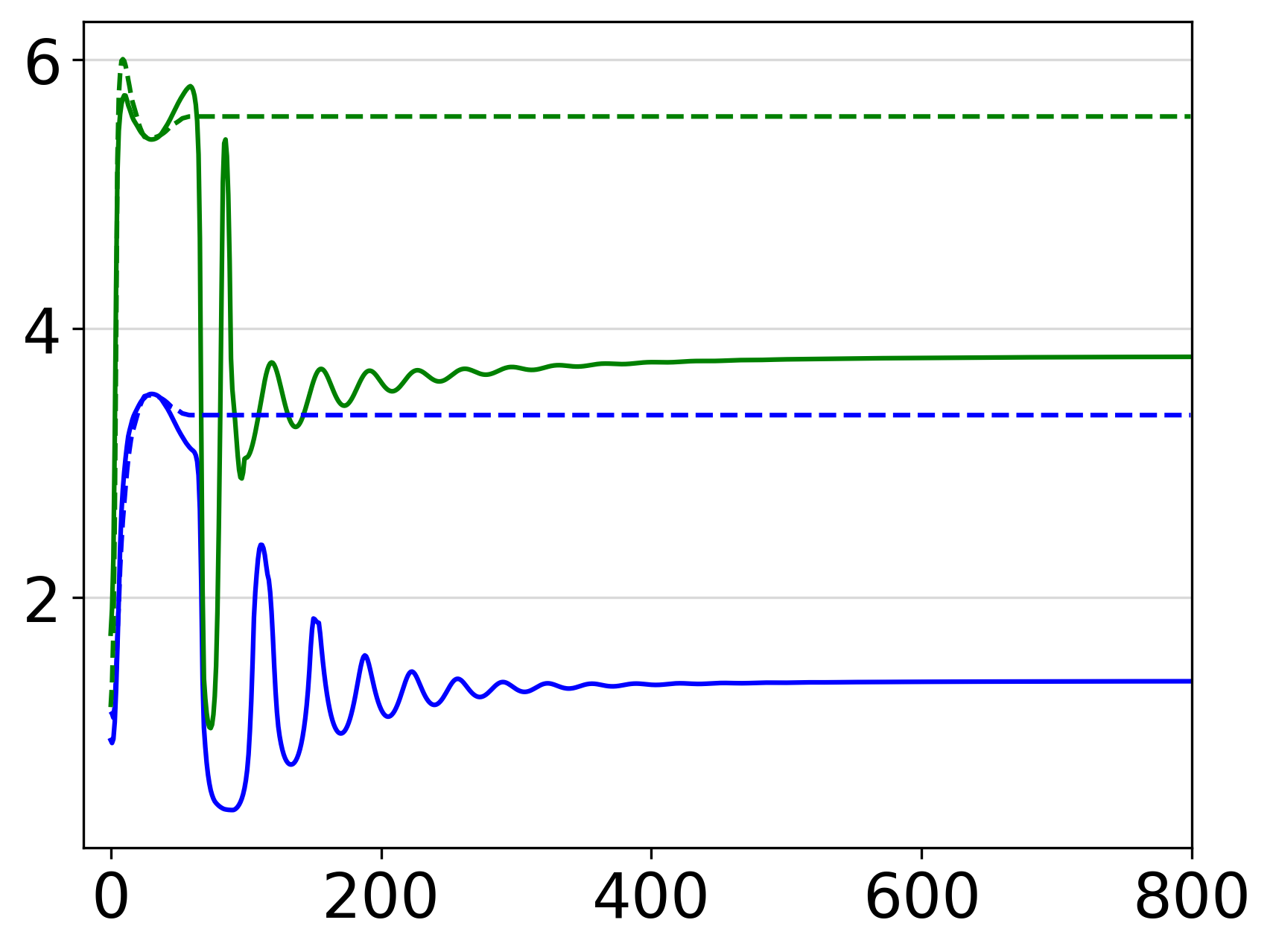}};
        \node[left = 0.4cm of img2,yshift=1.2cm, rotate=90]{utility values};
        \node[below = 0cm of img2]{iteration};
    \end{tikzpicture}
    \caption{Convergence performance of a non-resilient method ( $V_r^\pi(\rho)$: \ref{legend:red}, $V_{g_1}^\pi(\rho)$: \ref{legend:bluedash}, $V_{g_2}^\pi(\rho)$: \ref{legend:greendash}) and Algorithm~\ref{alg: resilient OPG} (ResOPG-PD, $V_r^\pi(\rho)$: \ref{legend:redsolid}, $V_{g_1}^\pi(\rho)$: \ref{legend:blue}, $V_{g_2}^\pi(\rho)$: \ref{legend:green}), with a cost functions $h(\xi) = \alpha \norm{\xi}^2$ for $\alpha = 0.08$, and stepsize $\eta = 0.05$.}
    \label{fig:comparealg}
\end{figure}
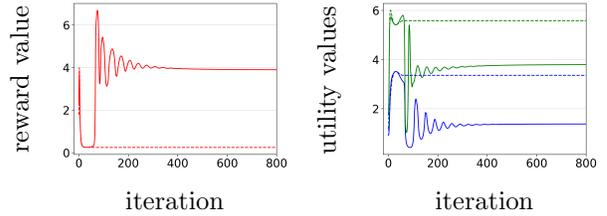

\begin{figure}[h]
    \centering
        \begin{tikzpicture}
        \node[](img1) at(0,0) {\includegraphics[width = 0.3\textwidth]{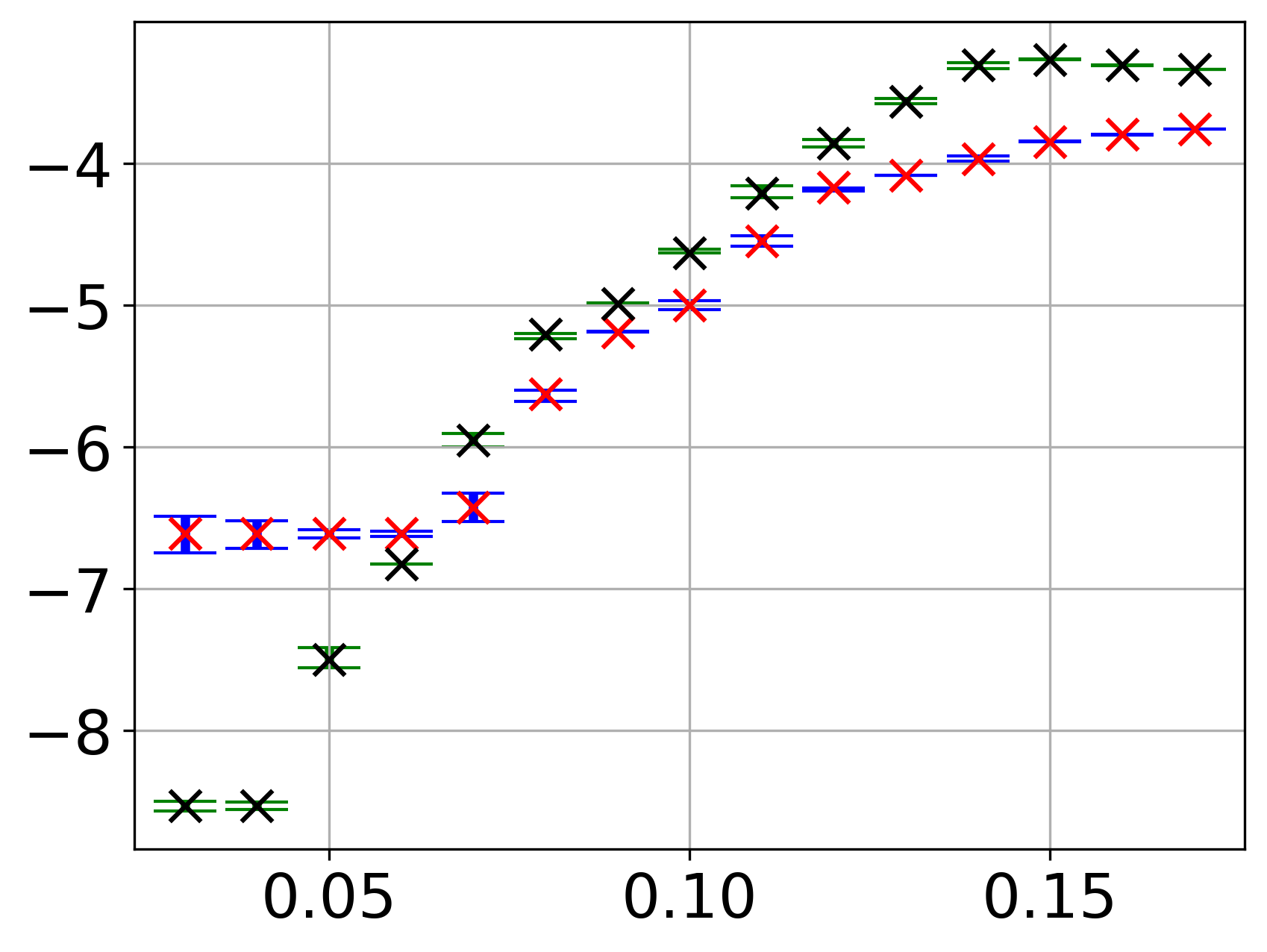}};
        \node[below = 0cm of img1]{relaxation cost};
         \node[left = 0.4cm of img1,yshift=1.0cm, rotate=90]{relaxation};
    \end{tikzpicture}
    \caption{Constraint specifications under different relaxation costs for Algorithm~\ref{alg: resilient PG} (ResPG-PD, $\xi_1$:~\ref{legend:errorbar}, $\xi_2$:~\ref{legend:greenerrorbar}\,) and Algorithm~\ref{alg: resilient OPG} (ResOPG-PD, $\xi_1$:~\ref{legend:x},  $\xi_2$:~\ref{legend:blackx}\,). The relaxation cost function is $h(\xi) = \alpha \norm{\xi}^2$.
    The horizontal axis is the value of $\alpha$ and the vertical axes are relaxations $\xi_1$ and $\xi_2$. 
    The height of~\ref{legend:errorbar} is the oscillation magnitude of ResPG-PD. We run ResPG-PD for $5000$ iterations with stepsize $\eta = 0.01$, and ResOPG-PD for $2000$ iterations with stepsize $\eta = 0.05$.   
    The initial distribution $\rho$ is uniform.
    }
    \label{fig:ResE3}
\end{figure}

Third, 
we generalize the previous monitoring problem to a larger state/action space in Figure~\ref{fig:Monitoring_nonres}, where in a given time a robot has to stay in  blue/green areas for a minimum amount of time while maximizing the time in red area in Figure~\ref{fig:Monitoring_nonres}; see Appendix~\ref{app: experiments} for the detail. We also choose \emph{infeasible} initial constraints. Figure~\ref{fig:Monitoring_nonres}(a) shows that applying a non-resilient method to this infeasible problem yields a policy that does not monitor $S_0$, where the non-resilient method is ResOPG-PD without relaxation update~\cite{ding2023last}. ResOPG-PD’s policy in Figure~\ref{fig:Monitoring_nonres}(b) balances three areas.
Figure~\ref{fig:comparealg} shows that ResOPG-PD gets higher reward value than the non-resilient method by modifying the constraints. However, the non-resilient method does not balance the reward and the constraints.
Figure~\ref{fig:ResE3} shows two relaxation curves for ResPG-PD and ResOPG-PD against different relaxation cost functions.


\section{CONCLUDING REMARKS}\label{sec:conclusion}

To specify the constraint specifications, we have presented an approach by making trade-off between the marginal decrease in the optimal reward value function that results from relaxation and the marginal increase in relaxation cost. We show the existence of such a resilient equilibrium under some mild regularity conditions, and provide a tractable constrained policy optimization that takes this equilibrium as an optimal solution. We provide two constrained policy search algorithms to search for such a resilient equilibrium with convergence guarantees on the optimality gap and constraint violation, generating nearly optimal policy and constraint specification. A series of computational experiments have demonstrated that our resilient constrained policy search methods effectively sustain the trade-off between the reward maximization and the constraint satisfaction even if the problem is infeasible.

For future work, our approach readily enhances existing constrained RL algorithms with resilience for solving constrained MDPs in online/offline learning settings. It is also of our interest to study tighter convergence analysis and sample-based algorithms. For application, it is important to investigate other practical constrained RL problems.     

\newpage
\clearpage
\bibliography{dd-bib.bib}
\bibliographystyle{apalike}

\newpage
\appendix
\onecolumn
\aistatstitle{Resilient Constrained Reinforcement Learning:\\[0.1cm]
Supplementary Materials}

\section{Proofs in Section~\ref{sec:CMDPs}}

We state proofs for claims made in Section~\ref{sec:CMDPs}. 

\subsection{Proof of Lemma~\ref{lem:primal_function}}\label{app:CMDPs}

By the definition, $V^\star(\xi)$ is real-valued for $\xi\in\Xi$. The monotonic property of $V^\star(\xi)$ over $\xi\in\Xi$ is straightforward from Problem~\eqref{eq:CMDP_relaxed}. We next prove the concavity.  

    
    We first check the convexity of the domain $\Xi$. For any $\xi$, $\xi'\in\Xi$, there exist $\pi$, $\pi'\in \Pi$ such that $V_{g_i}^\pi(\rho)\geq \xi_i$ and $V_{g_i}^{\pi'}(\rho)\geq \xi_i'$ for $i = 1,\ldots,m$. Let the occupancy measures associated with $\pi$, $\pi'$ be $q$, $q'$, respectively. Thus, $\langle g_i, q\rangle \geq \xi_i$ and $\langle g_i, q'\rangle \geq \xi_i$ for $i=1,\ldots,m$, which implies $\langle g_i, \alpha q+(1-\alpha)q'\rangle \geq \alpha\xi_i+(1-\alpha)\xi_i'$, i.e., the policy induced by $\alpha q+(1-\alpha)q'$ meets the constraint. Therefore, $\alpha \xi+(1-\alpha)\xi' \in \Xi$.

    We next show the convexity of $V^\star(\xi)$ over $\xi\in\Xi$. 
    We re-formulate Problem~\eqref{eq:CMDP_relaxed} in terms of occupancy measure,  \begin{equation}\label{eq:CMDP_occupancy}
    \begin{array}{rl}
         \displaystyle \maximize_{q\,\in\,\mathcal{Q}} & \langle r, q\rangle
         \\
         \subject&  
         \langle u_i, q\rangle \; \geq \; \xi_i \; \text{ for all } \; i = 1,\ldots,m
    \end{array}
\end{equation}
where $q$ is the occupancy measure that lives in a polytope $\mathcal{Q}$ specified by Bellman flow equations. Instead of policy $\pi$, we work with the occupancy measure $q$ in Problem~\eqref{eq:CMDP_occupancy}. The primal function $V^\star(\xi)$, $\xi\in\Xi$ does not change, because of the one-to-one correspondence between $\pi$ and $q$. The rest is straightforward from convex analysis. We define the function,
\[
F(q, \xi) 
\; = \;
\begin{cases}
    \langle r, q\rangle \;\;\;\;\text{ if } \; \langle g_i, q\rangle \geq \xi_i \, \text{ for } i=1,\ldots,m;
    \\
    -\infty \;\;\;\;\;\text{ otherwise. }
\end{cases}
\]
and its domain,
\[
    \text{dom}(F)  
    \; = \; 
    \left\{\,
    (q, \xi) \,\vert\,
    q \in \mathcal{Q}, \xi\in\mathbb{R}_\gamma^m, \langle g_i, q\rangle \geq \xi_i \, \text{ for } i=1,\ldots,m
    \,\right\}.
\]
Because of linearity, $\text{dom}(F)$ is a convex set, and $F(q,\xi)$ is a concave function in its domain. We notice that $V^\star(\xi) = \sup_{q\,\in\,\mathcal{Q}} F(q,\xi)$ and the convex domain $\Xi$. Therefore, $V^\star(\xi)$ is a concave function.

\subsection{Proof of Equation~\eqref{eq:conjugate_primal_function}}\label{app:dual_function} 

\[
\begin{array}{rcl}
     D(\lambda) & = & \displaystyle
     \sup_{\pi\,\in\,\Pi} \; V_{r+\lambda^\top g}^\pi(\rho)
     \\[0.2cm]
     & = & \displaystyle
     \sup_{  \{ (\xi,\pi) \,\vert\,\pi\,\in\,\Pi,\, \xi \,\in\,\Xi, V_{g_i}^\pi(\rho)\,\geq\, \xi_i,\, i = 1,\ldots\,m \}} \; V_{r+\lambda^\top g}^\pi(\rho) 
     \\[0.2cm]
     & = & \displaystyle
     \sup_{  \{ (\xi,\pi) \,\vert\,\pi\,\in\,\Pi,\, \xi \,\in\,\Xi, V_{g_i}^\pi(\rho)\,\geq\, \xi_i,\, i = 1,\ldots\,m \}} \; \left\{ V_{r}^\pi(\rho)+ \lambda^\top \xi \right\}
     \\[0.2cm]
     & = & \displaystyle
     \sup_{ \xi \,\in\,\Xi}
     \sup_{  \{\pi\,\in\,\Pi,\, V_{g_i}^\pi(\rho)\,\geq\, \xi_i,\, i = 1,\ldots\,m \}} \; \left\{ V_{r}^\pi(\rho)+ \lambda^\top \xi \right\}
     \\[0.2cm]
     & = & \displaystyle
     \sup_{ \xi \,\in\,\Xi}  \;  \left\{ \lambda^\top \xi - (- V^\star(\xi)) \right\}.
\end{array}
\]

\subsection{Proof of Lemma~\ref{lem:geometric_multiplier}}\label{app:geometric_multiplier}

    There are two directions.
    
    (ii)$\implies$(i): From  the geometric multiplier $\lambda$ for Problem~\eqref{eq:CMDP_relaxed}, 
        \[
        \begin{array}{rcl}
             V^\star(\xi)
             &  =  & \displaystyle \sup_{\pi\,\in\,\Pi} \left\{ V_{r+\lambda^\top g}^\pi(\rho) - \lambda^\top \xi\right\}
             \\[0.2cm]
            & = & \displaystyle
            \sup_{ \xi' \,\in\,\Xi}  \;  \left\{ \lambda^\top \xi' - (- V^\star(\xi')) \right\} - \lambda^\top \xi
        \end{array}
        \]
        where the second equality is due to~\eqref{eq:conjugate_primal_function}. Therefore,
        \begin{equation}\label{eq:subgradient}
         -V^\star(\xi') 
        \;\geq\; 
        -V^\star(\xi)+
        \lambda^\top (\xi' -\xi)
        \text{ for all } \xi'\in\mathbb{R}^m
        \end{equation}
        which shows that $\lambda$ is a subgradient of $-V^\star(\xi)$ at $\xi\in \Xi$.

        (i)$\implies$(ii): Assume that \eqref{eq:subgradient} holds for some $\lambda$. Since $V^\star(\xi)$ is monotonically non-increasing with respect to the coordinates of $\xi$, thus $\lambda\geq 0$. Otherwise, $ V^\star(\xi') +
        \lambda^\top (\xi' -\xi)$ would be unbounded below which makes~\eqref{eq:subgradient} invalid. From~\eqref{eq:subgradient}, we have
        \[
        V^\star(\xi)
        \;\geq\; \sup_{\xi'\,\in\,\Xi} \left\{
        V^\star(\xi') + \lambda^\top (\xi'-\xi)
        \right\}
        \;=\; D(\lambda) - \lambda^\top\xi
        \]
        where the second equality is due to~\eqref{eq:conjugate_primal_function}. We notice that $D(\lambda)-\lambda^\top\xi$ is the dual function for Problem~\eqref{eq:CMDP_relaxed} and the weak duality $V^\star(\xi)\leq D(\lambda)-\lambda^\top\xi$. Therefore, $\lambda$ is a geometric multiplier for Problem~\eqref{eq:CMDP_relaxed}.

\section{Proofs in Section~\ref{sec:resilient CRL}}

We state proofs for claims made in Section~\ref{sec:resilient CRL}. 

\subsection{Proof of Lemma~\ref{lem:resilient_eq}}\label{app:resilient_eq}

To show the existence, it is equivalent to show existence of subgradients for the function $\textbf{0} \in \partial (-V^\star(\bar\xi) + h(\bar\xi))$ for some $\bar\xi\in\Xi$.
        We introduce the set $\bar\Xi$,
        \begin{equation}\label{eq:primal_function_V+h}
        \bar\Xi 
        \; \DefinedAs \;
        \argmin_{\xi\,\in\,\Xi} 
        \; 
        \left\{\,-V^\star(\xi) \,+\, h(\xi) \,\right\}.
        \end{equation}
        We notice that $\Xi$ is an effective domain for $-V^\star(\xi)+h(\xi)$. Because of the concavity of $V^\star(\xi)$ in Lemma~\ref{lem:primal_function} and the convexity of $h(\xi)$, $-V^\star(\xi)+h(\xi)$ is a convex function on $\Xi$. Thus, the set $\bar\Xi$ is nonempty and $\textbf{0} \in \partial (-V^\star(\xi') + h(\xi'))$ for any $\xi'\in\bar\Xi$ according to the first-order optimality~\cite[Theorem~8.2]{recht-wright19}. Therefore, the existence is proved by simply taking a resilient equilibrium $\xi^\star = \bar\xi\in\bar\Xi$.

    From the further hypothesis on $h$, the function $-V^\star(\xi)+h(\xi)$ is strictly convex. Thus, the minimizer in Problem~\eqref{eq:primal_function_V+h} is unique or $\bar\Xi$ is a singleton. Therefore, the uniqueness holds.

\subsection{Proof of Lemma~\ref{lem:resilient_mono}}\label{app:resilient_mono}

    By the concavity of $V^\star$ in Lemma~\ref{lem:primal_function}, if $p \in \partial V^\star(\xi)$ and $p'\in \partial V^\star(\xi')$ for $\xi$, $\xi'\in\Xi$,  then,
    \[
    \begin{array}{rcl}
         V^\star(\xi') 
         & \leq &
         V^\star(\xi) + \langle p, \xi' - \xi\rangle
         \\[0.2cm]
         V^\star(\xi) 
         & \leq &
         V^\star(\xi') + \langle p', \xi - \xi'\rangle.
    \end{array}
    \]
    Thus,
    \[
        \langle p'-p, \xi - \xi'\rangle
        \;\geq \; 0
    \]
    which implies the second inequality.
    Similarly, the convexity of $h$ yields
    \[
    \langle
    \xi-\xi', \nabla h(\xi')-\nabla h(\xi)
    \rangle 
    \; \leq \; 0
    \]
    which implies the first inequality.

\subsection{Proof of Theorem~\ref{thm:resilient_geometric_multiplier}}\label{app:resilient_geometric_multiplier}

By the geometric multiplier $\lambda\geq0$,
    \[
    \begin{array}{rcl}
         V^\star(\bar\xi) & = & \displaystyle
         \sup_{\pi\,\in\,\Pi} \;  \left\{ V_{r+\lambda^\top g}^\pi(\rho) - \lambda^\top\bar\xi \right\}
         \\[0.2cm]
         & \geq & V_{r+\lambda^\top g}^{\pi^\star(\xi)}(\rho) - \lambda^\top\bar\xi 
         \\[0.2cm]
         & = & V_{r}^{\pi^\star(\xi)}(\rho) + \lambda^\top( V_{g}^{\pi^\star(\xi)}(\rho)- \bar\xi) 
         \\[0.2cm]
         & \geq & V^\star(\xi) +\lambda^\top(\xi-\bar\xi)
    \end{array}
    \]
    where we set $\pi = \pi^\star(\xi)$ in the first inequality, and the second inequality is due to that $V^\star(\xi)\DefinedAs V_r^{\pi^\star(\xi)}(\rho)$ and $V_{g_i}^{\pi^\star(\xi)}(\rho)\geq \xi_i$ for all $i=1,\ldots,m$.

    Therefore, $-V^\star(\xi) \geq -V^\star(\bar\xi) + \lambda^\top(\xi-\bar\xi)$ for all $\xi\in \Xi$, i.e., $\lambda$ is a subgradient of $-V^\star(\xi)$ at $\bar\xi$. By the assumption, $\nabla h(\bar\xi)$ is a subgradient of $V^\star(\bar\xi)$, which proves that $\bar\xi$ is a resilient equilibrium.

\subsection{Proof of Corollary~\ref{cor:resilient_geometric_multiplier}}\label{app:resilient_geometric_multiplier_cor}

It is straightforward to verify that $\lambda^\star(\bar\xi)$ is a geometric multiplier,
    \[
    V^\star(\bar\xi) 
    \; = \;
    D^\star(\bar\xi)
    \;  = \; 
    D(\lambda^\star(\bar\xi);\bar\xi)
    \; = \;
    \sup_{\pi\,\in\,\Pi} \left\{
        V_{r+(\lambda^\star(\bar\xi))^\top g}^{\pi}(\rho) - (\lambda^\star(\bar\xi))^\top\bar\xi \right\}.
    \]
    By Theorem~\ref{thm:resilient_geometric_multiplier} and $\nabla h(\bar\xi)+\lambda^\star(\bar\xi)=0$, $\bar\xi$ is a resilient equilibrium. 

\subsection{Proof of Lemma~\ref{lem:regularized_optimal_policy}}\label{app:regularized_optimal_policy}

It is equivalent to show that 
    \[
    \bar\xi^\star \; \in \; \argmax_{\xi\,\in\,\Xi} \; \left\{V^\star(\xi) - h(\xi)\right\}
    \]
    because of the concavity of $V^\star(\xi)-h(\xi)$ over $\xi\in\Xi$, and that the first-order optimality condition $\mathbf{0}\in \partial \left(-V^\star(\bar\xi^\star)+h(\bar\xi^\star)\right)$ equals to the resilient equilibrium's condition. 

    By the optimality of $(\bar\pi^\star,\bar\xi^\star)$,
    \[
    V_r^{\bar\pi^\star}(\rho) - h(\bar\xi^\star)
    \; \geq \;
    V_r^{\pi^\star(\xi)}(\rho) - h(\xi)
    \; = \; 
    V^\star(\xi) - h(\xi)
    \]
    where the right hand side of inequality particularly uses a pair $(\pi^\star(\xi),\xi)$ in which $\pi^\star(\xi)$ is an optimal policy of Problem~\eqref{eq:CMDP_relaxed} for some fixed $\xi\in\Xi$, and the equality is clear from $V^\star(\xi) \DefinedAs V_r^{\pi^\star(\xi)}(\rho)$. For the left hand side of the inequality, application of $\pi^\star(\bar\xi^\star)$ leads to $V^\star(\bar\xi^\star) \DefinedAs V_r^{\pi^\star(\bar\xi^\star)}(\rho) \geq V_r^{\bar\pi^\star}(\rho)$ and thereby,
    \[
    V^\star(\bar\xi^\star) - h(\bar\xi^\star)
    \; \geq \;
    V^\star(\xi) 
    - h(\xi)  \; \text{  for all } \xi\in\Xi.
    \]
    Therefore, $\textbf{0}$ is a subgradient of $-V^\star(\xi)+h(\xi)$ at $\bar\xi^\star\in\Xi$.

\subsection{Proof of Theorem~\ref{thm:strong_duality}}\label{app:strong_duality}

By the weak duality, $V_h^\star \leq D_h^\star$. The rest is to show that $V_h^\star \geq D_h^\star$. To proceed, we first make a few observations. It is easy to show the convexity of the set $\mathcal{Z}$,
    \[
    \mathcal{Z}
    \; \DefinedAs \;
    \left\{
    z \in \mathbb{R}^{m+1}
    \,\vert\,
    \exists (\pi,\xi)\in\Pi\times\Xi, V_r^{\pi}(\rho) - h(\xi) \geq z_0 \text{ and } V_{g_i}^{\pi}(\rho) - \xi_i\geq z_i \text{ for all } i = 1,\ldots,m 
    \right\}.
    \]
    By the optimality of $(\bar\pi^\star,\bar\xi^\star)$, $V_h^\star = V_r^{\bar\pi^\star}(\rho) - h(\bar\xi^\star)$ and $V_{g_i}^{\bar\pi^\star}(\rho)-\bar\xi_i^\star \geq 0$ for all $i=1,\ldots,m$, i.e., $(V_h^\star, \textbf{0})\in \mathcal{Z}$. Hence, $\mathcal{Z}$ is non-empty. In fact, $\mathcal{Z}$ is a convex set. Assume $z$, $z'\in\mathcal{Z}$ and $\alpha \in [0,1]$. Thus, there exist $\pi$, $\pi'\in\Pi$ and $\xi$, $\xi'\in\Xi$ such that 
    \[
    V_r^{\pi}(\rho) - h(\xi) \; \geq\; z_0
    \; \text{ and } \;
    V_r^{\pi'}(\rho) - h(\xi') \; \geq \; z_0'
    \]
    \[
    V_{g_i}^{\pi}(\rho) - \xi_i \; \geq\; z_i
    \; \text{ and } \;
    V_{g_i}^{\pi'}(\rho) - \xi'_i \;\geq\; z_i' \;\text{ for all } i=1,\ldots,m.
    \]
    Let the occupancy measures associated with $\pi$, $\pi'$ be $q$, $q'$, respectively. Thus, $\langle r, q\rangle -h(\xi)\geq z_0$ and $\langle r, q'\rangle - h(\xi')\geq z_0'$, which implies $\langle r, \alpha q+(1-\gamma)q'\rangle - h(\alpha \xi+(1-\alpha)\xi') \geq \alpha z_0 +(1-\alpha)z_0'$. Meanwhile, $\langle g_i, q\rangle - \xi_i \geq z_i$ and $\langle g_i, q'\rangle - \xi_i \geq z_i'$ for $i=1,\ldots,m$, which implies $\langle g_i, \alpha q+(1-\alpha)q'\rangle -( \alpha\xi_i+(1-\alpha)\xi_i')\geq \alpha z_i + (1-\alpha)z_i'$, i.e., the policy induced by $\alpha q+(1-\alpha)q'$ meets the constraint in $\mathcal{Z}$. Therefore, $\alpha z+(1-\alpha)z'\in\mathcal{Z}$. 

    To show that $V_h^\star\geq D_h^\star$, it is sufficient to prove that there exists $\lambda\in\Lambda$ such that
    \begin{equation}\label{eq:existence}
    V_h^\star
    \;  \geq \;
    D_h(\lambda) 
    \;  \DefinedAs \;
    \sup_{\pi\,\in\,\Pi}
    \left\{
        V_{r+\lambda^\top g}^\pi(\rho) - h(\xi) - \lambda^\top \xi
    \right\}.
    \end{equation}
    We notice that $(V_h^\star,\textbf{0})\in\partial{Z}$, where $\partial{Z}$ is the boundary set of $\mathcal{Z}$. If not, then
    $(V_h^\star,\textbf{0}) \in \text{int}(\mathcal{Z})$, i.e., there exists a small ball around $(V_h^\star,\textbf{0})$ inside $\mathcal{Z}$, which contradicts the optimality of $V_h^\star$. Since $\mathcal{Z}$ is a convex set and $(V_h^\star,\textbf{0})\in\partial{Z}$, by the supporting hyperplane theorem, there exists $\hat\lambda\DefinedAs (\hat\lambda_0, \hat\lambda_1,\ldots, \hat\lambda_m)\in\mathbb{R}^{m+1}$ such that
    \begin{equation}\label{eq:hyperplane}
    \big[ V_h^\star\;\;\textbf{0}^\top \big]
    \hat\lambda \;  \geq \; z^\top \hat\lambda \; \text{ for all } z\in\mathcal{Z}.
    \end{equation}
    We can show that $\hat\lambda\geq 0$ and $\hat\lambda_0>0$ by contradiction. Assume $\hat\lambda_i<0$ for some $i$. Since $\mathcal{Z}$ is unbounded from below, we can always select a very negative $z_i$ that fails~\eqref{eq:hyperplane}. Hence, $\hat\lambda_i\geq 0$ for all $i  = 0,1,\ldots,m$. On the other hand, assume $\hat\lambda_0=0$. Thus, \eqref{eq:hyperplane} reduces to $0  \geq z^\top \hat\lambda$ for all $z\in\mathcal{Z}$. Non-negativity of $\hat\lambda$ implies that there exists $i$ from $1$ to $m$ such that $z_i\leq 0$, which contradicts the strict feasibility that demands a positive $z\in \mathcal{Z}$. Therefore, we can denote $\lambda^\dagger\DefinedAs\hat\lambda/\hat\lambda_0$ and 
    \[
    (\pi^\dagger, \xi^\dagger) 
    \; \DefinedAs \;
    \argmax_{\pi\,\in\,\Pi,\, \xi\,\in\,\Xi}
    \left\{
        V_{r+(\lambda^\dagger)^\top g}^\pi(\rho) - h(\xi) -(\lambda^\dagger)^\top \xi
    \right\}.
    \]
    We notice that $(V_r^{\pi^\dagger}(\rho) - h(\xi^\dagger), V_{g}^{\pi^\dagger}(\rho) - \xi^\dagger) \in \mathcal{Z}$, and $\lambda^\dagger_0 = 1$. By the dual function and~\eqref{eq:hyperplane},
    \[
    D_h(\lambda^\dagger) 
    \; = \; 
    \big\langle(V_r^{\pi^\dagger}(\rho) - h(\xi^\dagger), V_{g}^{\pi^\dagger}(\rho) - \xi^\dagger),  \lambda^\dagger\big\rangle
    \; \leq\; V_h^\star
    \]
    which proves the existence~\eqref{eq:existence}.

\subsection{Proof of Corollary~\ref{cor:bounded_dual}}\label{app:bounded_dual}

We denote the level set of the dual function by $\Lambda_{a} \DefinedAs  \{ \lambda \in \mathbb{R}_+^m \,\vert\, D_h(\lambda) \leq a\}$ for $a\in\mathbb{R}$. For any $\lambda \in \Lambda_a$,
    \[
        a 
        \; \geq \;  
        D_h(\lambda)
        \; = \; 
        V_r^{\bar\pi}(\rho) - h(\bar\xi) 
        + \lambda^\top (V_g^{\bar\pi}(\rho)-\bar\xi)
        \; \geq \; 
        V_r^{\bar\pi}(\rho) - h(\bar\xi) 
        + c \, \lambda^\top \mathbf{1}
    \]
    where $(\bar\pi,\bar\xi)$ is a Slater point in Assumption~\ref{as:feasibility_regularized}. Taking $a = D_h^\star$ or  $ V_r^\star(\rho) - h(\bar\xi^\star)$ leads to $\Lambda_a = \Lambda^\star $ and
    \[
       \sum_{i\,=\,1}^m \lambda_i \;\leq\; 
       \frac{V_r^\star(\rho)  - h(\bar\xi^\star)  - (V_r^{\bar\pi} - h(\bar\xi))}{c} 
       \; \text{  for all }\lambda \in \Lambda_a
    \]
    which implies the bound on $\bar\lambda_i^\star$ for $i = 1,\ldots,m$.

\section{Proofs in Section~\ref{sec: resilient CPL}}

We state proofs for claims made in Section~\ref{sec: resilient CPL}. 

\subsection{Proof of Theorem~\ref{thm: average-value convergence}}\label{app: average-value convergence}

\begin{lemma}\label{lem:optimality NPG}
    In Algorithm~\ref{alg: resilient PG}, for any $\pi\in \Delta(A)$ and $\xi \in \Xi$
    \[
    \begin{array}{rcl}
         &  & 
         \!\!\!\!  \!\!\!\!  \!\! \displaystyle\eta \langle Q_{r+\lambda_t^\top g}^{\pi_t}(s,\cdot), (\pi- \pi_{t+1})(\cdot\,\vert\,s)\rangle + \frac{1}{2} \norm{ \pi_{t+1}(\cdot\,\vert\,s) - \pi_{t}(\cdot\,\vert\,s)}^2
         \\[0.2cm]
         & \leq & 
         \displaystyle 
         \frac{1}{2} \norm{\pi(\cdot\,\vert\,s) - \pi_t(\cdot\,\vert\,s)}^2 - \frac{1}{2}\norm{\pi(\cdot\,\vert\,s) - \pi_{t+1}(\cdot\,\vert\,s)}^2
    \end{array}
    \]
    and
    \[
    \begin{array}{rcl}
         &  & 
         \!\!\!\!  \!\!\!\!  \!\! \displaystyle
         \eta \langle -\nabla h(\xi_t) - \lambda_t, \xi  - \xi_{t+1} \rangle  + \frac{1}{2} \norm{\xi_{t+1} - \xi_t}^2
         \\[0.2cm]
         & \leq & \displaystyle
         \frac{1}{2}\norm{\xi - \xi_{t}}^2 - \frac{1}{2}\norm{\xi-\xi_{t+1}}^2.
    \end{array}
    \]
\end{lemma}
\begin{proof}
    By the optimality of $\pi_{t+1}$,
    \[
    \langle \eta Q_{r+\lambda_t^\top g}^{\pi_t}(s,\cdot) - (\pi_{t+1} - \pi_t)(\cdot\,\vert\,s), (\pi- \pi_{t+1})(\cdot\,\vert\,s)\rangle 
    \; \leq \; 0\; \text{ for any } \pi.
    \]
    Direct application of the equality $\frac{1}{2} \norm{\pi(\cdot\,\vert\,s) - \pi_t(\cdot\,\vert\,s)}^2 = \frac{1}{2} \norm{\pi_{t+1}(\cdot\,\vert\,s) - \pi_t(\cdot\,\vert\,s)}^2 + \langle \pi_{t+1}(\cdot\,\vert\,s)-\pi_{t}(\cdot\,\vert\,s), \pi(\cdot\,\vert\,s)- \pi_{t+1}(\cdot\,\vert\,s) \rangle+ \frac{1}{2} \norm{\pi(\cdot\,\vert\,s) - \pi_{t+1}(\cdot\,\vert\,s)}^2$ leads to the first inequality.
    Similarly, the optimality of $\xi_{t+1}$ shows that
    \[
        \displaystyle
        \langle \eta (-h(\xi_t) - \lambda_t) - (\xi_{t+1}- \xi_t),  \xi - \xi_{t+1}\rangle 
        \;\leq\; 0.
    \]
    In combination of the equality $\frac{1}{2}\norm{\xi - \xi_t}^2 = \frac{1}{2}\norm{\xi_{t+1} - \xi_t}^2 + \langle \xi_{t+1}-\xi_t, \xi-\xi_{t+1}\rangle + \frac{1}{2}\norm{\xi - \xi_{t+1}}^2$, we conclude the second inequality.
\end{proof}

\begin{lemma}
    In Algorithm~\ref{alg: resilient PG}, for any $s$ and $t$,
    \[
    \begin{array}{rcl}
         \displaystyle
        V_{r}^{\pi_{t+1}}(s) - V_r^{\pi_t}(s) + \lambda_t^\top (V_{g}^{\pi_{t+1}}(s) - V_g^{\pi_t}(s))
        & \geq &
        \displaystyle 
        \frac{1}{\eta(1-\gamma)} \left(  \norm{\pi_{t+1}(\cdot\,\vert\,s) - \pi_{t}(\cdot\,\vert\,s)}^2 \right).
    \end{array}
    \]
\end{lemma}
\begin{proof}
    By the performance difference lemma,
    \[
    \begin{array}{rcl}
         &&
         \!\!\!\!  \!\!\!\!  \!\! \displaystyle
         V_{r}^{\pi_{t+1}}(s) - V_r^{\pi_t}(s) + \lambda_t^\top (V_{g}^{\pi_{t+1}}(s) - V_g^{\pi_t}(s))
         \\[0.2cm]
         & = & 
         \displaystyle
         \frac{1}{1-\gamma} \mathbb{E}_{s'\,\sim\,d_{s}^{\pi_{t+1}}} 
         \left[ 
         \langle Q_{r+\lambda_t^\top g} (s',\cdot), (\pi_{t+1}-\pi_t)(\cdot\,\vert\,s')\rangle 
         \right].
    \end{array}
    \]
    Application of the first inequality in Lemma~\ref{lem:optimality NPG} with $\pi=\pi_t$ leads to our desired inequality.
\end{proof}

\begin{lemma}\label{lem: average bound PG}
    In Algorithm~\ref{alg: resilient PG}, for any $T>0$,
    \[
    \begin{array}{rcl}
         &&
         \!\!\!\!  \!\!\!\!  \!\! \displaystyle
         \frac{1}{T} \sum_{t \, = \,0}^{T-1} (V_r^\star(\rho) -h(\bar\xi^\star) - (V_r^{\pi_t}(\rho) -h(\xi_t))) +  \frac{1}{T} \sum_{t \, = \,0}^{T-1} \lambda_t^\top (V_g^\star(\rho) -\bar\xi^\star - (V_g^{\pi_t}(\rho) -\xi_t))
         \\[0.2cm]
         & \leq  & \displaystyle
         \frac{1}{(1-\gamma)^2 T}
         + \frac{1}{\eta (1-\gamma)T} + \frac{4\eta m}{(1-\gamma)^2} + \frac{\eta(L_h+1)^2 m}{(1-\gamma)^2}.
    \end{array}
    \]
\end{lemma}
\begin{proof}
    By the performance difference lemma,
    \begin{equation}\label{eq:PDL PG}
    \begin{array}{rcl}
         &&
         \!\!\!\!  \!\!\!\!  \!\! \displaystyle
          V_r^\star(s) -h(\bar\xi^\star) - (V_r^{\pi_t}(s) -h(\xi_t)) +  \lambda_t^\top (V_g^\star(s) -\bar\xi^\star - (V_g^{\pi_t}(s) -\xi_t))
          \\[0.2cm]
          & = &  \displaystyle
          \frac{1}{1-\gamma} \mathbb{E}_{s'\,\sim\,d_s^\star} \left[
          \langle Q_{r+\lambda_t^\top g}^{\pi_t}(s',\cdot), (\bar\pi^\star-\pi_{t})(\cdot\,\vert\,s')
          \right]
           -(h(\bar\xi^\star) - h(\xi_t)) + \lambda_t^\top (-\bar\xi^\star+ \xi_t)
           \\[0.2cm]
           & = & \displaystyle
          \frac{1}{1-\gamma} \mathbb{E}_{s'\,\sim\,d_s^\star} \left[
          \langle Q_{r+\lambda_t^\top g}^{\pi_t}(s',\cdot), (\bar\pi^\star-\pi_{t+1})(\cdot\,\vert\,s')
          \right]
          \\[0.2cm]
           &  & \displaystyle+
           \,\frac{1}{1-\gamma} \mathbb{E}_{s'\,\sim\,d_s^\star} \left[
          \langle Q_{r+\lambda_t^\top g}^{\pi_t}(s',\cdot), (\pi_{t+1}-\bar\pi^\star)(\cdot\,\vert\,s')
          \right]
          \\[0.2cm]
           &  & \displaystyle
           -\,(h(\bar\xi^\star) - h(\xi_t)) + \lambda_t^\top (-\bar\xi^\star+ \xi_t)
           \\[0.2cm]
           & \leq & \displaystyle
           \frac{1}{2\eta(1-\gamma)} \mathbb{E}_{s'\,\sim\,d_s^\star} \left[ \norm{\bar\pi^\star(\cdot\,\vert\,s) - \pi_t(\cdot\,\vert\,s)}^2 -\norm{\bar\pi^\star(\cdot\,\vert\,s) - \pi_{t+1}(\cdot\,\vert\,s)}^2 \right] 
           \\[0.2cm]
           &  & \displaystyle
           +\,\frac{1}{1-\gamma} \mathbb{E}_{s'\,\sim\,d_s^\star} \left[
          \langle Q_{r+\lambda_t^\top g}^{\pi_t}(s',\cdot), (\pi_{t+1}-\bar\pi^\star)(\cdot\,\vert\,s')
          \right]
          \\[0.2cm]
           &  & \displaystyle
           -\,(h(\bar\xi^\star) - h(\xi_t)) + \lambda_t^\top (-\bar\xi^\star+ \xi_t)
    \end{array}
    \end{equation}
    where the inequality is due to the first inequality in Lemma~\ref{lem:optimality NPG} with $\pi = \bar\pi^\star$. To further bound the inequality above, we first notice that $\langle Q_{r+\lambda_t^\top g}^{\pi_t}(s,\cdot), (\pi_{t+1}-\pi_t)(\cdot\,\vert\,s) \rangle \geq 0$ for any $s$, when we set $\pi= \pi_t$ in the first inequality in Lemma~\ref{lem:optimality NPG}. Thus, 
    \[
    \begin{array}{rcl}
         &&
         \!\!\!\!  \!\!\!\!  \!\! \displaystyle
        \mathbb{E}_{s'\,\sim\,d_\rho^\star}
        \left[\langle Q_{r+\lambda_t^\top g}^{\pi_t}(s',\cdot), (\pi_{t+1}-\pi_t)(\cdot\,\vert\,s') \rangle\right]
        \\[0.2cm]
        & = & \displaystyle
        \sum_{s'}\frac{d_\rho^\star(s')}{d_{d_\rho^\star}^{\pi_{t+1}}(s')}d_{d_\rho^\star}^{\pi_{t+1}}(s')
        \langle Q_{r+\lambda_t^\top g}^{\pi_t}(s',\cdot), (\pi_{t+1}-\pi_t)(\cdot\,\vert\,s') \rangle
         \\[0.2cm]
        & \leq & \displaystyle
        \frac{1}{1-\gamma}\sum_{s'}d_{d_\rho^\star}^{\pi_{t+1}}(s')
        \langle Q_{r+\lambda_t^\top g}^{\pi_t}(s',\cdot), (\pi_{t+1}-\pi_t)(\cdot\,\vert\,s') \rangle
        \\[0.2cm]
        & = & \displaystyle 
        (V_r^{\pi_{t+1}}(d_\rho^\star)-V_r^{\pi_{t}}(d_\rho^\star))
        + \lambda_t^\top (V_g^{\pi_{t+1}}(d_\rho^\star)-V_g^{\pi_{t}}(d_\rho^\star))
    \end{array}
    \]
    where the inequality is due to that $d_{d_\rho^\star}^{\pi_{t+1}}\geq (1-\gamma)d_\rho^\star$. Hence, we further bound~\eqref{eq:PDL PG} as
    \[
    \begin{array}{rcl}
         &&
         \!\!\!\!  \!\!\!\!  \!\! \displaystyle
          V_r^\star(\rho) -h(\bar\xi^\star) - (V_r^{\pi_t}(\rho) -h(\xi_t)) +  \lambda_t^\top(V_g^\star(\rho) -\bar\xi^\star - (V_g^{\pi_t}(\rho) -\xi_t))
          \\[0.2cm]
           & \leq & \displaystyle
           \frac{1}{2\eta(1-\gamma)} \mathbb{E}_{s'\,\sim\,d_\rho^\star} \left[ \norm{\bar\pi^\star(\cdot\,\vert\,s) - \pi_t(\cdot\,\vert\,s)}^2 -\norm{\bar\pi^\star(\cdot\,\vert\,s) - \pi_{t+1}(\cdot\,\vert\,s)}^2 \right] 
           \\[0.2cm]
           &  & \displaystyle
           +\,\frac{1}{1-\gamma} \left( (V_r^{\pi_{t+1}}(d_\rho^\star)-V_r^{\pi_{t}}(d_\rho^\star))
        + \lambda_t^\top (V_g^{\pi_{t+1}}(d_\rho^\star)-V_g^{\pi_{t}}(d_\rho^\star)) \right)
          \\[0.2cm]
           &  & \displaystyle
           -\,(h(\bar\xi^\star) - h(\xi_t)) + \lambda_t^\top (-\bar\xi^\star+ \xi_t).
    \end{array}
    \]
    By the convexity of $h$, $h(\bar\xi^\star) \geq h(\xi_t) + \langle \nabla h(\xi_t), \bar\xi^\star - \xi_t \rangle$. Thus, 
    \[
    \begin{array}{rcl}
         && \!\!\!\!  \!\!\!\!  \!\!
         \displaystyle
        - (h(\bar\xi^\star) - h(\xi_t)) + \lambda_t^\top (-\bar\xi^\star+\xi_t)
        \\[0.2cm]
        & \leq & \langle -\nabla h(\xi_t) - \lambda_t, \bar\xi^\star - \xi_t \rangle
        \\[0.2cm]
        & = & \langle -\nabla h(\xi_t) - \lambda_t, \bar\xi^\star - \xi_{t+1} \rangle+ \langle -\nabla h(\xi_t) - \lambda_t, \xi_{t+1} - \xi_{t} \rangle
        \\[0.2cm]
        & \leq & \displaystyle
        \frac{1}{2\eta}\norm{\bar\xi^\star - \xi_t}^2 -\frac{1}{2\eta}\norm{\bar\xi^\star - \xi_{t+1}}^2 
        + \langle -\nabla h(\xi_t) - \lambda_t, \xi_{t+1} - \xi_{t} \rangle
    \end{array}
    \]
    where the last inequality is due to the second inequality in Lemma~\ref{lem:optimality NPG} with $\xi = \bar\xi^\star$. Therefore,~\eqref{eq:PDL PG} reduces to
    \[
    \begin{array}{rcl}
         &&
         \!\!\!\!  \!\!\!\!  \!\! \displaystyle
          V_r^\star(\rho) -h(\bar\xi^\star) - (V_r^{\pi_t}(\rho) -h(\xi_t)) +  \lambda_t^\top(V_g^\star(\rho) - \bar\xi^\star - (V_g^{\pi_t}(\rho) -\xi_t))
          \\[0.2cm]
           & \leq & \displaystyle
           \frac{1}{2\eta(1-\gamma)} \mathbb{E}_{s'\,\sim\,d_\rho^\star} \left[ \norm{\bar\pi^\star(\cdot\,\vert\,s) - \pi_t(\cdot\,\vert\,s)}^2 -\norm{\bar\pi^\star(\cdot\,\vert\,s) - \pi_{t+1}(\cdot\,\vert\,s)}^2 \right] 
           \\[0.2cm]
           &  & \displaystyle
           +\,\frac{1}{1-\gamma} \left( (V_r^{\pi_{t+1}}(d_\rho^\star)-V_r^{\pi_{t}}(d_\rho^\star))
        + \lambda_t^\top (V_g^{\pi_{t+1}}(d_\rho^\star)-V_g^{\pi_{t}}(d_\rho^\star)) \right)
          \\[0.2cm]
           &  & \displaystyle
           +\,\frac{1}{2\eta}\norm{\bar\xi^\star - \xi_t}^2 -\frac{1}{2\eta}\norm{\bar\xi^\star - \xi_{t+1}}^2 
        + \langle -\nabla h(\xi_t) - \lambda_t, \xi_{t+1} - \xi_{t} \rangle.
    \end{array}
    \]
    Summing up the inequality above from $t=0$ to $t=T-1$ and dividing it by $T$ yield, 
    \begin{equation}\label{eq:PDL PG final}
    \begin{array}{rcl}
         &&
         \!\!\!\!  \!\!\!\!  \!\! \displaystyle
         \frac{1}{T}\sum_{t\,=\,0}^{T-1}\left(
          V_r^\star(\rho) -h(\bar\xi^\star) - (V_r^{\pi_t}(\rho) -h(\xi_t)) \right) + \frac{1}{T}\sum_{t\,=\,0}^{T-1}  \lambda_t^\top(V_g^\star(\rho) - \bar\xi^\star - (V_g^{\pi_t}(\rho) -\xi_t))
          \\[0.2cm]
           & \leq & \displaystyle
           \frac{1}{2\eta(1-\gamma) T} \mathbb{E}_{s'\,\sim\,d_\rho^\star} \left[ \norm{\bar\pi^\star(\cdot\,\vert\,s) - \pi_0(\cdot\,\vert\,s)}^2 -\norm{\bar\pi^\star(\cdot\,\vert\,s) - \pi_{T}(\cdot\,\vert\,s)}^2 \right] 
           \\[0.2cm]
           &  & \displaystyle
           +\,\frac{1}{(1-\gamma) T} (V_r^{\pi_{T}}(d_\rho^\star)-V_r^{\pi_{0}}(d_\rho^\star))
        + 
        \frac{1}{1-\gamma}\frac{1}{T}\sum_{t\,=\,0}^{T-1}\lambda_t^\top (V_g^{\pi_{t+1}}(d_\rho^\star)-V_g^{\pi_{t}}(d_\rho^\star))
          \\[0.2cm]
           &  & \displaystyle
           +\,\frac{1}{2\eta T}\norm{\bar\xi^\star - \xi_0}^2 -\frac{1}{2\eta T}\norm{\bar\xi^\star - \xi_{T}}^2 
        + \frac{1}{T}\sum_{t\,=\,0}^{T-1}\langle -\nabla h(\xi_t) - \lambda_t, \xi_{t+1} - \xi_{t} \rangle.
    \end{array}
    \end{equation}

    We notice that $\lambda_0 = 0$, $\lambda_T = \sum_{t\,=\,0}^{T-1}(\lambda_{t+1} - \lambda_t)$, and $|V_{g_i}^{\pi_t}(\rho)-\xi_{i,t}| \leq \frac{2}{1-\gamma}$ for $i=1,\ldots,m$. From the $\lambda$-update in~\eqref{eq:policy_gradient}, we have $|\lambda_{i,t}- \lambda_{i,t+1}|\leq \frac{2\eta}{1-\gamma}$ and $|\lambda_{i,T}|\leq \frac{2\eta T}{1-\gamma}$. Thus,
    \[
    \begin{array}{rcl}
         &&
         \!\!\!\!  \!\!\!\!  \!\! \displaystyle
    \sum_{t\,=\,0}^{T-1}\lambda_t^\top (V_g^{\pi_{t+1}}(d_\rho^\star)-V_g^{\pi_{t}}(d_\rho^\star))
    \\[0.2cm]
    & = & \displaystyle
    \sum_{t\,=\,0}^{T-1}(\lambda_{t+1}^\top  V_g^{\pi_{t+1}}(d_\rho^\star)-\lambda_t^\top V_g^{\pi_{t}}(d_\rho^\star)) + \sum_{t\,=\,0}^{T-1}(\lambda_{t}-\lambda_{t+1})^\top  V_g^{\pi_{t+1}}(d_\rho^\star)
    \\[0.2cm]
    & \leq & \displaystyle
    \lambda_{T}^\top  V_g^{\pi_{T}}(d_\rho^\star) + \sum_{t\,=\,0}^{T-1} \sum_{i\,=\, 1}^m|\lambda_{i,t}-\lambda_{i, t+1}| V_{g_i}^{\pi_{t+1}}(d_\rho^\star)
    \\[0.2cm]
    & \leq & \displaystyle
    \frac{4 \eta m T}{(1-\gamma)^2}.
    \end{array}
    \]
    Meanwhile, from the $\xi$-update in~\eqref{eq:policy_gradient}, we have $|\xi_{i,t+1} - \xi_{i,t}| \leq \eta | - \nabla h(\xi_t) - \lambda_t|_i  \leq \frac{\eta (L_h+1)}{1-\gamma}$. Thus,
    \[
        \displaystyle
        \sum_{t\,=\,0}^{T-1}\langle -\nabla h(\xi_t) - \lambda_t, \xi_{t+1} - \xi_{t} \rangle
        \leq  \frac{\eta  (L_h+1)^2 mT}{(1-\gamma)^2}.
    \]
    Finally, we combine these inequalities above with~\eqref{eq:PDL PG final} to get our desired inequality.
\end{proof}

\begin{proof}
    To show the first inequality, we notice that $\lambda_0 = 0$, $\lambda_{i,T}^2 = \sum_{t\,=\,0}^{T-1}(\lambda_{i,t+1}^2 - \lambda_{i,t}^2)$, and $|V_{g_i}^{\pi_t}(\rho)-\xi_{i,t}| \leq \frac{2}{1-\gamma}$ for $i=1,\ldots,m$. From the $\lambda$-update in~\eqref{eq:policy_gradient}, we have $|\lambda_{i,t}- \lambda_{i,t+1}|\leq \frac{2\eta}{1-\gamma}$ and $|\lambda_{i,T}|\leq \frac{2\eta T}{1-\gamma}$. Thus,
    \[
    \begin{array}{rcl}
        \displaystyle
    \lambda_{i,T}^2
    & = & \displaystyle
    \sum_{t\,=\,0}^{T-1} 
    (\lambda_{i,t+1})^2 - (\lambda_{i,t})^2
    \\[0.2cm]
    & = & \displaystyle
    \sum_{t\,=\,0}^{T-1} 
     -2 \eta \lambda_{i,t} ( V_{g_i}^{\pi_t}(\rho)-\xi_{i,t}) +  \eta^2 ( V_{g_i}^{\pi_t}(\rho)-\xi_{i,t})^2
    \\[0.2cm]
    & \leq & \displaystyle
    2 \eta\sum_{t\,=\,0}^{T-1} 
     \lambda_{i,t}(( V_{g_i}^{\star}(\rho)-\bar\xi_{i}^\star)-( V_{g_i}^{\pi_t}(\rho)-\xi_{i,t}))  +  \eta^2 ( V_{g_i}^{\pi_t}(\rho)-\xi_{i,t})^2
     \\[0.2cm]
    & \leq & \displaystyle
    2 \eta\sum_{t\,=\,0}^{T-1} 
     \lambda_{i,t}(( V_{g_i}^{\star}(\rho)-\bar\xi_{i}^\star)-( V_{g_i}^{\pi_t}(\rho)-\xi_{i,t}))  +  \frac{4 \eta^2 T}{(1-\gamma)^2}
    \end{array}
    \]
    where the first inequality is due to the feasibility $V_{g_i}^\star (\rho) \geq \bar\xi_i^\star$ for $i =1,\ldots,m$. Thus,
    \[
    \displaystyle
    -\frac{1}{T}\sum_{t\,=\,0}^{T-1} 
     \lambda_{t}^\top (( V_{g}^{\star}(\rho)-\bar\xi^\star)-( V_{g}^{\pi_t}(\rho)-\xi_{t}))  
     \; \leq \; \frac{2 \eta m }{(1-\gamma)^2}
    \]
    which can be added to the inequality in Lemma~\ref{lem: average bound PG} from both sides,
    \[
    \begin{array}{rcl}
         &&
         \!\!\!\!  \!\!\!\!  \!\! \displaystyle
         \frac{1}{T} \sum_{t \, = \,0}^{T-1} (V_r^\star(\rho) -h(\bar\xi^\star) - (V_r^{\pi_t}(\rho) -h(\xi_t))) 
         \\[0.2cm]
         & \leq  & \displaystyle
         \frac{1}{(1-\gamma)^2 T}
         + \frac{1}{\eta (1-\gamma)T} + \frac{4\eta m}{(1-\gamma)^2} + \frac{\eta(L_h+1)^2 m}{(1-\gamma)^2} + \frac{2\eta m}{(1-\gamma)^2}.
    \end{array}
    \]
    Hence, we obtain the first inequality by taking $\eta = \frac{1}{\sqrt{T}}$. 

    We next prove the second inequality. From the $\lambda$-update in~\eqref{eq:policy_gradient}, for any $\lambda \in \Lambda \DefinedAs \{\lambda\in\mathbb{R}^m\,\vert\, 0 \leq \lambda_i\leq C_h, i=1,\ldots,m\}$,
    \[
    (\lambda_{i,t+1}-\lambda_i)^2
    \;\leq\;
    (\lambda_{i,t}-\lambda_i)^2 - 2\eta (V_{g_i}^{\pi_t}(\rho)-\xi_{i,t})(\lambda_{i,t}-\lambda_i)+ \eta^2(V_{g_i}^{\pi_t}(\rho)-\xi_{i,t})^2
    \]
    which combines with $|V_{g_i}^{\pi_t}(\rho)-\xi_{i,t}|\leq \frac{2}{1-\gamma}$ to give us, 
    \begin{equation}\label{eq:dual add}
    \begin{array}{rcl}
         &&
         \!\!\!\!  \!\!\!\!  \!\! \displaystyle
     \frac{1}{T}\sum_{t\,=\,0}^{T-1}
    \sum_{i\,=\,1}^m(V_{g_i}^{\pi_t}(\rho)-\xi_{i,t})(\lambda_{i,t}-\lambda_i)
    \\[0.2cm]
    & \leq & \displaystyle
    \frac{1}{2\eta T}\sum_{t\,=\,0}^{T-1} \sum_{i\,=\,1}^m \left(
    (\lambda_{i,t}-\lambda_i)^2 -(\lambda_{i,t+1}-\lambda_i)^2
    \right)
    + \frac{2 \eta m}{\eta (1-\gamma)^2}
    \\[0.2cm]
    & \leq & \displaystyle
    \frac{1}{2\eta T} \sum_{i\,=\,1}^m 
    \lambda_i^2
    + \frac{2 \eta m}{ (1-\gamma)^2}.
    \end{array}
    \end{equation}
    Notice that $V_g^\star(\rho) \geq \bar\xi^\star$. If we add~\eqref{eq:dual add} to the inequality in Lemma~\ref{lem: average bound PG}, then
     \[
    \begin{array}{rcl}
         &&
         \!\!\!\!  \!\!\!\!  \!\! \displaystyle
         \frac{1}{T} \sum_{t \, = \,0}^{T-1} (V_r^\star(\rho) -h(\bar\xi^\star) - (V_r^{\pi_t}(\rho) -h(\xi_t))) +  \frac{1}{T} \sum_{t \, = \,0}^{T-1} \lambda^\top ( - (V_g^{\pi_t}(\rho) -\xi_t))
         \\[0.2cm]
         & \leq  & \displaystyle
         \frac{1}{(1-\gamma)^2 T}
         + \frac{1}{\eta (1-\gamma)T} + \frac{4\eta m}{(1-\gamma)^2} + \frac{\eta(L_h+1)^2 m}{(1-\gamma)^2} + \frac{1}{2\eta T} \norm{\lambda}^2 + \frac{2\eta m}{ (1-\gamma)^2}
    \end{array}
    \]
    where the RHS can be upper bounded by, if we take $\eta = \frac{1}{\sqrt{T}}$,
    \[
    \displaystyle
    \frac{2+(6+(L_h+1)^2) m}{(1-\gamma)^2\sqrt{T}}
    +
    \frac{1}{2\sqrt{T}}\norm{\lambda}^2.
    \]
    We next apply a fact in constrained convex optimization.
    We notice that $V_r^{\pi_t}(\rho)$ and $V_g^{\pi_t}$ are linear in the occupancy measure induced by $\pi_t$. By the convexity of the occupancy measure set, $\frac{1}{T}\sum_{t\,=\,0}^{T-1}V_r^{\pi_t}(\rho)$ and  $\frac{1}{T}\sum_{t\,=\,0}^{T-1}V_g^{\pi_t}(\rho)$ are linear in an occupancy measure induced by some policy~$\pi'$ and we denote them as $V_r^{\pi'}(\rho)$ and $V_g^{\pi'}(\rho)$. 
    Since $h$ is convex, there exists $\xi'$ such that $\frac{1}{T}\sum_{t\,=\,0}^{T-1} h(\xi_t) \geq h(\xi')$.
    If we choose $\lambda_i = 2C_h$ if $V_{g_i}^{\pi_t}(\rho)\leq \xi_i$ and $\lambda_i = 0$ otherwise, then 
    \[
    \begin{array}{rcl}
         &&
         \!\!\!\!  \!\!\!\!  \!\! \displaystyle
         V_r^\star(\rho) -h(\bar\xi^\star) - (V_r^{\pi'}(\rho) -h(\xi')) +  2C_h \sum_{i\,=\,1}^m \left[ \frac{1}{T}\sum_{t \, = \,0}^{T-1} \xi_{i,t}- V_{g_i}^{\pi'}(\rho) \right]_+
         \\[0.2cm]
         & \leq  & \displaystyle
         \displaystyle
    \frac{2+(6+(L_h+1)^2) m}{(1-\gamma)^2\sqrt{T}}
    +
    \frac{ 2 m C_h^2}{\sqrt{T}}.
    \end{array}
    \]
     
    Due to $2C_h \geq 2\lambda^\star$ and the strong duality,
    application of Lemma~\ref{lem:convex_violation_bound} leads to  
    \[
    \displaystyle
     \sum_{i\,=\,1}^m \left[ \frac{1}{T} \sum_{t \, = \,0}^{T-1}\xi_{i,t}- V_{g_i}^{\pi'}(\rho) \right]_+
    \; \leq\; 
    \frac{2+(6+(L_h+1)^2) m}{(1-\gamma)^2 C_h \sqrt{T}}
    +
    \frac{2m C_h}{\sqrt{T}}
    \]
    which shows the second inequality by replacing $V_{g_i}^{\pi'}(\rho)$ by $\frac{1}{T}\sum_{t\,=\,0}^{T-1}V_{g_i}^{\pi_t}(\rho)$.
\end{proof}

\subsection{Resilient Optimistic Policy Gradient Primal-Dual Method}\label{app: resilient OPG}

\begin{algorithm*}[h]
	\caption{Resilient optimistic policy gradient primal-dual (ResOPG-PD) method }
	\label{alg: resilient OPG}
	\begin{algorithmic}[1]
		\STATE
		\textbf{Parameters:} $\eta>0$. \\
		\textbf{Initialization}: Let $\pi_0 (a\,\vert\,s) = \hat\pi_0 (a\,\vert\,s) ={1}/{A}$ for $s\in\calS$, $a\in\calA$, and $\xi_0 = \hat\xi_0 = 0$, and $\lambda_0 = \hat\lambda_0 = 0$.
		\FOR{step $t=1,\ldots,T$} 
		\STATE Primal-dual update 
                    \begin{subequations}
                        \label{eq:optimistic_policy_gradient}
                \begin{equation}\label{eq:optimistic_policy_gradient_primal}
    \begin{array}{rcl}
         \pi_t(\cdot\,\vert\,s)
         &  =  & 
         \displaystyle\argmax_{\pi(\cdot\,\vert\,s)\,\in\,\Pi}\;
         \left\{
                \sum_{a} \pi(a\,\vert\,s) Q_{r+\lambda_{t-1}^\top g}^{\pi_{t-1}}(s,a)
                - 
                \frac{1}{2\eta} \norm{\pi(\cdot\,\vert\,s) - \hat\pi_t(\cdot\,\vert\,s)}^2
         \right\}
         \\[0.4cm]
         \xi_t
         &  =  &
         \displaystyle
         \argmax_{\xi\,\in\,\Xi} \;
         \left\{
                \xi^\top\left( -\nabla h(\xi_{t-1}) - \lambda_{t-1}\right) -
                \frac{1}{2\eta} \norm{\xi - \hat\xi_t}^2
         \right\}
         \\[0.4cm]
         \hat\pi_{t+1}(\cdot\,\vert\,s)
         &  =  & 
         \displaystyle\argmax_{\pi(\cdot\,\vert\,s)\,\in\,\Pi}\;
         \left\{
                \sum_{a} \pi(a\,\vert\,s) Q_{r+\lambda_{t}^\top g}^{\pi_{t}}(s,a)
                - 
                \frac{1}{2\eta} \norm{\pi(\cdot\,\vert\,s) - \hat\pi_t(\cdot\,\vert\,s)}^2
         \right\}
         \\[0.4cm]
         \hat\xi_{t+1}
         &  =  &
         \displaystyle
         \argmax_{\xi\,\in\,\Xi} \;
         \left\{
                \xi^\top\left( -\nabla h(\xi_{t}) - \lambda_{t}\right) -
                \frac{1}{2\eta} \norm{\xi - \hat\xi_t}^2
         \right\}
    \end{array}
\end{equation}
\begin{equation}\label{eq:optimistic_policy_gradient_dual}
    \begin{array}{rcl}
         \lambda_{t}
         &  =  &
         \displaystyle
         \argmax_{\lambda\,\in\,\Lambda} \;
         \left\{
                \lambda^\top\left( V_{g}^{\pi_{t-1}}(\rho)- \xi_{t-1}\right) +
                \frac{1}{2\eta} \norm{\lambda - \hat\lambda_t}^2
         \right\}
         \\[0.4cm]
         \hat\lambda_{t+1}
         &  =  &
         \displaystyle
         \argmax_{\lambda\,\in\,\Lambda} \;
         \left\{
                \lambda^\top\left( V_{g}^{\pi_t}(\rho)- \xi_{t}\right) +
                \frac{1}{2\eta} \norm{\lambda - \hat\lambda_t}^2
         \right\} 
    \end{array}
\end{equation}
\end{subequations}
		\ENDFOR
	\end{algorithmic}
\end{algorithm*}

\subsection{Proof of Theorem~\ref{thm: last-iterate convergence}}\label{app: last-iterate convergence}

We define
\[
    \begin{array}{rcl}
         \Theta_{t+1} & \DefinedAs &
         \displaystyle
         \frac{1}{2(1-\gamma)} \sum_{s}d_\rho^{\pi^\star}(s) \norm{\mathcal{P}_{\Pi^\star}( \hat\pi_{t+1}(\cdot\,\vert\,s)) - \hat\pi_{t+1}(\cdot\,\vert\,s)}^2
         +
         \frac{1}{2}\norm{\mathcal{P}_{\Xi^\star}(\hat\xi_{t+1})-\hat\xi_{t+1}}^2 
         +
         \frac{1}{2}  
            \norm{\mathcal{P}_{\Lambda^\star}(\hat\lambda_{t+1})-\hat\lambda_{t+1}}^2 
        \\[0.2cm]
         & &  \!\!\!\!  \!\!\!\!  \!\!\displaystyle
         +\,\frac{1}{4(1-\gamma)}  \sum_{s}d_\rho^{\pi^\star}(s) 
            \norm{\hat\pi_{t+1}(\cdot\,\vert\,s) - \pi_t(\cdot\,\vert\,s)}^2 + \frac{1}{4}
            \norm{\hat\xi_{t+1}-\xi_t}^2 +  \frac{1}{4}\norm{\hat\lambda_{t+1}-\lambda_t}^2
    \end{array}
\]
\[
    \begin{array}{rcl}
    \zeta_{t} & \DefinedAs &
     \displaystyle
         \,\frac{1}{2(1-\gamma)} \sum_{s}d_\rho^{\pi^\star}(s) 
             \norm{\hat\pi_{t+1}(\cdot\,\vert\,s) - \pi_t(\cdot\,\vert\,s)}^2
            +
        \frac{1}{2}
            \norm{\hat\xi_{t+1}-\xi_t}^2 +
         \frac{1}{2}  
            \norm{\hat\lambda_{t+1}-\lambda_t}^2
         \\[0.2cm]
         & &  \displaystyle
         +\,\frac{1}{2(1-\gamma)} \sum_{s}d_\rho^{\pi^\star}(s) 
             \norm{\pi_{t}(\cdot\,\vert\,s) - \hat\pi_t(\cdot\,\vert\,s)}^2
            +
        \frac{1}{2}
            \norm{\xi_{t}-\hat\xi_t}^2 +
         \frac{1}{2}  
            \norm{\lambda_{t}
            -\hat\lambda_t}^2
    \end{array}
\]
and 
\[
    \iota 
    \; \DefinedAs \;
    \max\left(
        \frac{|A|}{(1-\gamma)^2}+1, L_h^2+1, 2 |A|^2 \frac{\gamma(1+m (C_h)^2) \kappa_\rho}{(1-\gamma)^4\rho_{\min}} + \frac{|A| m\kappa_\rho^2}{(1-\gamma)^5}
    \right)
    \]
\[
    \eta_{\max}
    \;\DefinedAs \; 
            \min\left( \frac{1}{4\sqrt{|A|}}, \frac{1}{2(L_h+1)}, \frac{1}{5 \sqrt{m|A|}\kappa_\rho}, \frac{\rho_{\min}}{4 \gamma\sqrt{m}C_h |A|}, \frac{1}{2\sqrt{2\iota}}, \frac{4 \max(\frac{\kappa_\rho}{1-\gamma}, 1)}{\sqrt{ \Theta_1 C_{\rho,\gamma,\sigma}(1-\gamma)}}
    \right).
\]

\begin{lemma}\label{lem:non-increasing}
    In Algorithm~\ref{alg: resilient OPG}, for~\eqref{eq:optimistic_policy_gradient} with $\eta \leq \frac{1}{2\sqrt{2\iota}}$, 
    \[
    \Theta_{t+1}
    \; \leq \;
    \Theta_t - \frac{1}{2}\zeta_t.
    \]
\end{lemma}
\begin{proof}
    For any $(\pi^\star,\xi^\star,\lambda^\star) \in \Pi^\star\times\Xi^\star\times\Lambda^\star$, 
    \begin{equation}\label{eq:decomposition}
    \begin{array}{rcl}
         &  & 
         \!\!\!\!  \!\!\!\!  \!\!
         V_{r+\lambda_t^\top g}^{\pi^\star}(\rho) - h(\xi^\star) - \lambda_t^\top \xi^\star
        - \left(
        V_{r+(\lambda^\star)^\top g}^{\pi_t}(\rho) - h(\xi_t) - (\lambda^\star)^\top \xi_t
        \right)
         \\[0.2cm]
         & = & \underbrace{\left( V_{r+\lambda_t^\top g}^{\pi^\star}(\rho) - V_{r+\lambda_t^\top g}^{\pi_t}(\rho) \right)}_{\normalfont (a)} 
         \,+\, \underbrace{\left( -h(\xi^\star)-\lambda_t^\top\xi^\star + h(\xi_t) +\lambda_t^\top \xi_t \right)}_{\normalfont (b)}
         \\[0.2cm]
         &  &
         +\,
         \underbrace{\left(
             V_{r+\lambda_t^\top g}^{\pi_t}(\rho) -
             V_{r+(\lambda^\star)^\top g}^{\pi_t}(\rho)
             -
             \lambda_t^\top\xi_t + (\lambda^\star)^\top \xi_t
         \right)}_{\normalfont (c)}.
    \end{array}
    \end{equation}
    We next analyze three terms $(a)$, $(b)$, and $(c)$, separately.

    For $(a)$, we have
    \[
    \begin{array}{rcl}
         &  & \!\!\!\! \!\!\!\! \!\!
         V_{r+\lambda_t^\top g}^{\pi^\star}(\rho) - V_{r+\lambda_t^\top g}^{\pi_t}(\rho)
         \\[0.2cm]
         & = & \displaystyle \frac{1}{1-\gamma}
         \sum_{s,a} d_\rho^{\pi^\star}(s)
         (\pi^\star(a\,\vert\,s) - \pi_t(a\,\vert\,s)) Q_{r+\lambda_t^\top g}^{\pi_t}(s,a)
         \\[0.2cm]
         & = & \displaystyle\frac{1}{1-\gamma}
         \sum_{s,a} d_\rho^{\pi^\star}(s)
         (\pi^\star(a\,\vert\,s) - \hat\pi_{t+1}(a\,\vert\,s)) Q_{r+\lambda_t^\top g}^{\pi_t}(s,a)
         +\frac{1}{1-\gamma}\sum_{s,a} d_\rho^{\pi^\star}(s)
         (\hat\pi_{t+1}(a\,\vert\,s) - \pi_{t}(a\,\vert\,s)) Q_{r+\lambda_{t-1}^\top g}^{\pi_{t-1}}(s,a)
         \\[0.2cm]
         & & \displaystyle +\frac{1}{1-\gamma}\sum_{s,a} d_\rho^{\pi^\star}(s)
         (\hat\pi_{t+1}(a\,\vert\,s) - \pi_{t}(a\,\vert\,s)) \big(Q_{r+\lambda_{t}^\top g}^{\pi_{t}}(s,a)- Q_{r+\lambda_{t-1}^\top g}^{\pi_{t-1}}(s,a)\big)
         \\[0.2cm]
         & \leq & \displaystyle
         \frac{1}{2\eta(1-\gamma)} \sum_{s}d_\rho^{\pi^\star}(s)
         \left(
            \norm{\pi^\star(\cdot\,\vert\,s) - \hat\pi_t(\cdot\,\vert\,s)}^2
            -
        \norm{\pi^\star(\cdot\,\vert\,s) - \hat\pi_{t+1}(\cdot\,\vert\,s)}^2
        - \norm{\hat\pi_{t+1}(\cdot\,\vert\,s) - \hat\pi_t(\cdot\,\vert\,s)}^2
         \right)
         \\[0.2cm]
         & &  \displaystyle
         +\,\frac{1}{2\eta(1-\gamma)} \sum_{s}d_\rho^{\pi^\star}(s) 
         \left(
            \norm{\hat\pi_{t+1}(\cdot\,\vert\,s) - \hat\pi_t(\cdot\,\vert\,s)}^2
            - \norm{\hat\pi_{t+1}(\cdot\,\vert\,s) - \pi_t(\cdot\,\vert\,s)}^2
            -
            \norm{\pi_{t}(\cdot\,\vert\,s) - \hat\pi_t(\cdot\,\vert\,s)}^2
         \right)
         \\[0.2cm]
         & & \displaystyle +\frac{1}{1-\gamma}\sum_{s,a} d_\rho^{\pi^\star}(s)
         (\hat\pi_{t+1}(a\,\vert\,s) - \pi_{t}(a\,\vert\,s)) \big(Q_{r+\lambda_{t}^\top g}^{\pi_{t}}(s,a)- Q_{r+\lambda_{t-1}^\top g}^{\pi_{t-1}}(s,a)\big)
         \\[0.2cm]
         & \leq & \displaystyle
         \frac{1}{2\eta(1-\gamma)} \sum_{s}d_\rho^{\pi^\star}(s)
         \left(
            \norm{\pi^\star(\cdot\,\vert\,s) - \hat\pi_t(\cdot\,\vert\,s)}^2
            -
        \norm{\pi^\star(\cdot\,\vert\,s) - \hat\pi_{t+1}(\cdot\,\vert\,s)}^2
        - \norm{\hat\pi_{t+1}(\cdot\,\vert\,s) - \hat\pi_t(\cdot\,\vert\,s)}^2
         \right)
         \\[0.2cm]
         & &  \displaystyle
         +\,\frac{1}{2\eta(1-\gamma)} \sum_{s}d_\rho^{\pi^\star}(s) 
         \left(
            \norm{\hat\pi_{t+1}(\cdot\,\vert\,s) - \hat\pi_t(\cdot\,\vert\,s)}^2
            - \norm{\hat\pi_{t+1}(\cdot\,\vert\,s) - \pi_t(\cdot\,\vert\,s)}^2
            -
            \norm{\pi_{t}(\cdot\,\vert\,s) - \hat\pi_t(\cdot\,\vert\,s)}^2
         \right)
         \\[0.2cm]
         &  & +\,
         \displaystyle
         \frac{4\eta |A|}{(1-\gamma)^3} 
         \left( 
         \norm{\lambda_t -\hat\lambda_{t}}^2 + \norm{\hat\lambda_t -\lambda_{t-1}}^2
         \right)
         \\[0.2cm]
         & &
         \displaystyle
         +\,
         8\eta|A|^2 \frac{\gamma (1+m(C_h)^2) \kappa_\rho}{(1-\gamma)^5 \rho_{\min}} 
         \sum_s d_\rho^{\pi^\star}(s)\left(
         \norm{\pi_t(\cdot\,\vert\,s) - \hat\pi_{t}(\cdot\,\vert\,s)}^2
         +
        \norm{\hat\pi_t(\cdot\,\vert\,s) - \pi_{t-1}(\cdot\,\vert\,s)}^2
         \right)
    \end{array}
    \]
    where the first inequality is due to the optimality of $\hat\pi_{t+1}$ and $\pi_t$ that results from Lemma~\ref{lem:three-point}, and the second inequality is due to that 
    \[
    \begin{array}{rcl}
         &  & \!\!\!\!  \!\!\!\!  \!\!
         \displaystyle
         \sum_{s,a} d_\rho^{\pi^\star}(s)
         (\hat\pi_{t+1}(a\,\vert\,s) - \pi_{t}(a\,\vert\,s)) \big(Q_{r+\lambda_{t}^\top g}^{\pi_{t}}(s,a)- Q_{r+\lambda_{t-1}^\top g}^{\pi_{t-1}}(s,a)\big)  
         \\
         & \leq & \displaystyle
         \eta \sum_s d_\rho^{\pi^\star}(s) \norm{Q_{r+\lambda_{t}^\top g}^{\pi_{t}}(s,\cdot)- Q_{r+\lambda_{t-1}^\top g}^{\pi_{t-1}}(s,\cdot)}^2
         \\
         & \leq & \displaystyle
         2\eta \sum_s d_\rho^{\pi^\star}(s) \norm{ (\lambda_t  - \lambda_{t-1})^\top Q_{ g}^{\pi_{t}}(s,\cdot)}^2
         +
         2\eta \sum_s d_\rho^{\pi^\star}(s) \norm{ Q_{r+\lambda_{t-1}^\top g}^{\pi_{t}}(s,\cdot) - Q_{r+\lambda_{t-1}^\top g}^{\pi_{t-1}}(s,\cdot)}^2
         \\[0.2cm]
         & \leq &
         \displaystyle
         \frac{2\eta |A|}{(1-\gamma)^2} 
         \norm{\lambda_t -\lambda_{t-1}}^2
         +
         2\eta|A| \sum_s d_\rho^{\pi^\star}(s) \norm{ Q_{r+\lambda_{t-1}^\top g}^{\pi_{t}}(s,\cdot) - Q_{r+\lambda_{t-1}^\top g}^{\pi_{t-1}}(s,\cdot)}_\infty^2
        \\[0.2cm]
         & \leq &
         \displaystyle
         \frac{2\eta |A|}{(1-\gamma)^2} 
         \norm{\lambda_t -\lambda_{t-1}}^2
         +
         4\eta|A| \sum_s d_\rho^{\pi^\star}(s) \norm{ Q_{r}^{\pi_{t}}(s,\cdot) - Q_{r}^{\pi_{t-1}}(s,\cdot)}_\infty^2
         \\ [0.2cm]
         &  & \displaystyle
         +\,
         4\eta|A| \sum_s d_\rho^{\pi^\star}(s) 
         \norm{\lambda_{t-1}}^2
         m
         \max_i
         \norm{ Q_{g_i}^{\pi_{t}}(s,\cdot) - Q_{g_i}^{\pi_{t-1}}(s,\cdot)}_\infty^2
         \\[0.2cm]
         & \leq &
         \displaystyle
         \frac{2\eta |A|}{(1-\gamma)^2} 
         \norm{\lambda_t -\lambda_{t-1}}^2
         +
         4\eta|A| \frac{\gamma (1+m(C_h)^2)}{(1-\gamma)^2} 
         \max_s \norm{\pi_t(\cdot\,\vert\,s) - \pi_{t-1}(\cdot\,\vert\,s)}_1^2
         \\[0.2cm]
         & \leq &
         \displaystyle
         \frac{4\eta |A|}{(1-\gamma)^2} 
         \left( 
         \norm{\lambda_t -\hat\lambda_{t}}^2 + \norm{\hat\lambda_t -\lambda_{t-1}}^2
         \right)
         \\[0.2cm]
         & &
         \displaystyle
         +\,
         8\eta|A|^2 \frac{\gamma (1+m(C_h)^2) \kappa_\rho}{(1-\gamma)^4 \rho_{\min}} 
         \sum_s d_\rho^{\pi^\star}(s)\left(
         \norm{\pi_t(\cdot\,\vert\,s) - \hat\pi_{t}(\cdot\,\vert\,s)}^2
         +
        \norm{\hat\pi_t(\cdot\,\vert\,s) - \pi_{t-1}(\cdot\,\vert\,s)}^2
         \right)
    \end{array}
    \]
    where the first inequality is due to Lemma~\ref{lem:non-expansive}, we use the inequality $\norm{x+y}^2 \leq 2 \norm{x}^2 + 2\norm{y}^2$ in the second inequality, the third inequality is due to $\norm{x}\leq \sqrt{d}\norm{x}_\infty$ for any $x\in\mathbb{R}^d$, and the fourth inequality is due to the inequalities $\norm{x+y}^2 \leq 2 \norm{x}^2 + 2\norm{y}^2$, $\norm{x}\leq \sqrt{m}\norm{x}_\infty$ for any $x\in\mathbb{R}^m$, and
    $\norm{\sum_i x_i}_\infty\leq\sum_i\norm{x_i}_\infty$, application of Lemma~\ref{lem:policy-value-difference} leads to the fifth inequality together with $\norm{\lambda_{t-1}}\leq C_h$, and the last inequality is due to  $\norm{x}_1\leq \sqrt{d}\norm{x}$ for any $x\in\mathbb{R}^d$, $\norm{x+y}^2\leq 2 \norm{x}^2 + 2 \norm{y}^2$, and the property of $\kappa_\rho$,
    \[
    \frac{\kappa_\rho}{1-\gamma} d_\rho^{\pi^\star}(s)
    \; \geq \;
    d_\rho^\pi(s) \;\geq \; (1-\gamma)\rho_{\min}.
    \]

    For (b), we have
    \[
    \begin{array}{rcl}
         &  & \!\!\!\! \!\!\!\! \!\!
         - h(\xi^\star) - \lambda_t^\top \xi^\star + h(\xi_t) + \lambda_t^\top \xi_t
        \\[0.2cm]
        & \leq & \displaystyle
        (-\nabla h(\xi_t)-\lambda_t)^\top (\xi^\star - \xi_t)
        \\[0.2cm]
        & = & \displaystyle
        (-\nabla h(\xi_t)-\lambda_t)^\top (\xi^\star - \hat\xi_{t+1})
        + (-\nabla h(\xi_{t-1})-\lambda_{t-1})^\top (\hat\xi_{t+1}- \xi_t)
        \\[0.2cm]
        &  & \displaystyle
        + \, (-\nabla h(\xi_{t})-\lambda_{t}  + \nabla h(\xi_{t-1})+\lambda_{t-1})^\top (\hat\xi_{t+1}- \xi_t)
        \\[0.2cm]
        & \leq & \displaystyle
        \frac{1}{2\eta}\left(
            \norm{\xi^\star-\hat\xi_t}^2 - \norm{\xi^\star-\hat\xi_{t+1}}^2 -
            \norm{\hat\xi_{t+1} -\hat\xi_t}^2
        \right)
        + \frac{1}{2\eta}\left(
            \norm{\hat\xi_{t+1}-\hat\xi_t}^2-
            \norm{\hat\xi_{t+1}-\xi_t}^2-
            \norm{\xi_t-\hat\xi_t}^2
        \right)
        \\[0.2cm]
        &  & \displaystyle
        + \, (-\nabla h(\xi_{t})-\lambda_{t}  + \nabla h(\xi_{t-1})+\lambda_{t-1})^\top (\hat\xi_{t+1}- \xi_t)
        \\[0.2cm]
        & \leq & \displaystyle
        \frac{1}{2\eta}\left(
            \norm{\xi^\star-\hat\xi_t}^2 - \norm{\xi^\star-\hat\xi_{t+1}}^2 -
            \norm{\hat\xi_{t+1} -\hat\xi_t}^2
        \right)
        + \frac{1}{2\eta}\left(
            \norm{\hat\xi_{t+1}-\hat\xi_t}^2-
            \norm{\hat\xi_{t+1}-\xi_t}^2-
            \norm{\xi_t-\hat\xi_t}^2
        \right)
        \\[0.2cm]
        &  & \displaystyle
        + \, 4\eta \left( \norm{\lambda_t - \hat\lambda_{t}}^2 
        + \norm{\hat\lambda_t - \lambda_{t-1}}^2\right) + 4 \eta L_h^2 \left(\norm{\xi_{t} - \hat\xi_{t}}^2
        +  \norm{\hat\xi_{t} - \xi_{t-1}}^2\right)
    \end{array}
    \]
    where the first inequality is due to the convexity of $h$: $h(\xi^\star)\geq h(\xi_t) + \langle \nabla h(\xi_t), \xi^\star - \xi_t\rangle$, the second inequality is due to the optimality of $\hat\xi_{t+1}$ and $\xi_t$ that results from Lemma~\ref{lem:three-point}, and the last inequality is due to that
    \[
    \begin{array}{rcl}
         &  & \!\!\!\! \!\!\!\! \!\!
            (-\nabla h(\xi_{t})-\lambda_{t}  + \nabla h(\xi_{t-1})+\lambda_{t-1})^\top (\hat\xi_{t+1}- \xi_t)
        \\[0.2cm]
        & \leq & \eta \norm{-\nabla h(\xi_{t})-\lambda_{t}  + \nabla h(\xi_{t-1})+\lambda_{t-1}}^2
        \\[0.2cm]
        & \leq & 2\eta \norm{\lambda_t - \lambda_{t-1}}^2 + 2 \eta  \norm{-\nabla h(\xi_{t})+ \nabla h(\xi_{t-1})}^2
        \\[0.2cm]
        & \leq & 4\eta \left( \norm{\lambda_t - \hat\lambda_{t}}^2 
        + \norm{\hat\lambda_t - \lambda_{t-1}}^2\right) + 4 \eta L_h^2 \left(\norm{\xi_{t} - \hat\xi_{t}}^2
        +  \norm{\hat\xi_{t} - \xi_{t-1}}^2\right)
    \end{array}
    \]
    where the first inequality is due to Lemma~\ref{lem:non-expansive}, the second inequality is due to $\norm{x+y}^2 \leq 2\norm{x}^2+2\norm{y}^2$, the third inequality is due to the Lipschitz continuous gradient of $h$: $\norm{\nabla h(\xi) - \nabla h(\xi')}\leq L_h \norm{\xi-\xi'}$ for any $\xi$, $\xi'\in\Xi$, and $\norm{x+y}^2\leq 2\norm{x}^2+2\norm{y}^2$.

    For (c), we have
    \[
    \begin{array}{rcl}
         &  & \!\!\!\! \!\!\!\! \!\!
         V_{r+\lambda_t^\top g}^{\pi_t}(\rho) - V_{r+(\lambda^\star)^\top g}^{\pi_t}(\rho) - \lambda_t^\top \xi_t + (\lambda^\star)^\top \xi_t
         \\[0.2cm]
         & = & \displaystyle
         (\lambda_t - \lambda^\star)^\top \left( V_g^{\pi_t}(\rho) - \xi_t\right)
         \\[0.2cm]
         & = & \displaystyle
         (\lambda_t - \hat\lambda_{t+1})^\top \left( V_g^{\pi_{t-1}}(\rho) - \xi_{t-1}\right)
         + (\lambda_t - \hat\lambda_{t+1})^\top \left(
         V_g^{\pi_t}(\rho) - \xi_t -  V_g^{\pi_{t-1}}(\rho) + \xi_{t-1}
         \right)
        \\[0.2cm]
         &  & \displaystyle
         +\,  (\hat\lambda_{t+1} - \lambda^\star)^\top \left( V_g^{\pi_t}(\rho) - \xi_t \right)
         \\[0.2cm]
         & \leq & \displaystyle
         \frac{1}{2\eta} \left(
            \norm{\hat\lambda_{t+1}-\hat\lambda_t}^2   - 
            \norm{\hat\lambda_{t+1}-\lambda_t}^2 - \norm{\lambda_{t}-\hat\lambda_t}^2
         \right) 
         \\[0.2cm]
         &  & \displaystyle + (\lambda_t - \hat\lambda_{t+1})^\top \left(
         V_g^{\pi_t}(\rho) - \xi_t -  V_g^{\pi_{t-1}}(\rho) + \xi_{t-1}
         \right)
         \\[0.2cm]
         &  & \displaystyle +
         \frac{1}{2\eta} \left(
            \norm{\lambda^\star - \hat\lambda_t}^2 - \norm{\lambda^\star - \hat\lambda_{t+1}}^2 - \norm{\hat\lambda_{t+1}-\hat\lambda_t}^2
         \right)
         \\[0.2cm]
         & \leq & \displaystyle
         \frac{1}{2\eta} \left(
            \norm{\hat\lambda_{t+1}-\hat\lambda_t}^2   - 
            \norm{\hat\lambda_{t+1}-\lambda_t}^2 - \norm{\lambda_{t}-\hat\lambda_t}^2
         \right) 
          +
         \frac{1}{2\eta} \left(
            \norm{\lambda^\star - \hat\lambda_t}^2 - \norm{\lambda^\star - \hat\lambda_{t+1}}^2 - \norm{\hat\lambda_{t+1}-\hat\lambda_t}^2
         \right)
         \\[0.2cm]
         &  & \displaystyle +\,
         \frac{4\eta |A| m \kappa_\rho^2}{(1-\gamma)^6}
          \sum_s
         d_\rho^{\pi^\star}(s) \left(\norm{ \pi_t(\cdot\,\vert\,s) - \hat\pi_{t}(\cdot\,\vert\,s)}^2
         + \norm{ \hat\pi_t(\cdot\,\vert\,s) - \pi_{t-1}(\cdot\,\vert\,s)}^2\right)
         \\[0.2cm]
         &  & \displaystyle
         +\, 4\eta \left(
         \norm{
         \xi_t -  \hat\xi_{t}
         }^2 + \norm{
         \hat\xi_t -  \xi_{t-1}
         }^2
         \right)
    \end{array}
    \]
    where the first inequality is due to the optimality of $\lambda_{t}$ and $\hat\lambda_{t+1}$, and the the second inequality is due to that
    \[
    \begin{array}{rcl}
         &  & \!\!\!\! \!\!\!\! \!\!
         (\lambda_t - \hat\lambda_{t+1})^\top \left(
         V_g^{\pi_t}(\rho) - \xi_t -  V_g^{\pi_{t-1}}(\rho) + \xi_{t-1}
         \right) 
         \\[0.2cm]
         & \leq & \eta \norm{
         V_g^{\pi_t}(\rho) - \xi_t -  V_g^{\pi_{t-1}}(\rho) + \xi_{t-1}
         }^2
         \\[0.2cm]
         & \leq & \displaystyle
         2\eta m \left( 
         \frac{\kappa_\rho}{(1-\gamma)^3} \sum_s d_\rho^{\pi^\star}(s) \norm{ \pi_t(\cdot\,\vert\,s) - \pi_{t-1}(\cdot\,\vert\,s)}_1
         \right)^2 + 2\eta \norm{
         \xi_t -  \xi_{t-1}
         }^2
         \\[0.2cm]
         & \leq & \displaystyle
         \frac{2\eta m \kappa_\rho^2}{(1-\gamma)^6}
          \sum_s \left( 
         \sqrt{d_\rho^{\pi^\star}(s)} \norm{ \pi_t(\cdot\,\vert\,s) - \pi_{t-1}(\cdot\,\vert\,s)}_1
         \right)^2 + 2\eta \norm{
         \xi_t -  \xi_{t-1}
         }^2
         \\[0.2cm]
         & \leq & \displaystyle
         \frac{2\eta |A| m \kappa_\rho^2}{(1-\gamma)^6}
          \sum_s
         d_\rho^{\pi^\star}(s) \norm{ \pi_t(\cdot\,\vert\,s) - \pi_{t-1}(\cdot\,\vert\,s)}^2
         + 2\eta \norm{
         \xi_t -  \xi_{t-1}
         }^2
         \\[0.2cm]
         & \leq & \displaystyle
         \frac{4\eta |A| m \kappa_\rho^2}{(1-\gamma)^6}
          \sum_s
         d_\rho^{\pi^\star}(s) \left(\norm{ \pi_t(\cdot\,\vert\,s) - \hat\pi_{t}(\cdot\,\vert\,s)}^2
         + \norm{ \hat\pi_t(\cdot\,\vert\,s) - \pi_{t-1}(\cdot\,\vert\,s)}^2\right)
         \\[0.2cm]
         &  & \displaystyle
         +\, 4\eta \left(
         \norm{
         \xi_t -  \hat\xi_{t}
         }^2 + \norm{
         \hat\xi_t -  \xi_{t-1}
         }^2
         \right)
    \end{array}
    \]
    where the first inequality is due to Lemma~\ref{lem:non-expansive}, the second inequality is due to $\norm{x+y}^2 \leq 2\norm{x}^2+2\norm{y}^2$, $\norm{x} \leq \norm{x}_1$, and Lemma~\ref{lem:policy-value-difference}, the third inequality is due to Cauchy-Schwarz inequality, and the last inequality is due to 
    $\norm{x}_1\leq \sqrt{d}\norm{x}$ for any $x\in\mathbb{R}^d$ and $\norm{x+y}^2\leq 2 \norm{x}^2+2\norm{y}^2$.

    By applying the above upper bounds on (a), (b), and (c) to~\eqref{eq:decomposition}, it is ready to have
    \[
    \begin{array}{rcl}
         &  & 
         \!\!\!\!  \!\!\!\!  \!\!
         V_{r+\lambda_t^\top g}^{\pi^\star}(\rho) - h(\xi^\star) - \lambda_t^\top \xi^\star
        - \left(
        V_{r+(\lambda^\star)^\top g}^{\pi_t}(\rho) - h(\xi_t) - (\lambda^\star)^\top \xi_t
        \right)
         \\[0.2cm]
         & \leq & \displaystyle
         \frac{1}{2\eta(1-\gamma)} \sum_{s}d_\rho^{\pi^\star}(s)
         \left(
            \norm{\pi^\star(\cdot\,\vert\,s) - \hat\pi_t(\cdot\,\vert\,s)}^2
            -
        \norm{\pi^\star(\cdot\,\vert\,s) - \hat\pi_{t+1}(\cdot\,\vert\,s)}^2
        - \norm{\hat\pi_{t+1}(\cdot\,\vert\,s) - \hat\pi_t(\cdot\,\vert\,s)}^2
         \right)
         \\[0.2cm]
         & &  \displaystyle
         +\,\frac{1}{2\eta(1-\gamma)} \sum_{s}d_\rho^{\pi^\star}(s) 
         \left(
            \norm{\hat\pi_{t+1}(\cdot\,\vert\,s) - \hat\pi_t(\cdot\,\vert\,s)}^2
            - \norm{\hat\pi_{t+1}(\cdot\,\vert\,s) - \pi_t(\cdot\,\vert\,s)}^2
            -
            \norm{\pi_{t}(\cdot\,\vert\,s) - \hat\pi_t(\cdot\,\vert\,s)}^2
         \right)
         \\[0.2cm]
         &  & +\,
         \displaystyle
         \frac{4\eta |A|}{(1-\gamma)^3} 
         \left( 
         \norm{\lambda_t -\hat\lambda_{t}}^2 + \norm{\hat\lambda_t -\lambda_{t-1}}^2
         \right)
         \\[0.2cm]
         & &
         \displaystyle
         +\,
         8\eta|A|^2 \frac{\gamma (1+m(C_h)^2) \kappa_\rho}{(1-\gamma)^5 \rho_{\min}} 
         \sum_s d_\rho^{\pi^\star}(s)\left(
         \norm{\pi_t(\cdot\,\vert\,s) - \hat\pi_{t}(\cdot\,\vert\,s)}^2
         +
        \norm{\hat\pi_t(\cdot\,\vert\,s) - \pi_{t-1}(\cdot\,\vert\,s)}^2
         \right)
         \\[0.2cm]
        &  & \displaystyle +\,
        \frac{1}{2\eta}\left(
            \norm{\xi^\star-\hat\xi_t}^2 - \norm{\xi^\star-\hat\xi_{t+1}}^2 -
            \norm{\hat\xi_{t+1} -\hat\xi_t}^2
        \right)
        + \frac{1}{2\eta}\left(
            \norm{\hat\xi_{t+1}-\hat\xi_t}^2-
            \norm{\hat\xi_{t+1}-\xi_t}^2-
            \norm{\xi_t-\hat\xi_t}^2
        \right)
        \\[0.2cm]
        &  & \displaystyle
        + \, 4\eta \left( \norm{\lambda_t - \hat\lambda_{t}}^2 
        + \norm{\hat\lambda_t - \lambda_{t-1}}^2\right) + 4 \eta L_h^2 \left(\norm{\xi_{t} - \hat\xi_{t}}^2
        +  \norm{\hat\xi_{t} - \xi_{t-1}}^2\right)
        \\[0.2cm]
         & & \displaystyle +\,
         \frac{1}{2\eta} \left(
            \norm{\hat\lambda_{t+1}-\hat\lambda_t}^2   - 
            \norm{\hat\lambda_{t+1}-\lambda_t}^2 - \norm{\lambda_{t}-\hat\lambda_t}^2
         \right) 
          +
         \frac{1}{2\eta} \left(
            \norm{\lambda^\star - \hat\lambda_t}^2 - \norm{\lambda^\star - \hat\lambda_{t+1}}^2 - \norm{\hat\lambda_{t+1}-\hat\lambda_t}^2
         \right)
         \\[0.2cm]
         &  & \displaystyle +\,
         \frac{4\eta |A| m \kappa_\rho^2}{(1-\gamma)^6}
          \sum_s
         d_\rho^{\pi^\star}(s) \left(\norm{ \pi_t(\cdot\,\vert\,s) - \hat\pi_{t}(\cdot\,\vert\,s)}^2
         + \norm{ \hat\pi_t(\cdot\,\vert\,s) - \pi_{t-1}(\cdot\,\vert\,s)}^2\right)
         \\[0.2cm]
         &  & \displaystyle
         +\, 4\eta \left(
         \norm{
         \xi_t -  \hat\xi_{t}
         }^2 + \norm{
         \hat\xi_t -  \xi_{t-1}
         }^2
         \right).
    \end{array}
    \]
    
    We notice that $V_{r+\lambda_t^\top g}^{\pi^\star}(\rho) - h(\xi^\star) - \lambda_t^\top \xi^\star
        - \left(
        V_{r+(\lambda^\star)^\top g}^{\pi_t}(\rho) - h(\xi_t) - (\lambda^\star)^\top \xi_t
        \right) \geq 0$. Thus,
        \[
    \begin{array}{rcl}
         &  & 
         \!\!\!\!  \!\!\!\!  \!\!
         \displaystyle
         \frac{1}{2(1-\gamma)} \sum_{s}d_\rho^{\pi^\star}(s) \norm{\pi^\star(\cdot\,\vert\,s) - \hat\pi_{t+1}(\cdot\,\vert\,s)}^2
         +
         \frac{1}{2}\norm{\xi^\star-\hat\xi_{t+1}}^2 
         +
         \frac{1}{2}  
            \norm{\lambda^\star-\hat\lambda_{t+1}}^2 
         \\[0.2cm]
         & \leq & \displaystyle
         \frac{1}{2(1-\gamma)} \sum_{s}d_\rho^{\pi^\star}(s)
        \norm{\pi^\star(\cdot\,\vert\,s) - \hat\pi_t(\cdot\,\vert\,s)}^2
            +
            \frac{1}{2}
            \norm{\xi^\star-\hat\xi_t}^2 
            + 
         \frac{1}{2}
            \norm{\lambda^\star - \hat\lambda_t}^2 
         \\[0.2cm]
         &  & \displaystyle -\,
         \frac{1}{2(1-\gamma)} \sum_{s}d_\rho^{\pi^\star}(s)\norm{\hat\pi_{t+1}(\cdot\,\vert\,s) - \hat\pi_t(\cdot\,\vert\,s)}^2
         -
         \frac{1}{2}
            \norm{\hat\xi_{t+1} -\hat\xi_t}^2
        - 
         \frac{1}{2}  \norm{\hat\lambda_{t+1}-\hat\lambda_t}^2
         \\[0.2cm]
         & &  \displaystyle
         +\,\frac{1}{2(1-\gamma)} \sum_{s}d_\rho^{\pi^\star}(s) 
         \left(
            \norm{\hat\pi_{t+1}(\cdot\,\vert\,s) - \hat\pi_t(\cdot\,\vert\,s)}^2
            - \norm{\hat\pi_{t+1}(\cdot\,\vert\,s) - \pi_t(\cdot\,\vert\,s)}^2
            -
            \norm{\pi_{t}(\cdot\,\vert\,s) - \hat\pi_t(\cdot\,\vert\,s)}^2
         \right)
         \\[0.2cm]
        &  & \displaystyle +\,
        \frac{1}{2}\left(
            \norm{\hat\xi_{t+1}-\hat\xi_t}^2-
            \norm{\hat\xi_{t+1}-\xi_t}^2-
            \norm{\xi_t-\hat\xi_t}^2
        \right)
        \\[0.2cm]
         & & \displaystyle +\,
         \frac{1}{2} \left(
            \norm{\hat\lambda_{t+1}-\hat\lambda_t}^2   - 
            \norm{\hat\lambda_{t+1}-\lambda_t}^2 - \norm{\lambda_{t}-\hat\lambda_t}^2
         \right) 
         \\[0.2cm]
         & &
         \displaystyle
         +\,
         4\eta^2 \iota \frac{1}{1-\gamma}
         \sum_s d_\rho^{\pi^\star}(s)\left(
         \norm{\pi_t(\cdot\,\vert\,s) - \hat\pi_{t}(\cdot\,\vert\,s)}^2
         +
        \norm{\hat\pi_t(\cdot\,\vert\,s) - \pi_{t-1}(\cdot\,\vert\,s)}^2
         \right)
        \\[0.2cm]
        &  & \displaystyle +\, 4 \eta^2 \iota \left(\norm{\xi_{t} - \hat\xi_{t}}^2
        +  \norm{\hat\xi_{t} - \xi_{t-1}}^2\right)
        \\[0.2cm]
         &  & +\,
         \displaystyle
         4\eta^2 \iota
         \left( 
         \norm{\lambda_t -\hat\lambda_{t}}^2 + \norm{\hat\lambda_t -\lambda_{t-1}}^2
         \right)
    \end{array}
    \]
    where we use the defintion of $\iota$.
    After some re-combination, we have 
    \[
    \begin{array}{rcl}
         &  & 
         \!\!\!\!  \!\!\!\!  \!\!
         \displaystyle
         \frac{1}{2(1-\gamma)} \sum_{s}d_\rho^{\pi^\star}(s) \norm{\pi^\star(\cdot\,\vert\,s) - \hat\pi_{t+1}(\cdot\,\vert\,s)}^2
         +
         \frac{1}{2}\norm{\xi^\star-\hat\xi_{t+1}}^2 
         +
         \frac{1}{2}  
            \norm{\lambda^\star-\hat\lambda_{t+1}}^2 
         \\[0.2cm]
         & \leq & \displaystyle
         \frac{1}{2(1-\gamma)} \sum_{s}d_\rho^{\pi^\star}(s)
        \norm{\pi^\star(\cdot\,\vert\,s) - \hat\pi_t(\cdot\,\vert\,s)}^2
            +
            \frac{1}{2}
            \norm{\xi^\star-\hat\xi_t}^2 
            + 
         \frac{1}{2}
            \norm{\lambda^\star - \hat\lambda_t}^2 
         \\[0.2cm]
         & &  \displaystyle
         -\,\frac{1}{2(1-\gamma)} \sum_{s}d_\rho^{\pi^\star}(s) 
             \norm{\hat\pi_{t+1}(\cdot\,\vert\,s) - \pi_t(\cdot\,\vert\,s)}^2
            -
        \frac{1}{2}
            \norm{\hat\xi_{t+1}-\xi_t}^2 -
         \frac{1}{2}  
            \norm{\hat\lambda_{t+1}-\lambda_t}^2
         \\[0.2cm]
         & &  \displaystyle
         -\,\left(\frac{1}{2} -4\eta^2\iota\right) \left( \frac{1}{1-\gamma}\sum_{s}d_\rho^{\pi^\star}(s) 
            \norm{\pi_{t}(\cdot\,\vert\,s) - \hat\pi_t(\cdot\,\vert\,s)}^2 +
            \norm{\xi_t-\hat\xi_t}^2 +  \norm{\lambda_{t}-\hat\lambda_t}^2 \right)
         \\[0.2cm]
         & &
         \displaystyle
         +\,
         4\eta^2 \iota  \frac{1}{1-\gamma}
         \sum_s d_\rho^{\pi^\star}(s)
        \norm{\hat\pi_t(\cdot\,\vert\,s) - \pi_{t-1}(\cdot\,\vert\,s)}^2 + 4 \eta^2 \iota \norm{\hat\xi_{t} - \xi_{t-1}}^2 +
         \displaystyle
         4\eta^2 \iota
          \norm{\hat\lambda_t -\lambda_{t-1}}^2.
    \end{array}
    \]
    
    By taking 
    \[
    \pi^\star(\cdot\,\vert\,s) 
    \; = \;
    \mathcal{P}_{\Pi^\star}( \hat\pi_{t}(\cdot\,\vert\,s)),\; 
    \xi^\star \; = \;
    \mathcal{P}_{\Xi^\star}(\hat\xi_{t}), \;\text{ and }\;
    \lambda^\star \; = \;
    \mathcal{P}_{\Lambda^\star}(\hat\lambda_{t})
    \]
    and using 
    the non-expansivenss of projection operators $\mathcal{P}_{\Pi^\star}$, $\mathcal{P}_{\Xi^\star}$, and $\mathcal{P}_{\Lambda^\star}$,
    we obtain the following inequality,
    \[
    \begin{array}{rcl}
         &  & 
         \!\!\!\!  \!\!\!\!  \!\!
         \displaystyle
         \frac{1}{2(1-\gamma)} \sum_{s}d_\rho^{\pi^\star}(s) \norm{\mathcal{P}_{\Pi^\star}( \hat\pi_{t+1}(\cdot\,\vert\,s)) - \hat\pi_{t+1}(\cdot\,\vert\,s)}^2
         +
         \frac{1}{2}\norm{\mathcal{P}_{\Xi^\star}(\hat\xi_{t+1})-\hat\xi_{t+1}}^2 
         +
         \frac{1}{2}  
            \norm{\mathcal{P}_{\Lambda^\star}(\hat\lambda_{t+1})-\hat\lambda_{t+1}}^2 
         \\[0.2cm]
         & \leq & \displaystyle
         \frac{1}{2(1-\gamma)} \sum_{s}d_\rho^{\pi^\star}(s)
        \norm{\mathcal{P}_{\Pi^\star}( \hat\pi_{t}(\cdot\,\vert\,s)) - \hat\pi_t(\cdot\,\vert\,s)}^2
            +
            \frac{1}{2}
            \norm{\mathcal{P}_{\Xi^\star}(\hat\xi_{t})-\hat\xi_t}^2 
            + 
         \frac{1}{2}
            \norm{\mathcal{P}_{\Lambda^\star}(\hat\lambda_{t}) - \hat\lambda_t}^2 
         \\[0.2cm]
         & &  \displaystyle
         -\,\frac{1}{2(1-\gamma)} \sum_{s}d_\rho^{\pi^\star}(s) 
             \norm{\hat\pi_{t+1}(\cdot\,\vert\,s) - \pi_t(\cdot\,\vert\,s)}^2
            -
        \frac{1}{2}
            \norm{\hat\xi_{t+1}-\xi_t}^2 -
         \frac{1}{2}  
            \norm{\hat\lambda_{t+1}-\lambda_t}^2
         \\[0.2cm]
         & &  \displaystyle
         -\,\left(\frac{1}{2} -4\eta^2\iota\right) \left( \frac{1}{1-\gamma}\sum_{s}d_\rho^{\pi^\star}(s) 
            \norm{\pi_{t}(\cdot\,\vert\,s) - \hat\pi_t(\cdot\,\vert\,s)}^2 +
            \norm{\xi_t-\hat\xi_t}^2 +  \norm{\lambda_{t}-\hat\lambda_t}^2 \right)
         \\[0.2cm]
         & &
         \displaystyle
         +\,
         4\eta^2 \iota \frac{1}{1-\gamma}
         \sum_s d_\rho^{\pi^\star}(s)
        \norm{\hat\pi_t(\cdot\,\vert\,s) - \pi_{t-1}(\cdot\,\vert\,s)}^2 + 4 \eta^2 \iota \norm{\hat\xi_{t} - \xi_{t-1}}^2 +
         \displaystyle
         4\eta^2 \iota
          \norm{\hat\lambda_t -\lambda_{t-1}}^2.
    \end{array}
    \]
    If we choose $\eta>0$ such that $\frac{1}{2} - 4\eta^2 \iota\geq \frac{1}{4}$ and do some re-arrangement, we have
    \[
    \begin{array}{rcl}
         &  & 
         \!\!\!\!  \!\!\!\!  \!\!
         \displaystyle
         \frac{1}{2(1-\gamma)} \sum_{s}d_\rho^{\pi^\star}(s) \norm{\mathcal{P}_{\Pi^\star}( \hat\pi_{t+1}(\cdot\,\vert\,s)) - \hat\pi_{t+1}(\cdot\,\vert\,s)}^2
         +
         \frac{1}{2}\norm{\mathcal{P}_{\Xi^\star}(\hat\xi_{t+1})-\hat\xi_{t+1}}^2 
         +
         \frac{1}{2}  
            \norm{\mathcal{P}_{\Lambda^\star}(\hat\lambda_{t+1})-\hat\lambda_{t+1}}^2 
        \\[0.2cm]
         & &  \!\!\!\!  \!\!\!\!  \!\!\displaystyle
         +\,\frac{1}{4(1-\gamma)}  \sum_{s}d_\rho^{\pi^\star}(s) 
            \norm{\hat\pi_{t+1}(\cdot\,\vert\,s) - \pi_t(\cdot\,\vert\,s)}^2 + \frac{1}{4}
            \norm{\hat\xi_{t+1}-\xi_t}^2 +  \frac{1}{4}\norm{\hat\lambda_{t+1}-\lambda_t}^2 
         \\[0.2cm]
         & \leq & \displaystyle
         \frac{1}{2(1-\gamma)} \sum_{s}d_\rho^{\pi^\star}(s)
        \norm{\mathcal{P}_{\Pi^\star}( \hat\pi_{t}(\cdot\,\vert\,s)) - \hat\pi_t(\cdot\,\vert\,s)}^2
            +
            \frac{1}{2}
            \norm{\mathcal{P}_{\Xi^\star}(\hat\xi_{t})-\hat\xi_t}^2 
            + 
         \frac{1}{2}
            \norm{\mathcal{P}_{\Lambda^\star}(\hat\lambda_{t}) - \hat\lambda_t}^2 
        \\[0.2cm]
         & &  \displaystyle
         -\,\frac{1}{4(1-\gamma)} \sum_{s}d_\rho^{\pi^\star}(s) 
             \norm{\hat\pi_{t+1}(\cdot\,\vert\,s) - \pi_t(\cdot\,\vert\,s)}^2
            -
        \frac{1}{4}
            \norm{\hat\xi_{t+1}-\xi_t}^2 -
         \frac{1}{4}  
            \norm{\hat\lambda_{t+1}-\lambda_t}^2
         \\[0.2cm]
         & &  \displaystyle
         -\,\frac{1}{4(1-\gamma)} \sum_{s}d_\rho^{\pi^\star}(s) 
             \norm{\pi_{t}(\cdot\,\vert\,s) - \hat\pi_t(\cdot\,\vert\,s)}^2
            -
        \frac{1}{4}
            \norm{\xi_{t}-\hat\xi_t}^2 -
         \frac{1}{4}  
            \norm{\lambda_{t}-\hat\lambda_t}^2
         \\[0.2cm]
         & &
         \displaystyle
         +\,
         \frac{1}{4(1-\gamma)}
         \sum_s d_\rho^{\pi^\star}(s)
        \norm{\hat\pi_t(\cdot\,\vert\,s) - \pi_{t-1}(\cdot\,\vert\,s)}^2 + \frac{1}{4} \norm{\hat\xi_{t} - \xi_{t-1}}^2 +
         \displaystyle
         \frac{1}{4}
          \norm{\hat\lambda_t -\lambda_{t-1}}^2.
    \end{array}
    \]
    Finally, our desired inequality is obtained by using notation $\Theta$ and $\zeta$.
\end{proof}

\begin{lemma}\label{lem:lower_duality_gap}
    In Algorithm~\ref{alg: resilient OPG}, for~\eqref{eq:optimistic_policy_gradient} with $\eta \leq \min\left( \frac{1}{4\sqrt{|A|}}, \frac{1}{2(L_h+1)}, \frac{1}{5 \sqrt{m|A|}\kappa_\rho}, \frac{\rho_{\min}}{4 \gamma\sqrt{m}C_h |A|}
    \right)$, 
    \[
    \begin{array}{rcl}
     & &
     \!\!\!\!  \!\!\!\!  \!\!
     \displaystyle
         \frac{1}{1-\gamma}\sum_{s}d_\rho^{\pi^\star}(s) 
             \norm{\hat\pi_{t+1}(\cdot\,\vert\,s) - \pi_t(\cdot\,\vert\,s)}^2
            +
            \norm{\hat\xi_{t+1}-\xi_t}^2 
            +  
            \norm{\hat\lambda_{t+1}-\lambda_t}^2
         \\[0.2cm]
         & &  \displaystyle
         \!\!\!\!  \!\!\!\!  \!\!
         +\,\frac{1}{1-\gamma} \sum_{s}d_\rho^{\pi^\star}(s) 
             \norm{\pi_{t}(\cdot\,\vert\,s) - \hat\pi_t(\cdot\,\vert\,s)}^2
            +
            \norm{\xi_{t}-\hat\xi_t}^2 
            +  
            \norm{\lambda_{t}
            -\hat\lambda_t}^2
        \\[0.2cm]
         & \geq &  
         \displaystyle
         \frac{\eta^2}{9 \max\left( \frac{\kappa_\rho}{1-\gamma},1\right)^2} \frac{\left[ \left(
            V_{r+\hat\lambda_{t+1}^\top g}^{\pi}(\rho)
             - h(\xi) - \hat\lambda_{t+1}^\top \xi\right)
             - 
             \left(
            V_{r+\lambda^\top g}^{\hat\pi_{t+1}}(\rho)
             - h(\hat\xi_{t+1}) - \lambda^\top \hat\xi_{t+1}\right) \right]_+^2}{\left(
                \max_s \norm{\pi(\cdot\,\vert\,s)-\hat\pi_{t+1}(\cdot\,\vert\,s))}
                +
                \norm{\xi - \hat\xi_{t+1}}
                +
                \norm{\lambda - \hat\lambda_{t+1}}
            \right)^2}
    \end{array}
\]
for any $(\pi, \xi, \lambda) \neq (\hat\pi_{t+1}, \hat\xi_{t+1}, \hat\lambda_{t+1})$.
\end{lemma}
\begin{proof}
    By the optimality of $\hat\pi_{t+1}$,
    \[
    \left\langle
        Q_{r+\lambda_t^\top g}^{\pi_t}(s,\cdot) - \frac{1}{\eta}\left( \hat\pi_{t+1}(\cdot\,\vert\,s) - \hat\pi_t(\cdot\,\vert\,s) \right),  \hat\pi_{t+1}(\cdot\,\vert\,s) - \pi(\cdot\,\vert\,s)
    \right\rangle
    \; \geq \; 0\; \text{ for all }\pi \in\Pi.
    \]
    Hence,
    \begin{equation}\label{eq:primal_pi}
        \begin{array}{rcl}
             &  &
            \!\!\!\!  \!\!\!\!  \!\!
            \displaystyle
             \left\langle
                 \hat\pi_{t+1}(\cdot\,\vert\,s) - \hat\pi_t(\cdot\,\vert\,s),  \pi(\cdot\,\vert\,s) -\hat\pi_{t+1}(\cdot\,\vert\,s)
            \right\rangle
             \\[0.2cm]
             & \geq &
             \displaystyle \eta 
             \left\langle
            Q_{r+\lambda_t^\top g}^{\pi_t}(s,\cdot) ,   \pi(\cdot\,\vert\,s) -\hat\pi_{t+1}(\cdot\,\vert\,s) 
            \right\rangle
            \\[0.2cm]
             & = &
             \displaystyle \eta 
             \left\langle
            Q_{r+\hat\lambda_{t+1}^\top g}^{\pi_{t+1}}(s,\cdot) ,   \pi(\cdot\,\vert\,s) -\hat\pi_{t+1}(\cdot\,\vert\,s) 
            \right\rangle
            +
            \eta 
             \left\langle
            Q_{r+\lambda_{t}^\top g}^{\pi_{t}}(s,\cdot)-Q_{r+\lambda_{t}^\top g}^{\hat\pi_{t+1}}(s,\cdot) ,   \pi(\cdot\,\vert\,s) -\hat\pi_{t+1}(\cdot\,\vert\,s) 
            \right\rangle
            \\[0.2cm]
             &  &
             \displaystyle +\,\eta 
             \left\langle
            Q_{r+\lambda_{t}^\top g}^{\hat\pi_{t+1}}(s,\cdot)
            -Q_{r+{\hat\lambda}_{t+1}^\top g}^{\hat\pi_{t+1}}(s,\cdot),   \pi(\cdot\,\vert\,s) -\hat\pi_{t+1}(\cdot\,\vert\,s) 
            \right\rangle.
        \end{array}
    \end{equation}
    Similarly, the optimality of $\pi_{t+1}$, 
    \[
    \left\langle
        Q_{r+\lambda_t^\top g}^{\pi_t}(s,\cdot) - \frac{1}{\eta}\left( \pi_{t+1}(\cdot\,\vert\,s) - \hat\pi_{t+1}(\cdot\,\vert\,s) \right),  \pi_{t+1}(\cdot\,\vert\,s) - \pi(\cdot\,\vert\,s)
    \right\rangle
    \; \geq \; 0\; \text{ for all }\pi \in\Pi
    \]
    implies that
    \[
        \begin{array}{rcl}
             &  &
            \!\!\!\!  \!\!\!\!  \!\!
            \displaystyle
             \left\langle
                 \pi_{t+1}(\cdot\,\vert\,s) - \hat\pi_{t+1}(\cdot\,\vert\,s),  \pi(\cdot\,\vert\,s) -\pi_{t+1}(\cdot\,\vert\,s)
            \right\rangle
             \\[0.2cm]
             & \geq &
             \displaystyle \eta 
             \left\langle
            Q_{r+\lambda_t^\top g}^{\pi_t}(s,\cdot) ,   \pi(\cdot\,\vert\,s) -\pi_{t+1}(\cdot\,\vert\,s) 
            \right\rangle
            \\[0.2cm]
             & = &
             \displaystyle \eta 
             \left\langle
            Q_{r+\lambda_{t+1}^\top g}^{\pi_{t+1}}(s,\cdot) ,   \pi(\cdot\,\vert\,s) -\pi_{t+1}(\cdot\,\vert\,s) 
            \right\rangle
            +
            \eta 
             \left\langle
            Q_{r+\lambda_{t}^\top g}^{\pi_{t}}(s,\cdot)-Q_{r+\lambda_{t}^\top g}^{\pi_{t+1}}(s,\cdot) ,   \pi(\cdot\,\vert\,s) -\pi_{t+1}(\cdot\,\vert\,s) 
            \right\rangle
            \\[0.2cm]
             &  &
             \displaystyle +\,\eta 
             \left\langle
            Q_{r+\lambda_{t}^\top g}^{\pi_{t+1}}(s,\cdot)
            -Q_{r+\lambda_{t+1}^\top g}^{\pi_{t+1}}(s,\cdot),   \pi(\cdot\,\vert\,s) -\pi_{t+1}(\cdot\,\vert\,s) 
            \right\rangle.
        \end{array}
    \]

    By the optimality of $\hat\xi_{t+1}$,
    \[
    \left\langle
    -\nabla h(\xi_t) - \lambda_t  - \frac{1}{\eta}\left( \hat\xi_{t+1} - \hat\xi_t \right), \hat\xi_{t+1}-\xi
    \right\rangle
    \; \geq \; 0 \; \text{ for all }
    \xi \in \Xi.
    \]
    Hence,
    \begin{equation}\label{eq:primal_xi}
        \begin{array}{rcl}
             &  &
            \!\!\!\!  \!\!\!\!  \!\!
            \displaystyle
            \left\langle
            \hat\xi_{t+1} - \hat\xi_t,
            \xi - \hat\xi_{t+1}
            \right\rangle
            \\[0.2cm]
            & \geq &
            \eta \left\langle
            -\nabla h(\xi_t)-\lambda_t,
            \xi - \hat\xi_{t+1}
            \right\rangle
            \\[0.2cm]
            & = &
            \eta \left\langle
            -\nabla h(\hat\xi_{t+1})-\hat\lambda_{t+1},
            \xi - \hat\xi_{t+1}
            \right\rangle
            + \eta \left\langle
            -\nabla h(\xi_{t})-\lambda_{t} + \nabla h(\hat\xi_{t+1})+\lambda_t,
            \xi - \hat\xi_{t+1}
            \right\rangle
            \\[0.2cm]
            &  & +\,
             \eta \left\langle
            -\nabla h(\hat\xi_{t+1})-\lambda_{t} + \nabla h(\hat\xi_{t+1})+\hat\lambda_{t+1},
            \xi - \hat\xi_{t+1}
            \right\rangle.
        \end{array}
    \end{equation}
    Similarly, the optimality of $\xi_{t+1}$,
    \[
    \left\langle
    -\nabla h(\xi_t) - \lambda_t  - \frac{1}{\eta}\left( \xi_{t+1} - \hat\xi_{t+1} \right), \xi_{t+1}-\xi
    \right\rangle
    \; \geq \; 0 \; \text{ for all }
    \xi \in \Xi
    \]
    implies that
    \[
        \begin{array}{rcl}
             &  &
            \!\!\!\!  \!\!\!\!  \!\!
            \displaystyle
            \left\langle
            \xi_{t+1} - \hat\xi_{t+1},
            \xi - \xi_{t+1}
            \right\rangle
            \\[0.2cm]
            & \geq &
            \eta \left\langle
            -\nabla h(\xi_t)-\lambda_t,
            \xi - \xi_{t+1}
            \right\rangle
            \\[0.2cm]
            & = &
            \eta \left\langle
            -\nabla h(\xi_{t+1})-\lambda_{t+1},
            \xi - \xi_{t+1}
            \right\rangle
            + \eta \left\langle
            -\nabla h(\xi_{t})-\lambda_{t} + \nabla h(\xi_{t+1})+\lambda_t,
            \xi - \xi_{t+1}
            \right\rangle
            \\[0.2cm]
            &  & +\,
             \eta \left\langle
            -\nabla h(\xi_{t+1})-\lambda_{t} + \nabla h(\xi_{t+1})+\lambda_{t+1},
            \xi - \xi_{t+1}
            \right\rangle.
        \end{array}
    \]

    By the optimality of $\hat\lambda_{t+1}$,
    \[
    \left\langle
        V_g^{\pi_t}(\rho) - \xi_t +\frac{1}{\eta}\left( \hat\lambda_{t+1} - \hat\lambda_t \right), 
        \hat\lambda_{t+1} - \lambda
    \right\rangle
    \; \leq \; 0\;
    \text{ for all } \lambda \in \Lambda.
    \]
    Hence, 
    \begin{equation}\label{eq:dual_lambda}
    \begin{array}{rcl}
             &  &
            \!\!\!\!  \!\!\!\!  \!\!
            \displaystyle
            \left\langle
            \hat\lambda_{t+1} - \hat\lambda_{t},
            \lambda - \hat\lambda_{t+1}
            \right\rangle
            \\[0.2cm]
            & \geq &
            \eta \left\langle
            V_g^{\pi_t}(\rho) - \xi_t, \hat\lambda_{t+1} - \lambda
            \right\rangle
            \\[0.2cm]
            & = &
            \eta \left\langle
            V_g^{\hat\pi_{t+1}}(\rho) - \hat\xi_{t+1}, \hat\lambda_{t+1} - \lambda
            \right\rangle
            + 
            \eta \left\langle
            V_g^{\pi_{t}}(\rho) - \xi_{t} - V_g^{\hat\pi_{t+1}}(\rho) + \hat\xi_{t+1}
            , \hat\lambda_{t+1} - \lambda
            \right\rangle. 
    \end{array}
    \end{equation}
    Similarly, the optimality of $\lambda_{t+1}$,
    \[
    \left\langle
        V_g^{\pi_t}(\rho) - \xi_t +\frac{1}{\eta}\left( \lambda_{t+1} - \hat\lambda_{t+1} \right), 
        \lambda_{t+1} - \lambda
    \right\rangle
    \; \leq \; 0\;
    \text{ for all } \lambda \in \Lambda
    \]
    implies that
    \[
    \begin{array}{rcl}
             &  &
            \!\!\!\!  \!\!\!\!  \!\!
            \displaystyle
            \left\langle
            \lambda_{t+1} - \hat\lambda_{t+1},
            \lambda - \lambda_{t+1}
            \right\rangle
            \\[0.2cm]
            & \geq &
            \eta \left\langle
            V_g^{\pi_t}(\rho) - \xi_t, \lambda_{t+1} - \lambda
            \right\rangle
            \\[0.2cm]
            & = &
            \eta \left\langle
            V_g^{\pi_{t+1}}(\rho) - \xi_{t+1}, \lambda_{t+1} - \lambda
            \right\rangle
            + 
            \eta \left\langle
            V_g^{\pi_{t}}(\rho) - \xi_{t} - V_g^{\pi_{t+1}}(\rho) + \xi_{t+1}
            , \lambda_{t+1} - \lambda
            \right\rangle. 
    \end{array}
    \]

    Summing up the inequalities~\eqref{eq:primal_pi},~\eqref{eq:primal_xi}, and~\eqref{eq:dual_lambda} from both sides, with some state distribution $d_\rho^{\pi}$, yields,
    \[
    \begin{array}{rcl}
             &  &
            \!\!\!\!  \!\!\!\!  \!\!
            \displaystyle
            \frac{1}{1-\gamma}\sum_s d_\rho^{\pi}(s)
            \left\langle \hat\pi_{t+1}(\cdot\,\vert\,s) - \hat\pi_t(\cdot\,\vert\,s),
            \pi(\cdot\,\vert\,s) - \hat\pi_{t+1}(\cdot\,\vert\,s)\right\rangle
            \\[0.2cm]
            &  &
            \!\!\!\!  \!\!\!\!  \!\!
            \displaystyle
            + \left\langle \hat\xi_{t+1}-\hat\xi_t,
            \xi-\hat\xi_{t+1}
            \right\rangle
            + \left\langle
            \hat\lambda_{t+1} - \hat\lambda_t, 
            \lambda - \hat\lambda_{t+1}
            \right\rangle
            \\[0.2cm]
            & \geq &
             \displaystyle \frac{\eta}{1-\gamma} \sum_s d_\rho^{\pi}(s) 
             \left\langle
            Q_{r+\hat\lambda_{t+1}^\top g}^{\pi_{t+1}}(s,\cdot) ,   \pi(\cdot\,\vert\,s) -\hat\pi_{t+1}(\cdot\,\vert\,s) 
            \right\rangle
            \\[0.2cm]
            && \displaystyle+\,
            \frac{\eta}{1-\gamma} \sum_s d_\rho^{\pi}(s) 
             \left\langle
            Q_{r+\lambda_{t}^\top g}^{\pi_{t}}(s,\cdot)-Q_{r+\lambda_{t}^\top g}^{\hat\pi_{t+1}}(s,\cdot) ,   \pi(\cdot\,\vert\,s) -\hat\pi_{t+1}(\cdot\,\vert\,s) 
            \right\rangle
            \\[0.2cm]
             &  &
             \displaystyle +\,\frac{\eta }{1-\gamma}
             \sum_s d_\rho^{\pi}(s) \left\langle
            Q_{r+\lambda_{t}^\top g}^{\hat\pi_{t+1}}(s,\cdot)
            -Q_{r+{\hat\lambda}_{t+1}^\top g}^{\hat\pi_{t+1}}(s,\cdot),   \pi(\cdot\,\vert\,s) -\hat\pi_{t+1}(\cdot\,\vert\,s) 
            \right\rangle
            \\[0.2cm]
            &  & +\,
            \eta \left\langle
            -\nabla h(\hat\xi_{t+1})-\hat\lambda_{t+1},
            \xi - \hat\xi_{t+1}
            \right\rangle
            + \eta \left\langle
            -\nabla h(\xi_{t})-\lambda_{t} + \nabla h(\hat\xi_{t+1})+\lambda_t,
            \xi - \hat\xi_{t+1}
            \right\rangle
            \\[0.2cm]
            &  & +\,
             \eta \left\langle
            -\nabla h(\hat\xi_{t+1})-\lambda_{t} + \nabla h(\hat\xi_{t+1})+\hat\lambda_{t+1},
            \xi - \hat\xi_{t+1}
            \right\rangle
            \\[0.2cm]
            & & +\,
            \eta \left\langle
            V_g^{\hat\pi_{t+1}}(\rho) - \hat\xi_{t+1}, \hat\lambda_{t+1} - \lambda
            \right\rangle
            + 
            \eta \left\langle
            V_g^{\pi_{t}}(\rho) - \xi_{t} - V_g^{\hat\pi_{t+1}}(\rho) + \hat\xi_{t+1}
            , \hat\lambda_{t+1} - \lambda
            \right\rangle
            \\[0.2cm]
            & \geq &
             \displaystyle \eta  
             \left( 
            V_{r+\hat\lambda_{t+1}^\top g}^{\pi}(\rho)
            -
            V_{r+\lambda^\top g}^{\hat\pi_{t+1}}(\rho)
            \right)
            +
            \eta \left( h(\hat\xi_{t+1}) - h(\xi)\right)
            +
            \eta \left\langle
            -\hat\lambda_{t+1},
            \xi - \hat\xi_{t+1}
            \right\rangle
            +
            \eta \left\langle
            - \hat\xi_{t+1}, \hat\lambda_{t+1} - \lambda
            \right\rangle
            \\[0.2cm]
            && \displaystyle+\,
            \frac{\eta}{1-\gamma} \sum_s d_\rho^{\pi}(s) 
             \left\langle
            Q_{r+\lambda_{t}^\top g}^{\pi_{t}}(s,\cdot)-Q_{r+\lambda_{t}^\top g}^{\hat\pi_{t+1}}(s,\cdot) ,   \pi(\cdot\,\vert\,s) -\hat\pi_{t+1}(\cdot\,\vert\,s) 
            \right\rangle
            \\[0.2cm]
             &  &
             \displaystyle +\,\frac{\eta}{1-\gamma} 
             \sum_s d_\rho^{\pi}(s) \left\langle
            Q_{r+\lambda_{t}^\top g}^{\hat\pi_{t+1}}(s,\cdot)
            -Q_{r+{\hat\lambda}_{t+1}^\top g}^{\hat\pi_{t+1}}(s,\cdot),   \pi(\cdot\,\vert\,s) -\hat\pi_{t+1}(\cdot\,\vert\,s) 
            \right\rangle
            \\[0.2cm]
            &  & 
            +\, \eta \left\langle
            -\nabla h(\xi_{t})-\lambda_{t} + \nabla h(\hat\xi_{t+1})+\lambda_t,
            \xi - \hat\xi_{t+1}
            \right\rangle
            \\[0.2cm]
            &  & +\,
             \eta \left\langle
            -\nabla h(\hat\xi_{t+1})-\lambda_{t} + \nabla h(\hat\xi_{t+1})+\hat\lambda_{t+1},
            \xi - \hat\xi_{t+1}
            \right\rangle
            \\[0.2cm]
            & & 
            +\, 
            \eta \left\langle
            V_g^{\pi_{t}}(\rho) - \xi_{t} - V_g^{\hat\pi_{t+1}}(\rho) + \hat\xi_{t+1}
            , \hat\lambda_{t+1} - \lambda
            \right\rangle
            \\[0.2cm]
            & \geq &
             \displaystyle \eta  
             \left( 
            V_{r+\hat\lambda_{t+1}^\top g}^{\pi}(\rho)
            -
            V_{r+\lambda^\top g}^{\hat\pi_{t+1}}(\rho)
            \right)
            +
            \eta \left( h(\hat\xi_{t+1}) - h(\xi)\right)
            +
            \eta \left\langle
            -\hat\lambda_{t+1},
            \xi - \hat\xi_{t+1}
            \right\rangle
            +
            \eta \left\langle
            - \hat\xi_{t+1}, \hat\lambda_{t+1} - \lambda
            \right\rangle
            \\[0.2cm]
            && \displaystyle -\,
            \frac{\eta}{1-\gamma} \sum_s d_\rho^{\pi}(s) 
            \norm{
            Q_{r+\lambda_{t}^\top g}^{\pi_{t}}(s,\cdot)-Q_{r+\lambda_{t}^\top g}^{\hat\pi_{t+1}}(s,\cdot) }_\infty   \norm{\pi(\cdot\,\vert\,s) -\hat\pi_{t+1}(\cdot\,\vert\,s) }_1
            \\[0.2cm]
             &  &
             \displaystyle -\,\frac{\eta}{1-\gamma} 
             \sum_s d_\rho^{\pi}(s) 
             \norm{
            Q_{r+\lambda_{t}^\top g}^{\hat\pi_{t+1}}(s,\cdot)
            -Q_{r+{\hat\lambda}_{t+1}^\top g}^{\hat\pi_{t+1}}(s,\cdot)
            }_\infty
            \norm{\pi(\cdot\,\vert\,s) -\hat\pi_{t+1}(\cdot\,\vert\,s) 
            }_1
            \\[0.2cm]
            &  & 
            -  \eta \norm{
            \nabla h(\xi_{t})- \nabla h(\hat\xi_{t+1})}\norm{
            \xi - \hat\xi_{t+1}
            }
            -
             \eta \norm{\lambda_{t} - \hat\lambda_{t+1}}\norm{
            \xi - \hat\xi_{t+1}
            }
            \\[0.2cm]
            & & 
            -\, 
            \eta  \norm{
            V_g^{\pi_{t}}(\rho)  - V_g^{\hat\pi_{t+1}}(\rho) 
            }
            \norm{\hat\lambda_{t+1} - \lambda
            }
            -
            \eta  \norm{
             \xi_{t} -  \hat\xi_{t+1}
            }
            \norm{\hat\lambda_{t+1} - \lambda
            }
            \\[0.2cm]
            & \geq &
             \displaystyle \eta  
             \left( 
            V_{r+\hat\lambda_{t+1}^\top g}^{\pi}(\rho)
            -
            V_{r+\lambda^\top g}^{\hat\pi_{t+1}}(\rho)
            \right)
            +
            \eta \left( h(\hat\xi_{t+1}) - h(\xi)\right)
            +
            \eta \left\langle
            -\hat\lambda_{t+1},
            \xi - \hat\xi_{t+1}
            \right\rangle
            +
            \eta \left\langle
            - \hat\xi_{t+1}, \hat\lambda_{t+1} - \lambda
            \right\rangle
            \\[0.2cm]
            && \displaystyle -\,
            \frac{\gamma \eta \sqrt{m} C_h}{(1-\gamma)^4}
        \max_s  
        \norm{ 
        \pi_t(\cdot\,\vert\,s) - \hat\pi_{t+1}(\cdot\,\vert\,s)
        }_1 
            \sum_s d_\rho^{\pi}(s) 
          \norm{\pi(\cdot\,\vert\,s) -\hat\pi_{t+1}(\cdot\,\vert\,s) }_1
            \\[0.2cm]
             &  &
             \displaystyle -\,\frac{\eta}{(1-\gamma)^2} \norm{\lambda_t - \hat\lambda_{t+1}}
             \sum_s d_\rho^{\pi}(s) \norm{\pi(\cdot\,\vert\,s) -\hat\pi_{t+1}(\cdot\,\vert\,s) 
            }_1
            \\[0.2cm]
            &  & 
            -  \eta L_h \norm{
            \xi_{t} -\hat\xi_{t+1}}\norm{
            \xi - \hat\xi_{t+1}
            }
            -
             \eta \norm{\lambda_{t} - \hat\lambda_{t+1}}\norm{
            \xi - \hat\xi_{t+1}
            }
            \\[0.2cm]
            & & \displaystyle
            -\,  
            \frac{\eta \sqrt{m} \kappa_\rho}{(1-\gamma)^3}
            \norm{\hat\lambda_{t+1} - \lambda
            }\sum_s d_\rho^{\pi^\star}(s) \norm{  \pi_t(\cdot\,\vert\,s) - \hat\pi_{t+1}(\cdot\,\vert\,s)}_1
            -
            \eta  \norm{
             \xi_{t} -  \hat\xi_{t+1}
            }
            \norm{\hat\lambda_{t+1} - \lambda
            }
    \end{array}
    \]
    where the second inequality is due to the performance difference lemma and the convexity of $h(\xi)$,
    \[
    h(\xi) \; \geq \;
    h(\hat\xi_{t+1}) + \left\langle 
    \nabla h(\hat\xi_{t+1}), \xi - \hat\xi_{t+1}
    \right\rangle
    \]
    and the the last two inequalities is due to Cauchy–Schwarz inequality, and the inequalities in Lemma~\ref{lem:policy-value-difference} and the smoothness of $h$,
    \[
        \norm{
            Q_{r+\lambda_{t}^\top g}^{\pi_{t}}(s,\cdot)-Q_{r+\lambda_{t}^\top g}^{\hat\pi_{t+1}}(s,\cdot) }_\infty 
        \; \leq \;
        \frac{\gamma \sqrt{m} C_h}{(1-\gamma)^3}
        \max_s \; 
        \norm{ 
        \pi_t(\cdot\,\vert\,s) - \hat\pi_{t+1}(\cdot\,\vert\,s)
        }_1
    \]
    \[
    \norm{
            Q_{r+\lambda_{t}^\top g}^{\hat\pi_{t+1}}(s,\cdot)
            -Q_{r+{\hat\lambda}_{t+1}^\top g}^{\hat\pi_{t+1}}(s,\cdot)
            }_\infty
            \; \leq \;
        \frac{1}{1-\gamma} 
        \norm{ 
        \lambda_t - \hat\lambda_{t+1}
        }
    \]
    \[
    \norm{V_{g}^{\pi_t}(\rho) - V_g^{\hat\pi_{t+1}}(\rho)}
    \; \leq \;
    \frac{\sqrt{m} \kappa_\rho}{(1-\gamma)^3} \sum_s d_\rho^{\pi^\star}(s) \norm{  \pi_t(\cdot\,\vert\,s) - \hat\pi_{t+1}(\cdot\,\vert\,s)}_1 
    \]
    \[
    \norm{\nabla h(\xi_t) - \nabla h(\hat\xi_{t+1})}
    \; \leq \; L_h
    \norm{\xi_t - \hat\xi_{t+1}}.
    \]
    
    We further notice that
    \[
    \begin{array}{rcl}
             &  &
            \!\!\!\!  \!\!\!\!  \!\!
            \displaystyle
            \frac{1}{1-\gamma}\sum_s d_\rho^{\pi}(s)
            \left\langle \hat\pi_{t+1}(\cdot\,\vert\,s) - \hat\pi_t(\cdot\,\vert\,s),
            \pi(\cdot\,\vert\,s) - \hat\pi_{t+1}(\cdot\,\vert\,s)\right\rangle
            \\[0.2cm]
            &  &
            \!\!\!\!  \!\!\!\!  \!\!
            \displaystyle
            + \left\langle \hat\xi_{t+1}-\hat\xi_t,
            \xi-\hat\xi_{t+1}
            \right\rangle
            + \left\langle
            \hat\lambda_{t+1} - \hat\lambda_t, 
            \lambda - \hat\lambda_{t+1}
            \right\rangle
            \\[0.2cm]
            & \leq &
            \displaystyle
            \frac{1}{1-\gamma}\max_s \norm{ \pi(\cdot\,\vert\,s) - \hat\pi_t(\cdot\,\vert\,s)} \sum_s d_\rho^\pi(s) \norm{\hat\pi_{t+1}(\cdot\,\vert\,s) - \hat\pi_t(\cdot\,\vert\,s)}
            \\[0.2cm]
            &  &
            \displaystyle
            +\, \norm{ \hat\xi_{t+1}-\hat\xi_t}
            \norm{
            \xi-\hat\xi_{t+1}
            }
            + \norm{
            \hat\lambda_{t+1} - \hat\lambda_t}
            \norm{
            \lambda - \hat\lambda_{t+1}
            }.
    \end{array}
    \]
    Application of the inequality $ac  + bd\leq (a+b)(c+d)$ for $a\geq 0$, $b\geq 0$, $c\geq 0$, and $d\geq 0$ and $d_\rho^\pi(s) \leq \frac{\kappa_\rho}{1-\gamma} d_\rho^{\pi^\star}(s)$
    leads to
    \[
    \begin{array}{rcl}
             &  &
            \!\!\!\!  \!\!\!\!  \!\!
            \displaystyle
            \eta \left(
            V_{r+\hat\lambda_{t+1}^\top g}^{\pi}(\rho)
             - h(\xi) - \hat\lambda_{t+1}^\top \xi\right)
             - 
             \eta \left(
            V_{r+\lambda^\top g}^{\hat\pi_{t+1}}(\rho)
             - h(\hat\xi_{t+1}) - \lambda^\top \hat\xi_{t+1}\right) 
            \\[0.2cm]
            & \leq & \displaystyle
            \left(
                \max_s \norm{\pi(\cdot\,\vert\,s)-\hat\pi_{t+1}(\cdot\,\vert\,s))}
                +
                \norm{\xi - \hat\xi_{t+1}}
                +
                \norm{\lambda - \hat\lambda_{t+1}}
            \right)
            \\[0.2cm]
            &  & \displaystyle 
            \times\Bigg(
                \frac{1}{1-\gamma}\sum_s d_\rho^\pi(s) \norm{\hat\pi_{t+1}(\cdot\,\vert\,s) - \hat\pi_t(\cdot\,\vert\,s)}
                + \norm{\hat\xi_{t+1}-\hat\xi_t}
                + \norm{\hat\lambda_{t+1}-\hat\lambda_t}
            \\[0.2cm]
            && \displaystyle \;\;\;\; \;\;\;\; +\,
            \frac{\gamma \eta \sqrt{m} C_h |A|}{(1-\gamma)^4}
        \max_s  
        \norm{ 
        \pi_t(\cdot\,\vert\,s) - \hat\pi_{t+1}(\cdot\,\vert\,s)
        } 
            \\[0.2cm]
             &  &
             \displaystyle \;\;\;\; \;\;\;\; 
             +\,\frac{2\eta \sqrt{|A|}}{(1-\gamma)^2} \norm{\lambda_t - \hat\lambda_{t+1}}
            + \eta (L_h+1) \norm{
            \xi_{t} -\hat\xi_{t+1}}
            \\[0.2cm]
            & & \displaystyle
            \;\;\;\; \;\;\;\;
            +\,  
            \frac{\eta \sqrt{m |A|} \kappa_\rho}{(1-\gamma)^3}
            \sum_s d_\rho^{\pi^\star}(s) \norm{  \pi_t(\cdot\,\vert\,s) - \hat\pi_{t+1}(\cdot\,\vert\,s)} \Bigg)
        \\[0.2cm]
            & \leq & \displaystyle
            \left(
                \max_s \norm{\pi(\cdot\,\vert\,s)-\hat\pi_{t+1}(\cdot\,\vert\,s))}
                +
                \norm{\xi - \hat\xi_{t+1}}
                +
                \norm{\lambda - \hat\lambda_{t+1}}
            \right)
            \\[0.2cm]
            &  & \displaystyle 
            \times\Bigg( \frac{\kappa_\rho}{(1-\gamma)^2}
                \sum_s d_\rho^{\pi^\star}(s) \norm{\hat\pi_{t+1}(\cdot\,\vert\,s) - \hat\pi_t(\cdot\,\vert\,s)} 
            \\[0.2cm]
            &  & \displaystyle 
            \;\;\;\; \;\;\;\; 
                + \,\left(\frac{\eta \sqrt{m |A|} \kappa_\rho}{(1-\gamma)^2} +
            \frac{\gamma \eta \sqrt{m} C_h |A|}{(1-\gamma)^4 \rho_{\min}}\right)\frac{1}{1-\gamma}
            \sum_s d_\rho^{\pi^\star}(s) 
            \norm{ 
            \pi_t(\cdot\,\vert\,s) - \hat\pi_{t+1}(\cdot\,\vert\,s)
            } 
            \\[0.2cm]
             &  &
             \displaystyle \;\;\;\; \underbrace{\;\;\;\; 
            +\,\frac{2\eta \sqrt{|A|}}{(1-\gamma)^2} \norm{\lambda_t - \hat\lambda_{t+1}}
            + \norm{\hat\lambda_{t+1}-\hat\lambda_t}
            + \eta (L_h+1) \norm{
            \xi_{t} -\hat\xi_{t+1}}
            + \norm{\hat\xi_{t+1}-\hat\xi_t}\Bigg)}_{\displaystyle\text{$\DefinedAs\;$Diff}}.
    \end{array}
    \]
    
    If we take $\eta>0$ such that 
    \[
    \frac{2\eta \sqrt{|A|}}{(1-\gamma)^2}\;\leq\; \frac{1}{2},
    \;
    \eta (L_h+1) \;\leq\; \frac{1}{2},
    \;
    \frac{\eta \sqrt{m|A|}\kappa_\rho}{(1-\gamma)^2} \;\leq\; \frac{1}{4},
    \, \text{ and }
    \frac{\gamma \eta \sqrt{m} C_h|A|}{(1-\gamma)^4\rho_{\min}} \;\leq\; \frac{1}{4}
    \]
    then, 
    \[
    \begin{array}{rcl}
         \text{Diff}^2 
         & \leq & \displaystyle
         \max\left( 
         \frac{\kappa_\rho}{1-\gamma}, 1\right)^2 
         \Bigg( \frac{1}{1-\gamma}
            \sum_s d_\rho^{\pi^\star}(s) \left( \norm{\hat\pi_{t+1}(\cdot\,\vert\,s) - \hat\pi_t(\cdot\,\vert\,s)} + \frac{1}{2}\norm{\pi_t(\cdot\,\vert\,s)-\hat\pi_{t+1}(\cdot\,\vert\,s)}\right) 
         \\[0.2cm]
         & & \displaystyle 
         \;\;\;\;  \;\;\;\;
         \;\;\;\;  \;\;\;\;
         \;\;\;\;  \;\;\;\;
         \;\;\;\;  \;\;
         +\,
         \norm{\hat\lambda_{t+1}-\hat\lambda_t} + \frac{1}{2}\norm{\lambda_t - \hat\lambda_{t+1}}
         +
         \norm{\hat\xi_{t+1}-\hat\xi_t} + \frac{1}{2}\norm{\xi_t - \hat\xi_{t+1}}
         \Bigg)^2
        \\[0.2cm]
        & \leq & \displaystyle
         \max\left( 
         \frac{\kappa_\rho}{1-\gamma}, 1\right)^2 
         \Bigg(\frac{1}{1-\gamma}
            \sum_s d_\rho^{\pi^\star}(s) \left( \norm{\pi_{t}(\cdot\,\vert\,s) - \hat\pi_t(\cdot\,\vert\,s)} + \frac{3}{2}\norm{\pi_t(\cdot\,\vert\,s)-\hat\pi_{t+1}(\cdot\,\vert\,s)}\right) 
         \\[0.2cm]
         & & \displaystyle 
         \;\;\;\;  \;\;\;\;
         \;\;\;\;  \;\;\;\;
         \;\;\;\;  \;\;\;\;
         \;\;\;\;  \;\;
         +\,
         \norm{\lambda_{t}-\hat\lambda_t} + \frac{3}{2}\norm{\lambda_t - \hat\lambda_{t+1}}
         +
         \norm{\xi_{t}-\hat\xi_t} + \frac{3}{2}\norm{\xi_t - \hat\xi_{t+1}}
         \Bigg)^2
         \\[0.2cm]
        & \leq & \displaystyle 9
         \max\left( 
         \frac{\kappa_\rho}{1-\gamma}, 1\right)^2 
         \Bigg(\frac{1}{1-\gamma}
            \sum_s d_\rho^{\pi^\star}(s) \left( \norm{\pi_{t}(\cdot\,\vert\,s) - \hat\pi_t(\cdot\,\vert\,s)}^2 + \norm{\pi_t(\cdot\,\vert\,s)-\hat\pi_{t+1}(\cdot\,\vert\,s)}^2\right) 
         \\[0.2cm]
         & & \displaystyle 
         \;\;\;\;  \;\;\;\;
         \;\;\;\;  \;\;\;\;
         \;\;\;\;  \;\;\;\;
         \;\;\;\;  \;\;
         +\,
         \norm{\lambda_{t}-\hat\lambda_t}^2 + \norm{\lambda_t - \hat\lambda_{t+1}}^2
         +
         \norm{\xi_{t}-\hat\xi_t}^2 + \norm{\xi_t - \hat\xi_{t+1}}^2
         \Bigg)
    \end{array}
    \]
    where the second inequality is due to triangle inequality and the last inequality is due to relaxing the multiplier, $(x+y)^2\leq 2x^2+2y^2$, and  Jensen’s inequality.
\end{proof}

\begin{lemma}\label{lem:problem constant}
    In Algorithm~\ref{alg: resilient OPG},
    \[
    \begin{array}{rcl}
        && \!\!\!\!  \!\!\!\!  \!\!
        \displaystyle
         \sup_{\pi\,\in\,\Pi, \xi\,\in\,\Xi, \lambda\,\in\,\Lambda} \frac{ \left(
            V_{r+\hat\lambda_{t+1}^\top g}^{\pi}(\rho)
             - h(\xi) - \hat\lambda_{t+1}^\top \xi\right)
             - 
             \left(
            V_{r+\lambda^\top g}^{\hat\pi_{t+1}}(\rho)
             - h(\hat\xi_{t+1}) - \lambda^\top \hat\xi_{t+1}\right) }{
                \max_s \norm{\pi(\cdot\,\vert\,s)-\hat\pi_{t+1}(\cdot\,\vert\,s))}
                +
                \norm{\xi - \hat\xi_{t+1}}
                +
                \norm{\lambda - \hat\lambda_{t+1}}}
        \\[0.2cm]
             & \geq &
             \displaystyle
             C_{\rho,\gamma,\sigma}
             \left(
             \sum_s d_\rho^{\pi^\star}(s) \norm{\hat\pi_{t+1}(\cdot\,\vert\,s) - \mathcal{P}_{\Pi^\star}(\hat\pi_{t+1}(\cdot\,\vert\,s))}
             +
             \norm{\hat\xi_{t+1} - \mathcal{P}_{\Lambda^\star} (\hat\xi_{t+1}) }^2
             +
             \norm{\hat\lambda_{t+1} - \mathcal{P}_{\Lambda^\star}(\hat\lambda_{t+1})}\right)
    \end{array}
    \]    
    where $C_{\rho,\gamma,\sigma}>0$ is some problem-dependent constant.
\end{lemma}
\begin{proof}
    Using the saddle-point property of $(\pi^\star, \xi^\star)$, we first show that for all $\hat\pi_{t+1}$ and $\hat\xi_{t+1}$,
    \begin{equation}\label{eq:optimality of each}
    \begin{array}{rcl}
    \displaystyle 
    \min_{\lambda\,\in\,\Lambda} \;
    \left(
    V_{r+\lambda^\top g}^{\hat\pi_{t+1}} (\rho) - h(\hat\xi_{t+1}) - \lambda^\top \hat\xi_{t+1}
    \right)
    & \leq & \displaystyle
    \frac{1}{2} \min_{\lambda\,\in\,\Lambda} \;
    \left(
    V_{r+\lambda^\top g}^{\pi^\star} (\rho) - h(\hat\xi_{t+1}) - \lambda^\top \hat\xi_{t+1}
    \right) 
    \\[0.2cm]
    && \displaystyle
    +\,\frac{1}{2} \min_{\lambda\,\in\,\Lambda} \;
    \left(
    V_{r+\lambda^\top g}^{\hat\pi_{t+1}} (\rho) - h(\xi^\star) - \lambda^\top \xi^\star
    \right).
    \end{array}
    \end{equation}
    This can be proved using the linear program formulation of MDP, we can express the value function in terms of the occupancy measure, i.e., $V_{r+\lambda^\top g}^\pi(\rho) = \langle q^\pi, r+\lambda^\top g\rangle$, where $q^\pi$ is the occupancy measure that lives in a polytope $\mathcal{Q}$. Let the envelope function be $F(q^\pi, \xi) \DefinedAs \min_{\lambda\,\in\,\Lambda} ( \langle q^\pi, r+\lambda^\top g\rangle - h(\xi) -\lambda^\top \xi)$. We notice that the point-wise minimization is over a family of affine functions $\langle q^\pi, r+\lambda^\top g\rangle -\lambda^\top \xi$ in terms of $(q^\pi, \xi)$, and $h(\xi)$ is strongly convex. Thus, $F(q^\pi,\xi)$ is a concave function,
    \[
        F(q^\pi,\xi)
        \; \leq  \;
        F(q^\pi, \xi^\star)
        + \partial_\xi F(q^\pi, \xi^\star)^\top (\xi - \xi^\star)
        \; = \;
        F(q^\pi, \xi^\star)
        + \partial_\xi F(q^{\pi^\star}, \xi^\star)^\top (\xi - \xi^\star)
    \]
    \[
        F(q^\pi,\xi)
        \; \leq  \;
        F(q^{\pi^\star}, \xi)
        + \partial_q F(q^{\pi^\star}, \xi)^\top (q^{\pi} - q^{\pi^\star})
        \; =  \;
        F(q^{\pi^\star}, \xi)
        + \partial_q F(q^{\pi^\star}, \xi^\star)^\top (q^{\pi} - q^{\pi^\star})
    \]
    where $q^{\pi^\star}$ corresponds to $\pi^\star$ in the one-to-one way, and we notice that $\partial_\xi F(q^\pi, \xi^\star)  = \partial_\xi F(q^{\pi^\star}, \xi^\star)$ and $\partial_q F(q^{\pi^\star}, \xi^\star) = \partial_q F(q^{\pi^\star}, \xi)$ because of the decoupled structure of $q^\pi$ and $\xi$. 
    From the saddle-point property of $(q^{\pi^\star},\xi^\star)$, it also reaches the maximum of $F(q^\pi, \xi)$. Thus, by the optimality of $(q^{\pi^\star},\xi^\star)$,
    \[
    F(q^\pi,\xi) 
    \; \leq\;
    \frac{1}{2} 
    \left(
    F(q^\pi, \xi^\star)
    +
    F(q^{\pi^\star}, \xi) 
    \right) \; \text{ for all } q^\pi \in \mathcal{Q} \text{ and } \xi \in \Xi
    \]
    which proves~\eqref{eq:optimality of each} by taking $\pi = \hat\pi_{t+1}$ and $\xi = \hat\xi_{t+1}$. 
    
    Denote $V_h^\star \DefinedAs V_{r+(\lambda^\star)^\top g}^\star(\rho) - h(\xi^\star) - (\lambda^\star)^\top \xi^\star$ and \[
        D_{\max} 
        \; \DefinedAs \;
        \max_{\pi,\pi' \, \in\, \Pi, \, \xi, \xi'\, \in \, \Xi, \, \lambda, \lambda'\,\in\, \Lambda}
        \left(
            \max_s \norm{ \pi(\cdot\,\vert\,s) - \pi'(\cdot\,\vert\,s) }
            + 
            \norm{\xi - \xi'}
            +
            \norm{\lambda - \lambda'}
        \right).
    \]
    Thus,
    \[
        \begin{array}{rcl}
             &  &
             \!\!\!\!  \!\!\!\!  \!\!
            \displaystyle
         \sup_{\pi\,\in\,\Pi, \xi\,\in\,\Xi, \lambda\,\in\,\Lambda} \frac{ \left(
            V_{r+\hat\lambda_{t+1}^\top g}^{\pi}(\rho)
             - h(\xi) - \hat\lambda_{t+1}^\top \xi\right)
             - 
             \left(
            V_{r+\lambda^\top g}^{\hat\pi_{t+1}}(\rho)
             - h(\hat\xi_{t+1}) - \lambda^\top \hat\xi_{t+1}\right) }{
                \max_s \norm{\pi(\cdot\,\vert\,s)-\hat\pi_{t+1}(\cdot\,\vert\,s))}
                +
                \norm{\xi - \hat\xi_{t+1}}
                +
                \norm{\lambda - \hat\lambda_{t+1}}}
             \\[0.2cm]
             & \geq &
             \displaystyle
             \frac{1}{D_{\max}}
             \sup_{\pi\,\in\,\Pi, \xi\,\in\,\Xi, \lambda\,\in\,\Lambda} 
             \left(
             \left(
            V_{r+\hat\lambda_{t+1}^\top g}^{\pi}(\rho)
             - h(\xi) - \hat\lambda_{t+1}^\top \xi\right)
             - 
             \left(
            V_{r+\lambda^\top g}^{\hat\pi_{t+1}}(\rho)
             - h(\hat\xi_{t+1}) - \lambda^\top \hat\xi_{t+1}\right) \right)
             \\[0.2cm]
             & = &
             \displaystyle
             \frac{1}{2D_{\max}}
             \left(
             \sup_{\pi\,\in\,\Pi} 
             \left(
            V_{r+\hat\lambda_{t+1}^\top g}^{\pi}(\rho)
             - h(\xi^\star) - \hat\lambda_{t+1}^\top \xi^\star\right)
             - 
             \inf_{ \lambda\,\in\,\Lambda}
             \left(
            V_{r+\lambda^\top g}^{\hat\pi_{t+1}}(\rho)
             - h(\hat\xi_{t+1}) - \lambda^\top \hat\xi_{t+1}\right) \right)
             \\[0.2cm]
             &  &
             \displaystyle +\,
             \frac{1}{2D_{\max}}
             \left(
             \sup_{\xi\,\in\,\Xi} 
             \left(
            V_{r+\hat\lambda_{t+1}^\top g}^{\pi^\star}(\rho)
             - h(\xi) - \hat\lambda_{t+1}^\top \xi\right)
             - 
             \inf_{ \lambda\,\in\,\Lambda}
             \left(
            V_{r+\lambda^\top g}^{\hat\pi_{t+1}}(\rho)
             - h(\hat\xi_{t+1}) - \lambda^\top \hat\xi_{t+1}\right) \right)
             \\[0.2cm]
             & \geq &
             \displaystyle
             \frac{1}{2D_{\max}}
             \left(
             \sup_{\pi\,\in\,\Pi} 
             \left(
            V_{r+\hat\lambda_{t+1}^\top g}^{\pi}(\rho)
             - h(\xi^\star) - \hat\lambda_{t+1}^\top \xi^\star\right)
             - 
             \inf_{ \lambda\,\in\,\Lambda}
             \left(
            V_{r+\lambda^\top g}^{\hat\pi_{t+1}}(\rho)
             - h(\xi^\star) - \lambda^\top \xi^\star\right) \right)
             \\[0.2cm]
             &  &
             \displaystyle +\,
             \frac{1}{2D_{\max}}
             \left(
             \sup_{\xi\,\in\,\Xi} 
             \left(
            V_{r+\hat\lambda_{t+1}^\top g}^{\pi^\star}(\rho)
             - h(\xi) - \hat\lambda_{t+1}^\top \xi\right)
             - 
             \inf_{ \lambda\,\in\,\Lambda}
             \left(
            V_{r+\lambda^\top g}^{\pi^\star}(\rho)
             - h(\hat\xi_{t+1}) - \lambda^\top \hat\xi_{t+1}\right) \right)
        \end{array}
    \]
    where the first inequality is due to the domain's diameter, and the second inequality is due to~\eqref{eq:optimality of each}.

    Denote $V_h^\star \DefinedAs V_{r+(\lambda^\star)^\top g}^{\pi^\star} (\rho) - h(\xi^\star) - (\lambda^\star)^\top \xi^\star$. If we can prove that there exist constants $c_1>0$, $c_2>0$, and $c_3>0$ such that
    \begin{subequations}
    \begin{equation}\label{eq:max_pi}
    \max_{\pi\,\in\,\Pi} 
    \left(
    V_{r+\hat\lambda_{t+1}^\top g}^\pi (\rho) - h(\xi^\star) - \hat\lambda_{t+1}^\top \xi^\star
    \right)
    - V_h^\star
    \; \geq \; 
    c_1 \norm{\hat\lambda_{t+1} - \mathcal{P}_{\Lambda^\star}(\hat\lambda_{t+1})}
    \end{equation}
    \begin{equation}\label{eq:min_lambda}
    V_h^\star
    -
    \min_{\lambda\,\in\,\Lambda} 
    \left(
    V_{r+\lambda^\top g}^{\hat\pi_{t+1}} (\rho) - h(\xi^\star) - \lambda^\top \xi^\star
    \right)
    \; \geq \; 
    c_2\sum_s d_\rho^{\pi^\star}(s) \norm{\hat\pi_{t+1}(\cdot\,\vert\,s) - \mathcal{P}_{\Pi^\star}(\hat\pi_{t+1}(\cdot\,\vert\,s))}
    \end{equation}
    \begin{equation}\label{eq:max_min_xi_lambda}
    \sup_{\xi\,\in\,\Xi} 
             \left(
            V_{r+\hat\lambda_{t+1}^\top g}^{\pi^\star}(\rho)
             - h(\xi) - \hat\lambda_{t+1}^\top \xi\right)
             - 
             \inf_{ \lambda\,\in\,\Lambda}
             \left(
            V_{r+\lambda^\top g}^{\pi^\star}(\rho)
             - h(\hat\xi_{t+1}) - \lambda^\top \hat\xi_{t+1}\right)
             \; \geq \;
             c_3 \norm{\hat\xi_{t+1} - \mathcal{P}_{\Lambda^\star} (\hat\xi_{t+1}) }^2
    \end{equation}
    \end{subequations}
    then,
    \[
        \begin{array}{rcl}
             &  &
             \!\!\!\!  \!\!\!\!  \!\!
            \displaystyle
         \sup_{\pi\,\in\,\Pi, \xi\,\in\,\Xi, \lambda\,\in\,\Lambda} \frac{ \left(
            V_{r+\hat\lambda_{t+1}^\top g}^{\pi}(\rho)
             - h(\xi) - \hat\lambda_{t+1}^\top \xi\right)
             - 
             \left(
            V_{r+\lambda^\top g}^{\hat\pi_{t+1}}(\rho)
             - h(\hat\xi_{t+1}) - \lambda^\top \hat\xi_{t+1}\right) }{
                \max_s \norm{\pi(\cdot\,\vert\,s)-\hat\pi_{t+1}(\cdot\,\vert\,s))}
                +
                \norm{\xi - \hat\xi_{t+1}}
                +
                \norm{\lambda - \hat\lambda_{t+1}}}
             \\[0.2cm]
             & \geq &
             \displaystyle
             \frac{\min(c_1, c_2, c_3)}{2D_{\max}}
             \left(
             \sum_s d_\rho^{\pi^\star}(s) \norm{\hat\pi_{t+1}(\cdot\,\vert\,s) - \mathcal{P}_{\Pi^\star}(\hat\pi_{t+1}(\cdot\,\vert\,s))}
             +
             \norm{\hat\xi_{t+1} - \mathcal{P}_{\Lambda^\star} (\hat\xi_{t+1}) }^2
             +
             \norm{\hat\lambda_{t+1} - \mathcal{P}_{\Lambda^\star}(\hat\lambda_{t+1})}\right)
        \end{array}
    \]
    which proves the desired inequality by taking $C_{\rho,\gamma,\sigma} = \frac{\min(c_1, c_2, c_3)}{2D_{\max}}$.

    To complete the proof, we first show~\eqref{eq:max_pi} and~\eqref{eq:min_lambda} by resorting the following partial bilinear saddle-point problem,
    \[
    \maximize_{q^\pi \,\in \, \mathcal{Q}}\; \minimize_{\lambda\,\in\,\Lambda} \; \langle q^\pi, r+\lambda^\top g\rangle - h(\xi^\star) - \lambda^\top \xi^\star
    \; = \;
     \minimize_{\lambda\,\in\,\Lambda} \; \maximize_{q^\pi \,\in \, \mathcal{Q}}\;\langle q^\pi, r+\lambda^\top g\rangle - h(\xi^\star) - \lambda^\top \xi^\star.
    \]
    Hence, we can apply \cite[Lemma~10]{ding2023last} to obtain~\eqref{eq:max_pi} and~\eqref{eq:min_lambda} directly. To see~\eqref{eq:max_min_xi_lambda}, we consider the following partial saddle-point problem,
    \[
    \maximize_{\xi \,\in \, \Xi}\; \minimize_{\lambda\,\in\,\Lambda} \; \langle q^{\pi^\star}, r+\lambda^\top g\rangle - h(\xi) - \lambda^\top \xi
    \; = \;
     \minimize_{\lambda\,\in\,\Lambda} \; \maximize_{\xi \,\in \, \Xi}\;\langle q^{\pi^\star}, r+\lambda^\top g\rangle - h(\xi) - \lambda^\top \xi
    \]
    which has a strongly concave and linear saddle-point objective function.
    
    We notice that $V_h^\star = \langle q^{\pi^\star}, r+(\lambda^\star)^\top g\rangle - h(\xi^\star) - (\lambda^\star)^\top \xi^\star$.
    It is straightforward to verify that
    \[
    \begin{array}{rcl}
         &  & \displaystyle
         \!\!\!\!  \!\!\!\!  \!\!
         \max_{\xi\, \in\, \Xi} \left( 
         \langle q^{\pi^\star}, r+\lambda^\top g\rangle - h(\xi) - \lambda^\top \xi
         \right) -
         \min_{\lambda \, \in\, \Lambda} \left( 
         \langle q^{\pi^\star}, r+\lambda^\top g\rangle - h(\xi) - \lambda^\top \xi
         \right)
         \\[0.2cm]
         & \geq & \displaystyle
         \left( 
         \langle q^{\pi^\star}, r+\lambda^\top g\rangle - h(\xi^\star) - \lambda^\top \xi^\star
         \right) - V_h^\star
         \\[0.2cm]
         && \displaystyle +\,
         V_h^\star -
         \left( 
         \langle q^{\pi^\star}, r+(\lambda^\star)^\top g\rangle - h(\xi) - (\lambda^\star)^\top \xi
         \right)
         \\[0.2cm]
         & \geq & \displaystyle
         \frac{\sigma}{2}\norm{\xi-\xi^\star}^2
    \end{array}
    \]
    where we fix $\xi = \xi^\star$  and $\lambda = \lambda^\star$ for the first inequality and the second inequality is from the strong convexity and the optimality of $(\xi^\star,\lambda^\star)$,
    \[
    \langle q^{\pi^\star}, r+(\lambda^\star)^\top g\rangle - h(\xi) - (\lambda^\star)^\top \xi
    \; \leq \;
    V_h^\star + (\nabla h(\xi^\star)-\lambda^\star)^\top (\xi-\xi^\star)- \frac{\sigma}{2} \norm{\xi - \xi^\star}^2
    \]
    \[
    \langle q^{\pi^\star}, r+\lambda^\top g\rangle - h(\xi^\star) - \lambda^\top \xi^\star
    \; \geq \;V_h^\star + (\langle q^{\pi^\star}, g\rangle-\xi^\star)^\top (\lambda- \lambda^\star).
    \]
    
    Finally, we notice that replacing $\xi^\star$ by $\mathcal{P}_{\Xi^\star}(\xi)$ does not alter the argument above, which proves~\eqref{eq:max_min_xi_lambda} by taking $c_3 = \frac{\sigma}{2}$.
\end{proof}

\begin{proof}
    By the non-increasing sequence in Lemma~\ref{lem:non-increasing} and the definition of $\zeta_t$, we have
    \[
    \displaystyle
         \,\frac{1}{2(1-\gamma)} \sum_{s}d_\rho^{\pi^\star}(s) 
             \norm{\hat\pi_{t+1}(\cdot\,\vert\,s) - \pi_t(\cdot\,\vert\,s)}^2
            +
        \frac{1}{2}
            \norm{\hat\xi_{t+1}-\xi_t}^2 +
         \frac{1}{2}  
            \norm{\hat\lambda_{t+1}-\lambda_t}^2
    \; \leq \;
    \zeta_t 
    \; \leq \;
    2\Theta_t
    \; \leq \;
    2\Theta_1.
    \]
    
    Meanwhile, by the definition of $\zeta_t$ and Lemma~\ref{lem:lower_duality_gap}, we have 
    \[
    \begin{array}{rcl}
     \zeta_t & \geq &
     \displaystyle
         \frac{1}{4(1-\gamma)}\sum_{s}d_\rho^{\pi^\star}(s) 
             \norm{\hat\pi_{t+1}(\cdot\,\vert\,s) - \pi_t(\cdot\,\vert\,s)}^2
            +
            \frac{1}{4(1-\gamma)}\norm{\hat\xi_{t+1}-\xi_t}^2 
            +  
            \frac{1}{4}\norm{\hat\lambda_{t+1}-\lambda_t}^2
     \\[0.2cm]
     & &
     \displaystyle +\,
         \frac{1}{4(1-\gamma)}\sum_{s}d_\rho^{\pi^\star}(s) 
             \norm{\hat\pi_{t+1}(\cdot\,\vert\,s) - \pi_t(\cdot\,\vert\,s)}^2
            +
            \frac{1}{4}\norm{\hat\xi_{t+1}-\xi_t}^2 
            +  
            \frac{1}{4}\norm{\hat\lambda_{t+1}-\lambda_t}^2
         \\[0.2cm]
         & &  \displaystyle
         +\, \frac{1}{4(1-\gamma)}\sum_{s}d_\rho^{\pi^\star}(s) 
             \norm{\pi_{t}(\cdot\,\vert\,s) - \hat\pi_t(\cdot\,\vert\,s)}^2
            +
            \frac{1}{4}\norm{\xi_{t}-\hat\xi_t}^2 
            +  
            \frac{1}{4}\norm{\lambda_{t}
            -\hat\lambda_t}^2
        \\[0.2cm]
         & \geq &  
        \displaystyle
         \frac{1}{4(1-\gamma)}\sum_{s}d_\rho^{\pi^\star}(s) 
             \norm{\hat\pi_{t+1}(\cdot\,\vert\,s) - \pi_t(\cdot\,\vert\,s)}^2
            +
            \frac{1}{4}\norm{\hat\xi_{t+1}-\xi_t}^2 
            +  
            \frac{1}{4}\norm{\hat\lambda_{t+1}-\lambda_t}^2
            \\[0.2cm]
         &  &
         \displaystyle +\,
         \frac{\eta^2}{32 \max\left( \frac{\kappa_\rho}{1-\gamma},1\right)^2} \frac{\left[ \left(
            V_{r+\hat\lambda_{t+1}^\top g}^{\pi}(\rho)
             - h(\xi) - \hat\lambda_{t+1}^\top \xi\right)
             - 
             \left(
            V_{r+\lambda^\top g}^{\hat\pi_{t+1}}(\rho)
             - h(\hat\xi_{t+1}) - \lambda^\top \hat\xi_{t+1}\right) \right]_+^2}{\left(
                \max_s \norm{\pi(\cdot\,\vert\,s)-\hat\pi_{t+1}(\cdot\,\vert\,s))}
                +
                \norm{\xi - \hat\xi_{t+1}}
                +
                \norm{\lambda - \hat\lambda_{t+1}}
            \right)^2}
            \\[0.2cm]
             & \geq &
            \displaystyle
         \frac{1}{4(1-\gamma)}\sum_{s}d_\rho^{\pi^\star}(s) 
             \norm{\hat\pi_{t+1}(\cdot\,\vert\,s) - \pi_t(\cdot\,\vert\,s)}^2
            +
            \frac{1}{4}\norm{\hat\xi_{t+1}-\xi_t}^2 
            +  
            \frac{1}{4}\norm{\hat\lambda_{t+1}-\lambda_t}^2
            \\[0.2cm]
             && \displaystyle +\,
         \frac{\eta^2 C_{\rho,\gamma,\sigma}}{32 \max\left( \frac{\kappa_\rho}{1-\gamma},1\right)^2} 
             \Bigg(
             \sum_s d_\rho^{\pi^\star}(s) \norm{\hat\pi_{t+1}(\cdot\,\vert\,s) - \mathcal{P}_{\Pi^\star}(\hat\pi_{t+1}(\cdot\,\vert\,s))}^2
             \\[0.2cm]
             &  &
             \displaystyle
             \;\;\;\;  \;\;\;\;
             \;\;\;\;  \;\;\;\;
             \;\;\;\;  \;\;\;\;
             \;\;\;\;  \;\;\;\;
             +\,
             {\norm{\hat\xi_{t+1} - \mathcal{P}_{\Lambda^\star} (\hat\xi_{t+1}) }^4}
             +
             {\norm{\hat\lambda_{t+1} - \mathcal{P}_{\Lambda^\star}(\hat\lambda_{t+1})}^2}\Bigg)
            \\[0.2cm]
             & \gtrsim &
            \displaystyle
         \frac{1}{4(1-\gamma)}\sum_{s}d_\rho^{\pi^\star}(s) 
             \norm{\hat\pi_{t+1}(\cdot\,\vert\,s) - \pi_t(\cdot\,\vert\,s)}^4
            +
            \frac{1}{4}\norm{\hat\xi_{t+1}-\xi_t}^4
            +  
            \frac{1}{4}\norm{\hat\lambda_{t+1}-\lambda_t}^4
            \\[0.2cm]
             && \displaystyle +\,
         \frac{\eta^2 C_{\rho,\gamma,\sigma}}{32 \max\left( \frac{\kappa_\rho}{1-\gamma},1\right)^2} 
             \Bigg(
             \sum_s d_\rho^{\pi^\star}(s) \norm{\hat\pi_{t+1}(\cdot\,\vert\,s) - \mathcal{P}_{\Pi^\star}(\hat\pi_{t+1}(\cdot\,\vert\,s))}^4
             \\[0.2cm]
             &  &
             \displaystyle
             \;\;\;\;  \;\;\;\;
             \;\;\;\;  \;\;\;\;
             \;\;\;\;  \;\;\;\;
             \;\;\;\;  \;\;\;\;
             +\,
             \norm{\hat\xi_{t+1} - \mathcal{P}_{\Lambda^\star} (\hat\xi_{t+1}) }^4
             +
             \norm{\hat\lambda_{t+1} - \mathcal{P}_{\Lambda^\star}(\hat\lambda_{t+1})}^4\Bigg)
             \\[0.2cm]
             & \geq &
             \displaystyle \min\left( \frac{1}{4},
         \frac{\eta^2 C_{\rho,\gamma,\sigma}(1-\gamma)}{32 \max\left( \frac{\kappa_\rho}{1-\gamma},1\right)^2} \right)
             \Bigg( \frac{1}{1-\gamma}
                \sum_{s}d_\rho^{\pi^\star}(s) 
             \norm{\hat\pi_{t+1}(\cdot\,\vert\,s) - \pi_t(\cdot\,\vert\,s)}^2
            +
            \norm{\hat\xi_{t+1}-\xi_t}^2 
            +  
            \norm{\hat\lambda_{t+1}-\lambda_t}^2
             \\[0.2cm]
             &  &
             \displaystyle
             +\, \frac{1}{1-\gamma}
             \sum_s d_\rho^{\pi^\star}(s) \norm{\hat\pi_{t+1}(\cdot\,\vert\,s) - \mathcal{P}_{\Pi^\star}(\hat\pi_{t+1}(\cdot\,\vert\,s))}^2
             +
             \norm{\hat\xi_{t+1} - \mathcal{P}_{\Lambda^\star} (\hat\xi_{t+1}) }^2
             +
             \norm{\hat\lambda_{t+1} - \mathcal{P}_{\Lambda^\star}(\hat\lambda_{t+1})}^2\Bigg)^2
             \\[0.2cm]
             & \geq &
             \displaystyle \min\left( \frac{1}{4},
         \frac{\eta^2 C_{\rho,\gamma,\sigma}(1-\gamma)}{32 \max\left( \frac{\kappa_\rho}{1-\gamma},1\right)^2} \right)
         \Theta_{t+1}^2
    \end{array}
\]
where $\gtrsim$ means $\geq$ up to some normalization constants for $(\xi,\lambda)$-relevant terms that can be normalized to one due to the boundedness.

Denote $C_\eta \DefinedAs \min\left( \frac{1}{4},
         \frac{\eta^2 C_{\rho,\gamma,\sigma}(1-\gamma)}{32 \max\left( \frac{\kappa_\rho}{1-\gamma},1\right)^2} \right)$. Thus,
         \[
         \Theta_{t+1}
         \; \leq \;
         \Theta_t - C_\eta \Theta_{t+1}^2.
         \]
         By Lemma~\ref{lem:optimistic_last_rate}, we have
         \[
         \Theta_t \;=\; O\left(\frac{1}{t}\right)
         \]
         where the stepsize $\eta$ satisfies  
         \[
            \eta \;\leq\; 
            \min\left( \frac{1}{4\sqrt{|A|}}, \frac{1}{2(L_h+1)}, \frac{1}{5 \sqrt{m|A|}\kappa_\rho}, \frac{\rho_{\min}}{4 \gamma\sqrt{m}C_h |A|}, \frac{1}{2\sqrt{2\iota}}, \frac{4 \max(\frac{\kappa_\rho}{1-\gamma}, 1)}{\sqrt{ \Theta_1 C_{\rho,\gamma,\sigma}(1-\gamma)}}
    \right)  \;\DefinedAs \; \eta_{\max}.
         \]
\end{proof}

\subsection{Proof of Corollary~\ref{cor:policy and relaxation}}
\label{app:policy and relaxation}
    According to Theorem~\ref{thm: last-iterate convergence}, if $t = \Omega(1/\epsilon)$, then,
    \[
    \begin{array}{rcl}
         \displaystyle
         \frac{1}{2} \sum_{s}d_\rho^{\pi^\star}(s) \norm{\mathcal{P}_{\Pi^\star}( \hat\pi_{t}(\cdot\,\vert\,s)) - \hat\pi_{t}(\cdot\,\vert\,s)}^2
         & = & O(\epsilon)
         \\[0.2cm]
        \displaystyle
         \frac{1}{2}\norm{\mathcal{P}_{\Xi^\star}(\hat\xi_{t})-\hat\xi_{t}}^2 
         & = & O(\epsilon)
         \\[0.2cm]
         \displaystyle
         \frac{1}{2}  
            \norm{\mathcal{P}_{\Lambda^\star}(\hat\lambda_{t})-\hat\lambda_{t}}^2
        & = & O(\epsilon).
    \end{array}
\]
Let $\hat\pi_{t}^\star(\cdot\,\vert\,s)\DefinedAs\mathcal{P}_{\Pi^\star}( \hat\pi_{t}(\cdot\,\vert\,s))$, $\hat\xi_{t}^\star\DefinedAs\mathcal{P}_{\Xi^\star}(\hat\xi_{t})$, and $\hat\lambda_{t}^\star \DefinedAs\mathcal{P}_{\Lambda^\star}(\hat\lambda_{t})$. Because of the interchangeability of saddle points, $(\hat\pi_{t}^\star, \hat\xi_{t}^\star, \hat\lambda_{t}^\star)$ is a saddle point in $\Pi^\star\times\Xi^\star\times\Lambda^\star$. 

First, we have
\[
\begin{array}{rcl}
     V_r^{\star}(\rho) - h(\bar\xi^\star) - (V_r^{\hat\pi_t}(\rho) - h(\hat\xi_t))  & = &
     \displaystyle \frac{1}{1-\gamma}\sum_{s,a} d_\rho^{\hat{\pi}_{t+1}}(s)(\hat{\pi}_{t}^\star(a\,\vert\,s) - \hat\pi_t(a\,\vert\,s
     )) Q_r^{\hat\pi_t}(s,a) - h(\bar\xi^\star) + h(\hat\xi_t)
     \\[0.2cm]
     & \leq &  \displaystyle\frac{1}{(1-\gamma)^2}\sum_{s} d_\rho^{\hat{\pi}_{t+1}}(s)\norm{\hat{\pi}_{t}^\star(\cdot\,\vert\,s) - \hat\pi_t(\cdot\,\vert\,s
     ) }_1 - h(\bar\xi^\star) + h(\hat\xi_t) 
     \\[0.2cm]
     & \leq &  \displaystyle\frac{\sqrt{|A|}}{(1-\gamma)^2}\sum_{s} d_\rho^{\hat{\pi}_{t+1}}(s)\norm{\hat{\pi}_{t}^\star(\cdot\,\vert\,s) - \hat\pi_t(\cdot\,\vert\,s
     ) } - h(\bar\xi^\star) + h(\hat\xi_t) 
     \\[0.2cm]
     & \leq &  \displaystyle\frac{\sqrt{|A|}}{(1-\gamma)^2}\sqrt{\sum_{s} d_\rho^{\hat{\pi}_{t+1}}(s)\norm{\hat{\pi}_{t}^\star(\cdot\,\vert\,s) - \hat\pi_t(\cdot\,\vert\,s
     ) }^2} - h(\bar\xi^\star) + h(\hat\xi_t) 
\end{array}
\]
where we use Cauchy–Schwarz inequality in the first and third inequalities, and the second inequality is due to that $\norm{x}_1\leq \sqrt{d}\norm{x}_2$ for $x\in\mathbb{R}^d$. We note that $h$ is continuous. Thus, 
\[
    V_r^{\star}(\rho) - h(\bar\xi^\star) - (V_r^{\hat\pi_t}(\rho) - h(\hat\xi_t)) \; =\; O(\sqrt{\epsilon})  
\]
where $V_r^\star(\rho) = V_r^{\hat\pi_t^\star}(\rho)$, and $|h(\bar\xi^\star) - h(\hat\xi_t)| \leq \sqrt{\epsilon}$ for small $\epsilon$. 

Second, we have 
\[
    \begin{array}{rcl}
         \norm{\hat\xi_{t} - V_{g}^{\hat\pi_t}(\rho)}
        &\leq&
        \norm{V_{g}^{\hat\pi_t^\star}(\rho)- \hat\xi_{t}^\star+
    \hat\xi_{t} - V_{g}^{\hat\pi_t}(\rho)}
         \\[0.2cm]
         & = &   
    \norm{\hat\xi_{t} -\hat\xi_{t}^\star}+ \norm{V_{g}^{\hat\pi_t^\star}(\rho)- V_{g}^{\hat\pi_t}(\rho)}
    \\[0.2cm]
         & \leq & O(\sqrt{\epsilon})
    \end{array}
\]
where $|V_{g_i}^{\hat\pi_t^\star}(\rho)- V_{g_i}^{\hat\pi_t}(\rho)| = O(\sqrt{\epsilon})$ is similar as we did for the reward value function. 

Finally, we replace $\sqrt{\epsilon}$ by $\epsilon$ and combine big $O$ notation to complete the proof. 

\section{Experiment Setup and Additional Results}
\label{app: experiments}

We provide details of our experiments and additional results. We conduct this experiment on Google Colab in Jupyter Notebook.

\subsection{Resilient Policy Search}\label{sec:policysearch}

In this experiment, we consider a tabular constrained MDP with a randomly generated transition kernel, a discount factor $\gamma=0.9$, uniform rewards $r\in[0,1]$ and utilities $g\in[-1,1]$, and a uniform initial state distribution $\rho$. The relaxation cost function is $h(\xi) = \alpha\xi^2$, where parameter $\alpha$ balances the closeness to original constraints and the reward maximization objective. The initial constraint is $V_g^\pi(\rho)\geq c$, where $c$ is some constant that often makes the problem infeasible. For comparison, we solve a linear program in occupancy-measure space to find the maximal utility value, which is the minimum value of $c$ to make this problem infeasible. Then, we solve a quadratic program to find the maximal value $V_{r}^{\pi^\star}(\rho)-h(\xi)$ at an optimal policy $\pi^\star$. Throughout all experiments, the random seed is fixed, and the minimal $c$ to make the problem infeasible is $5.56$. 

\begin{figure}[tbh]
    \centering
        \begin{tikzpicture}
        \node[](img1) at(0,0) {\includegraphics[width = 0.33\textwidth]{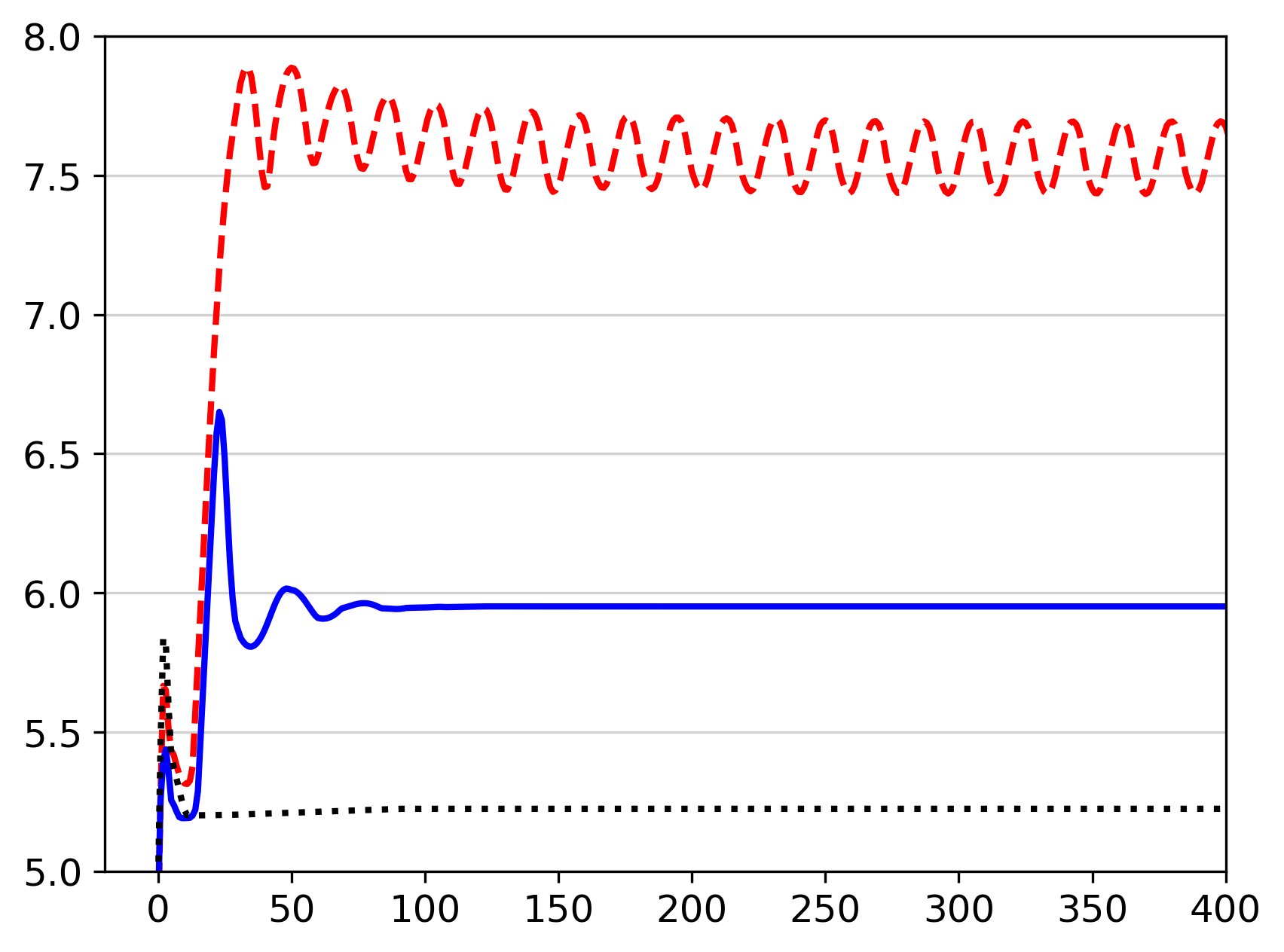}};
        \node[below = 0cm of img1]{iteration};
         \node[left = 0.4cm of img1,yshift=1.3cm, rotate=90]{reward value};
    \end{tikzpicture}
    \hspace{0.1\textwidth}
            \begin{tikzpicture}
        \node[](img1) at(0,0) {\includegraphics[width = 0.33\textwidth]{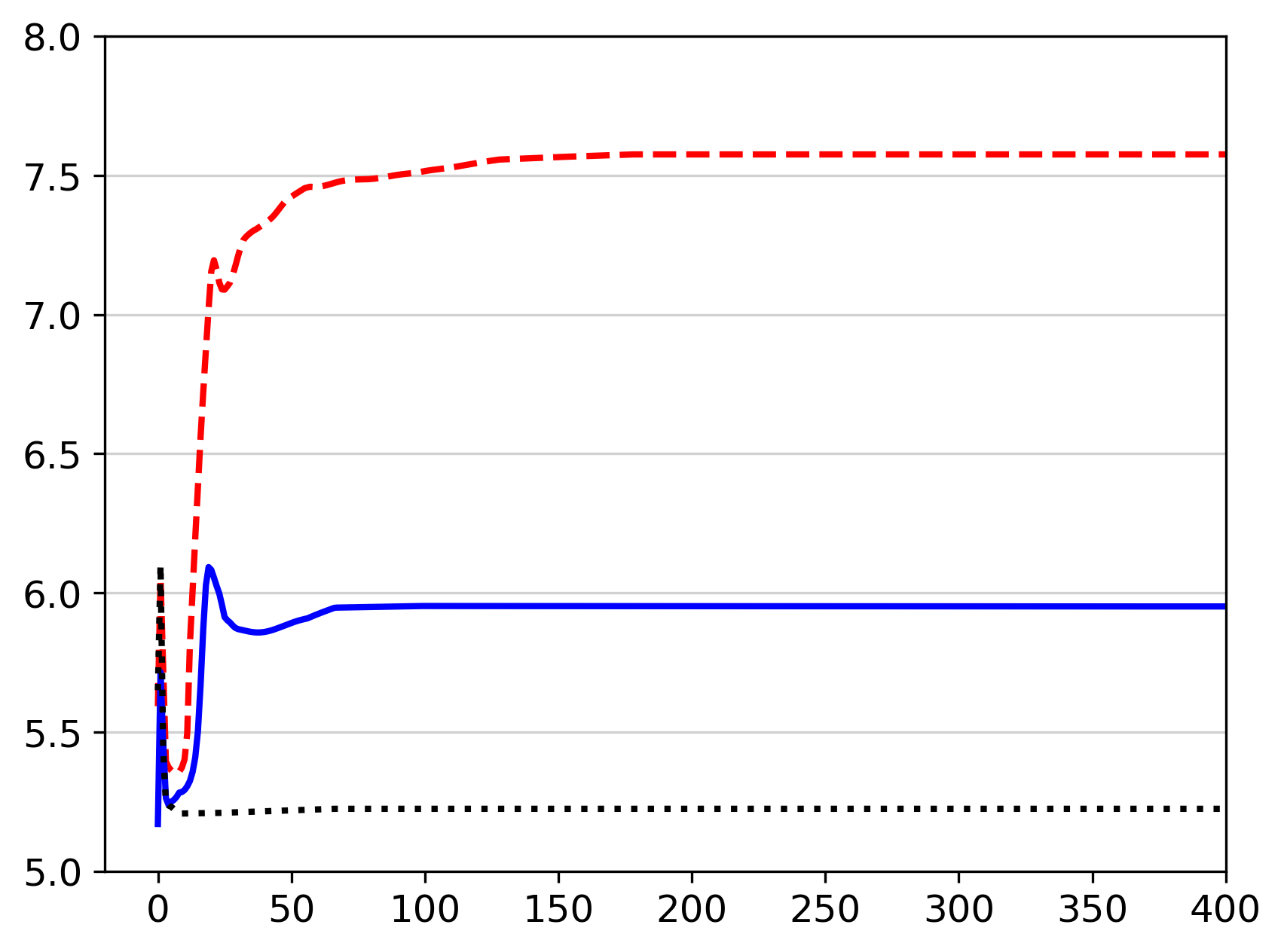}};
        \node[below = 0cm of img1]{iteration};
    \end{tikzpicture}
    \caption{Reward value convergence of ResPG-PD (Algorithm~\ref{alg: resilient PG}, left) and ResOPG-PD (Algorithm~\ref{alg: resilient OPG}, right), with different cost functions $h(\xi) = \alpha \xi^2$, where $\alpha = 0.03$ (\ref{legend:red}), $\alpha = 0.2$ (\ref{legend:blue}), $\alpha = 1$ (\ref{legend:black}), and stepsize $\eta=0.2$ in the policy search problem of Section~\ref{sec:policysearch}.
    }
    \label{fig:AppE1RewardvsIter}
    \end{figure}

\begin{figure}[tbh]
    \centering
        \begin{tikzpicture}
        \node[](img1) at(0,0) {\includegraphics[width = 0.33\textwidth]{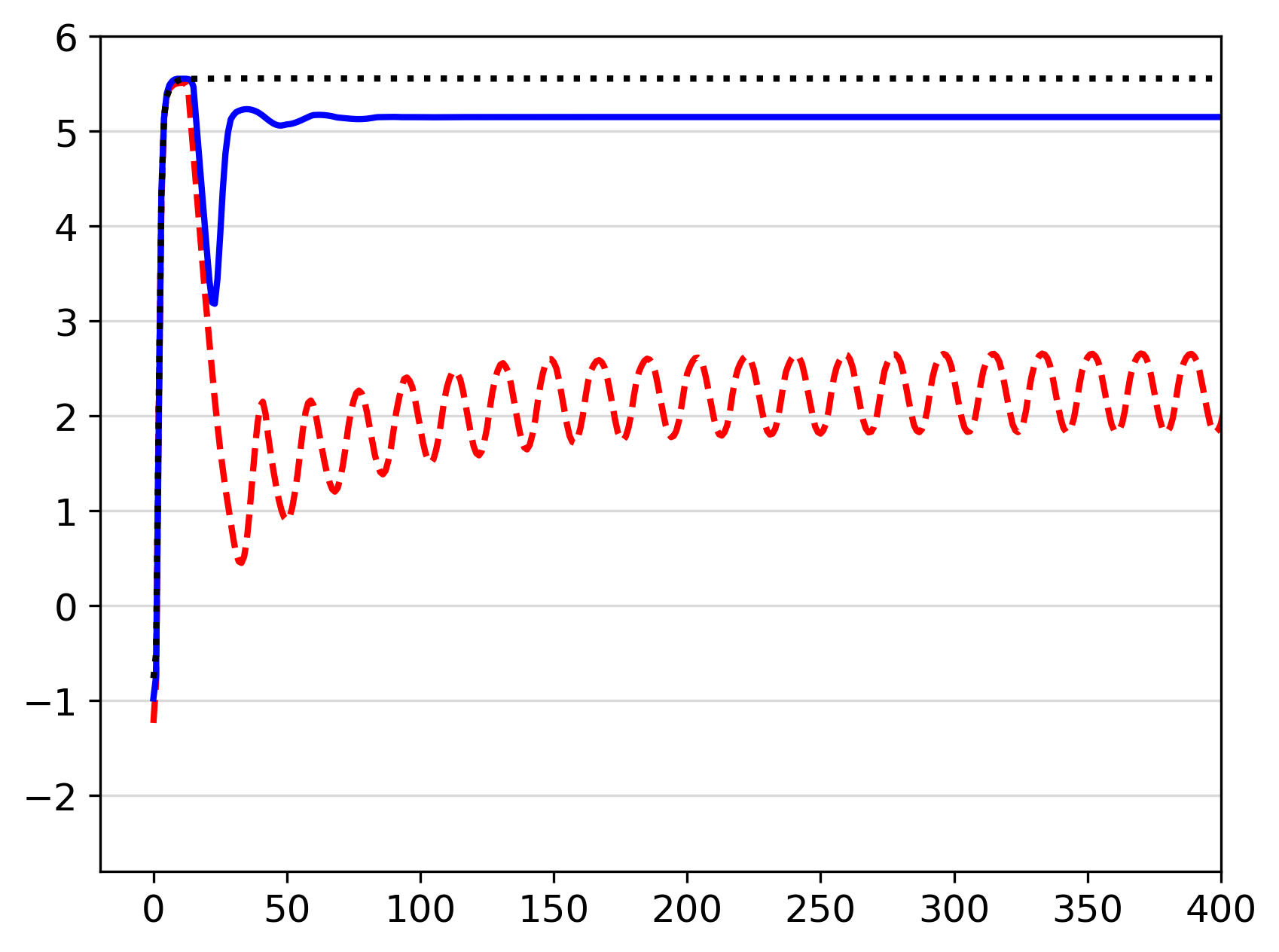}};
        \node[below = 0cm of img1]{iteration};
         \node[left = 0.4cm of img1,yshift=1.3cm, rotate=90]{utility value};
    \end{tikzpicture}
    \hspace{0.1\textwidth}
            \begin{tikzpicture}
        \node[](img1) at(0,0) {\includegraphics[width = 0.33\textwidth]{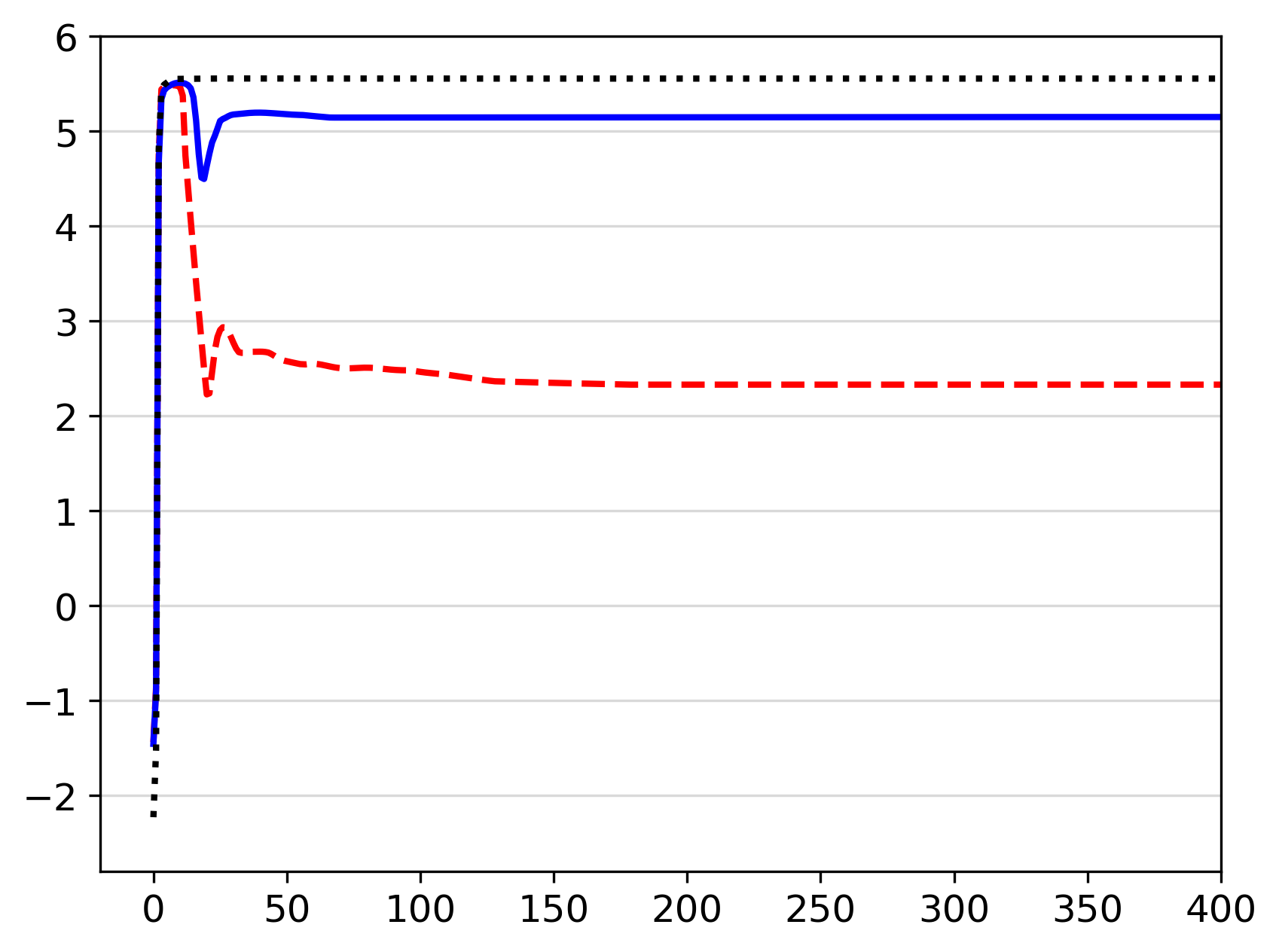}};
        \node[below = 0cm of img1]{iteration};
    \end{tikzpicture}
    \caption{Utility value convergence of ResPG-PD (Algorithm~\ref{alg: resilient PG}, left) and ResOPG-PD (Algorithm~\ref{alg: resilient OPG}, right), with different cost functions $h(\xi) = \alpha \xi^2$, where $\alpha = 0.03$ (\ref{legend:red}), $\alpha = 0.2$ (\ref{legend:blue}), $\alpha = 1$ (\ref{legend:black}), and stepsize $\eta=0.2$ in the policy search problem of Section~\ref{sec:policysearch}.
    }
    \label{fig:AppE1UtilityvsIter}
    \end{figure}

As an example, we take $c=8$. We report the reward and utility value convergence of our two methods: Algorithm~\ref{alg: resilient PG} (ResPG-PD) and Algorithm~\ref{alg: resilient OPG} (ResOPG-PD) in Figures~\ref{fig:AppE1RewardvsIter}--\ref{fig:AppE1UtilityvsIter}, which show similar convergence performance as the relaxation in Figure~\ref{fig:E1XivsIter}. We note, for large $\alpha$, ResPG-PD behaves similar as ResOPG-PD, which is mainly due to the strongly convex cost function $h$. However, when $\alpha$ is small, ResPG-PD starts to oscillate while ResOPG-PD still converges, which demonstrates the advantage of the last-iterate convergence in Theorem~\ref{thm: last-iterate convergence}.

\begin{figure}[tbh]
    \centering
        \begin{tikzpicture}
        \node[](img1) at(0,0) {\includegraphics[width = 0.33\textwidth]{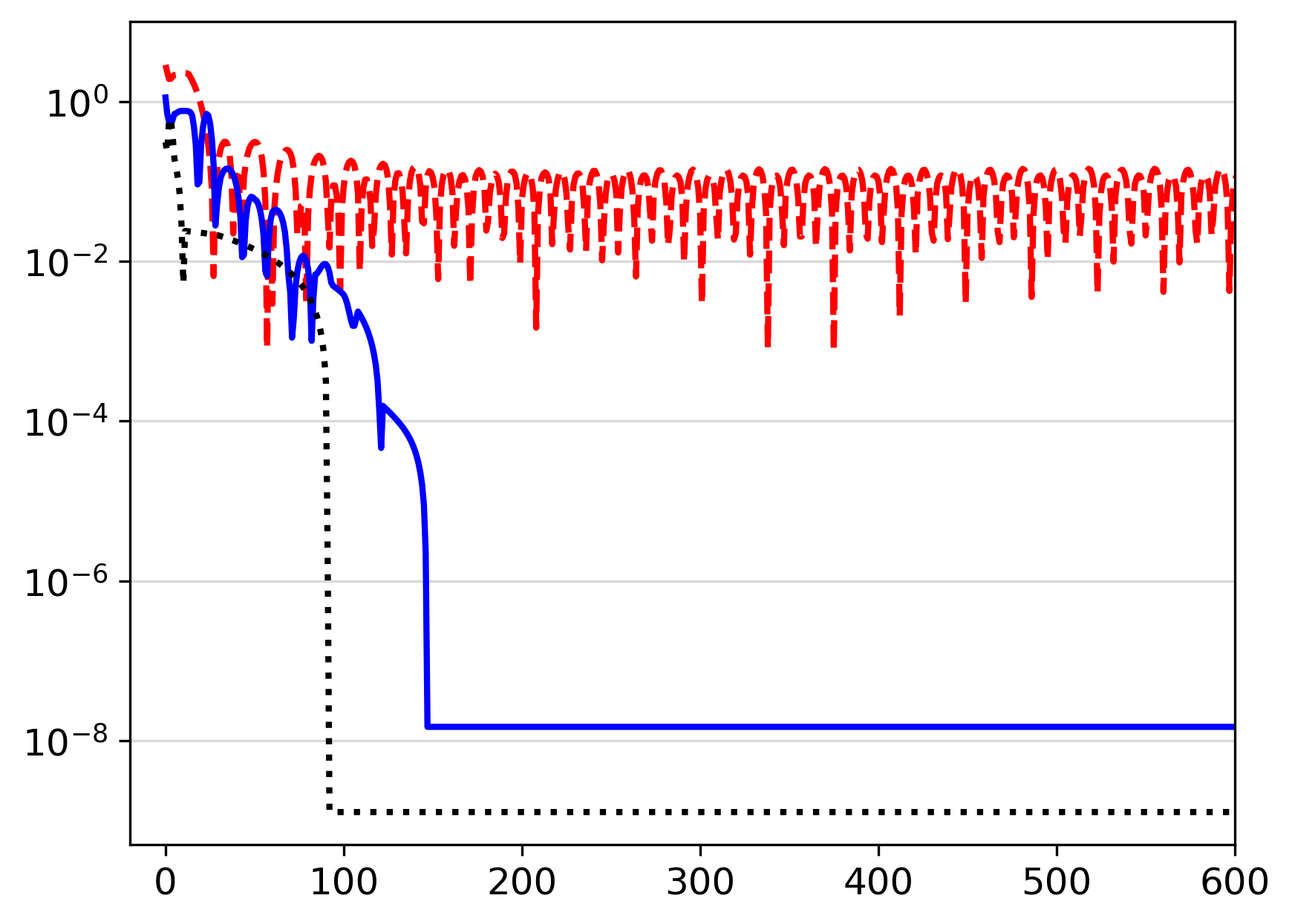}};
        \node[below = 0cm of img1]{iteration};
         \node[left = 0.4cm of img1,yshift=2.0cm, rotate=90]{reward optimality gap};
    \end{tikzpicture}
    \hspace{0.1\textwidth}
            \begin{tikzpicture}
        \node[](img1) at(0,0) {\includegraphics[width = 0.33\textwidth]{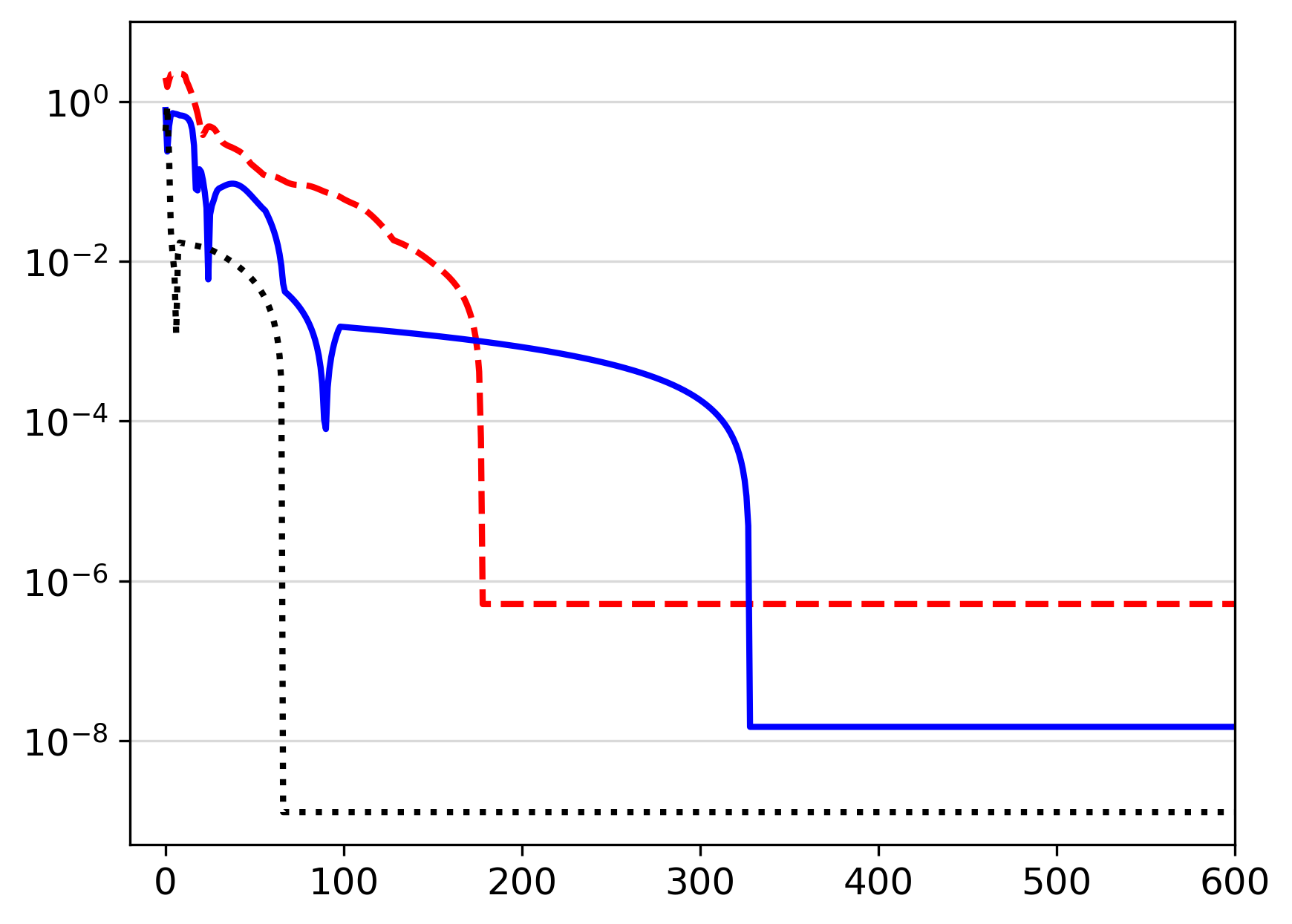}};
        \node[below = 0cm of img1]{iteration};
    \end{tikzpicture}
    \caption{Reward optimality gap of ResPG-PD (Algorithm~\ref{alg: resilient PG}, left) and ResOPG-PD (Algorithm~\ref{alg: resilient OPG}, right), with different cost functions $h(\xi) = \alpha {\xi}^2$, where $\alpha = 0.03$ (\ref{legend:red}), $\alpha = 0.2$ (\ref{legend:blue}), $\alpha = 1$ (\ref{legend:black}), and stepsize $\eta=0.2$ in the policy search problem of Section~\ref{sec:policysearch}.
    }
    \label{fig:AppE1RewardErrorvsIter}
    \end{figure}

\begin{figure}[tbh]
    \centering
        \begin{tikzpicture}
        \node[](img1) at(0,0) {\includegraphics[width = 0.33\textwidth]{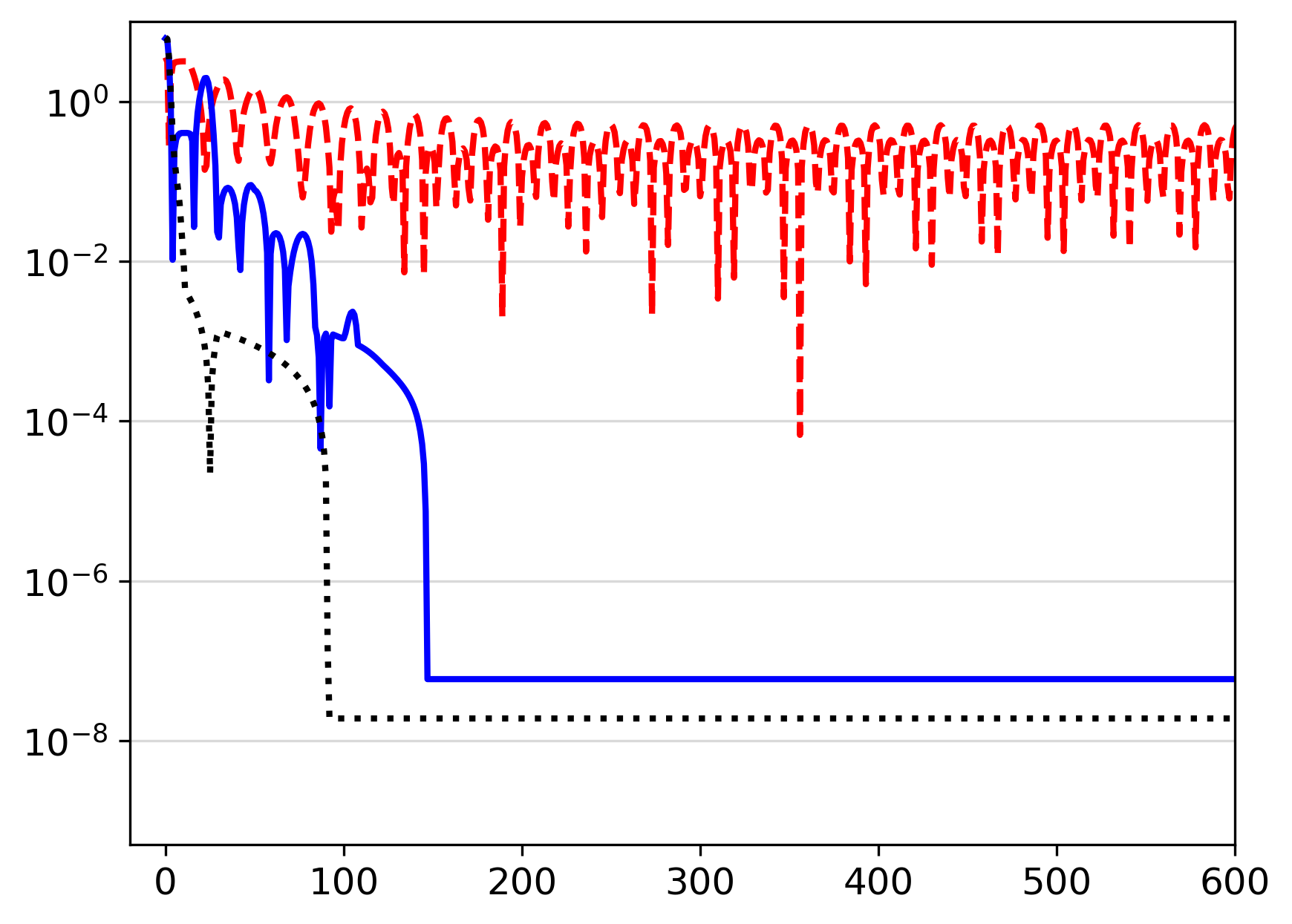}};
        \node[below = 0cm of img1]{iteration};
         \node[left = 0.4cm of img1,yshift=2.0cm, rotate=90]{utility optimality gap};
    \end{tikzpicture}
    \hspace{0.1\textwidth}
            \begin{tikzpicture}
        \node[](img1) at(0,0) {\includegraphics[width = 0.33\textwidth]{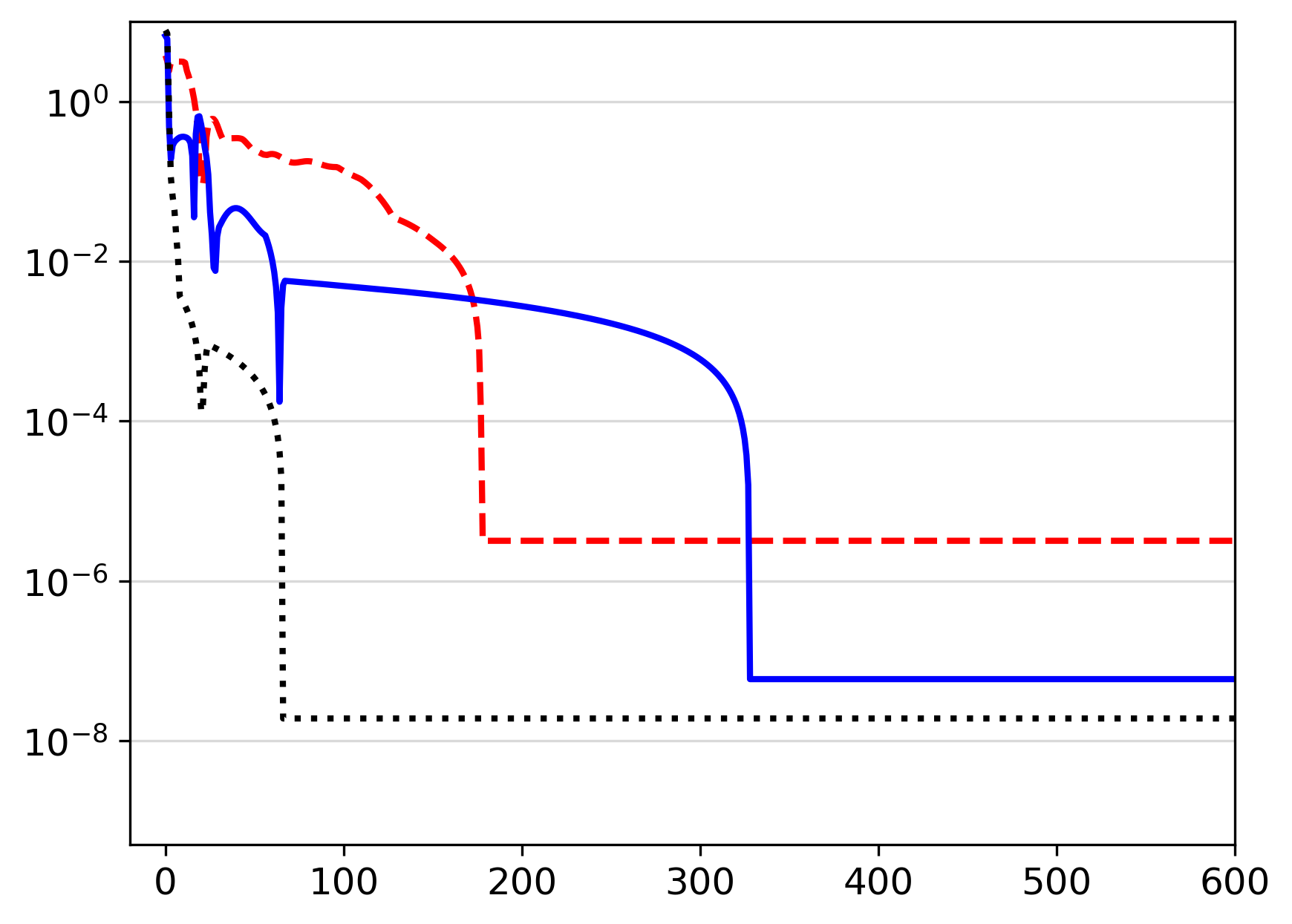}};
        \node[below = 0cm of img1]{iteration};
    \end{tikzpicture}
    \caption{Utility optimality gap of ResPG-PD (Algorithm~\ref{alg: resilient PG}, left) and ResOPG-PD (Algorithm~\ref{alg: resilient OPG}, right), with different cost functions $h(\xi) = \alpha {\xi}^2$, where $\alpha = 0.03$ (\ref{legend:red}), $\alpha = 0.2$ (\ref{legend:blue}), $\alpha = 1$ (\ref{legend:black}), and stepsize $\eta=0.2$ in the policy search problem of Section~\ref{sec:policysearch}.
    }
    \label{fig:AppE1UtilityErrorvsIter}
    \end{figure}

\begin{figure}[tbh]
    \centering
        \begin{tikzpicture}
        \node[](img1) at(0,0) {\includegraphics[width = 0.33\textwidth]{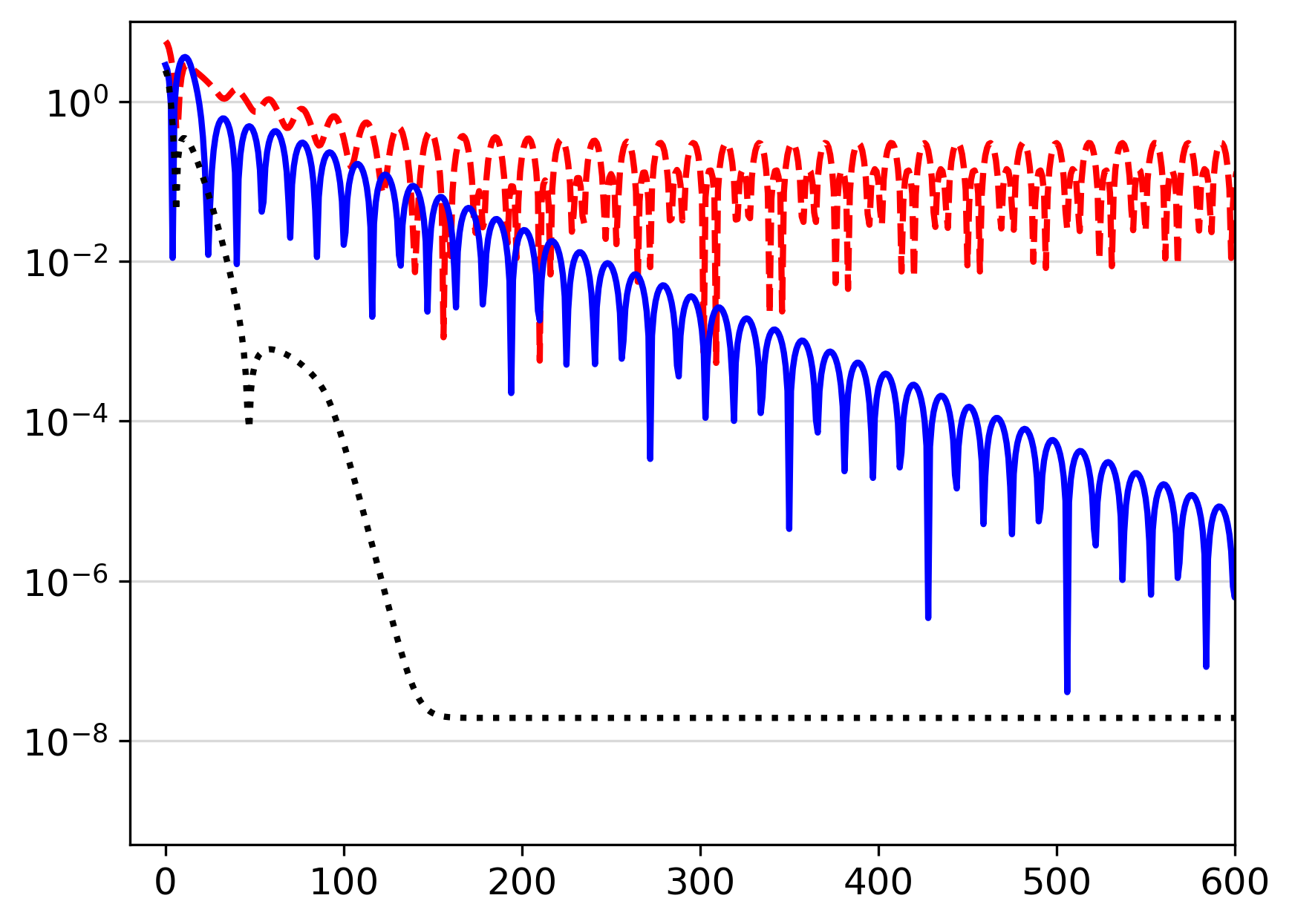}};
        \node[below = 0cm of img1]{iteration};
         \node[left = 0.4cm of img1,yshift=2.1cm, rotate=90]{relaxation optimality gap};
    \end{tikzpicture}
    \hspace{0.1\textwidth}
            \begin{tikzpicture}
        \node[](img1) at(0,0) {\includegraphics[width = 0.33\textwidth]{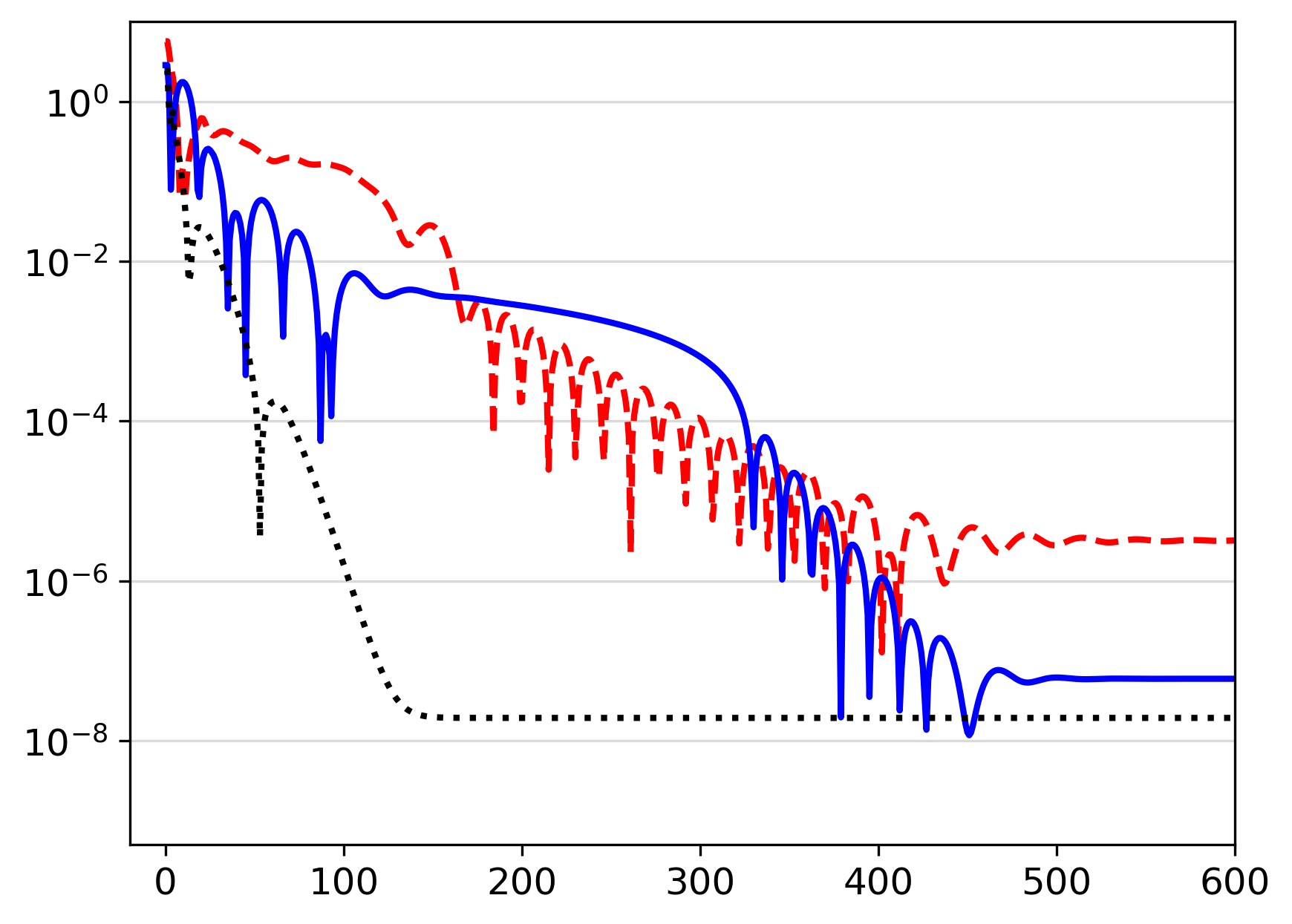}};
        \node[below = 0cm of img1]{iteration};
    \end{tikzpicture}
    \caption{Relaxation optimality gap of ResPG-PD (Algorithm~\ref{alg: resilient PG}, left) and ResOPG-PD (Algorithm~\ref{alg: resilient OPG}, right), with different cost functions $h(\xi) = \alpha {\xi}^2$, where $\alpha = 0.03$ (\ref{legend:red}), $\alpha = 0.2$ (\ref{legend:blue}), $\alpha = 1$ (\ref{legend:black}), and stepsize $\eta=0.2$ in the policy search problem of Section~\ref{sec:policysearch}.
    }
    \label{fig:AppE1XiErrorvsIter}
    \end{figure}

We then report the reward, utility, and relaxation gaps by comparing the iterates generated by the two methods with the optimal values from the quadratic program in Figures~\ref{fig:AppE1RewardErrorvsIter}--\ref{fig:AppE1XiErrorvsIter}. For three choices of $\alpha$, reward, utility, and relaxation of ResOPG-PD converge to the optimal ones, reaching low platforms due to the accuracy of quadratic program. ResPG-PD behaves similarly for large $\alpha$: $0.2$ or $1$, however, oscillates when small $\alpha$: $0.03$. 

\subsection{Resilient Monitoring: Small State Space}\label{sec:smallmonitor}

  





We consider a partial monitoring problem with three locations $S_0$, $S_1$, and $S_2$ in Figure~\ref{fig:MointoringProblem}. Each location represents a state of an agent.  
In each state, the choice of action determines the next state, $s_{t+1} = a_t$. In state $S_0$, possible actions are $\{S_1,S_2\}$. In state $S_i$ with $i\neq 0$, possible actions are $\{S_0,S_i\}$. We define the reward functions as
\[
    \begin{array}{rcl}
    r_i(s_t,a_t) 
    \; = \;
    b_i\mathbb{I}(s_t=S_i)\;\text{  for }\; i = 0, 1, 2
    \end{array}
\]
where $b_i$ is the reward for a single time step in a state $S_i$. To formulate a constrained MDP, we have: (i) two constraints $V_i^{\pi}\geq c_i$ for $i=1,2$, where the value functions $V_i^{\pi}$ are associated with the corresponding rewards $r_i(s_t,a_t)$; (ii) the objective $V_0^{\pi}$ is the value function associated with the reward $r_0(s_t,a_t)$. This problem has a uniform initial distribution $\rho$.  

As an example, we choose $b_0=b_1=1, b_2 = 1.2$, a discount factor $\gamma=0.9$ and initial constraints $V_{g_1}^{\pi}(\rho) = V_1^{\pi}(\rho)\geq 7$ and $V_{g_2}^{\pi}(\rho)=V_2^{\pi}(\rho)\geq 9$, which are infeasible. The relaxation controls the trade-off between the closeness to the original constraints and the gain of rewards. One extreme case is that the agent always moves to $S_0$ if not $S_0$. In this case, the reward reaches the maximum value and the relaxation is the maximum. Another extreme case is that the agent always stays in its original state $S_i$ if not in $S_0$. In this case, the agent spends as much time as possible in $S_1$, $S_2$ to gain utility, which makes relaxed constraints close to the original constraints and makes the reward value smallest. The relaxation cost function is $h(\xi) = \alpha\|\xi\|^2$ in which we set $\alpha=0.1$. We observe that our two methods can relax the two initial constraints to make this problem feasible, but keep away from the extreme cases.

\begin{figure}[tbh]
    \centering
        \begin{tikzpicture}
        \node[](img1) at(0,0) {\includegraphics[width = 0.33\textwidth]{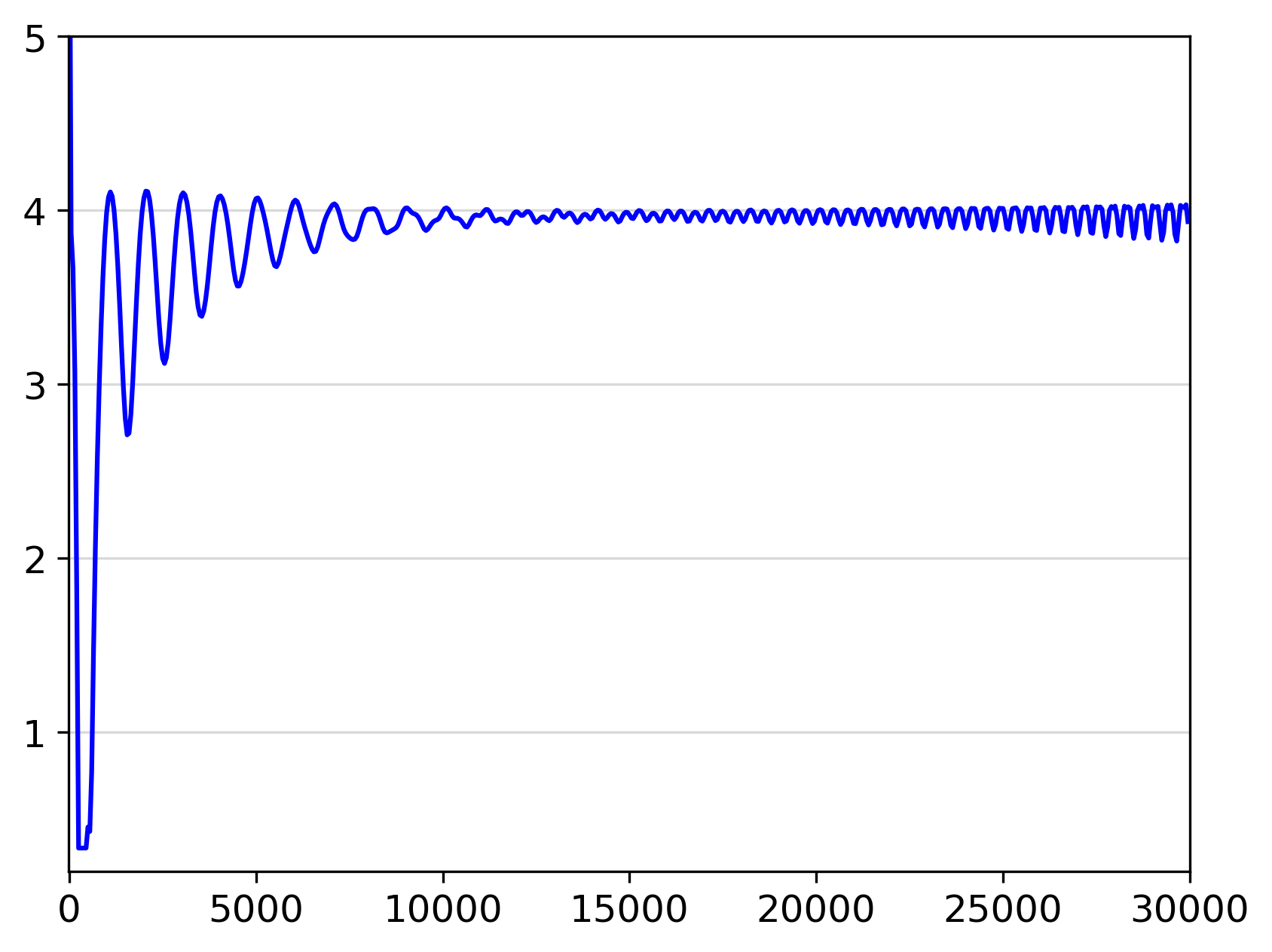}};
        \node[below = 0cm of img1]{iteration};
         \node[left = 0.4cm of img1,yshift=1.3cm, rotate=90]{reward value};
    \end{tikzpicture}
    \hspace{0.1\textwidth}
            \begin{tikzpicture}
        \node[](img1) at(0,0) {\includegraphics[width = 0.33\textwidth]{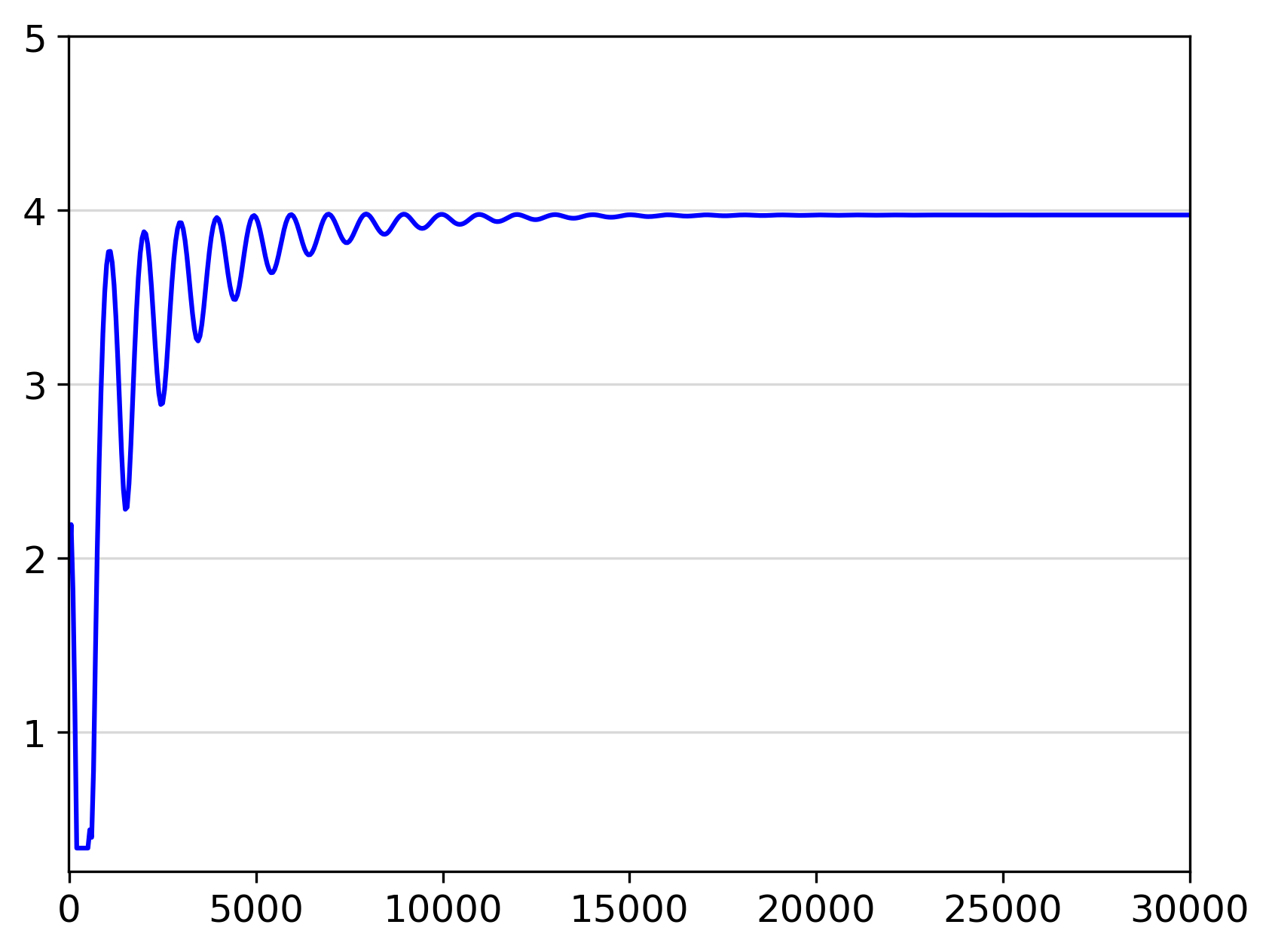}};
        \node[below = 0cm of img1]{iteration};
    \end{tikzpicture}
    \caption{Reward value convergence of ResPG-PD (Algorithm~\ref{alg: resilient PG}, left) and ResOPG-PD (Algorithm~\ref{alg: resilient OPG}, right), with a cost functions $h(\xi) = \alpha \norm{\xi}^2$, where $\alpha = 0.1$,  and stepsize $\eta=0.005$ in the monitoring problem of Section~\ref{sec:smallmonitor}.}
    \label{fig:AppE2RewardvsIter}
    \end{figure}
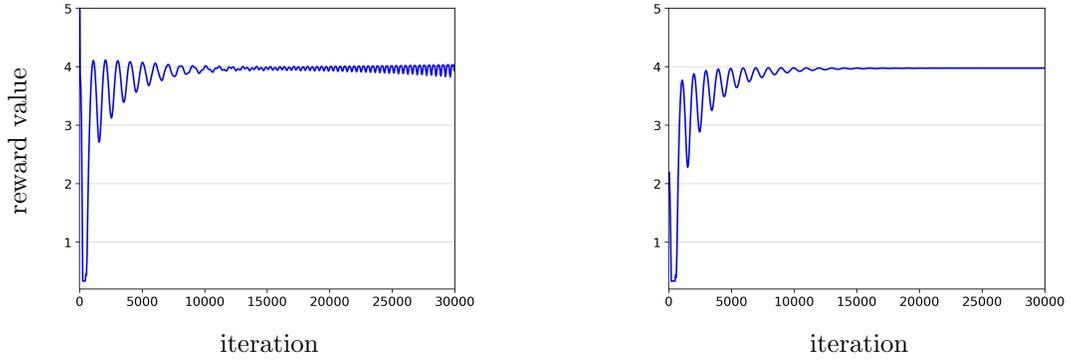

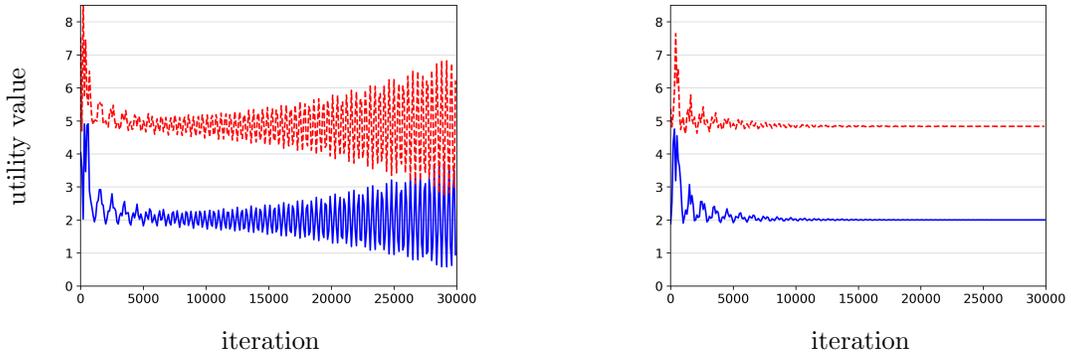
\begin{figure}[tbh]
    \centering
        \begin{tikzpicture}
        \node[](img1) at(0,0) {\includegraphics[width = 0.33\textwidth]{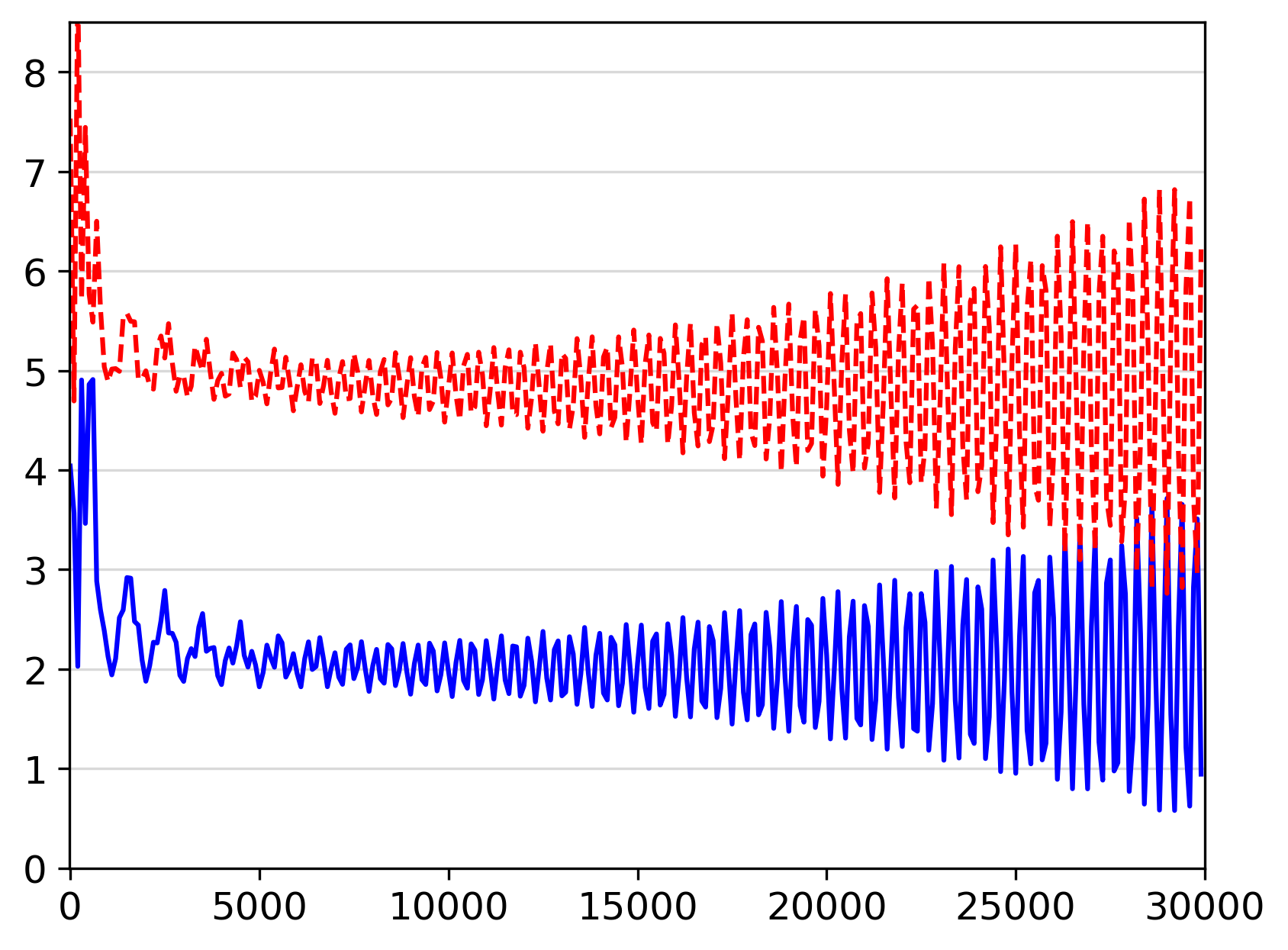}};
        \node[below = 0cm of img1]{iteration};
         \node[left = 0.4cm of img1,yshift=1.3cm, rotate=90]{utility value};
    \end{tikzpicture}
    \hspace{0.1\textwidth}
            \begin{tikzpicture}
        \node[](img1) at(0,0) {\includegraphics[width = 0.33\textwidth]{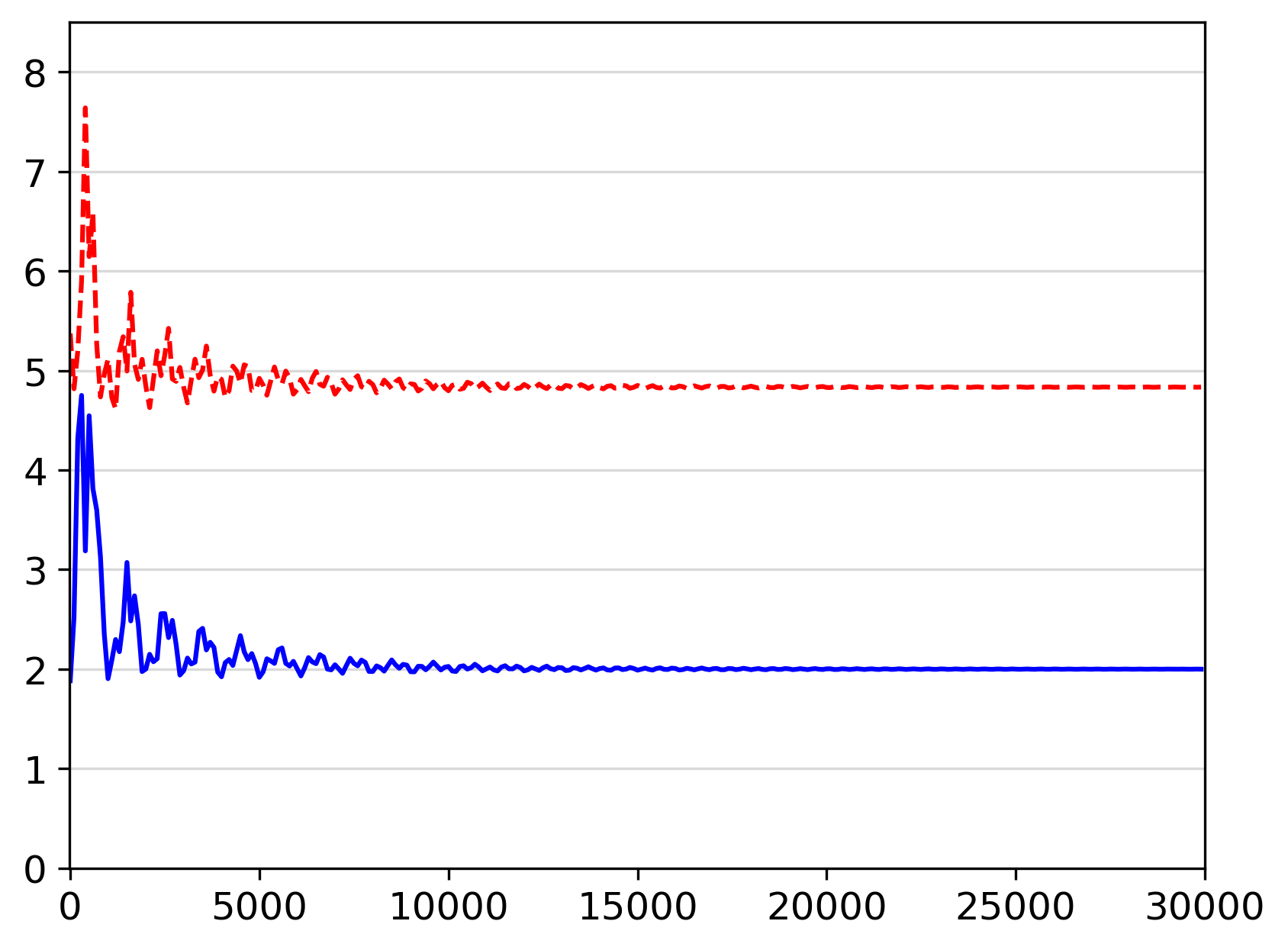}};
        \node[below = 0cm of img1]{iteration};
    \end{tikzpicture}
    \caption{Utility value convergence ($V_{g_1}^\pi(\rho)$:~\ref{legend:blue}, $V_{g_2}^\pi(\rho)$:~\ref{legend:red}) of ResPG-PD (Algorithm~\ref{alg: resilient PG}, left) and ResOPG-PD (Algorithm~\ref{alg: resilient OPG}, right), with a cost functions $h(\xi) = \alpha \norm{\xi}^2$ for $\alpha = 0.1$, and stepsize $\eta = 0.005$ in the monitoring problem of Section~\ref{sec:smallmonitor}. 
    }
    \label{fig:AppE2UtilityvsIter}
    \end{figure}

\begin{figure}[tbh]
    \centering
        \begin{tikzpicture}
        \node[](img1) at(0,0) {\includegraphics[width = 0.33\textwidth]{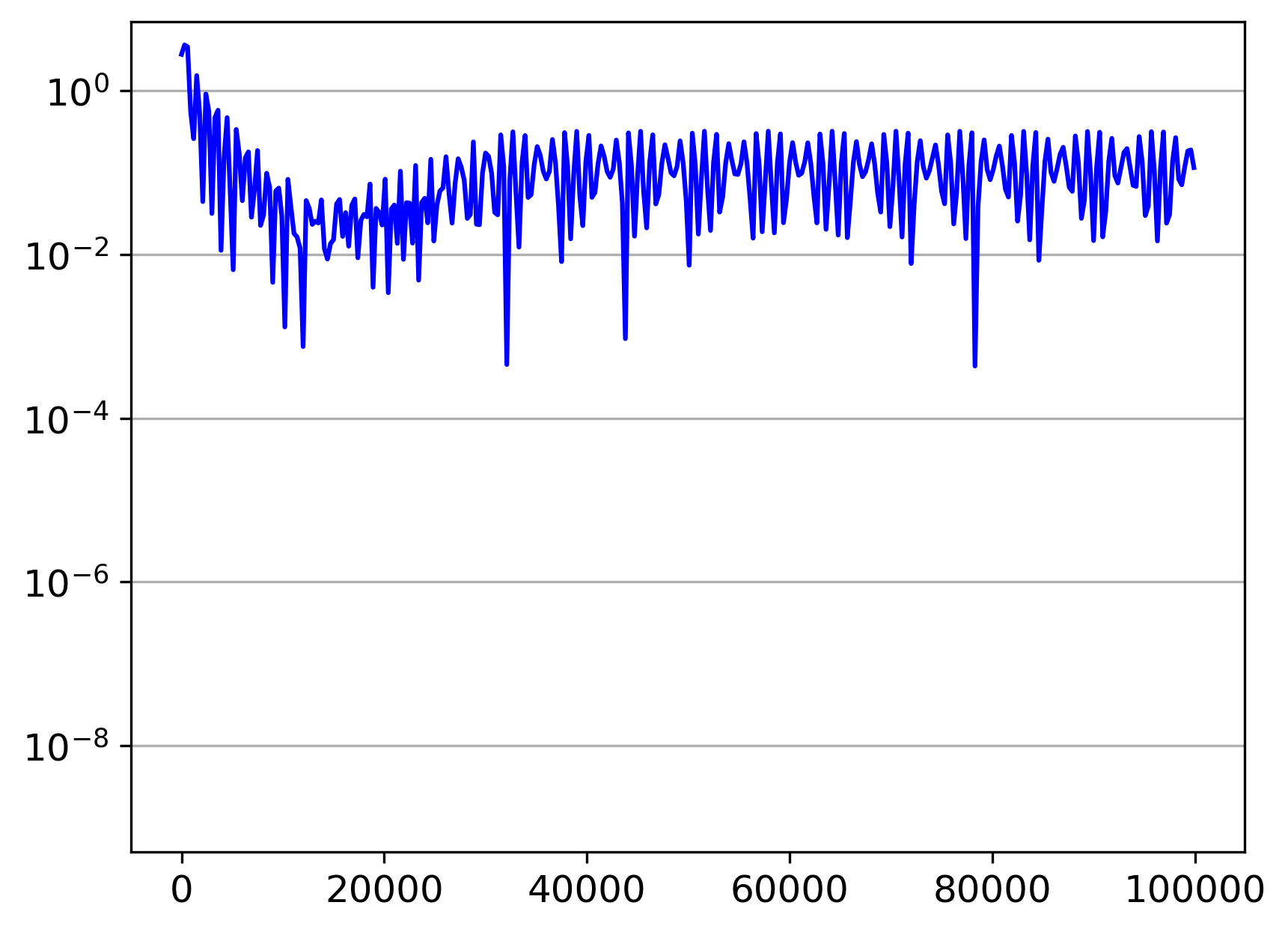}};
        \node[below = 0cm of img1]{iteration};
         \node[left = 0.4cm of img1,yshift=2.1cm, rotate=90]{reward optimality gap};
    \end{tikzpicture}
    \hspace{0.1\textwidth}
            \begin{tikzpicture}
        \node[](img1) at(0,0) {\includegraphics[width = 0.33\textwidth]{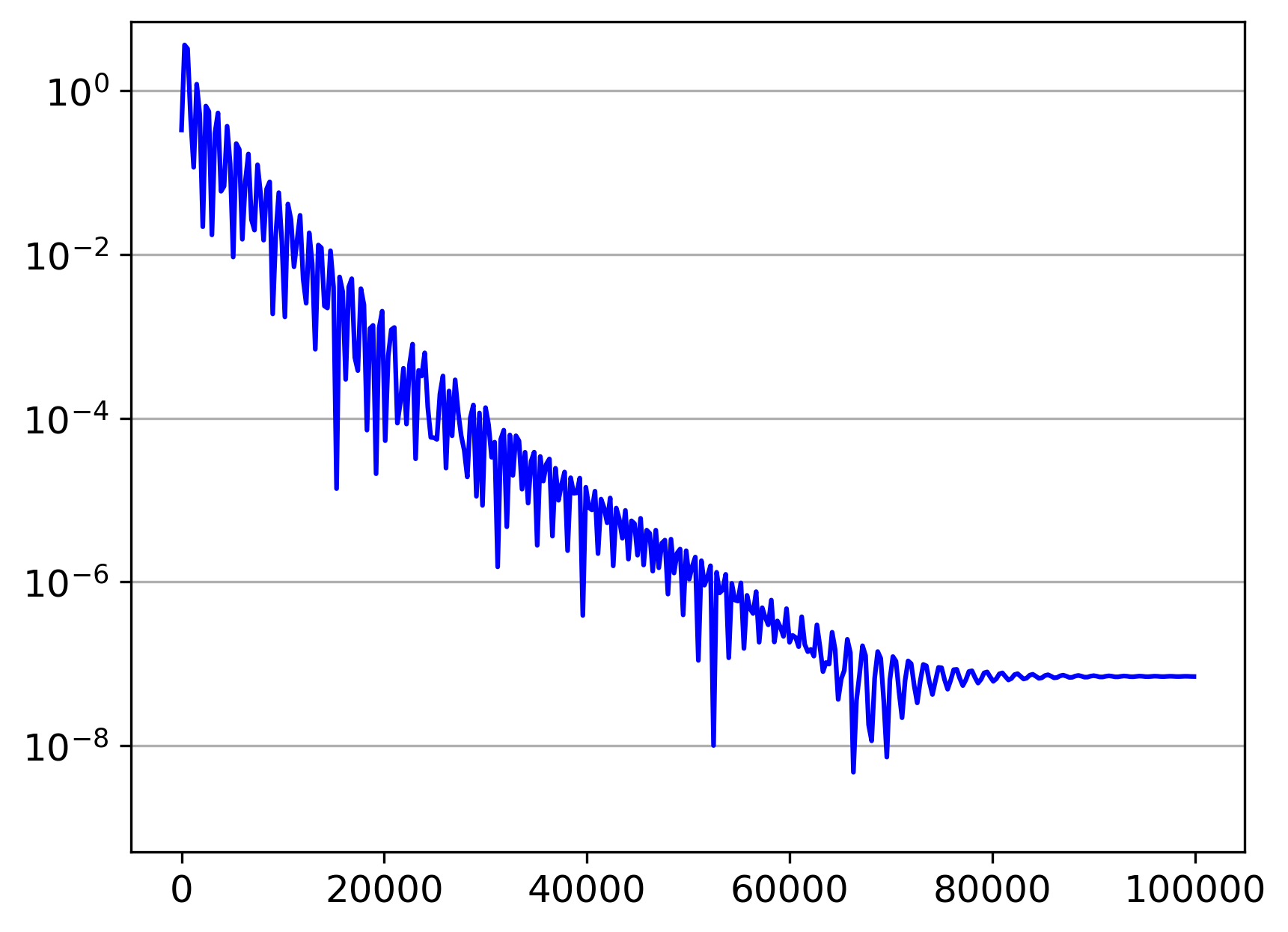}};
        \node[below = 0cm of img1]{iteration};
    \end{tikzpicture}
    \caption{Reward optimality gap of ResPG-PD (Algorithm~\ref{alg: resilient PG}, left) and ResOPG-PD (Algorithm~\ref{alg: resilient OPG}, right), with a cost functions $h(\xi) = \alpha \norm{\xi}^2$, where $\alpha = 0.1$ and stepsize $\eta=0.005$ in the monitoring problem of Section~\ref{sec:smallmonitor}.
    }
    \label{fig:AppE2RewardErrorvsIter}
    \end{figure}
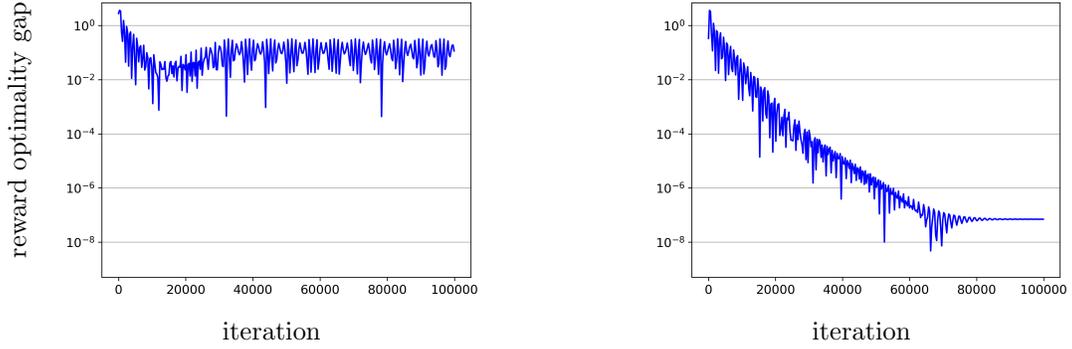

\begin{figure}[tbh]
    \centering
        \begin{tikzpicture}
        \node[](img1) at(0,0) {\includegraphics[width = 0.22\textwidth]{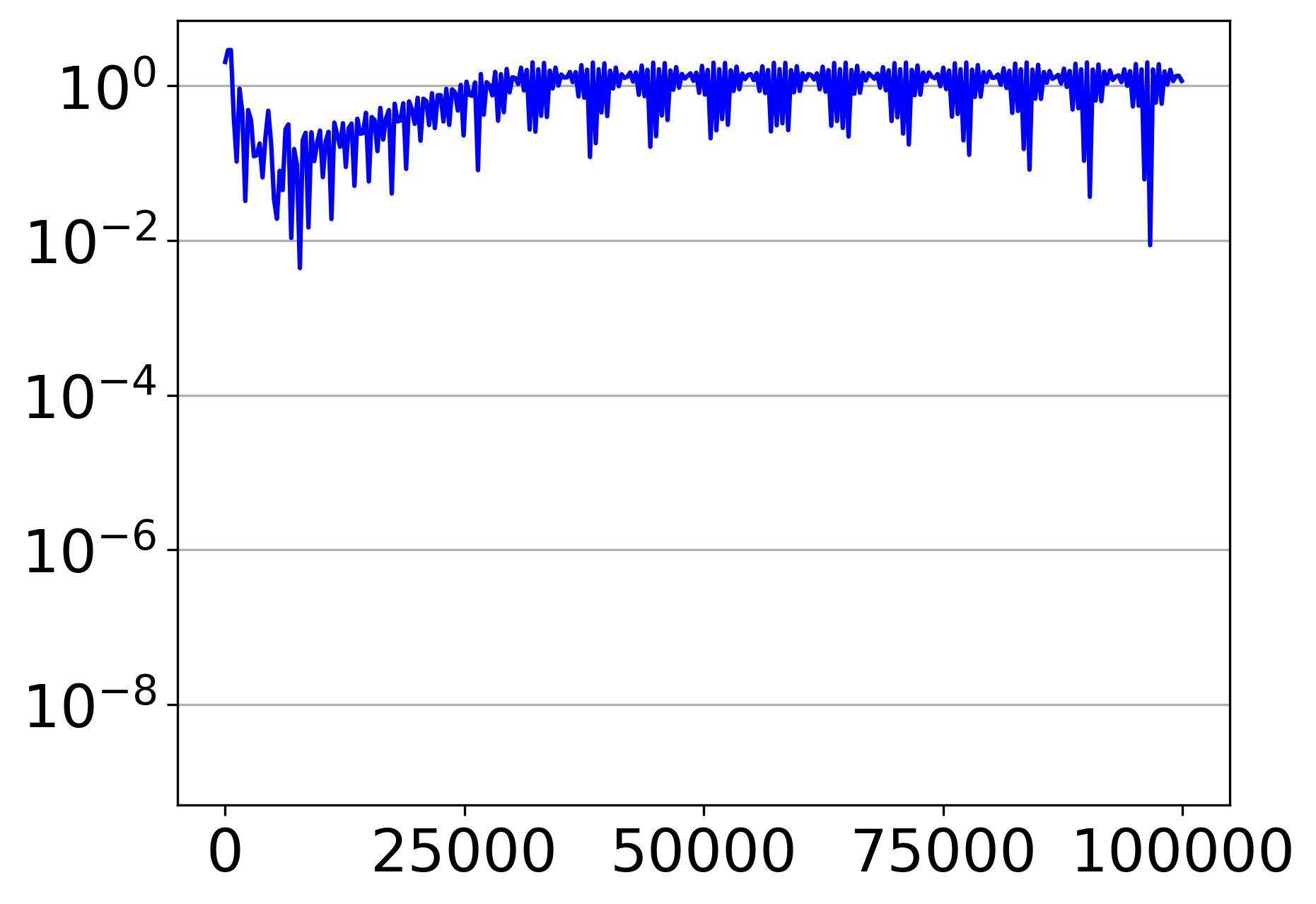}};
        \node[below = 0cm of img1]{iteration};
         \node[left = 0.4cm of img1,yshift=2cm, rotate=90]{utility optimality gap};
    \end{tikzpicture}
            \begin{tikzpicture}
        \node[](img1) at(0,0) {\includegraphics[width = 0.22\textwidth]{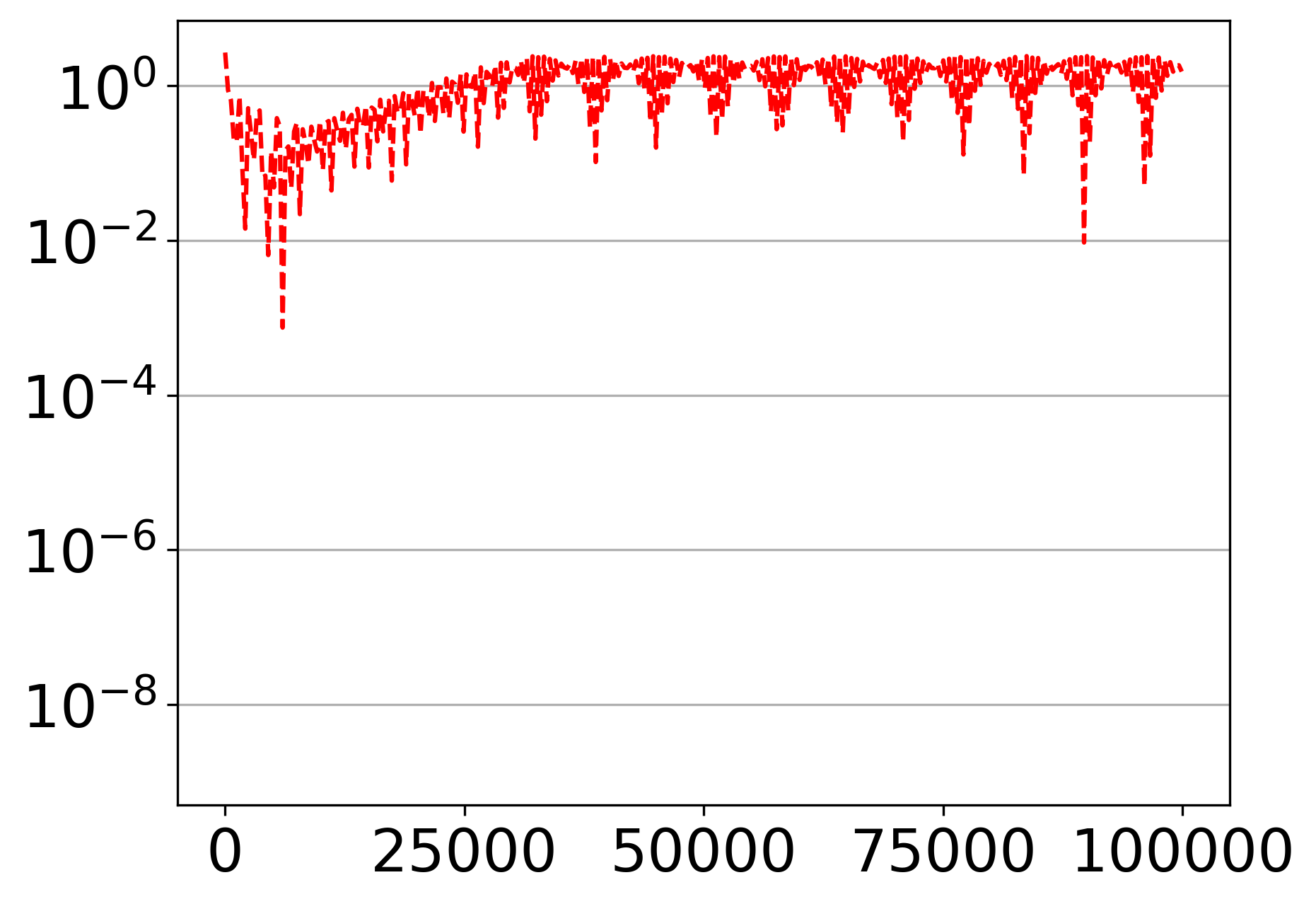}};
        \node[below = 0cm of img1]{iteration};
    \end{tikzpicture}
            \begin{tikzpicture}
        \node[](img1) at(0,0) {\includegraphics[width = 0.22\textwidth]{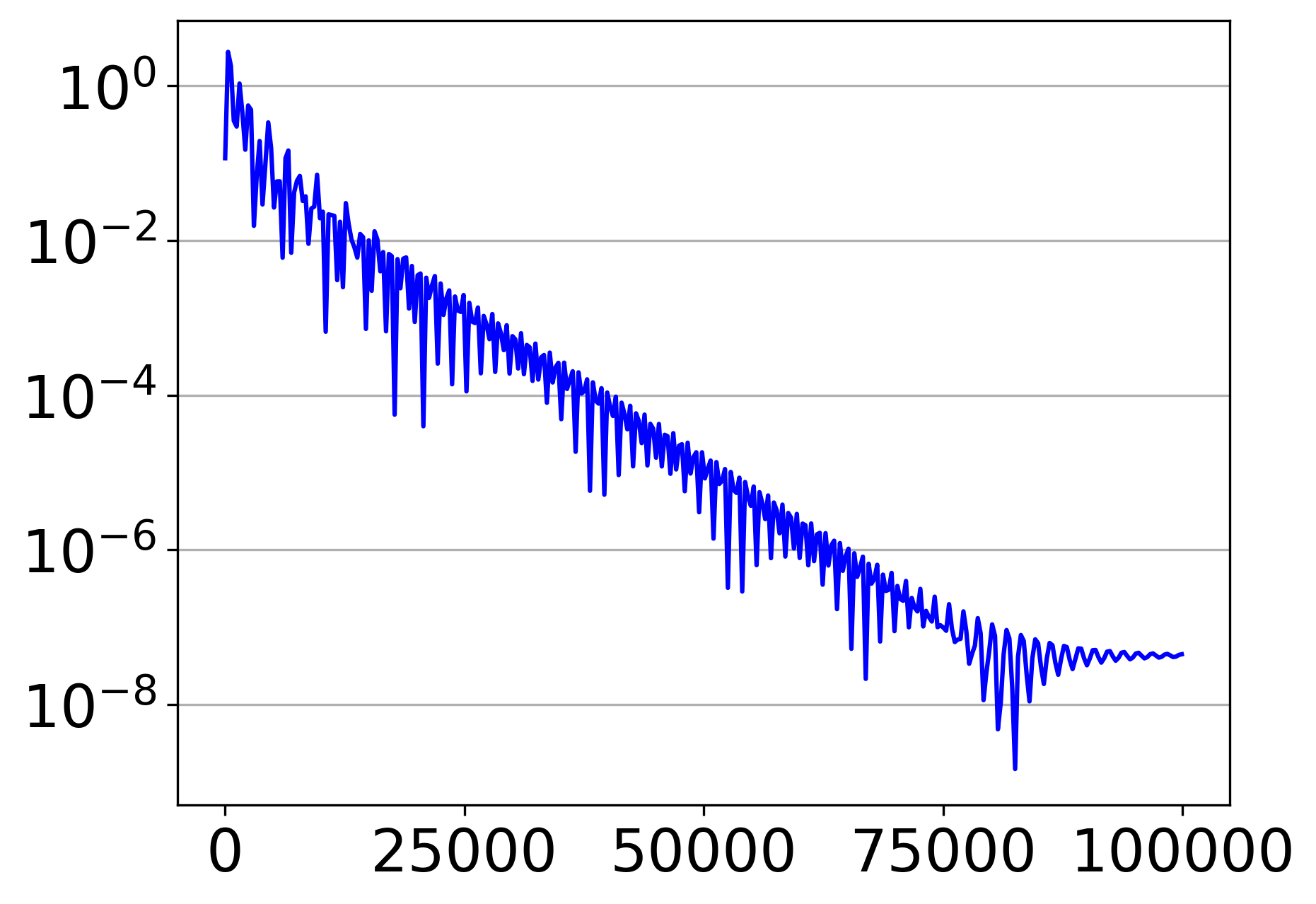}};
        \node[below = 0cm of img1]{iteration};
    \end{tikzpicture}
            \begin{tikzpicture}
        \node[](img1) at(0,0) {\includegraphics[width = 0.22\textwidth]{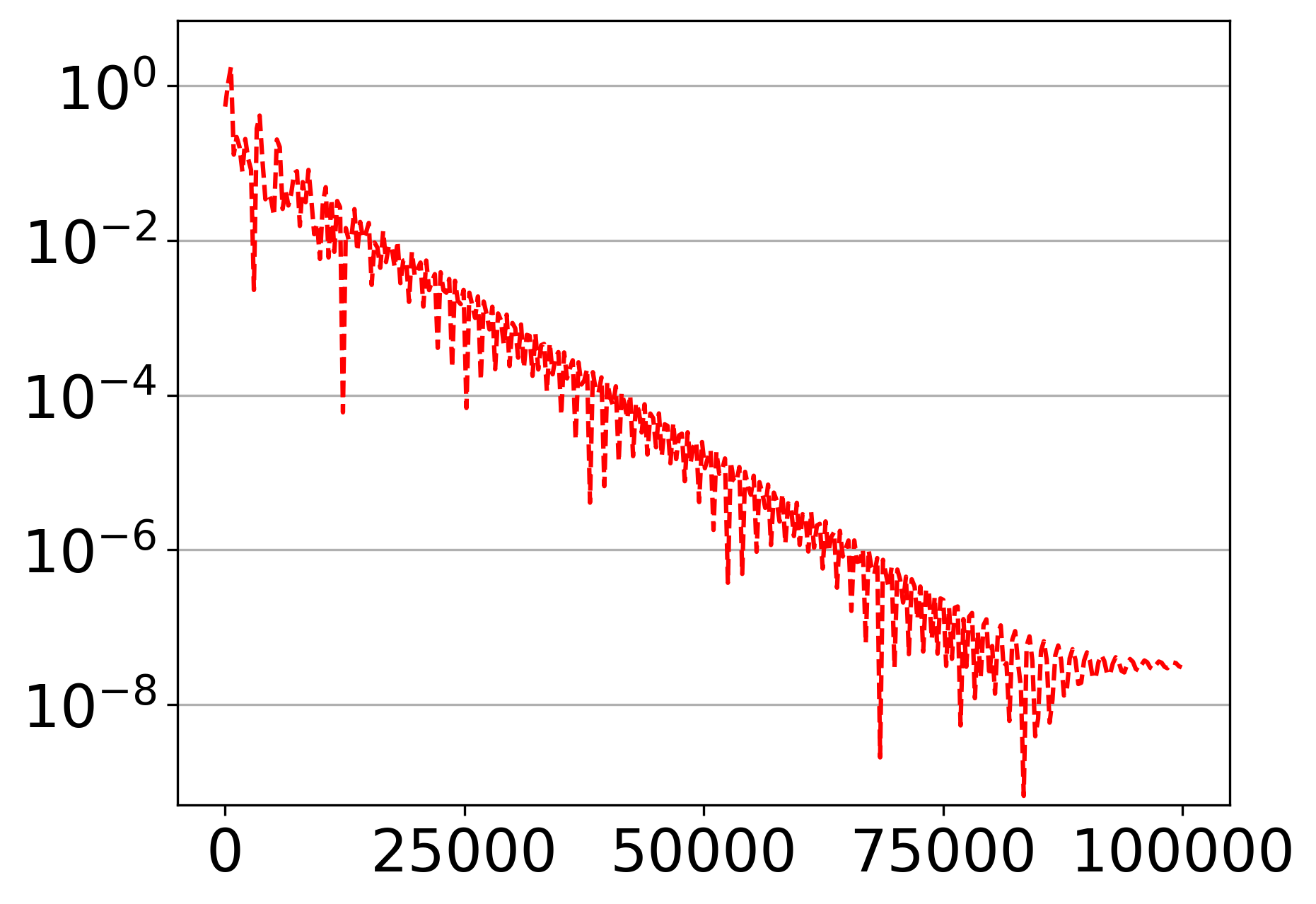}};
        \node[below = 0cm of img1]{iteration};
    \end{tikzpicture}
    \caption{Utility optimality gap ($V_{g_1}^\pi(\rho)$:~\ref{legend:blue}, $V_{g_2}^\pi(\rho)$:~\ref{legend:red}) of ResPG-PD (Algorithm~\ref{alg: resilient PG}, two figures on the left) and ResOPG-PD (Algorithm~\ref{alg: resilient OPG}, two figures on the right), with a cost functions $h(\xi) = \alpha \norm{\xi}^2$, where $\alpha = 0.1$ and stepsize $\eta=0.005$ in the monitoring problem of Section~\ref{sec:smallmonitor}.
    }
    \label{fig:AppE2UtilityErrorvsIter}
    \end{figure}
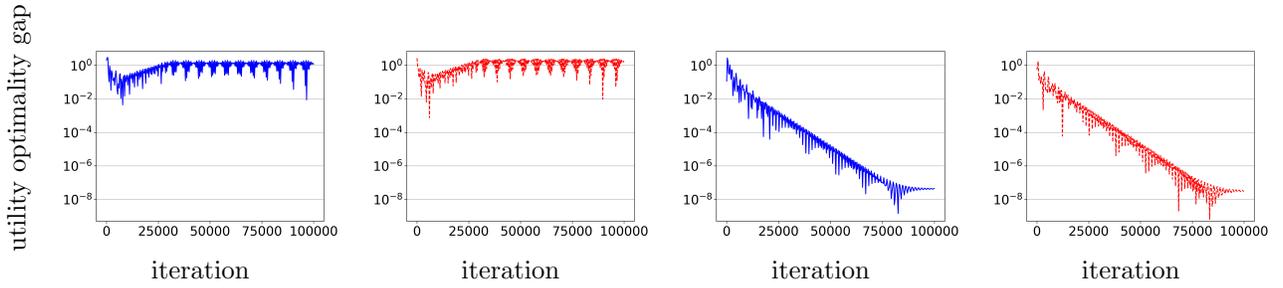

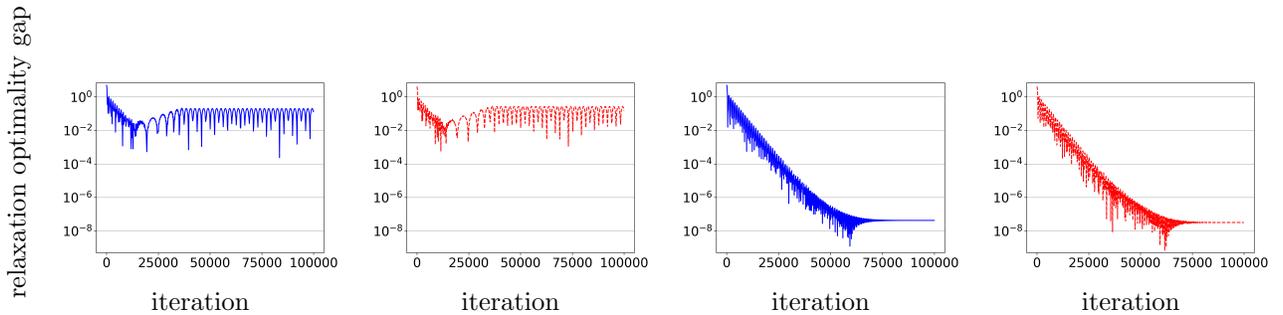
\begin{figure}[tbh]
    \centering
        \begin{tikzpicture}
        \node[](img1) at(0,0) {\includegraphics[width = 0.22\textwidth]{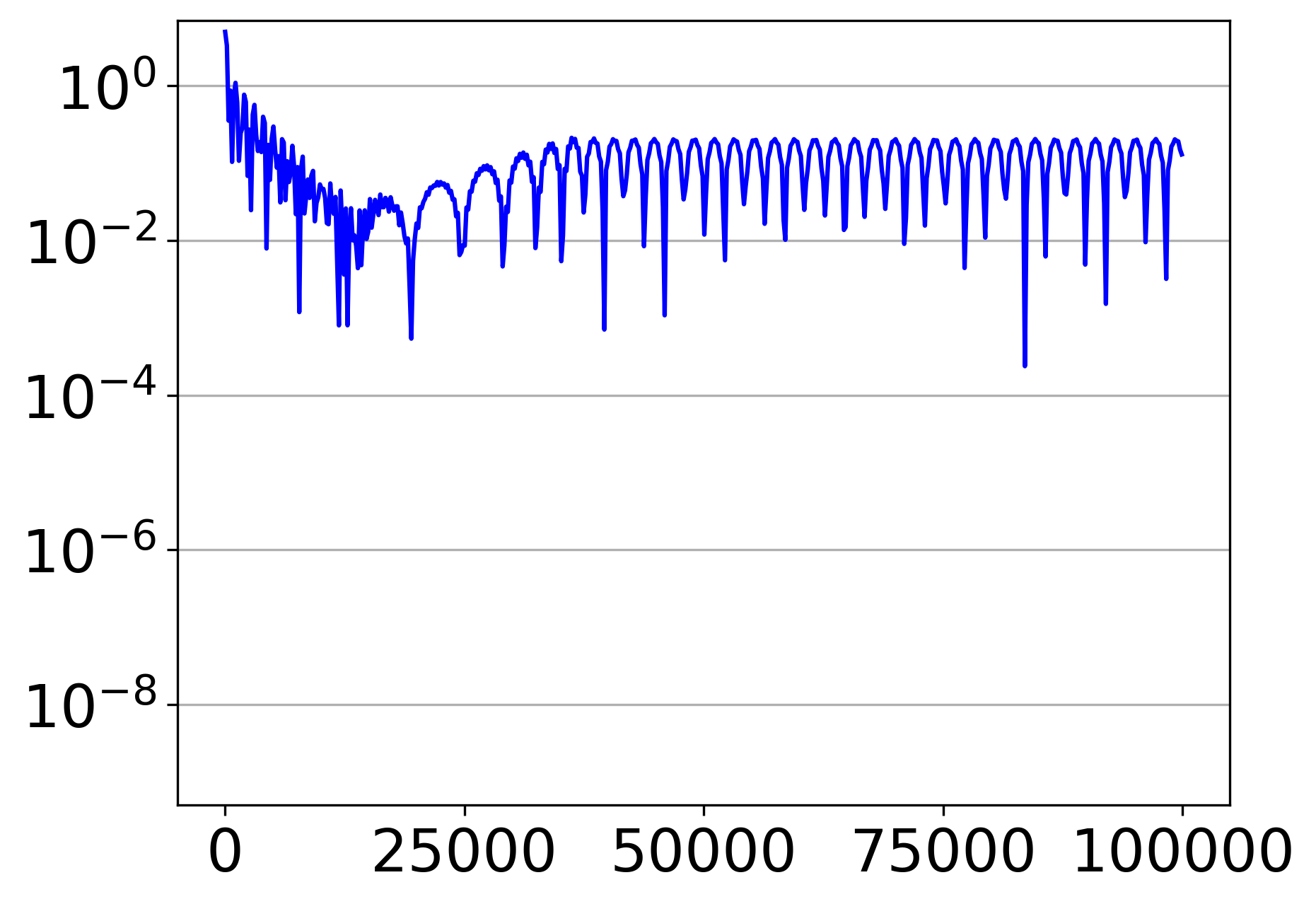}};
        \node[below = 0cm of img1]{iteration};
         \node[left = 0.4cm of img1,yshift=2.4cm, rotate=90]{relaxation optimality gap};
    \end{tikzpicture}
            \begin{tikzpicture}
        \node[](img1) at(0,0) {\includegraphics[width = 0.22\textwidth]{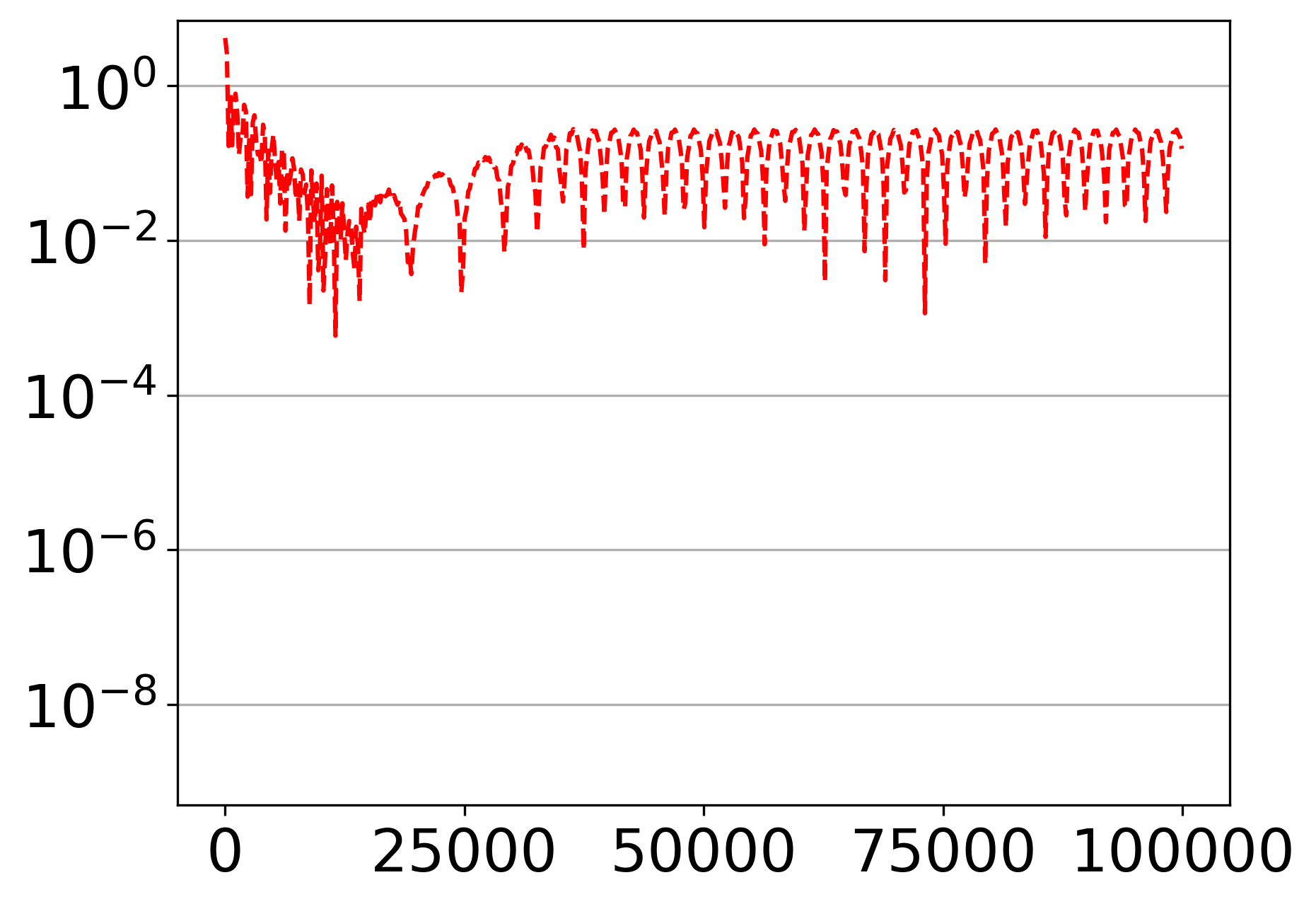}};
        \node[below = 0cm of img1]{iteration};
    \end{tikzpicture}
            \begin{tikzpicture}
        \node[](img1) at(0,0) {\includegraphics[width = 0.22\textwidth]{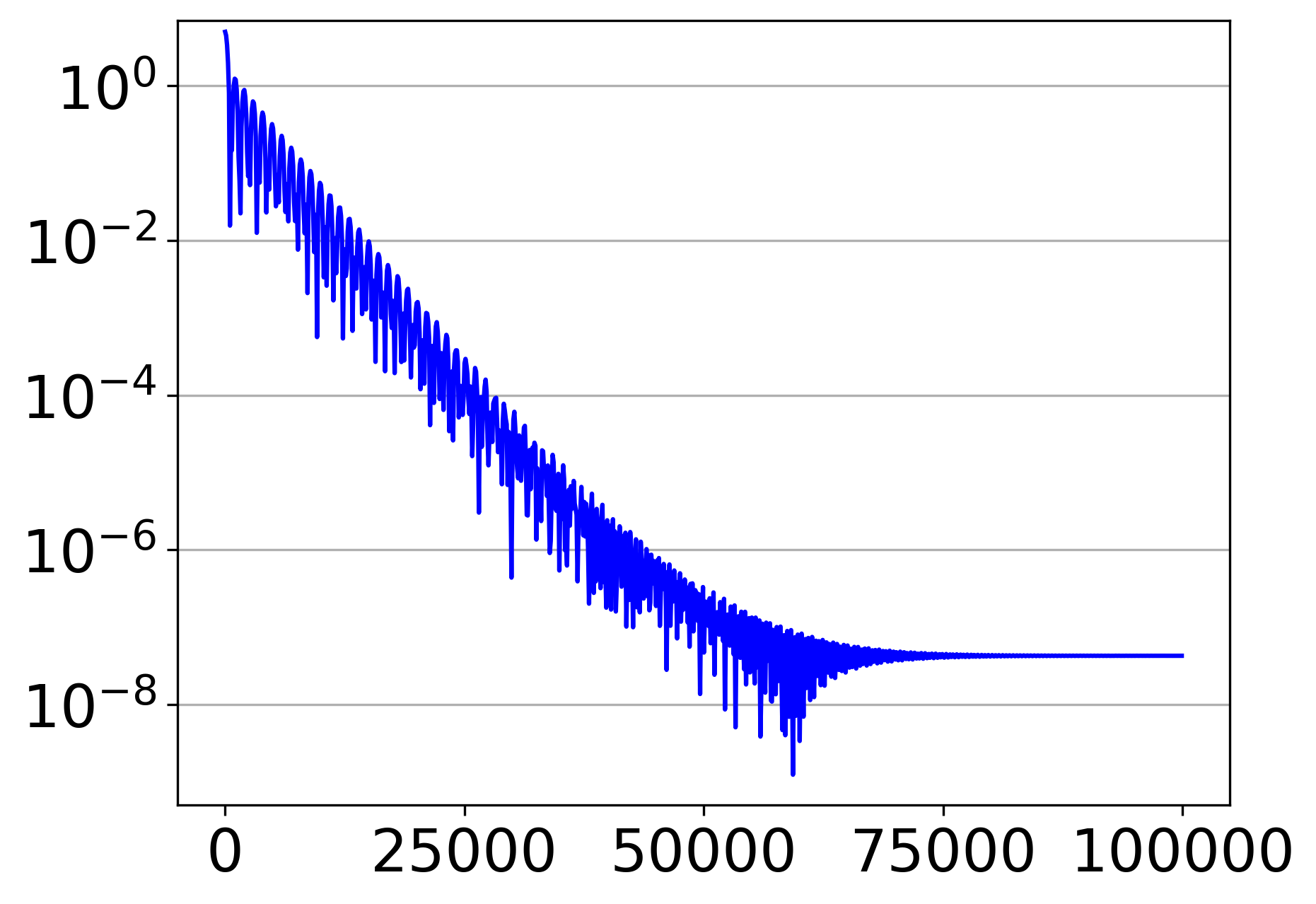}};
        \node[below = 0cm of img1]{iteration};
    \end{tikzpicture}
            \begin{tikzpicture}
        \node[](img1) at(0,0) {\includegraphics[width = 0.22\textwidth]{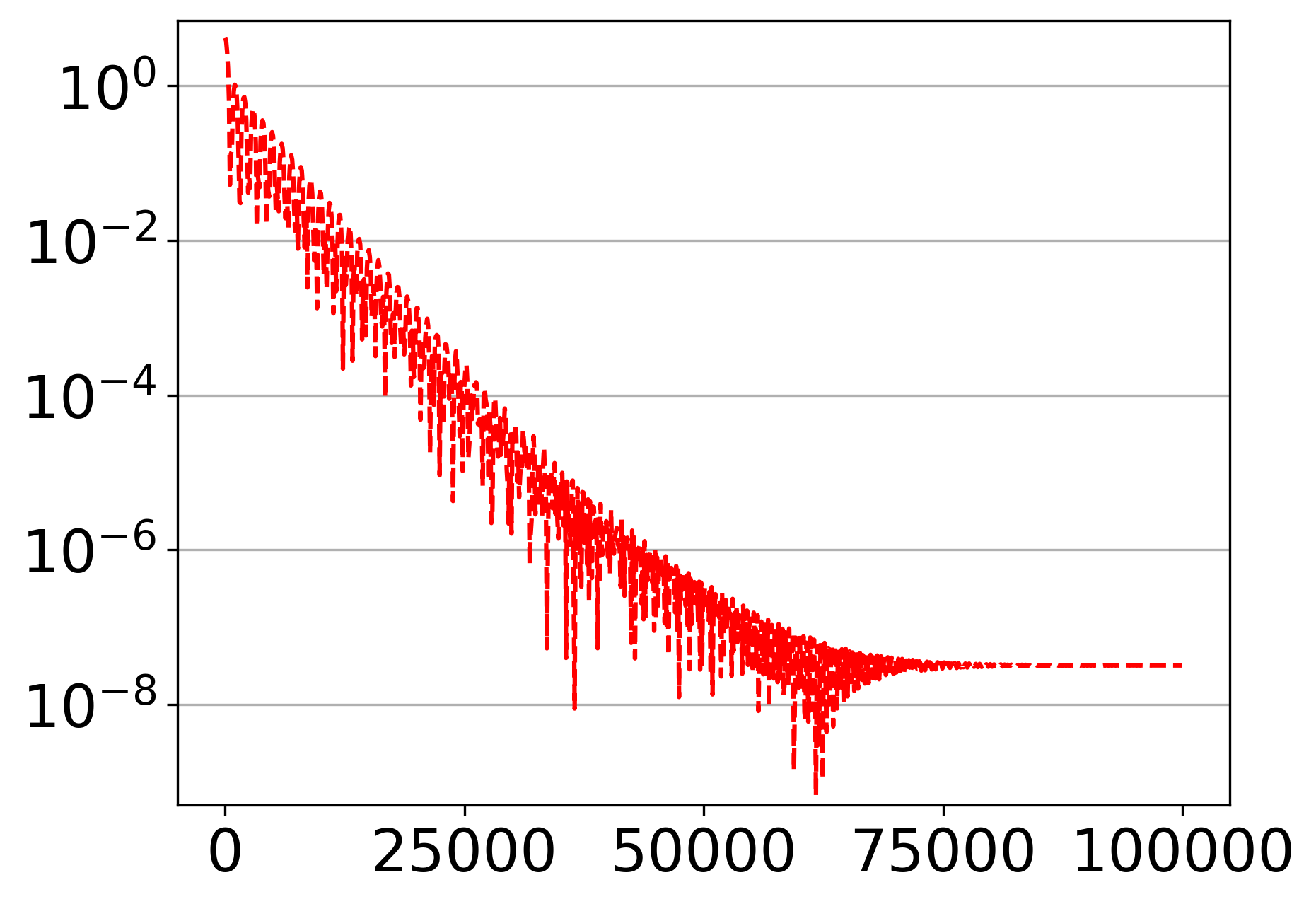}};
        \node[below = 0cm of img1]{iteration};
    \end{tikzpicture}
    \caption{Relaxation optimality gap ($\xi_1$: \ref{legend:blue}, $\xi_2$: \ref{legend:red}) of ResPG-PD (Algorithm~\ref{alg: resilient PG}, two figures on the left) and ResOPG-PD (Algorithm~\ref{alg: resilient OPG}, two figures on the right), with a cost functions $h(\xi) = \alpha \norm{\xi}^2$, where $\alpha = 0.1$ and stepsize $\eta=0.005$ in the monitoring problem of Section~\ref{sec:smallmonitor}.
    }
    \label{fig:AppE2XiErrorvsIter}
    \end{figure}

We recall Figure~\ref{fig:ResE2} that both methods can successfully relax the two constraints, which are initially infeasible, to be feasible for different cost functions. To show the reward and utility convergence performance of our two methods: Algorithm~\ref{alg: resilient PG} (ResPG-PD) and Algorithm~\ref{alg: resilient OPG} (ResOPG-PD), we report Figure~\ref{fig:AppE2RewardvsIter}--\ref{fig:AppE2UtilityvsIter} for $\alpha = 0.1$. We show the reward, utility and relaxation optimality gaps of the two methods in Figures~\ref{fig:AppE2RewardErrorvsIter}--\ref{fig:AppE2XiErrorvsIter}. As we expect, ResPG-PD often behave oscillating during training while ResOPG-PD can overcome oscillation, yielding a nearly-optimal policy in the last iterate.

\subsection{Resilient Monitoring: Large State Space}\label{sec:largemonitor}

We generalize the problem in Section~\ref{sec:smallmonitor} to a robot monitoring problem in a $10\times10$ grid as shown in Figure~\ref{fig:Monitoring_nonres},where each grid point is a state. In each state, four possible actions are: left, right, up, and down. The choice of the action and current state determines next state. The next state is the current state moving towards the action selected for one unit. If the next state falls outside the grid, the robot remains in the current state. Three areas $S_0$, $S_1$, $S_2$ in the grid represent three areas to be monitored. We define the reward functions as
\[
    \begin{array}{rcl}
    r_i(s_t,a_t) 
    \; = \;
    b_i\mathbb{I}(s_t\in S_i)\;\text{  for }\; i = 0, 1, 2
    \end{array}
\]
where $b_i$ is the reward for a single time step in an area $S_i$. We also have: (i) two constraints $V_i^{\pi}\geq c_i$ for $i=1,2$, where the value function $V_i^{\pi}$ is associated with the reward $r_i(s_t,a_t)$; (ii) the objective $V_0^{\pi}$ is the value function associated with the reward $r_0(s_t,a_t)$. The initial distribution $\rho$ is uniform.  

We use the same problem parameter setting as Section~\ref{sec:smallmonitor}, i.e., $b_0=b_1=1, b_2 = 1.2$, a discount factor $\gamma=0.9$ and initial constraints $V_{g_1}^{\pi}(\rho) = V_1^{\pi}(\rho)\geq 7$ and $V_{g_2}^{\pi}(\rho)=V_2^{\pi}(\rho)\geq 9$, which are infeasible. The relaxation cost function is $h(\xi) = \alpha\|\xi\|^2$ with $\alpha=0.08$.

\begin{figure}[tbh]
    \centering
        \begin{tikzpicture}
        \node[](img1) at(0,0) {\includegraphics[width = 0.33\textwidth]{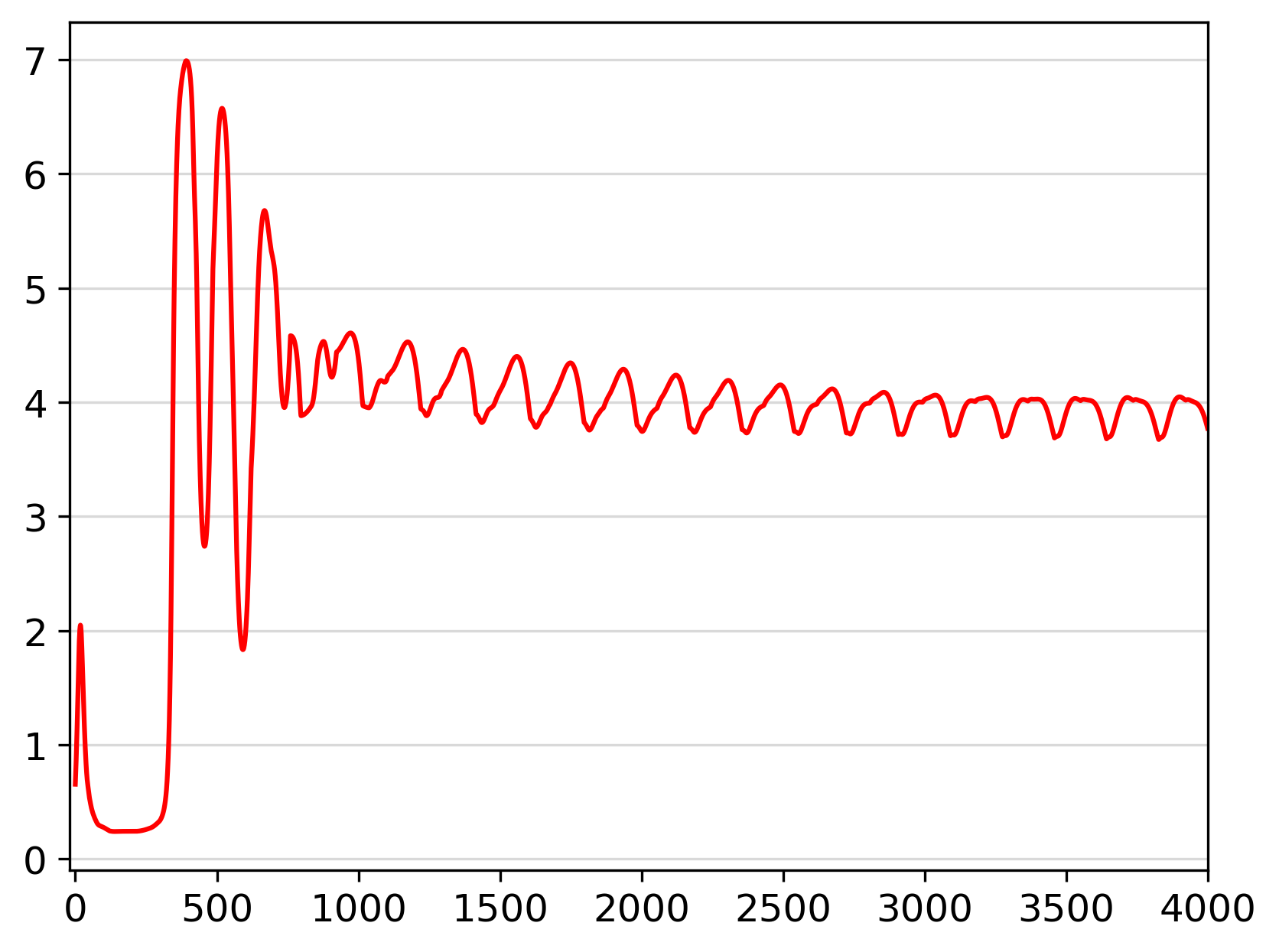}};
        \node[below = 0cm of img1]{iteration};
         \node[left = 0.4cm of img1,yshift=1.3cm, rotate=90]{reward value};
    \end{tikzpicture}
    \hspace{0.1\textwidth}
            \begin{tikzpicture}
        \node[](img1) at(0,0) {\includegraphics[width = 0.33\textwidth]{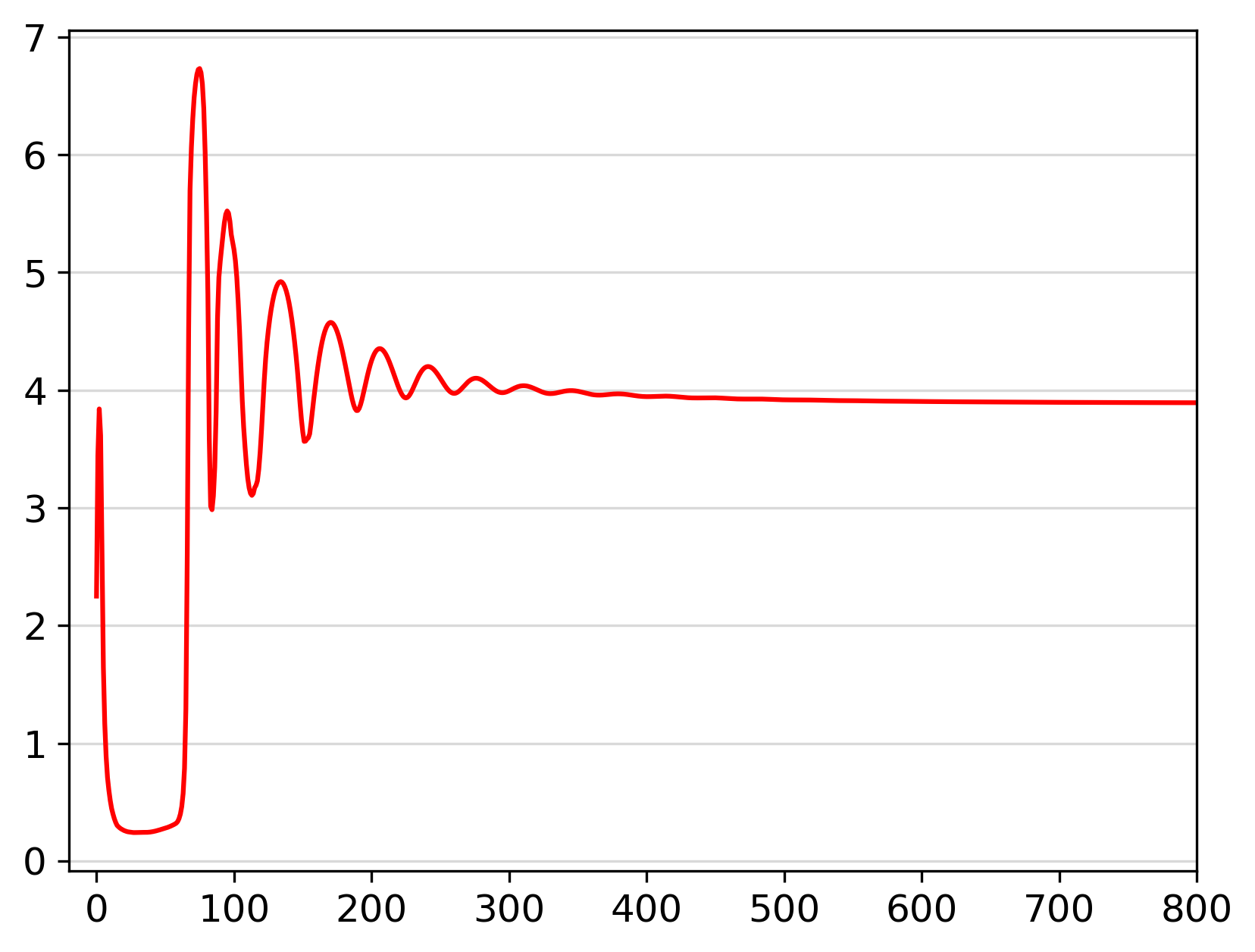}};
        \node[below = 0cm of img1]{iteration};
    \end{tikzpicture}
    \caption{Reward value convergence of ResPG-PD (Algorithm~\ref{alg: resilient PG}, left) and ResOPG-PD (Algorithm~\ref{alg: resilient OPG}, right), with a cost functions $h(\xi) = \alpha \norm{\xi}^2$, where $\alpha = 0.08$,  in the monitoring problem of Section~\ref{sec:largemonitor}. The stepsize for ResPG-PD is $\eta=0.01$ and the stepsize for ResOPG-PD is $\eta=0.05$.
    }
    \label{fig:AppE3RewardvsIter}
    \end{figure}

\begin{figure}[tbh]
    \centering
        \begin{tikzpicture}
        \node[](img1) at(0,0) {\includegraphics[width = 0.33\textwidth]{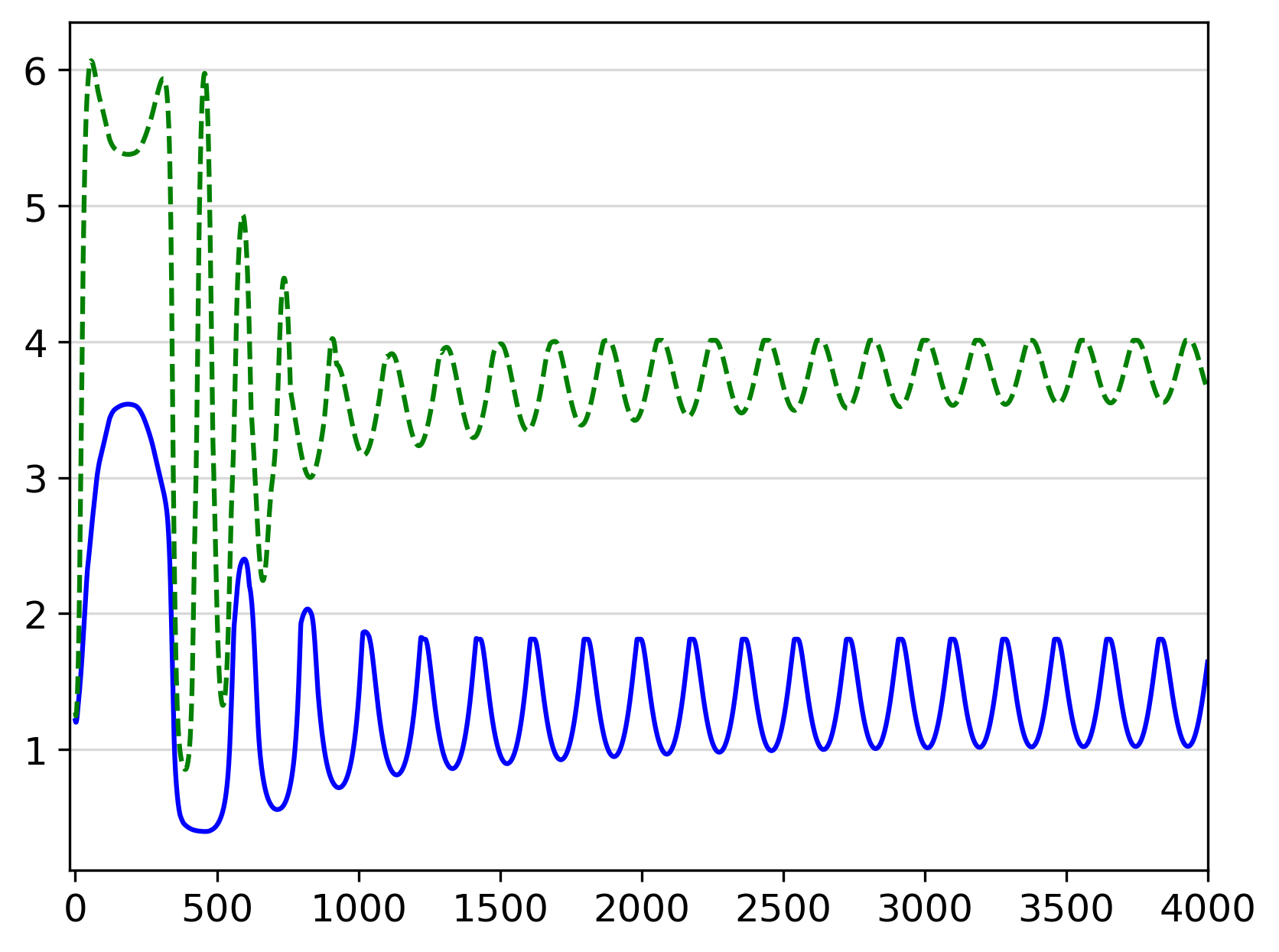}};
        \node[below = 0cm of img1]{iteration};
         \node[left = 0.4cm of img1,yshift=1.3cm, rotate=90]{utility value};
    \end{tikzpicture}
    \hspace{0.1\textwidth}
            \begin{tikzpicture}
        \node[](img1) at(0,0) {\includegraphics[width = 0.33\textwidth]{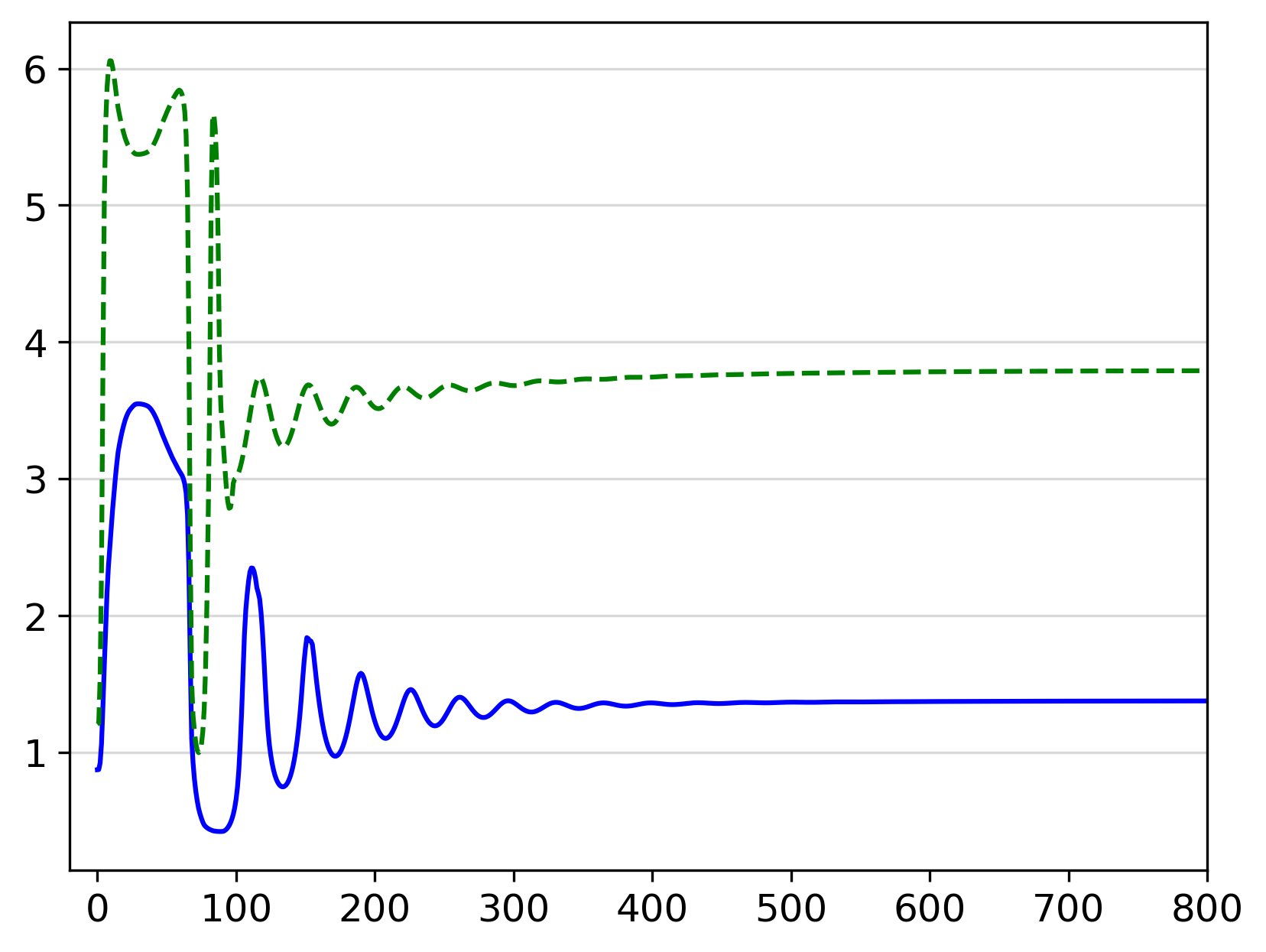}};
        \node[below = 0cm of img1]{iteration};
    \end{tikzpicture}
    \caption{Utility value convergence ($V_{g_1}^\pi(\rho)$:~\ref{legend:blue}, $V_{g_2}^\pi(\rho)$:~\ref{legend:greendash}) of ResPG-PD (Algorithm~\ref{alg: resilient PG}, left) and ResOPG-PD (Algorithm~\ref{alg: resilient OPG}, right), with a cost functions $h(\xi) = \alpha \norm{\xi}^2$ for $\alpha = 0.08$, in the monitoring problem of Section~\ref{sec:largemonitor}. The stepsize for ResPG-PD is $\eta=0.01$ and the stepsize for ResOPG-PD is $\eta=0.05$.
    }
    \label{fig:AppE3UtilityvsIter}
    \end{figure}
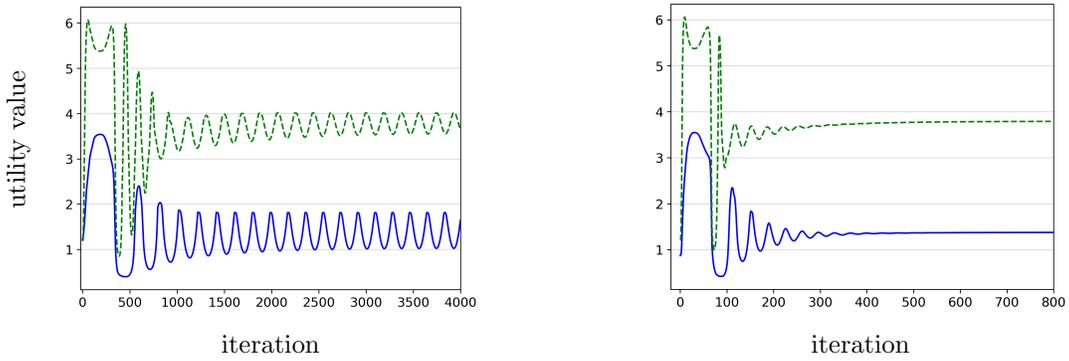

\begin{figure}[tbh]
    \centering
        \begin{tikzpicture}
        \node[](img1) at(0,0) {\includegraphics[width = 0.33\textwidth]{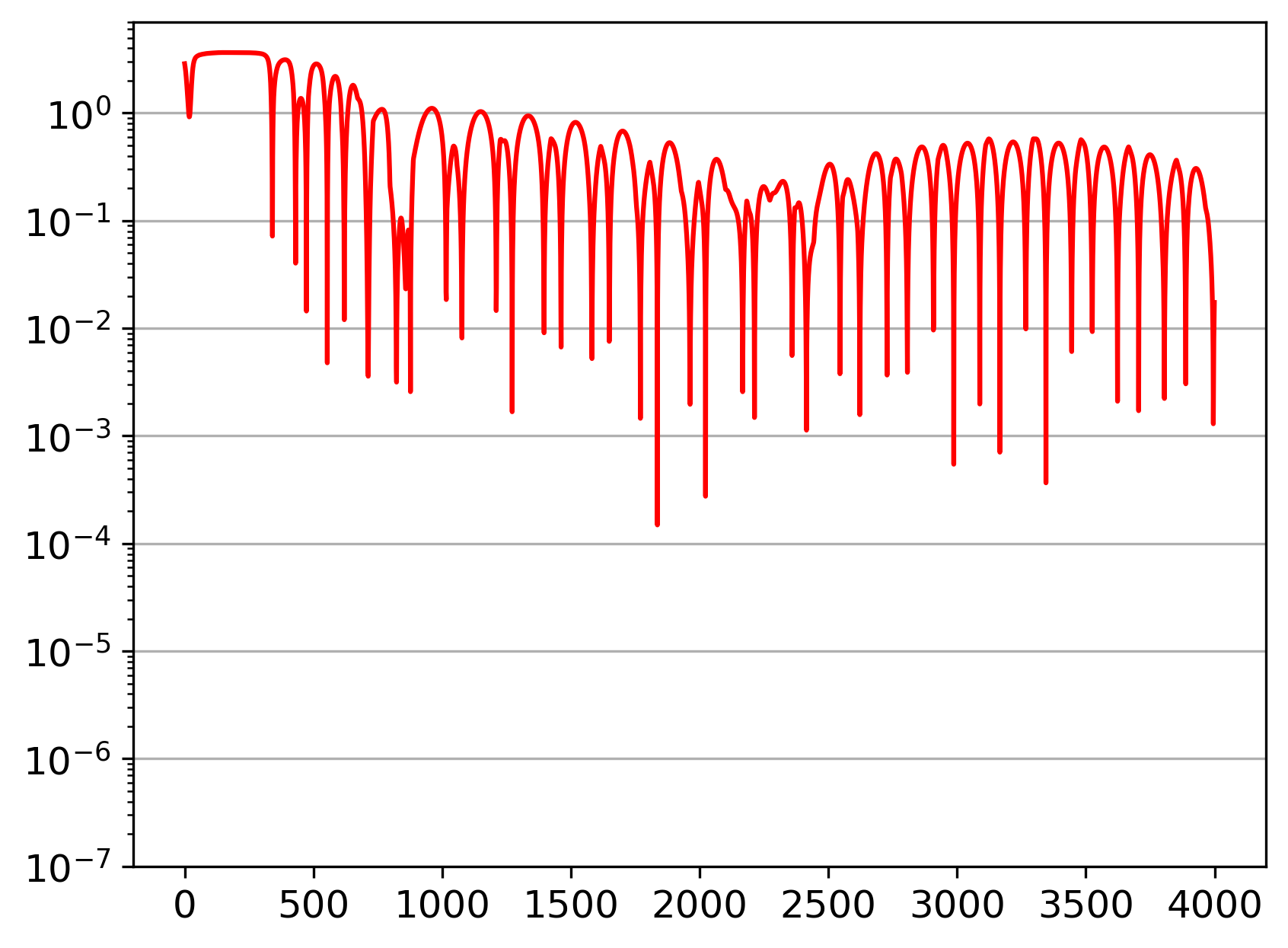}};
        \node[below = 0cm of img1]{iteration};
         \node[left = 0.4cm of img1,yshift=2.1cm, rotate=90]{reward optimality gap};
    \end{tikzpicture}
    \hspace{0.1\textwidth}
            \begin{tikzpicture}
        \node[](img1) at(0,0) {\includegraphics[width = 0.33\textwidth]{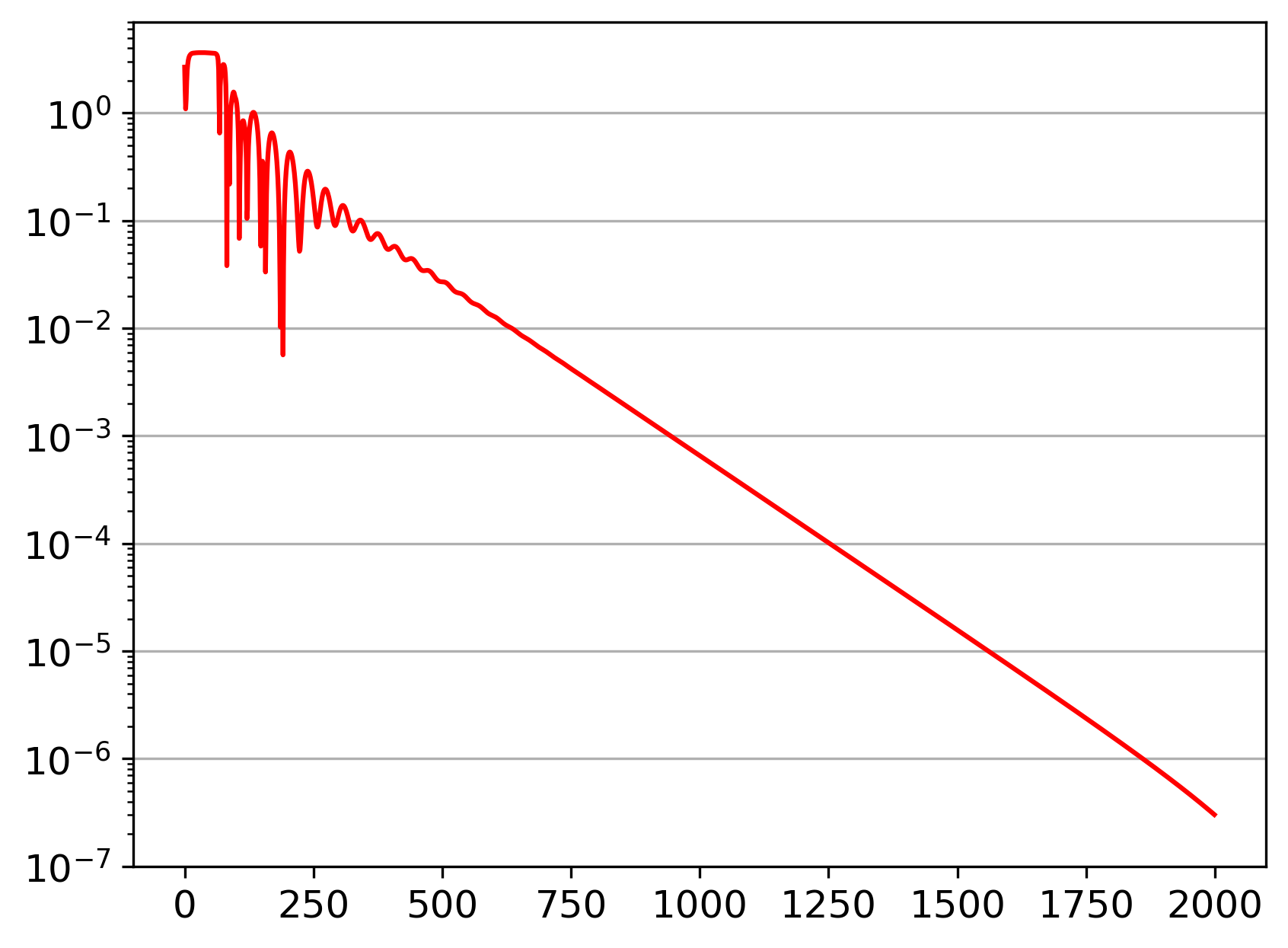}};
        \node[below = 0cm of img1]{iteration};
    \end{tikzpicture}
    \caption{Reward optimality gap of ResPG-PD (Algorithm~\ref{alg: resilient PG}, left) and ResOPG-PD (Algorithm~\ref{alg: resilient OPG}, right), with a cost functions $h(\xi) = \alpha \norm{\xi}^2$, where $\alpha = 0.08$  in the monitoring problem of Section~\ref{sec:largemonitor}. The stepsize for ResPG-PD is $\eta=0.01$ and the stepsize for ResOPG-PD is $\eta=0.05$.
    }
    \label{fig:AppE3RewardErrorvsIter}
    \end{figure}

\begin{figure}[tbh]
    \centering
        \begin{tikzpicture}
        \node[](img1) at(0,0) {\includegraphics[width = 0.22\textwidth]{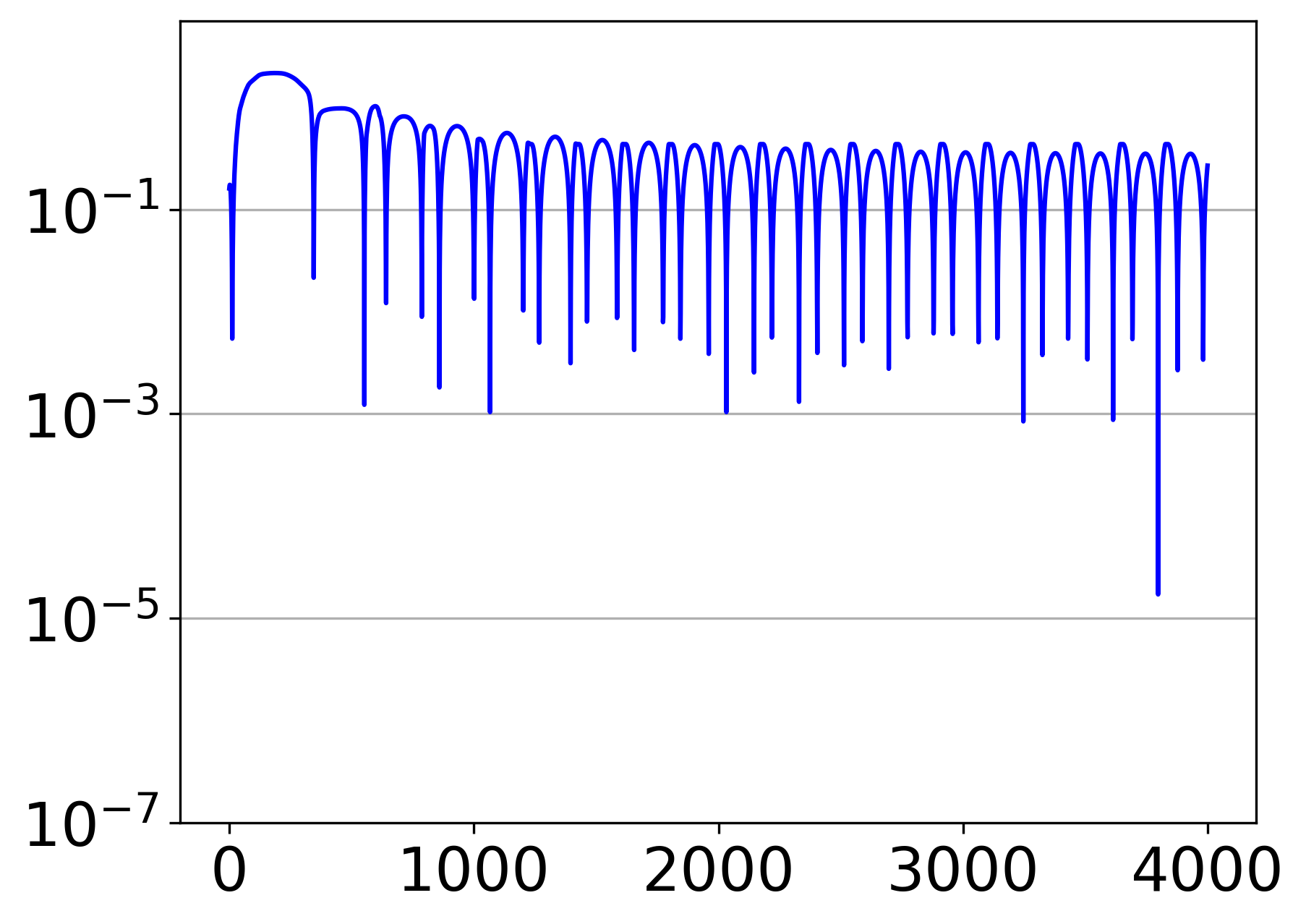}};
        \node[below = 0cm of img1]{iteration};
         \node[left = 0.4cm of img1,yshift=2cm, rotate=90]{utility optimality gap};
    \end{tikzpicture}
            \begin{tikzpicture}
        \node[](img1) at(0,0) {\includegraphics[width = 0.22\textwidth]{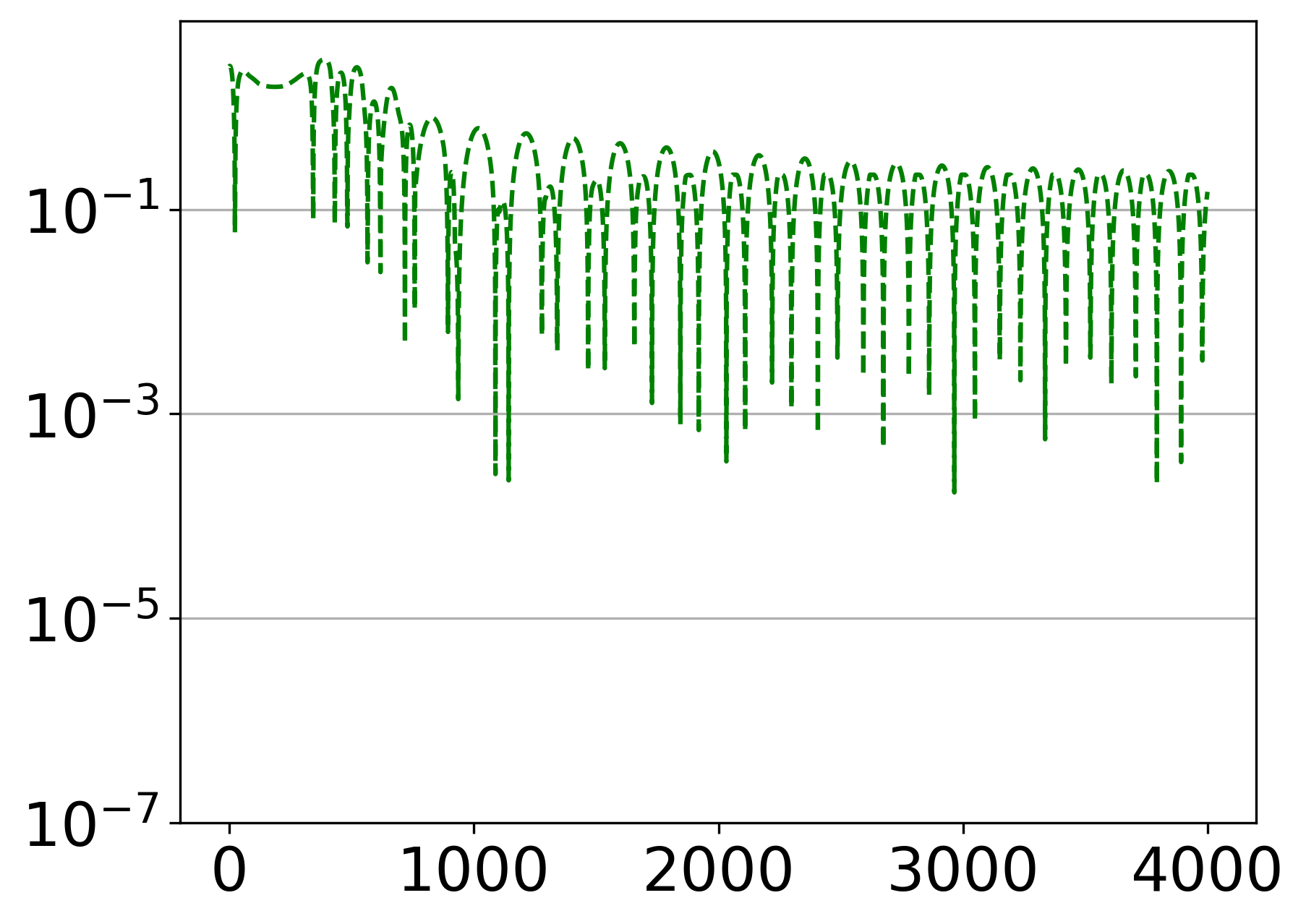}};
        \node[below = 0cm of img1]{iteration};
    \end{tikzpicture}
            \begin{tikzpicture}
        \node[](img1) at(0,0) {\includegraphics[width = 0.22\textwidth]{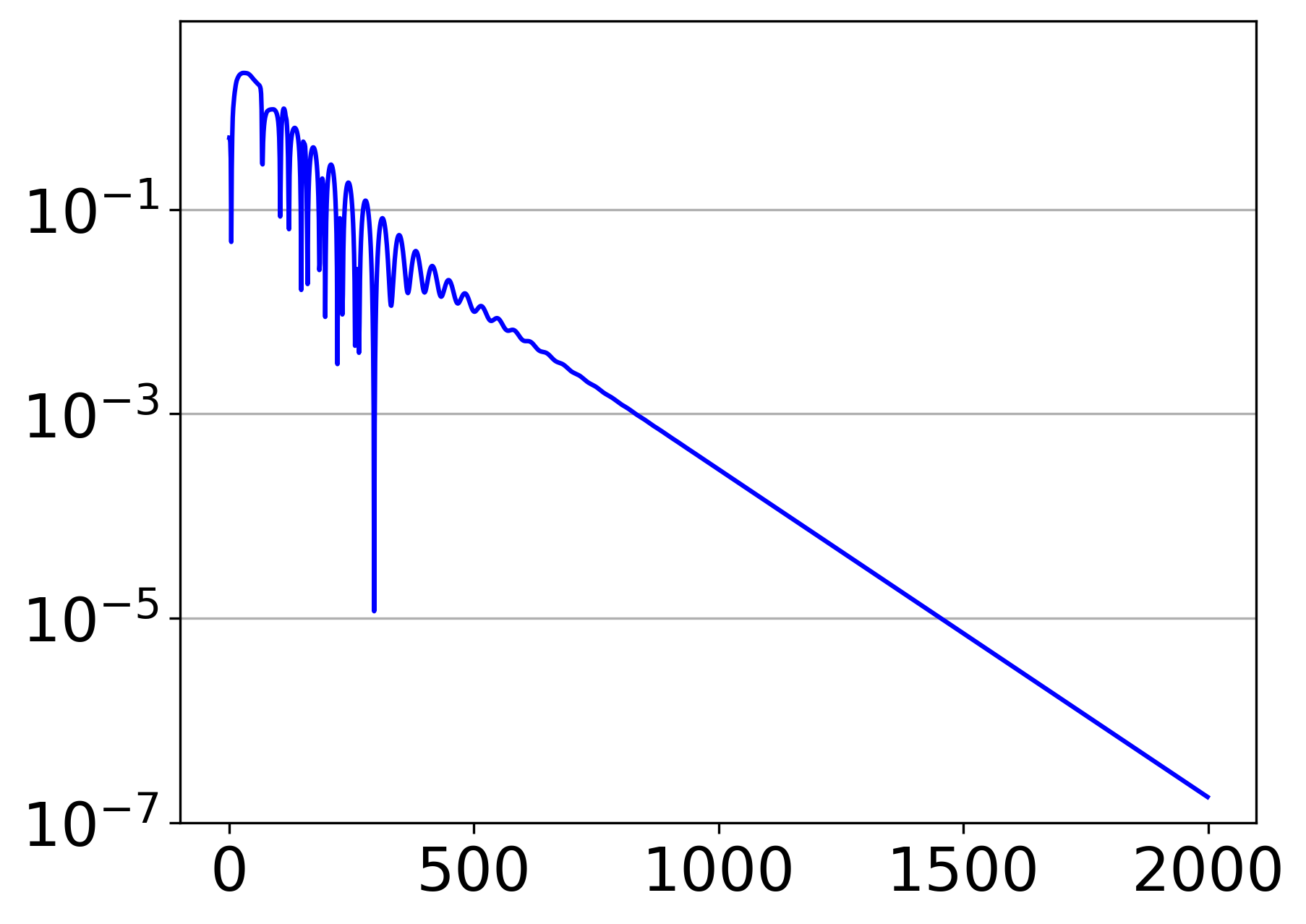}};
        \node[below = 0cm of img1]{iteration};
    \end{tikzpicture}
            \begin{tikzpicture}
        \node[](img1) at(0,0) {\includegraphics[width = 0.22\textwidth]{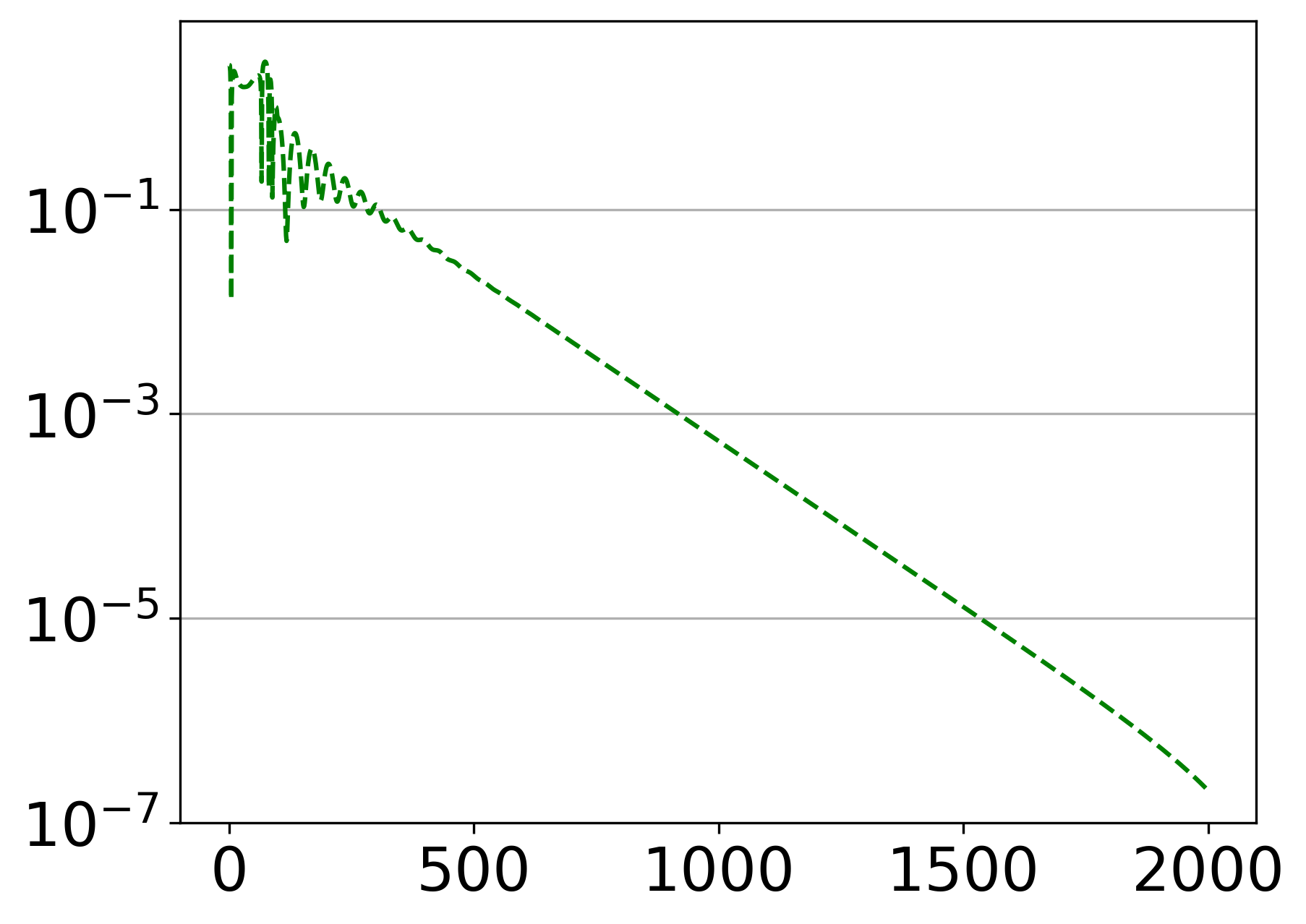}};
        \node[below = 0cm of img1]{iteration};
    \end{tikzpicture}
    \caption{Utility optimality gap ($V_{g_1}^\pi(\rho)$:~\ref{legend:blue}, $V_{g_2}^\pi(\rho)$:~\ref{legend:greendash}) of ResPG-PD (Algorithm~\ref{alg: resilient PG}, two figures on the left) and ResOPG-PD (Algorithm~\ref{alg: resilient OPG}, two figures on the right), with a cost functions $h(\xi) = \alpha \norm{\xi}^2$, where $\alpha = 0.08$  in the monitoring problem of Section~\ref{sec:largemonitor}. The stepsize for ResPG-PD is $\eta=0.01$ and the stepsize for ResOPG-PD is $\eta=0.05$.
    }
    \label{fig:AppE3UtilityErrorvsIter}
    \end{figure}
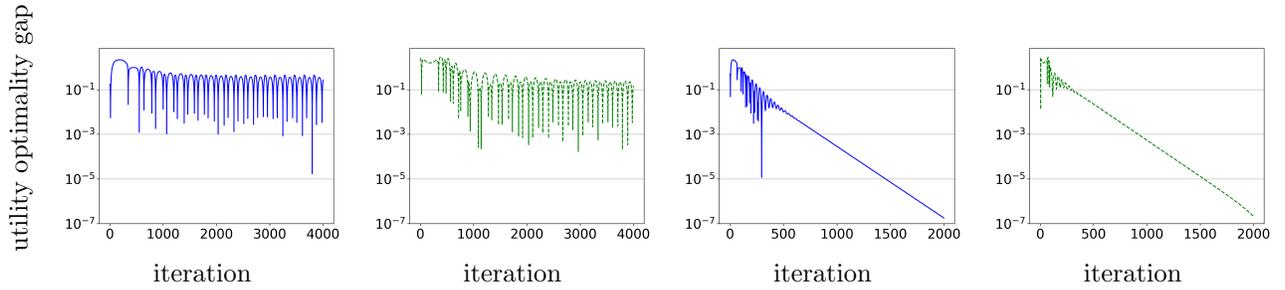

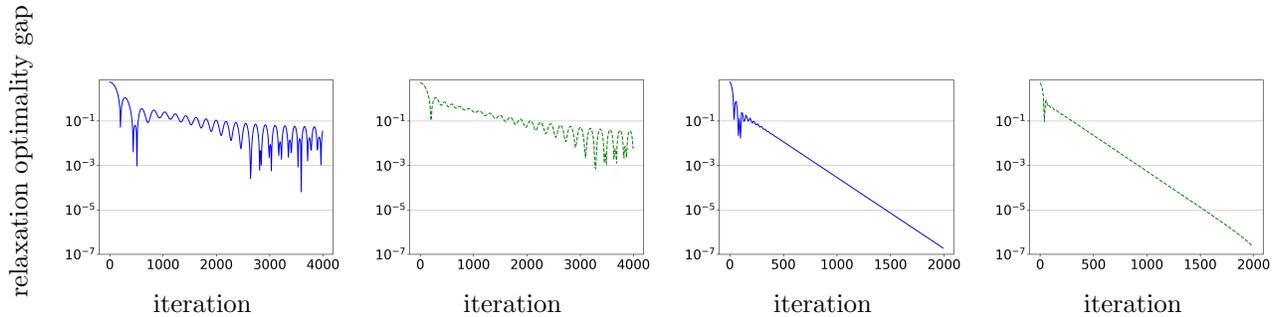
\begin{figure}[tbh]
    \centering
        \begin{tikzpicture}
        \node[](img1) at(0,0) {\includegraphics[width = 0.22\textwidth]{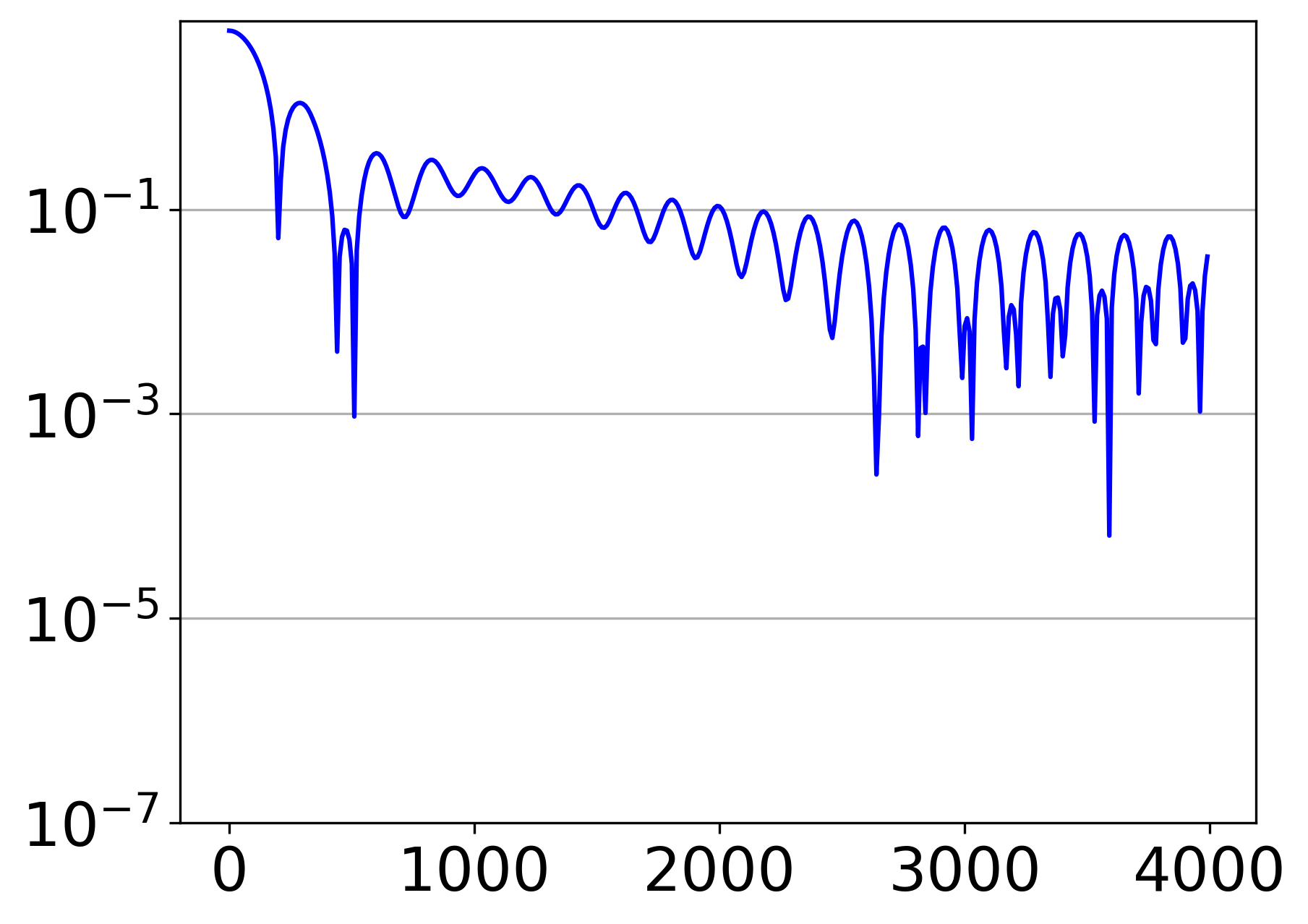}};
        \node[below = 0cm of img1]{iteration};
         \node[left = 0.4cm of img1,yshift=2.4cm, rotate=90]{relaxation optimality gap};
    \end{tikzpicture}
            \begin{tikzpicture}
        \node[](img1) at(0,0) {\includegraphics[width = 0.22\textwidth]{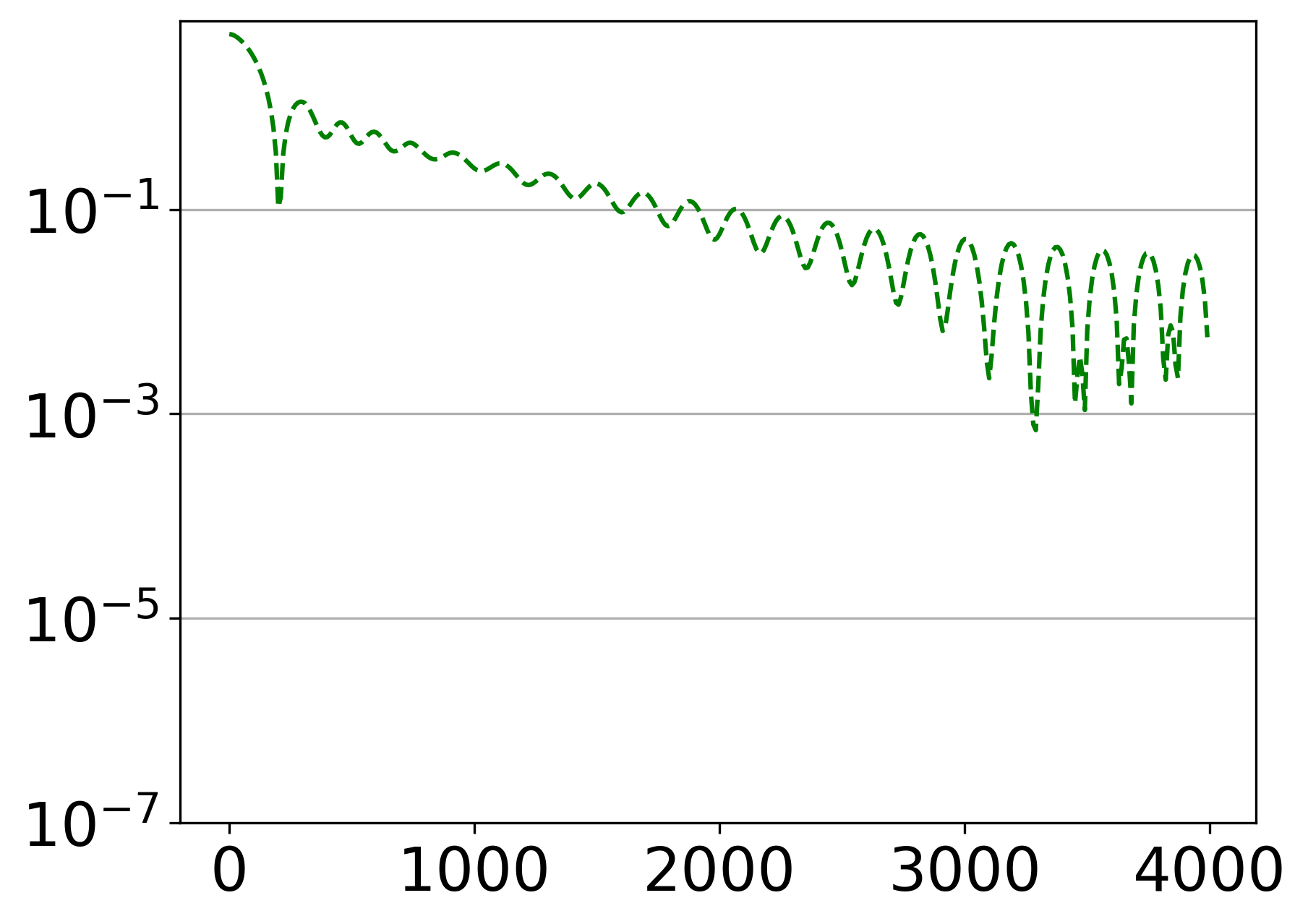}};
        \node[below = 0cm of img1]{iteration};
    \end{tikzpicture}
            \begin{tikzpicture}
        \node[](img1) at(0,0) {\includegraphics[width = 0.22\textwidth]{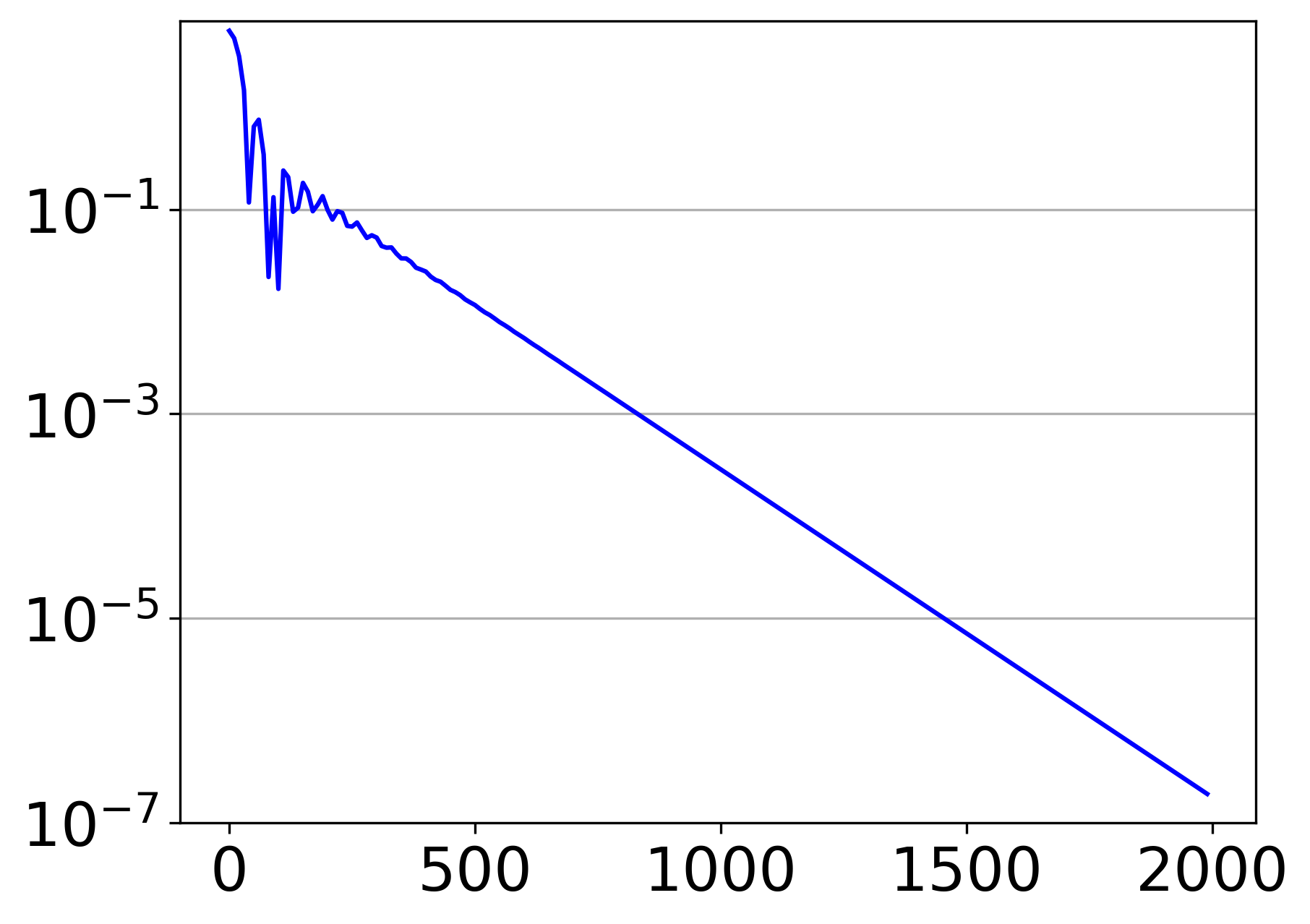}};
        \node[below = 0cm of img1]{iteration};
    \end{tikzpicture}
            \begin{tikzpicture}
        \node[](img1) at(0,0) {\includegraphics[width = 0.22\textwidth]{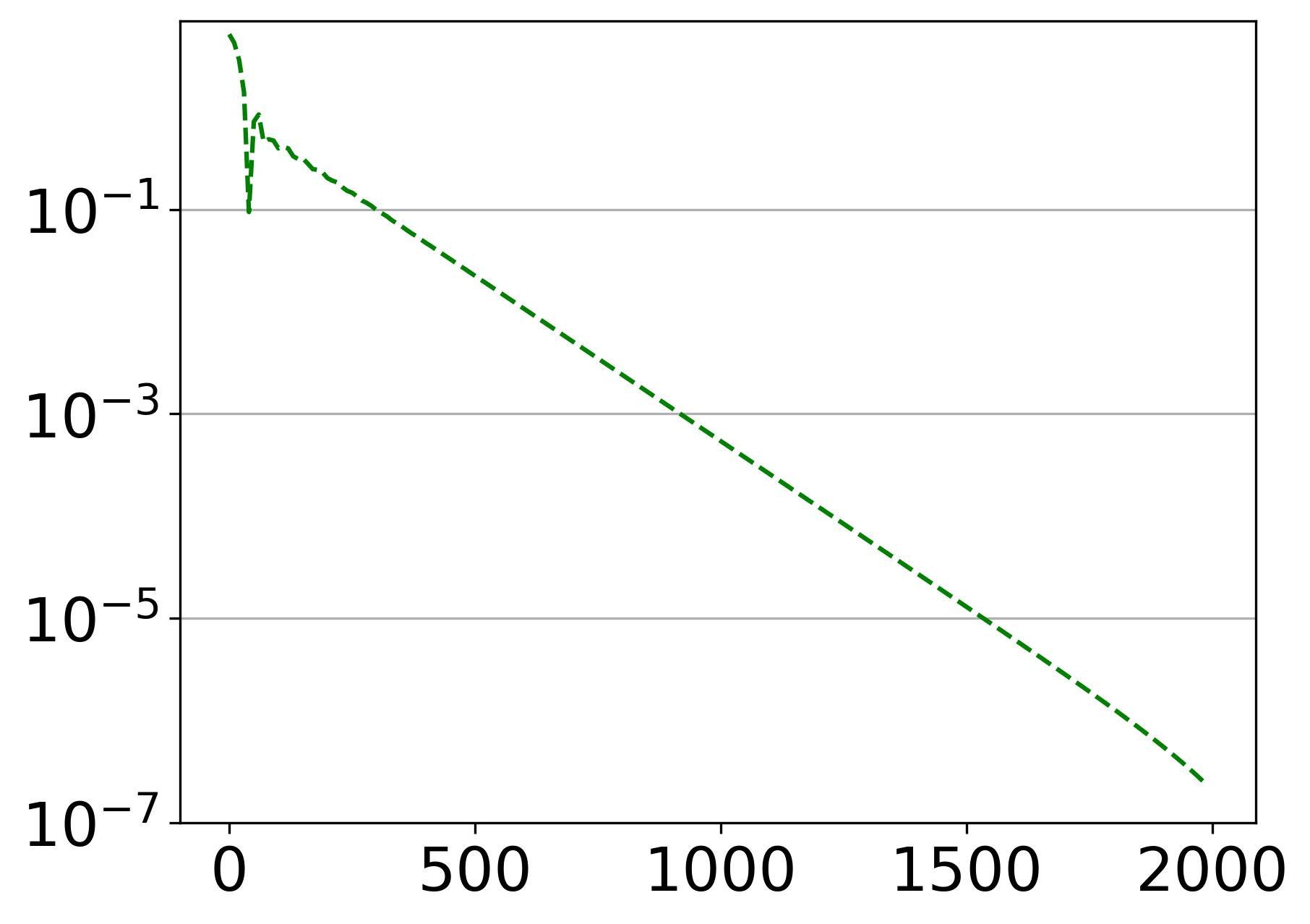}};
        \node[below = 0cm of img1]{iteration};
    \end{tikzpicture}
    \caption{Relaxation optimality gap ($\xi_1$: \ref{legend:blue}, $\xi_2$: \ref{legend:greendash}) of ResPG-PD (Algorithm~\ref{alg: resilient PG}, two figures on the left) and ResOPG-PD (Algorithm~\ref{alg: resilient OPG}, two figures on the right), with a cost functions $h(\xi) = \alpha \norm{\xi}^2$, where $\alpha = 0.08$  in the monitoring problem of Section~\ref{sec:largemonitor}. The stepsize for ResPG-PD is $\eta=0.01$ and the stepsize for ResOPG-PD is $\eta=0.05$.
    }
    \label{fig:AppE3XiErrorvsIter}
    \end{figure}

As previously shown in Figure~\ref{fig:ResE3}, both methods can relax the infeasible constraints. We then report reward and utility convergence performance for the two methods in Figure~\ref{fig:AppE3RewardvsIter}--\ref{fig:AppE3UtilityvsIter} and report the reward, utility and relaxation optimality gaps of two methods in Figures~\ref{fig:AppE3RewardErrorvsIter}--\ref{fig:AppE3XiErrorvsIter}. The behavior of two methods in a large state/action space is the same as in Section~\ref{sec:smallmonitor} as expected.

\section{Some Useful Lemmas}

\begin{lemma}\label{lem:convex_violation_bound}
    Let Assumption~\ref{as:feasibility_regularized} hold for Problem~\eqref{eq:CMDP_regularized}. For any $C\geq 2\bar\lambda^\star$, if there exists a policy $\pi\in\Pi$, $\xi\in\Xi$, and $\delta>0$ such that $V_r^\star(\rho)-h(\bar\xi^\star)-(V_r^{\pi}(\rho)-h(\xi)) + C\sum_{i\,=\,1}^m[\xi_i - V_{g_i}^{\pi}(\rho)]_+ \leq \delta$, then $\sum_{i\,=\,1}^m[\xi_i - V_{g_i}^{\pi}(\rho)]_+ \leq 2\delta/C$, where $[x]_+ = \max(x,0)$.
\end{lemma}
\begin{proof}
    For Problem~\eqref{eq:CMDP_regularized}, we introduce its value function as
    \[
    v(\tau) 
    \; = \;
    \maximize_{\pi\,\in\,\Pi,\xi\,\in\,\Xi} 
    \;
    \left\{\,
    V_r^\pi(\rho) - h(\xi)
    \,\vert\, V_{g_i}^\pi(\rho) \geq \xi_i +\tau_i, i = 1,\ldots,m
    \,\right\}
    \]
    We note that, $v(\tau)$ is a concave function, which is similar as Lemma~\ref{lem:primal_function}, and $v(0) = V_h^\star = V_r^\star(\rho) - h(\bar\xi^\star)$. The rest of proof is straightforward from~\cite[Lemma~4]{ding2022convergence}. 
\end{proof}

For any convex differentiable function $\psi$: $X\to\mathbb{R}$, the Bregman divergence of $x$, $x'\in X$ is given by $D_{\psi}(x',x) \DefinedAs \psi(x')  - \psi(x)  - \langle \nabla\psi(x), x'-x\rangle$. When $\psi$ is $\sigma$-strongly convex, $D_\psi(x',x) \geq \frac{\sigma}{2}\norm{x-x'}^2$ for any $x$, $x'\in X$. Specifically, when $\psi(x) = \frac{1}{2}\norm{x}^2$, $D_\psi(x',x) =  \frac{1}{2}\norm{x'-x}^2$.

\begin{lemma}\label{lem:three-point}
    Let $X$ be a convex set. If $x' = \argmin_{\bar x\,\in\,X} \langle \bar x, g\rangle + D_\psi(\bar x, x)$, then for any $x^\star\in X$,
    \[
    \langle x' - x^\star, g\rangle 
    \; \leq  \;
    D_\psi(x^\star,x) -D_\psi(x^\star,x')-D_\psi(x',x).
    \]
\end{lemma}
\begin{proof}
    See~\cite[Lemma~10]{wei2020linear}.
\end{proof}

\begin{lemma}\label{lem:non-expansive}
    Assume that $D_\psi(x,x')  \geq \frac{1}{2}\norm{x-x'}_p^2$ for some $\psi$ and $p\geq 1$. If 
    \[
    x_1 \; = \; \argmin_{\bar x\,\in\,X}\; \langle \bar x, g_1\rangle + D_\psi(\bar x, x)
    \;\text{ and }\;
    x_2 \; = \; \argmin_{\bar x\,\in\,X}\; \langle \bar x, g_2\rangle + D_\psi(\bar x, x)
    \]
    then,
    \[
    \norm{x_1-x_2}_p \; \leq  \; \norm{g_1-g_2}_q
    \]
    where $\frac{1}{p}+\frac{1}{q} = 1$.
\end{lemma}
\begin{proof}
    See~\cite[Lemma~11]{wei2020linear}.
\end{proof}

\begin{lemma}\label{lem:policy-value-difference}
    For any two policies $\pi$ and $\pi'$, we have
    \[
    \begin{array}{rcl}
         \norm{Q_r^\pi(\cdot,\cdot) - Q_r^{\pi'}(\cdot,\cdot)}_\infty
         & \leq &
         \displaystyle
         \frac{\gamma}{(1-\gamma)^2} \max_s \;\norm{\pi(\cdot\,\vert\,s) - \pi'(\cdot\,\vert\,s)}_1
         \\[0.3cm]
         |V_r^{\pi}(\rho) - V_r^{\pi'}(\rho)|& \leq &
         \displaystyle
         \frac{\kappa_\rho}{(1-\gamma)^3}
         \sum_{s} d_\rho^{\pi^\star}(s) \norm{\pi(\cdot\,\vert\,s) - \pi'(\cdot\,\vert\,s)}_1.
    \end{array}
    \]
\end{lemma}
\begin{proof}
    See~\cite[Lemma~11]{ding2023last}.
\end{proof}

\begin{lemma}\label{lem:optimistic_last_rate}
    Let a non-negative sequence $\{ B_t \}_{t\,\geq\,1}$ satisfy that for any $p>0$ and $q>0$,
    \[
    B_{t+1}\; \leq \; B_t - qB_{t+1}^{p+1}
    \; \text{ and } \;
    q(1+p) B_1^p \; \leq \;1.
    \]
    Then, $B_t \leq C\, t^{-1/p}$, where $C \DefinedAs \max (B_1, (1/(pq))^{1/p})$.
\end{lemma}
\begin{proof}
    See~\cite[Lemma~12]{wei2020linear}.
\end{proof}

    \clearpage

    \vfill
\clearpage

\end{document}